\documentclass[12pt]{article}
\usepackage{amsfonts}
\usepackage{amssymb}
\usepackage{amsmath}
\usepackage{txfonts}
\usepackage{tipa}
\usepackage{graphicx,color,overpic}
\usepackage{lscape}
\usepackage{CJK}
\usepackage{epstopdf}
\usepackage{booktabs}
\usepackage{multirow}
\usepackage{floatrow}
\newtheorem{theorem}{Theorem}  
\newtheorem{lemma}{Lemma}

\newtheorem{proposition}{Proposition}
\newtheorem{definition}{Definition}
\newtheorem{example}{Example}
\newtheorem{continue example}{Continue Example}
\newtheorem{remark}{Remark}

\newcommand{\ba}{\begin{array}}
\newcommand{\ea}{\end{array}}
\newcommand{\bt}{\begin{tabular}}
\newcommand{\et}{\end{tabular}}
\newcommand{\btb}{\begin{table}}
\newcommand{\etb}{\end{table}}
\newcommand{\bc}{\begin{center}}
\newcommand{\ec}{\end{center}}
\newcommand{\bea}{\begin{eqnarray}}
\newcommand{\eea}{\end{eqnarray}}
\newcommand{\Bea}{\begin{eqnarray*}}
\newcommand{\Eea}{\end{eqnarray*}}
\newcommand{\beq}{\begin{equation}}
\newcommand{\eeq}{\end{equation}}

\def \bfm#1{\mbox{\boldmath$#1$}}

   \def \b {{\bfm b}}
\def \B {{\bfm B}}  
\def \d {{\bfm d}}  
\def \E {{\bfm E}}

\def \N {{\bfm N}}  \def \0 {{\bfm 0}}

\def \x {{\bfm x}} 
\def \Y {{\bfm Y}} 
\def \Z {{\bfm Z}}

\def \one {{\bf 1}} \def \zero {{\bf 0}}

\def \ll {{\bfm l-1}}

\def \ca {{\mathcal A}} \def \cC {{\mathcal C}}
  
\def \cf {{\mathcal F}}
  
  \def \cn {{\mathcal N}}
 
\def \cu {{\mathcal U}}  

\def \cv {{\mathcal V}}
\def \ct {{\mathcal T}}

\def \aug {augmented }
\def \Aug {Augmented }
\def \caug {column augmented design}
\def \caugs {column augmented designs}

\newcommand{\revA}[1]{{\color{red} #1}}

\renewcommand{\revA}[1]{{#1}}

\begin{document}

\baselineskip 16 pt
\title{ {\Large \bf   Mixed-Level  Column \Aug Uniform Designs }}

\author{{ Feng Yang$^a$,\ \
Yong-Dao Zhou$^b$\thanks{Corresponding author. Email: ydzhou@nankai.edu.cn}},\ \ Aijun Zhang$^c$ \\
{\small \it $^a$  College of Mathematics, Sichuan University,
Chengdu
610064, China} \\
{\small \it $^b$ School of Statistics and Data Science \& LPMC,
Nankai
University, Tianjin 300071, China } \\
{\small \it $^c$ Department of Statistics and Actuarial Science, The
University of Hong Kong, China} }

\maketitle

\noindent \emph{Abstract}: 
Follow-up experimental designs are popularly used in industry. In many
follow-up designs, some additional factors with two or three levels may 
be added in the follow-up stage since they are quite important but may 
be neglected in the first stage. Such follow-up designs are called mixed-level 
column augmented designs. In this paper, based on the initial designs,
mixed-level column augmented uniform designs are proposed by using 
the uniformity criterion, wrap-around $L_2$-discrepancy (WD). The 
multi-stage augmented procedure which adds the additional design points
stage by stage is also investigated. We present the analytical expressions 
and the corresponding lower bounds of the WD of the column augmented 
designs. It is shown that the column \aug uniform designs are also the
optimal designs under the non-orthogonality criterion, $E(f_{NOD})$.
Furthermore, a construction algorithm for the column augmented uniform 
design is provided. Some examples show that the lower bounds
are tight and the construction algorithm is effective.

\vspace{2mm}


\section{Introduction} \label{sec intro}

In many industrial experiments such as semiconductor fabrication,
the cost of each run is very expensive; therefore, small run sizes
are preferred. Assume the number of runs in the initial design is
much smaller than the number of parameters in the full second-order
model, where the second-order model is commonly used to explore 
the nonlinear relationship between factors and responses. If analytical 
results based on the initial design satisfy the objective, experiment is
terminated. Otherwise, a follow-up design may be needed to collect 
further information from the experimental system. Sometimes, a 
two-stage design may still not satisfy the needs, then more design 
points may be added stage by stage until meeting certain requirements. 
Such experimental strategy is called multi-stage \aug design.

Usually, in the field of the follow-up design, the researchers only
considered adding some runs to the initial design but did not add
any factor; see \cite{LL03,LLY03,EQ15a} and \cite{QCO13}. In many 
cases, however, the experimenters may only select some factors of 
interest, and ignore other possibly important factors because of the
limited run size. For instance, Lai et al. \cite{LPT03} considered a real-world 
application of column augmented designs where there are seven factors
$\{x_1,\dots,x_7\}$ in the initial design. Afterwards, in the second 
stage, a new glycerol factor $x_8$ was added to examine the possibility 
of enhancing the response lovastatin production. It turns out that the 
addition of $x_8$ is able to improve production significantly. Another 
example is from an industrial production, where the factor of reaction 
pressure is not considered to be significant initially. \revA{Thus, this} factor was 
fixed as the standard atmosphere pressure (0.1 MPa), in order to limit 
the number of runs for saving cost. However, after analysis of the initial 
experiment, it shows that the reaction pressure may be an important 
factor and need further investigation in the follow-up stage. Usually,
the initial design and the follow-up portions are usually conducted at 
different time, with different equipment or under different operators. 
The researchers may want to test whether there is a system difference 
or not for these nuisance variables between the different stages. Then 
a blocking factor should be added in the follow-up designs to analyze 
the significance of the nuisance variables. Therefore, in the follow-up 
stages, the experimenter may consider not only some additional runs 
but also some additional factors.

\revA{Yang et al.}  \cite{YZZ17} considered augmenting not only some
 number of runs but also some two-level factors, in which the design structure 
 for the two-level factors \revA{is} very special such that they are completely
correlated. Compared to \cite{YZZ17}, this paper discusses more
complex situations with increased technical difficulty. Adding three-level 
factors \revA{is worth} to be considered since they can be used to explore 
the nonlinear relationship based on second-order models. We consider 
the following types of follow-up designs:

\begin{enumerate}
\item Augmenting several three-level factors to an initial design.
\item Augmenting one or two blocking factors in the follow-up steps  
for assessing the system deviation.
\item Constructing multi-stage augmented designs.
\end{enumerate}
\noindent 
The blocking factors can be also treated as common two-level
additional factors. Since the two-level, three-level and blocking
factors may be added in the next stage, such follow-up designs are
called mixed-level column augmented \revA{designs, which} will have a wide range of applications. In particular, the
resulting \caugs~do not have any two completely correlated factors,
which overcomes the shortcoming of that in \cite{YZZ17}.

As we know, for those expensive and time-consuming experiments, 
optimal designs are a widely used type of designs; see
\cite{W43,KW59,F72} and \cite{P93}. However, optimal designs need
to specify a priori underlying models. When the relationship
between the factors and the response is unknown, uniform design
proposed by \cite{F80} is an effective experimental method. The
main idea of uniform design is to scatter the design points uniformly 
on the experimental domain. Yue and Hickernel \cite{YH99} showed that 
uniform design 
is robust against model specification. Xie and Fang \cite{XF00} showed that uniform 
designs have the property of admissibility and minimaxity. Moreover, 
the number of runs in uniform designs is very flexible and can be 
chosen to be any integer. This appealing property often allows us to 
save experimental costs. Furthermore, the uniformity criterion
has close relationship with the generalized minimum aberration
criterion, which is widely used in orthogonal design theory; see
\cite{FM00} and \cite{ZX14}. Hence, the uniformity criterion is a
reasonable consideration for assessing the goodness of the column
augmented designs especially when the true model is unknown.
In this paper, we discuss the best mixed-level \caugs~under the 
uniformity criterion, wrap-around $L_2$-discrepancy (WD, Hickernell, 
1998), and call the resulting designs as mixed-level column \aug 
uniform designs.

Apart from the uniformity criterion, the orthogonality is also an important
assessing criterion in experimental designs. $E(f_{NOD})$ criterion is 
popularly employed for comparing different designs from the viewpoint 
of non-orthogonality; see \cite{FLL03b}, \cite{KM05} and \cite{SLL11}. 
The $E(f_{NOD})$-value of any design is nonnegative, and less
$E(f_{NOD})$-value means the designs have better orthogonality. We
will show that there exists some interesting relationship between
the uniformity criterion WD and the $E(f_{NOD})$ criterion of the \caugs.

The rest of this paper is organized as follows. Section \ref{sec wd}
gives the expressions of WD and the corresponding lower bounds of
\caug~with three-level additional factors and at most two additional 
blocking factors. Moreover, the multi-stage situation is also discussed. 
The $E(f_{NOD})$ criterion of the \caugs~is discussed in Section 
\ref{sec f_NOD}, where the connection between WD and $E(f_{NOD})$ 
is also established. Section \ref{sec algorithm} presents the construction 
algorithm for the column augmented uniform designs, and shows some
examples to demonstrate that the construction algorithm is powerful and 
the lower bounds in this paper are relatively tight. Some conclusions are 
summarized in Section \ref{sec conclusions}. Some proofs are in Appendix A
and the design matrices mentioned in Section \ref{sec algorithm} are
listed in \revA{Appendix B}.

\section{WD criterion for \caug } \label{sec wd}

According to the projection uniformity on one dimension, it is preferred to 
restrict the designs to be balanced, i.e., U-type designs. An asymmetric 
U-type design $U(n;q_1,\ldots, q_m)$ corresponds to a $n\times m$ matrix 
$X=(\x_1,\ldots,\x_m)$, such that each column $\x_i$ takes values from 
a set of $q_i$ integers, say $\{0,1,\ldots,q_i-1\}$, equally often. If some 
$q_i$'s are equal, we denote this asymmetrical U-type design by
$U(n;q_1^{m_1},\ldots,q_s^{m_s})$, where $m=\sum_{i=1}^sm_i$. Denote
all of the $U(n;q_1,\ldots, q_m)$ and $U(n,q_1^{m_1},\ldots,q_s^{m_s})$ 
by $\cu(n;q_1,\ldots, q_m)$ and $\cu(n,q_1^{m_1},\ldots,q_s^{m_s})$, 
respectively. For each $\d\in\cu(n;q_1,\ldots, q_m)$, the $n$ runs of $\d$
transform into $n$ points in $C^m=[0,1]^m$ by mapping $f: x_{ik} 
\longrightarrow (2x_{ik}+1)/(2q_k)$. There are different criteria for 
measuring the uniformity. Among them the WD has many good properties; 
for example, it is invariant under reordering the runs, relabeling
coordinates and coordinate shift. The squared WD-value of $\d\in
\cu(n;q_1,\ldots, q_m)$ is $\mbox{WD}(\d)=-\left(\frac{4}{3}\right)^m+\frac{1}
{n^2} \sum_{i=1}^{n}\sum_{j=1}^{n}\prod_{k=1}^m\left(\frac32-\left|u_{ik}-u_{jk}
\right|\left(1-\left|u_{ik}-u_{jk}\right|\right)\right),
$ where $u_{ik}=(2x_{ik}+1)/(2q_k).$

\subsection[CAUD with additional three-level factors]{Column augmented 
uniform designs with additional three-level factors}
\label{sec column augmented}

Let the initial design $\d_0 \in \cu(n;2^{m_1}3^{m_2})$, $m=m_1+m_2$. 
In the follow-up stage, one wants to add $n_1$ additional runs and $r$ 
additional three-level factors. It is reasonable to assume that each 
follow-up stage may not augment too many runs. Without loss of 
generality, let $n_1\le n$. Denote $\zero_{t\times k}$ and $\one_{t\times k}$ 
be the $t\times k$ matrix whose elements are zeros and ones, respectively.
A $t\times k$ design matrix $\d\in \{0,1,\ldots,q-1\}^{t\times k}$ means that 
each element of $\d$ is chosen from the set $\{0,1,\ldots,q-1\}$.

\begin{definition} \label{defi1}
A design $D_{3}=\left(\ba{cc} \d_0 &\zero_{n\times r} \\ \d_1 &\d_2 \ea\right)$ 
is a \emph{column augmented design with additional three-level factors}, 
if the initial design $\d_0 \in \cu(n;2^{m_1}3^{m_2})$ is augmented with 
$\d_1 \in \cu(n_1;2^{m_1}3^{m_2})$, and $\d_2 \in \{0,1,2\}^{n_1\times r}$.
Denote all such column augmented designs  by $\cC_{3}(n+n_1;2^{m_1}
3^{m_2}\bullet 3^{r})$.
\end{definition}

\begin{remark} Each of the $r$ additional three-level factors can be 
a common quantitative or qualitative factor, or even a blocking factor. 
If $r=0$, the \caugs~become the row \aug designs defined by \cite{YZZ17}.
\end{remark}

The explanation of Definition \ref{defi1} is as follows. In the initial design, 
the level of initial ignored factors $\{x_{m+1},\cdots,x_{m+r}\}$ is usually 
fixed. For example, if the temperature is ignored in the initial stage, the 
researchers often set it to the room temperature in the initial $n$ runs. 
It can be easily shown that the WD-values of designs are not changed 
for mixed-level designs when permuting the levels of each column.
Without loss of generality, assume all the levels of the additional factors 
$\{x_{m+1},\cdots,x_{m+r}\}$ be labeled as 0 in the initial design,
then the initial design can be represented by $(\ba{ccc} \d_0
&\zero_{n\times r} \ea)$ with the $m+r$ factors. The design matrix in 
the follow-up stage is $(\ba{ccc} \d_1 &\d_2 \ea)$, where $\d_1$ and 
$\d_2 $ are the design matrices in the second stage for the initial $m$ 
factors and the additional $r$ factors, respectively. When $n_1\le n$, 
to make the column augmented design as uniform as possible, the 
additional part $\d_1$ should be limited to be U-type, and each column 
of the $r$ columns, $\d_2$, should occur the same number of 1 and 2 
in the follow-up stage, because the design can be more uniform when 
the elements in each column are more balanced. Thereby, the number 
of adding runs $n_1$ meets the following requirements,

\bc {\small \tabcolsep=8pt \bt {c|l|l} \hline $n_1$ & the cases of $m_1$ 
and $m_2$ & the cases of $r$
\\\hline
multiple of 2 & $m_2=0$   &  $r\geq 0$ \\
\hline multiple of 3 & $m_1=0$   & $r=0$  \\ \hline
 multiple of 6  & $m_1=0$    & $r>0$  \\
         & $m_1>0,m_2>0$ & $r\geq 0$ \\
\hline \et} \ec 
Compared with other kinds of follow-up design, such as foldover and 
semifoldover designs which respectively limit $n_1=n$ and $n_1=n/2$, 
our requirement for the number of additional runs is more flexible. 
Moreover, it should be mentioned that the restriction of $\d_1$ to 
be a U-type design can be relaxed, i.e., one can augment any number 
of runs $n_1$ based on the initial design. In the rest of the paper, 
we will consider the cases when $\d_1$ is U-type.

According to Definition \ref{defi1}, there are many alternative column 
\aug designs for a given initial design $\d_0$. Under the uniformity 
criterion, one tends to add the follow-up part $(\ba{cc} \d_1 &\d_2\ea)$ 
such that the column \aug design $D_{3}$ with $n+n_1$ points is as 
uniform as possible.

\begin{definition} \label{defi2} A column augmented design $D_{3}^*\in 
\cC_{3}(n+n_1;2^{m_1}3^{m_2}\bullet 3^{r})$ is a \emph{column 
augmented uniform design}  if $D_{3}^*$ has the smallest WD-value 
among the design set.

\end{definition}
According to the expression of WD, we can get the squared WD-value 
of the column augmented design $D_{3} \in \cC_{3}(n+n_1;2^{m_1}
3^{m_2}\bullet 3^{r})$,
\begin{align}&\mbox{WD}(D_{3})=-\left(\frac{4}{3}\right)^{m+r}+
\frac{1}{(n+n_1)^2}\sum_{i=1}^{n+n_1}\sum_{j=1}^{n+n_1}\prod_{k=1}^{m+r}
\left(\frac32-\left|u_{ik}-u_{jk}\right|\left(1-\left|u_{ik}-u_{jk}\right|\right)\right),
\label{WD} \end{align} where $n+n_1$ is the total number of runs and
$m+r$ is the total number of factors for the column augmented design.

In practice, we often choose a uniform design for the initial design 
$\d_0$ when the relationship between factors and response is
unknown, and we may terminate the experiments based on the data
analysis of the first stage. If more runs and factors should be added 
after $\d_0$, one searches the follow-up part $(\ba{cc} \d_1 &\d_2\ea)$ 
such that the column \aug design $D_{3}=\left(\ba{cc} \d_0 &\zero_{n\times 
r}  \\ \d_1  &\d_2 \ea\right) $ is as uniform as possible. Moreover,
the initial design $\d_0$ is assumed to be known in the follow-up
stage. Therefore, it is necessary to derive the expression of the WD of
the column \aug design based on $\d_0$. According to the coincident
numbers between any two rows, which represents the number of places
where two rows of the design $\d$ take the same value, we can rewrite 
the expression of the WD-value in (\ref{WD}).

\begin{proposition}\label{th D3} Given an initial design $\d_0 \in \cu(n;
2^{m_1}3^{m_2})$, its column augmented designs $D_{3} \in
\cC_{3}(n+n_1;2^{m_1}3^{m_2}\bullet 3^{r})$ have \begin{align}
\mbox{WD}(D_{3})=& ~C(r)+\frac{n^2}{(n+n_1)^2}\left(\frac{3}{2}\right)^{r}
WD(\d_0)+\frac{1}{(n+n_1)^2}\left(\frac{5}{4}\right)^{m_1}\left(\frac{23}
{18}\right)^{m_2+r}\times\nonumber\\
&\left(\sum_{i=n+1}^{n+n_1}\sum_{j(\neq i)=n+1}^{n+n_1}\left(\frac{6}
{5}\right)^{\cf_{ij}}\left(\frac{27}{23}\right)^{\cv_{ij}+\ct_{ij}}+2\sum_{i=1}^{n}
\sum_{j=n+1}^{n+n_1} \left(\frac{6}{5}\right)^{\cf_{ij}}\left(\frac{27}{23}\right)
^{\cv_{ij}}\right), \label{WD D3}
\end{align}
where \begin{align}C(r)=-\left(\left(\frac{4}{3}\right)^{r}-\frac{n^2}{(n+n_1)^2}
\left(\frac{3}{2}\right)^{r}\right)\left(\frac{4}{3}\right)^{m}+\frac {n_1}
{(n+n_1)^2}\left(\frac{3}{2}\right)^{m+r},\label{cr}\end{align}
$\cf_{ij} =\#\left\{~k:u_{ik}=u_{jk}, k=1,2,\dots,m_1 \right\},~\cv_{ij} =
\#\left\{~k :u_{ik}=u_{jk}, k= m_1+1,\dots,m\right\},\\
~\ct_{ij} =\#\left\{~k :(u_{ik},u_{jk})\in\left\{\left(\frac12,\frac12\right),
\left(\frac56,\frac56\right)\right\},\right.$ $ k=m+1,\dots,\left.m+r,\right\}$, 
$\#\left\{S\right\}$ is the number of elements in the set $S$.
\end{proposition}

The WD-value of \caug~$D_3$ is a function of $WD(\d_0)$. The
equation (\ref{WD D3}) can be used for obtaining the lower bound of
WD-values of the column augmented designs, which can be served 
as a benchmark to judge that whether a design is uniform or not. If the
WD-value of a design reaches the lower bound, then this design must
have the smallest WD-value among the design space, i.e., it is a
uniform design.

\begin{theorem} \label{th LBW3} For a given initial design $\d_0 \in 
\cu(n;2^{m_1}3^{m_2})$, its column augmented designs $D_{3}\in 
\cC_{3}(n+n_1;2^{m_1}3^{m_2}\bullet 3^{r})$ have $ \mbox{WD}(D_3)
\geq LBW_3, $ where \begin{align} &LBW_3=  C(r)+\frac{n^2}
{(n+n_1)^2}\left(\frac{3}{2}\right)^{r}WD(\d_0)+\frac{1}{(n+n_1)^2}\left(\frac{5}
{4}\right)^{m_1}\left(\frac{23}{18}\right)^{m_2+r}\left(T_1+2T_2 \right),
\label{LBW3_3}
\end{align} 
and $C(r)$ refers to (\ref{cr}), {$a=\ln\left(\frac{6}{5}\right), ~b=\ln\left(\frac{27}
{23}\right)$, }$\varphi_1=\frac{am_1(n_1-2)}{2(n_1-1)}+\frac{bm_2(n_1-3)}
{3(n_1-1)}+\frac{br(n_1-2)}{2(n_1-1)}, ~\varphi_2=\frac{am_1}{2}+\frac{bm_2}{3},
~T_1=n_1(n_1-1)e^{\varphi_1}, ~T_2=nn_1e^{\varphi_2}.$ The lower
bound can be achieved if and only if $\cf_{ij}=\frac{m_1(n_1-2)}{2(n_1-1)},
\cv_{ij}+\ct_{ij}=\frac{m_2(n_1-3)}{3(n_1-1)}+\frac{r(n_1-2)}{2(n_1-1)}$,
$i=n+1,\dots, n+n_1, j(\neq i)=n+1,\dots, n+n_1$, and $\cf_{ij}=\frac{m_1}{2}$,
$\cv_{ij}=\frac{m_2}{3},i=1,\dots, n, j=n+1,\dots, n+n_1$.
\end{theorem}

The column augmented design $D_3$ for initial mixed-level design can
be reduced to that for initial symmetrical two-level or three-level designs, 
i.e., $\cC_3(n+n_1;2^{m_1}3^{m_2}\bullet 3^{r})$ becomes $\cC_3(n+n_1
;2^{m}\bullet 3^{r})$ when $m_2=0$, or $\cC_3(n+n_1;3^{m}\bullet 3^{r})$ 
when $m_1=0$, according to actual demands. For these cases, we can 
derive more accurate lower bounds for these cases.

\begin{theorem}\label{th_twolevel_lbw} Given an initial design $\d_0\in\cu(n;
2^{m})$, the WD-value of \caugs~$D_3\in \cC_3(n+n_1;2^{m}\bullet 3^{r})$ 
has the lower bound
\begin{align} LBW_{3} = & ~C(r)+\frac{n^2}{(n+n_1)^2}\left(\frac{3}{2}\right)^{r}
WD(\d_0) +\frac{1}{(n+n_1)^2}\left(\frac{5}{4}\right)^{m}\left(\frac{23}
{18}\right)^{r}\times\nonumber\\
 &\left(T_1' +2\left(p_1\left(\frac{6}{5}\right)^{w_1}+q_1\left(\frac{6}
 {5}\right)^{w_1+1}\right) \right),\label{lbw3,2}
\end{align}
and $C(r)$ is in (\ref{cr}), $\varphi_1'=\frac{am(n_1-2)}{2(n_1-1)}+
\frac{br(n_1-2)}{2(n_1-1)}, ~T_1'=n_1(n_1-1)e^{\varphi_1'}, ~w_{1}=
\left\lfloor \frac{m}{2}\right\rfloor, ~p_{1}w_{1}+q_{1}(w_{1}+1)=\frac{mnn_1}{2}, 
~p_{1}+q_{1}=nn_1.$ The lower bound can be reached if and only if 
$\cf_{ij}=\frac{m_1(n_1-2)}{2(n_1-1)}, \ct_{ij}=\frac{r(n_1-2)}{2(n_1-1)}$, 
$i=n+1,\dots, n+n_1, j(\neq i)=n+1,\dots, n+n_1$, and there are $p_1$ 
number of $\cf_{ij}$ take $w_1$, $q_1$ number of $\cf_{ij}$ take $w_1+1$, 
$i=1,\dots, n, j=n+1,\dots, n+n_1$. 
\end{theorem}

If the number of initial factors $m$ is even, according to the discussion 
after Lemma 2 in the Appendix A, the lower bound in Theorem 
\ref{th_twolevel_lbw} is equivalent to that in Theorem \ref{th LBW3}. 
However, if $m$ is odd, the lower bound in Theorem \ref{th_twolevel_lbw} 
is more tight. For instance, choose the design in Example 2 of \cite{EQ15c} 
as the initial design $\d_0$, and assume $n_1=4, r=1$; then both of the 
lower bounds in Theorem \ref{th LBW3} and Theorem \ref{th_twolevel_lbw} 
are 53.5134, because $m=14$ is even. However, consider the design in 
Example 1 of \cite{EQ15c} as the initial design $\d_0$, and let $n_1=2, 
r=1$, then the lower bounds in Theorem \ref{th LBW3} and Theorem
\ref{th_twolevel_lbw} are 3.4094 and 3.4307, respectively, since
$m=7$ is odd.

\begin{theorem}\label{th_threelevel_lbw_d0} Given an initial design
$\d_0\in\cu(n;3^{m})$, the lower bound of WD-value of \caugs~$D_3
\in \cC_3(n+n_1;3^{m}\bullet 3^{r})$ is
\begin{align}
LBW_{3} =& ~C(r)+\frac{n^2}{(n+n_1)^2}\left(\frac{3}{2}\right)^{r}WD(\d_0)
 +\frac{1}{(n+n_1)^2}\left(\frac{23}{18}\right)^{m+r}\times \nonumber\\
&\left(p_2\left(\frac{27}{23}\right)^{w_2}+q_2\left(\frac{27}{23}\right)^{w_2+1}+
2\left(p_3\left(\frac{27}{23}\right)^{w_3}+q_3\left(\frac{27}{23}\right)^{w_3+1}\right) \right),\label{lbw3,3}
 \end{align}
and $C(r)$ refers to (\ref{cr}), $w_{2}=\left\lfloor\left(\frac{m(n_1-3)}{3(n_1-1)}+
\frac{r(n_1-2)}{2(n_1-1)}\right)\right\rfloor,$$~p_{2}+q_{2}=n_1(n_1-1),
~p_{2}w_{2}+q_{2}(w_{2}+1)=\frac{mn_1(n_1-3)}{3}+\frac{rn_1(n_1-2)}{2};
$$w_{3}=\left\lfloor \frac{m}{3}\right\rfloor$ $ ~p_{3}+q_{3}=nn_1,
~p_{3}w_{3}+q_{3}(w_{3}+1)=\frac{mnn_1}{3}.$ The lower bound can be
achieved if and only if $p_2$ number of $\cv_{ij}+\ct_{ij}$ take $w_2$, 
$q_2$ number of $\cv_{ij}+\ct_{ij}$ take $w_2+1$ , $i=n+1,\dots, n+n_1, 
j(\neq i)=n+1,\dots, n+n_1$, and $p_3$ number of $\cv_{ij}$ take $w_3$, 
$q_3$ number of $\cv_{ij}$ take $w_3+1$,$i=1,\dots, n, j=n+1,\dots, n+n_1$. 
\end{theorem}

If $\frac{m(n_1-3)}{3(n_1-1)}+\frac{r(n_1-2)}{2(n_1-1)}$ and $m/3$
are integers simultaneously, the lower bound in Theorem 
\ref{th_threelevel_lbw_d0} is equivalent to that in Theorem \ref{th
LBW3}. Otherwise, the lower bound in Theorem \ref{th_threelevel_lbw_d0} 
is more tight. The proof of Theorem \ref{th_threelevel_lbw_d0} is analogous 
to that of Theorem \ref{th_twolevel_lbw} and we omit it.

\subsection{Column \aug uniform designs with one additional blocking factor} 
\label{sec block_factor}

This subsection discusses the column \aug designs with one additional 
blocking factor. Similarly, the levels of the blocking factor can be fixed 
as 0 in the first stage. For the level of the blocking factor in the second 
stage, we can often take 1. This is because the level of the blocking 
factor in the follow-up stage is often different from that in the first stage.
A design $D_{3b}=\left(\ba{ccc} \d_0 &\zero_{n\times r} &\zero_{n\times 1} \\
\d_1 & \d_2 & \one_{n_1\times 1} \ea\right)$ is called a \caug~with one 
additional blocking factor. Denote all such column \aug designs by 
$\cC_{3b}(n+n_1;2^{m_1}3^{m_2}\bullet 3^{r}\bullet \b)$. Furthermore, 
to judge the uniformity of this type of \caugs, the lower bound is presented.

\begin{theorem}\label{th D3b} Given an initial design $\d_0 \in \cu(n;2^{m_1}
3^{m_2})$, $ \cu(n;2^{m})$ or $\cu(n;3^{m})$, the lower bound of the 
\caug~$D_{3b}$ with one additional blocking factor is
\begin{align*} LBW_{3b} =&
~C(r+1)+\frac{n^2}{(n+n_1)^2}\left(\frac{3}{2}\right)^{r+1}WD(\d_0)
+\frac{1}{(n+n_1)^2} \revA{Z,}  \end{align*}
 where
$Z=\left(\frac{5}{4}\right)^{m_1}\left(\frac{23}{18}\right)^{m_2+r}
\left( \left(\frac{3}{2}\right)T_1+2\left(\frac{5}{4}\right)T_2
\right)$ when $D_{3b}\in \cC_{3b}(n+n_1;2^{m_1}3^{m_2}\bullet
3^{r}\bullet \b)$, $Z=\left(\frac{5}{4}\right)^{m}\left(\frac{23}{18}\right)^{r}
\left(\left(\frac{3}{2}\right) T_1' +2\left(\frac{5}{4}\right)\left(p_1\left(\frac{6}
{5}\right)^{w_1}+q_1\left(\frac{6}{5}\right)^{w_1+1}\right) \right)
$ when $D_{3b}\in\cC_{3b}(n+n_1;2^{m}\bullet 3^{r}\bullet \b)$, and
$Z=\left(\frac{23}{18}\right)^{m+r}\left(\left(\frac{3}{2}\right)\left(p_2
\left(\frac{27}{23}\right)^{w_2}+q_2\left(\frac{27}{23}\right)^{w_2+1}\right)
+2\left(\frac{5}{4}\right)\left(p_3\left(\frac{27}{23}\right)^{w_3}
+q_3\left(\frac{27}{23}\right)^{w_3+1}\right) \right)$ when $D_{3b}\in\cC_{3b}
(n+n_1;3^{m}\bullet 3^{r}\bullet \b)$, respectively. Here, all the parameters 
are defined in Theorems \ref{th LBW3}-\ref{th_threelevel_lbw_d0}.

\end{theorem}

For the \caug~with a blocking factor $D_{3b}=\left(\ba{ccc} \d_0 &\zero_
{n\times r} &\zero_{n\times 1} \\ \d_1 & \d_2 & \one_{n_1\times 1} \ea\right)$, 
denote $\b=(\ba{cc}\zero_{n\times 1}^T & \one_{n_1\times 1}^T\ea)^T$. 
One wants to search the best additional part $\d_1$ and $\d_2$ such that
$D_{3b}$ is uniform under WD. There is only a different positive coefficient 
between the lower bounds in Theorems \ref{th LBW3}-
\ref{th_threelevel_lbw_d0} and that in Theorem \ref{th D3b} for each term. 
Then, we have the following result directly, and omit its proof.

\begin{proposition}\label{co equivalency1} Given an initial design
$\d_0\in\cu(n;2^{m_1}3^{m_2})$, $\cu(n;2^{m})$ or $\cu(n;3^{m})$,
the \caug~$D_3^*$ achieves its lower bounds in Theorems \ref{th
LBW3}-\ref{th_threelevel_lbw_d0} if and only if the corresponding
column \aug design $D_{3b}^*=\left( D_3^* ~~ \b \right)$ achieves
its lower bound in Theorem \ref{th D3b}.
\end{proposition}

Proposition \ref{co equivalency1} means that one can construct
$D_{3b}$ through $D_3$, i.e., for constructing $D_{3b}^*$, one only
needs to construct the corresponding $D_3^*$, then add $\b$ into
$D_3^*$ to obtain the design $D_{3b}^*=(\ba{cc}D_3^* & \b\ea)$.

\subsection{Column \aug uniform designs with two additional blocking factors}
 \label{sec block2factor}

If one wants to add two blocking factors in the follow-up design
due to the practical requirements, as the same idea in Subsection
\ref{sec block_factor}, the level of the two blocking factors in
the initial design should be 0 and that in the second stage should
be 1, i.e., the design matrix of the two additional blocking factors is
$\E=(\ba{cc}\zero_{n\times 2}^T & \one_{n_1\times 2}^T\ea)^T$.
However, the two columns in $\E$ are fully correlated such that the
effects of them cannot be distinguished. If we wish to assess these 
blocking effects accurately, we have to adjust the structure of the two 
additional blocking factors to reduce the correlation. Usually, if an 
initial design has been done, the design matrix is fixed and cannot 
be altered, we have to consider changing the structure in the second-
stage. Naturally, replace $\E$ by $\B=\left(\ba{cc} \zero_{n\times 1}
&\zero_{n\times 1}  \\\one_{n_1\times 1}  & \ba{c} \zero_{\frac{n_1}{2}
\times 1} \\
\one_{\frac{n_1}{2}\times 1} \ea \ea\right)$,  denote
$D_{3B}=\left(\ba{cc} \ba{cc} \d_0 &\zero_{n\times r}  \\
\d_1 &\d_2 \ea & \B \ea\right)=(\ba{cc}D_3 &\B \ea) \in
\cC_{3B}(n+n_1;2^{m_1}3^{m_2}\bullet 3^{r}\bullet \B)$.
This solution sacrifices the level balance of the second additional
blocking factor to reduce the non-orthogonality between the two
additional blocking factors. For assessing the uniformity of $D_{3B}$, 
we have the following lower bound.

\begin{theorem} \label{theorem lbw 3b} Given an initial design $\d_0\in
\cu(n;2^{m_1}3^{m_2})$, any $D_{3B}$ has the following lower bound 
of WD-value, \begin{align*} LBW_{3B} =&~C(r+2)+\frac{n^2}{(n+n_1)^2} 
\left(\frac{3}{2}\right)^{r+2} WD(\d_0)+\frac{1}{(n+n_1)^2}\left(\frac{5}
{4}\right)^{m_1+1}\left(\frac{23}{18}\right)^{m_2+r}\times\nonumber\\
 &\left(\left(\frac{3}{2}\right)T_{1B}+2\left(\frac{5}{4}\right)T_{2B} \right), 
 \end{align*}
and the function $C(\cdot)$ is defined in (\ref{cr}), $\varphi_{1B}=
\frac{a(m_1+1)(n_1-2)}{2(n_1-1)}+\frac{bm_2(n_1-3)}{3(n_1-1)}+
\frac{br(n_1-2)}{2(n_1-1)},~\varphi_{2B}=\frac{a(m_1+1)}{2}+\frac{bm_2}{3},$
~$T_{1B}=n_1(n_1-1)e^{\varphi_{1B}}, ~T_{2B}=nn_1e^{\varphi_{2B}}.$

\end{theorem}

The proof is similar to Theorem \ref{th LBW3} and we omit it here. If 
the initial design is a symmetrical two-level $\d_0\in\cu(n;2^{m})$ 
or three-level $\d_0\in\cu(n;3^{m})$, a lower bound can be derived in 
the same way.

Another method to solve the high correlation between the two
additional blocking factors is replacing $\E$ by $\B'=(\ba{cc}
\zero_{n\times 2}^T & \d_{B'}^T\ea)^T$ where $\d_B'\in\cu(n_1;
2^{2})$, which sacrifices the level balance of two blocking factors. 
Denote the corresponding \caug~by $D_{3B'}=\left(\ba{cc} \ba{cc} 
\d_0 &\zero_{n\times r} \\ \d_1 &\d_2 \ea & \B'\ea\right)=(\ba{cc}
D_3 &\B' \ea) \in \cC_{3B'}(n+n_1;2^{m_1}3^{m_2}\bullet 3^{r}
\bullet \B')$. It can be easily seen that \caug~$D_{3B}$ has smaller 
a WD-value than $D_{3B'}$, because the former is more balanced 
than the latter.

\begin{remark} The additional blocking factors in Subsections 
\ref{sec block_factor} and \ref{sec block2factor} can also serve as 
additional common two-level factors. Therefore, the \caug~$D_{3b}$ 
and $D_{3B}$ are \caugs~with additional mixed-level factors.
\end{remark}

In fact, whether the additional two-level factors are used for blocking or 
not does not affect the structure of the design matrix; it only influences 
the modeling aspect.

\subsection{Multi-stage augmented designs}

If there is no additional factor in the multi-stage design and only
some rows are added in each stage, then the $l$-stage row \aug
design can defined as $D^{(l)}=(\d_0^T, \dots, \d_{l-1}^T )^T,$ and 
the additional number of runs of the $i$-th stage portion is $n_{i-1}$. 
Specially, $n_0=n$. For constructing three-stage row \aug design, 
one can take the first and the second stage design $D^{(2)}=( \d_0^T, 
\d_{1}^T)^T $ as the initial design, then add the third-stage portion 
as the follow-up portion. Next, take the first three stages portion as 
the initial design, and add the fourth-stage portion as the follow-up 
portion, and so on.

If there exists an additional blocking factor in the multi-stage
\caug, as similar as the discussion before, let the levels of the 
blocking factor in $l$ stages take $0,1,\dots,l-1$, respectively. 
For the initial design $\d_0$, the $l$-stage \caug ~with one blocking 
factor, with $r$ three-level additional factors, and with $r$ three-level 
additional factors and one blocking factor can be respectively defined 
as follows, $$D_{b}^{(l)}= \left(\ba{cc}
\d_0 &\zero_{n\times 1} \\
\d_1 & \one_{n_1\times 1} \\
 \vdots  & \vdots\\
\d_{l-1} & (\ll)_{n_{l-1}\times 1} \ea\right),D_{3}^{(l)}= \left(\ba{cc}
\d_0 & \zero_{n\times r}  \\
\d_1 & \d_1'  \\
 \vdots & \vdots \\
\d_{l-1} &\d_{l-1}' \ea\right),
D_{3b}^{(l)}= \left(\ba{ccc}
\d_0 & \zero_{n\times r} &\zero_{n\times 1} \\
\d_1 & \d_1' & \one_{n_1\times 1} \\
 \vdots & \vdots & \vdots\\
\d_{l-1} &\d_{l-1}' & (\ll)_{n_{l-1}\times 1} \ea\right). $$ It is assumed 
that the additional three-level factors may be considered in the second 
stage, after which no factor is added.

\begin{proposition}\label{th equi2} Given an initial design $\d_0\in\cu(n;
2^{m_1}3^{m_2})$, its $l$-stage \caug~$D_{b}^{(l)}$ achieves its 
lower bound if and only if the corresponding $l$-stage row \aug design 
$D^{(l)}$ achieves its lower bound, and we have $ \mbox{WD}(D^{(l)})\geq 
LBW^{(l)}, $ where \begin{align} \label{WD LBWl}
LBW^{(l)} = & -\frac{n_{l-1}^2+2N_{l-1}n_{l-1}}{(N_{l-1}+n_{l-1})^2}\left(
\frac{4}{3}\right)^{m} +\frac {n_{l-1}}{(N_{l-1}+n_{l-1})^2}\left(\frac{3}{2}
\right)^{m}+\frac{N_{l-1}^2}{(N_{l-1}+n_{l-1})^2}WD(D^{(l-1)})\nonumber\\
&+\frac{1}{(N_{l-1}+n_{l-1})^2}\left(\frac{5}{4}\right)^{m_1}\left(\frac{23}
{18}\right)^{m_2}\left(T_{l+2} +2T_{l+1} \right),
 \end{align} and $D^{(l-1)}=(\ba{ccc}\d_0^T & \dots & \d_{l-2}^T\ea)^T$ 
 is the $(l-1)$-stage row \aug design, $D^{(1)}=\d_0$; $N_{l-1}=n+n_1+
 \dots+n_{l-2}$, representing the total run number of the first $(l-1)$ stages, 
 $N_1=n$, $T_{l+1}=N_{l-1}n_{l-1}e^{\varphi_{l+1}}$,$\varphi_{l+1}
 =am_1/2+bm_2/3$, $T_{l+2}=n_{l-1}(n_{l-1}-1)e^{\varphi_{l+2}}$,
$\varphi_{l+2}=\frac{am_1(n_{l-1}-2)}{2(n_{l-1}-1)}+\frac{bm_2(n_{l-1}-3)}
{3(n_{l-1}-1)}$.
\end{proposition}

The proof is in the Appendix B. The result in Proposition \ref{th equi2} is a
generalization from two-stage to $l$-stage situation. Thus, one can 
arrange $l$-stage row \aug design $D^{(l)}$ stage by stage, and then 
construct $l$-stage column \aug design with one additional blocking 
factor $D_{b}^{(l)}$ through adding $(\zero_{n\times 1}^T, \dots,
(\ll)_{n_{l-1}\times 1}^T )^T$ to $D^{(l)}$.

Similarly, the multi-stage \caug~$D_{3b}^{(l)}$ reaches the lower bound 
if and only if the $D_{3}^{(l)}$ reaches the lower bound. However, its 
recursive lower bounds of multi-stage \caug~$D_{3}^{(l)}$ like 
(\ref{WD LBWl}) are too complicated, since the lower bounds closely 
depend on the total additional number of runs, $\overline{N_{l-1}}=
n_1+\dots + n_{l-1}$. When $\overline{N_{l-1}}\leq 2n$, each column 
for each of the follow-up stages $\d_1',\dots, \d_{l-1}'$ should take 
the same number of levels 1 and 2, for keeping the balance of the 
levels of the $r$ additional factors. When $\overline{N_{l-1}} > 2n$, 
let $h=max\{~i~|~\overline{N_{i-1}}\leq 2n~\}$. Each column of the
follow-up stages $\d_1',\dots, \d_{h-1}'$ should take the same
number of levels 1 and 2. While, each column of the follow-up stage
$\d_h',\dots, \d_{l-1}'$ should take values from $\{0,1,2\}$. Hence,
it has some technical difficulty for deriving the lower bounds of the 
$l$-stage \caug~$D_{3}^{(l)}$ when $l>2$. The Proposition
\ref{th equi2} is based on the initial mixed-level design $\d_0\in
\cu(n;2^{m_1}3^{m_2})$. Specially, under symmetrical two-level 
$\d_0\in\cu(n;2^{m})$ or three-level $\d_0\in\cu(n;3^{m})$ initial 
designs, similar results can be obtained easily.

\section{$E(f_{NOD})$ criterion for \caug} \label{sec f_NOD}
In this section, we study the \caugs~from the view of non-orthogonality, 
by using the $E(f_{NOD})$ criterion. For a design $\d \in\cu(n;q_1,
\ldots, q_m)$, define the non-orthogonality between the $k$-th 
and $l$-th columns of $\d$ as $f^{kl}_{NOD}(\d)=\sum\limits_{u=0}
^{q_k-1}\sum\limits_{v=0}^{q_l-1} \left(n_{uv}^{kl}(\d)-\frac{n}{q_kq_l}
\right)^2,$ where $n_{uv}^{kl}(\d)$ is the number of $(u,v)$-pairs in 
the $k$-th and $l$-th columns in $\d$, and $\frac{n}{q_kq_l}$ represents 
for the average frequency of all level-combinations in each pair of the
$k$-th and $l$-th columns of $\d$. Here, the subscript $NOD$ stands
for non-orthogonality of the design. The $k$-th and $l$-th columns
of $\d$ are orthogonal if and only if $f^{kl}_{NOD}(\d)=0$. Define
$E(f_{NOD}(\d))= \sum_{1\leq k<l\leq m}f^{kl}_{NOD}(\d)/\dbinom{m}{2}$, 
which measures the average non-orthogonality among the columns 
of design $\d$. Especially, $\d$ is an orthogonal design if and only if 
$E(f_{NOD}(\d))=0$. The smaller $E(f_{NOD})$ value implies that 
the design has better orthogonality. Therefore, one prefers a design 
with small $E(f_{NOD})$. A design is $E(f_{NOD})$-optimal if it has 
the smallest value of $E(f_{NOD})$ among the design space.

The lower bound plays a key role in detecting the $E(f_{NOD})$-optimal 
design. The lower bounds of $E(f_{NOD})$ of the \caug~are impacted by 
the non-orthogonality of the initial design, $E(f_{NOD}(\d_0))$, that is, 
the lower bounds of $E(f_{NOD})$ for \caug~contain the nonorthogonality 
information of the initial design $\d_0$.

Given an initial design $\d_0\in\cu(n;2^{m_1}3^{m_2})$, the lower
bounds of $E(f_{NOD})$ of \caugs~$D_3$, $D_{3b}$ and $D_{3B}$
respectively are $LBf_{3}$, $LBf_{3b}$, $LBf_{3B}$, which are given
in Lemma \ref{thlbf} in the Appendix A.

\begin{proposition}\label{co_wd_fnod}
For a given initial design $\d_0\in\cu(n;2^{m_1}3^{m_2})$, $\cu(n;
2^{m})$ or $\cu(n;3^{m})$, its \caug~$D_3$ reaches $LBf_3$ if
and only if the corresponding \caug~$D_{3b}$ reaches $LBf_{3b}$.
\end{proposition}

Proposition \ref{co_wd_fnod} is similar to Proposition \ref{co equivalency1}. 
The proof is straightforward, \revA{omitting it}. It suffices to 
discuss the property of $D_3$; the $D_{3b}$ has the similar property. 
Moreover, although the non-orthogonality criterion $E(f_{NOD})$ is 
defined by any two columns in designs, it can be shown that $E(f_{NOD})$ 
is a function of coincident numbers between any pairs of rows. Refer 
to the WD-value in (\ref{WD D3}), which is also related to coincident
numbers. Then we build some connections between WD and $E(f_{NOD})$
as follows.

\begin{theorem}\label{connection}
(1) For an initial design $\d_0\in\cu(n;2^{m_1}3^{m_2})$, $\cu(n;2^{m})$ 
or $\cu(n;3^{m})$ and $r\geq 0$, if $WD({D_3})$ or $WD({D_{3B}})$ 
reaches its lower bound $LBW_3$ or $LBW_{3B}$, then $E(f_{NOD}
(D_{3}))$ or $E(f_{NOD}(D_{3B}))$ reaches its lower bound $LBf_3$ 
or $LBf_{3B}$, respectively.

(2) For an initial design $\d_0\in\cu(n;3^{m})$ and $r\geq 0$, if $E(f_{NOD}
(D_3))$ reaches the lower bound $LBf_3$, then $WD({D_3})$ also reaches 
its lower bound $LBW_3$.

(3) For an initial design $\d_0\in\cu(n;2^{m})$  and $r=0$, if $E(f_{NOD}(D))$ 
or $E(f_{NOD}(D_{B}))$ reaches its lower bound then $WD({D})$ or 
$WD({D_{B}})$ also reaches its lower bound, respectively.

\end{theorem}

In Theorem \ref{connection}, we omit the discussion about $D_{3b}$,
because it has the same result as $D_3$ according to Propositions
\ref{co equivalency1} and \ref{co_wd_fnod}. From Theorem
\ref{connection}(1), column \aug uniform designs are $E(f_{NOD})$
-optimal, but the converse may not hold unless the condition in (2) 
or (3) is satisfied. It is known that if a \caug~has good uniformity, it 
often has good orthogonality. Thus, studying \caug~from the view of 
uniformity criterion is reasonable.

\section{Construction algorithm and illustrative examples } \label{sec algorithm}

In this section, the construction algorithm for the column augmented 
uniform design $D_3^*\in\cC_{3}(n+n_1;2^{m_1}3^{m_2}\bullet 3^r)$ 
is discussed and $D_{3B}^*$ can be constructed similarly.

For saving experimental runs, researchers prefer to choose a
uniform design $\d_0$ as the initial design, since they may stop the
experiment after the data analysis of the first stage. Then, the
searching procedure for $D_3^*$ is as follows,\\

\vspace{2mm}
Step 1) Search a uniform design $\d_0^*$ from the design space $\cu(n;
 2^{m_1}3^{m_2})$ as the initial design.\\
 Step 2) Search the additional design $\d_1$ and $\d_2$ to make 
$D_3^*=\left(\ba{cc} \d_0^* &\zero_{n\times r}  \\\d_1  & \d_2 
\ea\right)$ be as uniform as possible.

\vspace{2mm}

From the structure of column augmented design $D_3^*$, if the
initial design $\d_0^*$ is chosen as the empty set $\emptyset$,
$n_1=n$ and $r=0$, Step 2) reduces to constructing a uniform
design $\d_0^*$ with $n$ runs. Hence, we only give the detailed
construction algorithm for Step 2). This construction algorithm can
also be used for the multi-stage \aug designs. For example, for
constructing the third-stage portion under WD, one can return to
Step 1), denote the obtained two-stage \caug~$D_3^*$ as the new
initial design $\d_0^*$, and repeat the whole procedure again.

In the searching procedure for $D_3^*$, some stochastic algorithms
may be used. Among them, the threshold accepting (TA) algorithm is
widely used to construct uniform designs; see \cite{WF97},
\cite{FLW03a}, \cite{ZF13}, etc. The TA algorithm is a modified
simulated annealing (SA) method. SA method is a heuristics
optimization algorithm and very powerful for NP hard problem(Kirkpatrick et
al., 1983). The main idea of SA for an optimization problem is
to start with an initial solution, and iteratively update the current 
solution to its neighboring candidate in a stochastic way.
If the neighboring solution is better, accept it with probability
one; otherwise accept it with a small probability that tends to
zero as the number of iterations increases. In the SA algorithm, the
reason of accepting worse solution is to jump from the local
optimum. The TA algorithm accepts the worse neighboring
solution by some hard threshold instead of using probability, while
gradually shrinking the threshold towards zero. Dueck and Scheuer \cite{DS90}
pointed out that the TA algorithm is much simpler and has higher
convergence rate than the SA algorithm. Moreover, the TA algorithm
is also a global searching algorithm (Alth\"{o}fer and Koschnick, 1991). 
For further details on the TA algorithm for uniform designs construction,
readers can refer to \cite{FLW03a}.

Algorithm 1 uses the TA algorithm to search the optimal follow-up
part in Step 2). It can be seen that only the second-stage portion
$(\d_1~ \d_2)$ is updated in $D_3$ and the initial design $\d_0^*$
is fixed in the whole procedure, which is the main difference
between  Algorithm 1 and the TA algorithm used in the aforementioned 
references. It seems that Algorithm 1 is similar to the ``exchange 
algorithms'', which are the most popular method for constructing 
optimal designs; see \cite{F72} and \cite{MN95}. However, the 
exchange algorithms greedily pursue current optimality at every 
exchange step, such that they may drop into the local optimum, 
while Algorithm 1 is a global searching algorithm. Moreover, compared 
with the coordinate descent algorithm in \cite{W15}, our algorithm can 
process the cases of multimodal and discrete function optimization.

\noindent \rule{\textwidth}{0.4mm} {\bf Algorithm 1}~~Construction
of
the column augmented uniform design $D_3^*$ \\
\rule{\textwidth}{0.4mm}

\noindent 1. \textbf{Input:} The initial uniform design $\d_0^*$,
the parameters  $n_1$  and   $r$.   \\
2. Initialize $I,J$ and the sequence of thresholds $T_i, i=1,2,\dots,I$ \\
3. Generate starting design $D_3=\left(\ba{cc} \d_0^*
&\zero_{n\times r}
\\\d_1   & \d_2
\ea\right)$,  $\d_1\in \cu(n_1;2^{m_1}3^{m_2})$ and {$\d_2 \in
\{1,2\}^{n_1\times r}$},\\
4. \textbf{for }$i=1:I$\\
5. $~~~$ \textbf{for }$j=1:J$\\
6. $~~~~~~~$ \textbf{if} $\mbox{WD}(D_3)=LBW_3$, go to line 14; otherwise\\
7. $~~~~~~~~~~~$ Generate $\d_{1_{new}}\in \cn(\d_1)$ and
$\d_{2_{new}}\in \cn(\d_2)$,  obtain
$D_{3_{new}}=\left(\ba{cc}\d_0^* &\zero_{n\times r}
\\\d_{1_{new}} & \d_{2_{new}}
\ea\right)$ \\
8. $~~~~~~~~~~~$ \textbf{if} $\mbox{WD}(D_{3_{new}})-\mbox{WD}(D_3)\leq T_i$ \\
9. $~~~~~~~~~~~~~~~$ update $D_3=D_{3_{new}}$, $\d_1=\d_{1_{new}}$
and $\d_2=\d_{2_{new}}$ \\
10. $~~~~~~~~~~~$\textbf{ end if}\\
11. $~~~~~~~$\textbf{end if}\\
12. $~~~$ \textbf{end for}\\
13. \textbf{end for}\\
14. \textbf{Output: }$D_3^*=D_3$  \\
\rule{\textwidth}{0.4mm}

From the discussion of the Subsection \ref{sec column augmented},
there are some limits for the run size of the follow-up portion $n_1$. 
In Line 2, the iteration times $I$, $J$ and the thresholds $T_i, i=1,2,
\dots,I$ can be selected in the same way as \cite{ZF13}. We choose 
the thresholds $T_1>T_2> \dots > T_I=0$. The positive threshold 
$T_i$ aims to prevent falling into the local optimum. The threshold 
decreases for the algorithm to avoid an endless loop. In Line 3, the 
U-type designs $\d_1$ and $\d_2$ can be randomly generated. Usually, 
we also can adopt the good lattice point method combined with the level
transformation to generate $\d_1$ and $\d_2$, as in \cite{ZF13}. In
Line 6, if the WD-value of \caug~reaches the lower bound, the
search procedure ends. In Line 7, $\d_{1_{new}}$ and $\d_{2_{new}}$
are randomly chosen from the neighborhood of $\d_1$ and $\d_2$,
respectively, denoted by $\cn(\d_1)$ and $\cn(\d_2)$, where the
neighborhood of $\d$, $\cn(\d)$, means all of the designs that are
obtained by exchanging two randomly chosen elements in a random
column of $\d$. Due to the randomness of Algorithm 1, when the
column augmented design does not  reach the lower bounds, we can
repeat the algorithm several times and select the best one as the
final design. According to our experience, the threshold $T_i$ can
be chosen as follows. At first, randomly generate $P$ designs in the
neighborhood of the current $D_3$, denoted by $D_{3_{new},1},
\dots,D_{3_{new},P}$, and compute the $P$ differences $\triangle 
WD_i=WD(D_{3_{new},i})-WD(D_{3}),~i=1,\dots, P.$ The threshold 
$T_1$ can be defined as the 5th-percentile of $F$. $F$
is the empirical distribution of the difference $ \triangle WD$
which is larger than 0. Then, let $T_i=(I-i)/(I-1)*T_1,~i=1,\dots,P.$

Next, we give some examples to show the usefulness and
effectiveness of Algorithm 1. For comparison, we define the
efficiencies of designs as follows, 
\begin{align*}
f_3=\frac{LBW_3}{\mbox{WD}(D_3)},
f^{(l)}=\frac{LBW^{(l)}}{\mbox{WD}(D^{(l)})},
f_{3B}=\frac{LBW_{3B}}{\mbox{WD}(D_{3B})},
ff_{3}=\frac{LBf_{3}}{E(f_{NOD}(D_{3}))}.
\end{align*}

\noindent When the efficiencies equal 1, the \caugs~are uniform, and
when the efficiencies are close to 1, the \caugs~are nearly uniform. 
Usually, in practice, in one sequential stage, the number of augmented 
factors $r$ may be small and the additional runs $n_1$ may also be 
not too large. Thus, we choose $n_1\leq n, r\leq 4$ in the following examples.
The design matrices of initial designs in Example 3 and
all the optimal \caugs~are listed in Appendix B.

\begin{example}\label{example1} \emph{Consider the initial design $\d_0^*
\in \cu(12;2^53^{7})$. The initial $\d_0^*$ is a nearly uniform design since
$WD(\d_0)=11.3240$ and $LBW=11.2206$,}
\end{example}\vspace{-6mm}
$$ \d_0^*=\left[\ba{cccccccccccc}
 1 1  0  0  1  1  1  0  2  0  0  2\\
0  1  0  1  0  2  2  1  1  0  0  0\\
0  0  0  0  0  1  2  0  0  0  2  1\\
0  1  1  1  1  1  0  2  0  2  0  1\\
0  0  0  1  1  0  1  2  2  1  2  2\\
1  0  1  1  1  2  2  0  1  1  0  1\\
1  0  1  1  0  1  0  1  2  0  1  2\\
1  0  0  0  1  2  0  1  0  2  2  0\\
1  1  1  1  0  0  1  0  0  2  2  0\\
0  1  1  0  1  0  2  1  2  1  1  0\\
1  1  0  0  0  0  0  2  1  1  1  1\\
0  0  1  0  0  2  1  2  1  2  1  2
 \ea\right]. $$
According to Subsection \ref{sec column augmented}, $n_1$ should be
a multiple of 6. We search the \caug ~$D_3^*$ and ~$D_{3B}^*$
under some $n_1$ and $r$, and the results are shown in Table 1. We
can see that the column augmented designs $D_3^*$ ($D_{3B}^*$) are
nearly uniform since $f_3$ ($f_{3B}$) are close to 1 in all cases. For the 
right portion of Table 1, there are $r$ additional three-level factors and 
2 additional two-level factors in $D_{3B}^*$. It can be seen that the 
design efficiency increases with the increase of $r$, when $n_1$ is fixed.

\begin{table} \caption{The efficiencies of optimal column augmented designs
based on $\d_0^*\in \cu(12;2^53^{7})$}\label{table effi1}
\vspace{-2mm}
 \bc {\small \tabcolsep=2pt \bt {cc|ccc|cc|ccc}\hline $n_1$ & $r$
& $LBW_3$ &~~$WD(D_3^*)$~~ & $f_3$ & $n_1$ & $r$ & $LBW_{3B}$ 
&~~$WD(D_{3B}^*)$~~ & $f_{3B}$ \\\hline
6 & 1 & 16.0690 & 16.3116 & 0.9851 & 6 & 0 & 27.4964 & 27.9462 & 0.9839 \\
6 & 2 & 25.6762 & 25.9766 & 0.9884 & 6 & 1 & 43.2330 & 43.7927 & 0.9872\\
6 & 3 & 40.8687 & 41.2569 & 0.9906 & 6 & 2 & 67.8435 & 68.5665 & 0.9895\\
6 & 4 & 64.7699 & 65.3419 & 0.9912 & 12 & 0 & 23.6793 & 24.5695 & 0.9638\\
12 & 1 & 13.9002 & 14.4545 & 0.9617 &12 & 1 & 35.4254 & 36.5552 & 0.9691 \\
12& 2 & 20.9598 & 21.6406 & 0.9685 & 12 & 2 & 53.3404 & 54.7688 & 0.9739\\
 \cline{6-10}
12 & 3 & 31.8361 & 32.6828 & 0.9741 \\
12 & 4 & 48.5969 & 49.7171 & 0.9775 \\
\cline{1-5}
 \et} \ec
 \end{table}
Moreover, consider the multi-stage \aug design based on the
initial design $\d_0^*$ in Example \ref{example1}. Assume one adds 6
runs in each follow-up stage and there is no additional factor. We
search the optimal additional parts $\d_1,\d_2,\d_3$ are as follows,
respectively. 
 \bc {\small 
 \tabcolsep=4pt \bt {c|c|c} \hline
 $\d_1$ & $\d_2$ & $\d_3$ \\\hline
011002002122 & 000010200202 & 011101202211\\
100100212201 & 101012102011 & 010012022001\\
101001120100 & 110001222220 & 101110220020\\
110112220012 & 111111011122 & 000011111110\\
000111001210 & 010102110111 & 100102101122\\
011010111021 & 001100021000 & 111000010202\\\hline
 \et} \ec
The efficiencies of $D^{(1)}=\d_0^*,$ $ D^{(2)}= (\ba{cc} \d_0^{*T}
& \d_1^T  \ea)^T$, $D^{(3)}= (\ba{ccc} \d_0^{*T}  & \d_1^T  &
\d_2^T  \ea)^T $ and $D^{(4)}= (\ba{cccc} \d_0^{*T} & \d_1^T  &
\d_2^T & \d_3^T \ea)^T $
are 0.9909, 0.9803, 0.9683 and 0.9634, respectively.

\begin{example} \label{example2} 
\emph{Consider the initial design $\d_0^*\in \cu(6;3^{10})$, which
is a uniform design under WD since
$WD(\d_0^*)=LBW_0=5.1774$,}\end{example} \vspace{-4mm}\Bea
 \d_0^*=\left[\ba{cccccccccccc}
0 1 0 1 0 1 0 0 2 2\\
1 0 1 0 0 0 0 1 0 0\\
1 0 0 1 1 2 1 2 1 1\\
2 1 1 2 2 1 2 1 1 1\\
0 2 2 2 1 0 2 2 2 0\\
2 2 2 0 2 2 1 0 0 2
 \ea\right]. \label{d06}
 \Eea
The result is shown in Table \ref{table effi2}. It is known that the column 
augmented designs $D_3^*$ are uniform since $f_3=1$ in all the cases, 
and the designs $D_3^*$ are also $E(f_{NOD})$-optimal since $ff_3=1$. 
The result also verifies the connection in Theorem \ref{connection}. 
Moreover, the column augmented designs $D_{3B}^*$ are nearly 
uniform since $f_{3B}$ are very close to 1.

\begin{table}[!t] \caption{The efficiencies of optimal column augmented designs
based on $\d_0^*\in \cu(6;3^{10})$} \label{table effi2} \bc {\small
\tabcolsep=1pt \bt {cc|ccc|ccc|cc|ccc}\hline $n_1$ & $r$ & $LBW_3$
&~~$WD(D_3^*)$~~ & $f_3$ & $LBf_3$ & $E(f_{NOD}(D_3^*))$ & $ff_3$ &
 $n_1$ & $r$ & $LBW_{3B}$ &~~$WD(D_{3B}^*)$~~ & $f_{3B}$ \\\hline
6 & 1 & 5.7673 & 5.7673 & 1  & 3.8545 & 3.8545 & 1
& 6 & 0 & 10.4016 & 10.5390 & 0.9870 \\
6 & 2 & 9.1444 & 9.1444 & 1 & 4.2727 & 4.2727 & 1
&6 & 1 & 16.1739 & 16.3443  & 0.9896 \\
6 & 3 & 14.5244 & 14.5244 & 1 & 4.9231  & 4.9231  & 1 &6 & 2 &
25.2053 & 25.4247 & 0.9914 \\\cline{9-13} 6 & 4 & 23.0480 &
23.0480 & 1 &  5.6923 & 5.6923 & 1 \\\cline{1-8}
 \et} \ec
\end{table}

\begin{table}[!t] \label{table effi3} \caption{The efficiencies of optimal column
 augmented designs based on $\d_0^*\in \cu(72;2^43^{45})$}
 \bc {\small \tabcolsep=2pt \bt {cc|ccc|cc|ccc}\hline $n_1$ & $r$
& $LBW_3/10^7$ &~~$WD(D_3^*)/10^7$~~ & $f_3$ & $n_1$ & $r$ &
$LBW_{3B}/10^7$ &~~$WD(D_{3B}^*)/10^7$~~ & $f_{3B}$ \\\hline

24 & 1 & 0.8103 & 0.8166 & 0.9923
& 24 & 0 & 1.2260  & 1.2352  & 0.9926 \\
24 & 2 & 1.2173  & 1.2251  & 0.9936
& 24 & 1 & 1.8412  & 1.8529  & 0.9936\\
24 & 3 & 1.8298  & 1.8398  & 0.9946
& 24 & 2 & 2.7663  & 2.7817  & 0.9945\\
24 & 4 & 2.7516  & 2.7646  & 0.9953
& 48 & 0 & 1.0122  & 1.0311  & 0.9817\\
48 & 1 & 0.6715 & 0.6845 & 0.9810
&48 & 1 & 1.5055  & 1.5296  & 0.9842 \\
48 & 2 & 0.9984 & 1.0151  & 0.9835
& 48 & 2 & 2.2432  &  2.2740  &  0.9864\\
48 & 3 & 1.4873  & 1.5099  & 0.9850
&72 & 0 & 0.8753 &  0.9010 &  0.9714 \\
48 & 4 & 2.2193  & 2.2482  & 0.9871
&72 & 1 & 1.2932  &  1.3277  &  0.9740 \\
72 & 1 & 0.5812 & 0.5992 & 0.9699 &72 & 2 & 1.9154  & 1.9597
& 0.9774 \\\cline{6-10}
72 & 2 & 0.8581 & 0.8817 & 0.9733\\
72 & 3 & 1.2705  & 1.3016  & 0.9761 \\
72 & 4 & 1.8852  & 1.9254  & 0.9792 \\
\cline{1-5}
 \et} \ec
 \end{table}

\begin{example} \label{example3} 
 \emph{Consider the initial design $\d_0^*\in \cu(72;2^43^{45})$.\vspace{3mm} \\
To save space, we only list the result of $n_1=24,48$ and $72$ in
Table 3. When $n_1=72$, the total run size of the \caug~is $144$. It
shows that the efficiency is also close to 1 in each case for the
designs with larger run sizes.} \end{example}

Examples \ref{example1}-\ref{example3} present column augmented
(nearly) uniform designs $D_3^*$ ($D_{3B}^*$) and also show that the
lower bounds in this paper are fairly tight because the efficiencies are 
very close to 1, and even can reach 1 in some situations.

\section{Conclusion and discussion} \label{sec conclusions}
The mixed-level \caugs~proposed in this paper are suitable for the
situation that the initial stage ignores a few significant factors.
Column augmented designs with additional three-level factors and 
one or two blocking factors are discussed. In fact, the additional
blocking factor can be used as the ordinary two-level factor. That
is, we consider adding mixed-level factors to the initial designs
in the follow-up designs. For generalizing them, one could discuss 
the situation for adding $r_2$ two-level blocking factors with $r_2>2$. 
Although this is similar in principle, the discussion would be more 
complex. Moreover, one or two additional blocking factors often meet
the actual requirement. Hence, we do not discuss such a case. It is
interesting to find that the column \aug uniform designs are also
$E(f_{NOD})$-optimal. Furthermore, three examples show that the
construction algorithm is efficient and the lower bounds are tight
in many cases.

\section*{Acknowledgements}
The authors would like to thank the associate editor and the two
referees for their valuable suggestions, thank Prof. Robert Mee and
Dr. Siyuan Yang for their help and suggestions. This work is
supported by National Natural Science Foundation of China
(11471229,11871288), China Scholarship Council and Graduate Student
Innovation \& Practice Ability Enhancement Project of Sichuan
University.

\section*{Appendix A}

\noindent For simplicity, denote
$\bigtriangleup_{ijk}=\frac32-\left|u_{ik}-
u_{jk}\right|\left(1-\left|u_{ik}-u_{jk}\right|\right)$ throughout the Appendix A.

\vspace{4mm} \noindent \textbf{Proof of Proposition 1}. It is easily
obtained that $\bigtriangleup_{ijk}=3/2$ when $u_{ik}=u_{jk},
k=1,2,\dots,m+r$; $\bigtriangleup_{ijk}=5/4$ when $u_{ik}\neq
u_{jk}, k=1,2,\dots,m_1$; and $\bigtriangleup_{ijk}= {23}/{18}$ when
$u_{ik}\neq u_{jk}, k=m_1+1,\dots,m+r$. Then, the formula
(\ref{WD}) can be rewritten as
\begin{align}
&\mbox{WD}(D_{3})=-\left(\frac{4}{3}\right)^{m+r}+\frac{1}{(n+n_1)^2}
\sum_{i=1}^{n+n_1}\sum_{j=1}^{n+n_1}\prod_{k=1}^{m+r}\bigtriangleup_{ijk}\nonumber\\
=&-\left(\frac{4}{3}\right)^{m+r}+\frac{1}{(n+n_1)^2}\left(\sum_{i=1}^{n}\sum_{j=1}^{n}
\prod_{k=1}^{m+r}\bigtriangleup_{ijk}+\sum_{i=n+1}^{n+n_1}\sum_{j=n+1}^{n+n_1}
\prod_{k=1}^{m+r}\bigtriangleup_{ijk}\right. \left.
+2\sum_{i=1}^{n}\sum_{j=n+1}^{n+n_1}
\prod_{k=1}^{m+r}\bigtriangleup_{ijk}\right) \nonumber\\
=&-\left(\frac{4}{3}\right)^{m+r}+\frac{1}{(n+n_1)^2}\left(\frac{3}{2}\right)^{r}
\sum_{i=1}^{n}\sum_{j=1}^{n}\prod_{k=1}^{m}\bigtriangleup_{ijk}
+\frac{1}{(n+n_1)^2}
\left(\frac{5}{4}\right)^{m_1}\left(\frac{23}{18}\right)^{m_2+r}
\times\nonumber\\
&~~~~\left(\sum_{i=n+1}^{n+n_1}\sum_{j=n+1}^{n+n_1}
\left(\frac{6}{5}\right)^{\cf_{ij}}\left(\frac{27}{23}\right)^{\cv_{ij}+\ct_{ij}}
+2\sum_{i=1}^{n}\sum_{j=n+1}^{n+n_1}
\left(\frac{6}{5}\right)^{\cf_{ij}}\left(\frac{27}{23}\right)^{\cv_{ij}}
\right).\label{proof11}
\end{align}

 \noindent According to the expression of WD,
\begin{align}\label{proof13}
\sum_{i=1}^{n}\sum_{j=1}^{n}\prod_{k=1}^m\bigtriangleup_{ijk}
=&n^2\left(WD(\d_0)+\left(\frac{4}{3}\right)^m\right). \end{align}

\noindent Then,  (\ref{WD D3})  can be obtained by substituting
(\ref{proof13}) into (\ref{proof11}).

\vspace{4mm}
For obtaining the lower bounds of WD, we give the two lemmas first.
The proof of Lemma 1 is easy and we omit it.
\begin{lemma}\label{le1}For any column augmented design
 $D_{3}\in \cC_{3}(n+n_1;2^{m_1}3^{m_2}\bullet 3^{r})$,\\
\noindent
(1)$\sum\limits_{i=1}^n\sum\limits_{j(\neq i)=1}^n\cf_{ij} =\frac {m_1n(n-2)}{2}$,
(2)$\sum\limits_{i=n+1}^{n+n_1}\sum\limits_{j(\neq i)=n+1}^{n+n_1}\cf_{ij} 
    =\frac {m_1n_1(n_1-2)}{2}$,
(3)$\sum\limits_{i=1}^n\sum\limits_{j=n+1}^{n+n_1}\cf_{ij} =\frac
     {m_1nn_1}{2}$,\\
(4)$\sum\limits_{i=1}^n\sum\limits_{j(\neq i)=1}^n\cv_{ij} =\frac {m_2n(n-3)}{3}$,
(5)$\sum\limits_{i=n+1}^{n+n_1}\sum\limits_{j(\neq i)=n+1}^{n+n_1}\cv_{ij} 
    =\frac {m_2n_1(n_1-3)}{3}$,
(6)$\sum\limits_{i=1}^n\sum\limits_{j=n+1}^{n+n_1}\cv_{ij} =\frac{m_2nn_1}{3}$, \\
(7)$\sum\limits_{i=n+1}^{n+n_1}\sum\limits_{j(\neq i)=n+1}^{n+n_1}\ct_{ij} 
     =\frac {rn_1(n_1-2)}{2}$.
\end{lemma}

\begin{lemma} \label{le2}
(1) (Elsawah and Qin \cite{EQ15b}) When $\sum_{i=1}^nx_i=c$ and
$x_1,\cdots,x_n$ are nonnegative, then for any positive integer $t$,
we have \bea\label{lemma20} \sum_{i=1}^n\left(x_i\right)^t \geq
n\varphi^t, \eea where $\varphi=c/n$, the equality holds if and only
if $x_1=\dots=x_n=\varphi.$

(2) (Chatterjee et al. \cite{CLQ12}) Let $x_1,\cdots,x_n$ be the $n$
nonnegative integers and
 $\sum_{i=1}^nx_i=c$, then, any positive integer $t$,
\bea\label{lemma21} \sum_{i=1}^n\left(x_i\right)^t \geq
pw^t+q\left(w+1\right)^t, \eea  where $p$, $q$ are integers such
that $w=\left\lfloor\frac{c}{n}\right\rfloor$,$~p+q=n,
~pw+q(w+1)=c$. The equality holds if and only if among the $n$
number of $x_i$, $p$ of them take $w$, $q$ of them take $w+1$.
\end{lemma}

Lemma \ref{le2}(1) is suitable for nonnegative non-integer situation, thus 
there is less restriction for $x_i$ in Lemma \ref{le2}(1) than Lemma 
\ref{le2}(2).  When $x_i$ meets the assumptions in Lemma \ref{le2}(2), 
inequalities (\ref{lemma20}) and (\ref{lemma21}) are simultaneously valid. 
While, we know that the right of (\ref{lemma21}) is larger than the right 
of (\ref{lemma20}) according to the property of convex function. Hence, 
(\ref{lemma21}) is more compact than (\ref{lemma20}). In other words, 
if $\sum_{i=1}^n\left(x_i\right)^t$ achieves the lower bound of 
(\ref{lemma20}), then it is sure to reach the lower bound of 
(\ref{lemma21}), but not vice versa. If $x_i$ are nonnegative integers, 
we should use (\ref{lemma21}) rather than (\ref{lemma20}).

\vspace{4mm} \noindent \textbf{Proof of Theorem 1}. The portion of
the fourth term in  (\ref{WD D3}),
\begin{align}  &\sum_{i=n+1}^{n+n_1}\sum_{j(\neq
i)=n+1}^{n+n_1}\left(\frac{6}{5}\right)^{\cf_{ij}}\left(\frac{27}{23}\right)^{\cv_{ij}+\ct_{ij}}
=\sum_{i=n+1}^{n+n_1}\sum_{j(\neq i)=n+1}^{n+n_1}
e^{a\cf_{ij}+b(\cv_{ij}+\ct_{ij})} \nonumber\\
& = \sum_{i=n+1}^{n+n_1}\sum_{j(\neq
i)=n+1}^{n+n_1}\sum_{t=0}^{\infty}\frac{\left(
a\cf_{ij}+b(\cv_{ij}+\ct_{ij})\right)^t}{t!} \nonumber\\
&=\sum_{t=0}^{\infty}\frac{1}{t!}\left(\sum_{i=n+1}^{n+n_1}\sum_{j(\neq
i)=n+1}^{n+n_1}\left(a\cf_{ij}+b(\cv_{ij}+\ct_{ij})\right)^t \right)\nonumber\\
&\geq \sum_{t=0}^{\infty}\frac{1}{t!}n_1(n_1-1)\varphi_1^t
 =n_1(n_1-1)e^{\varphi_1} =T_1. \label{proof21_2}
\end{align}

\noindent The last inequality holds resulting from Lemma \ref{le1}
and \ref{le2}(1).  Similarly, we have
 \begin{align}\label{proof26}
\sum_{i=1}^{n}\sum_{j=n+1}^{n+n_1}\left(\frac{6}{5}\right)^{\cf_{ij}}\left(\frac{27}{23}\right)^{\cv_{ij}
} \geq   nn_1e^{\varphi_2}=T_2. \end{align} The lower bound in
 (\ref{LBW3_3})  can be obtained by substituting (\ref{proof21_2}) and
(\ref{proof26}) into  (\ref{WD D3}).
 Moreover, the
equality in (\ref{proof21_2}) holds if and only if
$a\cf_{ij}+b(\cv_{ij}+\ct_{ij})=\varphi_1$, $i=n+1,\dots, n+n_1,
j(\neq i)=n+1,\dots, n+n_1$ and the equality in (\ref{proof26})
holds if and only if $a\cf_{ij}+b\cv_{ij}=\varphi_2$, $i=1,\dots, n,
j=n+1,\dots, n+n_1$.

\vspace{4mm} \noindent \textbf{Proof of Theorem 2}. When
$\d_0\in\cu(n;2^{m})$, let $m_1=m, m_2=0$ in (\ref{WD D3}), then
the WD-value of column augmented design $D_3\in
\cC_3(n+n_1;2^{m}\bullet 3^{r})$ is,  \begin{align} \mbox{WD}(D_3)=&
C(r) +\frac{n^2}{(n+n_1)^2}\left(\frac{3}{2}\right)^{r}WD(\d_0)
  +\frac{1}{(n+n_1)^2}\left(\frac{5}{4}\right)^{m}\left(\frac{23}{18}\right)^{r} \times \nonumber\\
   &\left(
\sum_{i=n+1}^{n+n_1}\sum_{j(\neq
i)=n+1}^{n+n_1}\left(\frac{6}{5}\right)^{\cf_{ij}}\left(\frac{27}{23}\right)^{\ct_{ij}}
+2\sum_{i=1}^{n}\sum_{j=n+1}^{n+n_1}
\left(\frac{6}{5}\right)^{\cf_{ij}} \right). \label{WD D3_2}
\end{align} The part of the fourth term in (\ref{WD D3_2}),
\begin{align}
&\sum_{i=1}^{n}\sum_{j=n+1}^{n+n_1}\left(\frac{6}{5}\right)^{\cf_{ij}}
=\sum_{i=1}^{n}\sum_{j=n+1}^{n+n_1} e^{a\cf_{ij}}  =
\sum_{i=1}^{n}\sum_{j=n+1}^{n+n_1}\sum_{t=0}^{\infty}\frac{\left(
a\cf_{ij}\right)^t}{t!} =\sum_{t=0}^{\infty}\frac{1}{t!}\left(\sum_{i=1}^{n}\sum_{j=n+1}^{n+n_1}\left(a\cf_{ij}\right)^t \right) \nonumber\\
&\geq  \sum_{t=0}^{\infty}\frac{1}{t!}\left[p_{1}
(aw_{1})^t+q_{1}(a(w_{1}+1))^t \right] \label{proof4_1} =
p_{1}\left(\frac{6}{5}\right)^{w_{1}}
+q_{1}\left(\frac{6}{5}\right)^{w_{1}+1}.\end{align}
The (\ref{proof4_1}) holds because of Lemma \ref{le1} and
\ref{le2}(2). The proof is complete.

\vspace{4mm} \noindent \textbf{Proof of Theorem 4}. For
$\d_0\in\cu(n;2^{m_1}3^{m_2})$, we have
\begin{align}\label{proof51}
&\mbox{WD}(D_{3b})
=-\left(\frac{4}{3}\right)^{m+r+1}+\frac{1}{(n+n_1)^2}
\sum_{i=1}^{n+n_1}\sum_{j=1}^{n+n_1}\prod_{k=1}^{m+r+1}\bigtriangleup_{ijk}\nonumber\\
&=-\left(\frac{4}{3}\right)^{m+r+1}+\frac{1}{(n+n_1)^2}\left(\left(\frac{3}{2}\right)^{r+1}\sum_{i=1}^{n}\sum_{j=1}^{n}
\prod_{k=1}^{m}\bigtriangleup_{ijk}+\sum_{i=n+1}^{n+n_1}\sum_{j=n+1}^{n+n_1}
\prod_{k=1}^{m+r+1}\bigtriangleup_{ijk}\right.\nonumber\\
&~~~~\left.+2\left(\frac{5}{4}\right)\left(\frac{23}{18}\right)^{r}\sum_{i=1}^{n}\sum_{j=n+1}^{n+n_1}
\prod_{k=1}^{m}\bigtriangleup_{ijk}\right)\nonumber\\
&=
 C(r+1) +\frac{n^2}{(n+n_1)^2}\left(\frac{3}{2}\right)^{r+1}WD(\d_0) +\frac{1}{(n+n_1)^2}\left(\frac{5}{4}\right)^{m_1}\left(\frac{23}{18}\right)^{m_2+r}\times\nonumber\\
& ~~~~\left(
\left(\frac{3}{2}\right)\sum_{i=n+1}^{n+n_1}\sum_{j(\neq
i)=n+1}^{n+n_1}\left(\frac{6}{5}\right)^{\cf_{ij}}\left(\frac{27}{23}\right)^{\cv_{ij}+\ct_{ij}}
+2\left(\frac{5}{4}\right)\sum_{i=1}^{n}\sum_{j=n+1}^{n+n_1}
\left(\frac{6}{5}\right)^{\cf_{ij}}\left(\frac{27}{23}\right)^{\cv_{ij}
}  \right).
\end{align}
There is only a different positive coefficient in each term between
(\ref{proof51}) and  (\ref{WD D3}). For $\d_0\in\cu(n;2^{m})$ and
$\d_0\in\cu(n;3^{m})$, yield the similar result, then Theorem 4
can be easily obtained.

\vspace{4mm} {\noindent \textbf{Proof of Proposition 3}. Without
loss of generality, we only consider the three-stage \caug, more
stages are similar.
\begin{align}\label{proof equi2 2}
&\mbox{WD}(D^{(3)})=-\left(\frac{4}{3}\right)^{m}+\frac{1}{(n+n_1+n_2)^2}
\sum_{i=1}^{n+n_1+n_2}\sum_{j=1}^{n+n_1+n_2}\prod_{k=1}^{m}\bigtriangleup_{ijk}\nonumber\\
&=-\left(\frac{4}{3}\right)^{m}+\frac{1}{(n+n_1+n_2)^2}\left(\sum_{i=1}^{n}\sum_{j=1}^{n}
\prod_{k=1}^{m}\bigtriangleup_{ijk}+\sum_{i=n+1}^{n+n_1}\sum_{j=n+1}^{n+n_1}
\prod_{k=1}^{m}\bigtriangleup_{ijk}
+\sum_{i=n+n_1+1}^{n+n_1+n_2}\sum_{j=n+n_1+1}^{n+n_1+n_2}
\prod_{k=1}^{m}\bigtriangleup_{ijk}\right.\nonumber\\
&~~~~\left. +2\sum_{i=1}^{n}\sum_{j=n+1}^{n+n_1}
\prod_{k=1}^{m}\bigtriangleup_{ijk}
+2\sum_{i=1}^{n}\sum_{j=n+n_1+1}^{n+n_1+n_2}
\prod_{k=1}^{m}\bigtriangleup_{ijk}+2\sum_{i=n+1}^{n+n_1}\sum_{j=n+n_1+1}^{n+n_1+n_2}
\prod_{k=1}^{m}\bigtriangleup_{ijk}\right), \end{align} and
\begin{align}
\mbox{WD}(D_{b}^{(3)})&=-\left(\frac{4}{3}\right)^{m+1}+\frac{1}{(n+n_1+n_2)^2}
\sum_{i=1}^{n+n_1+n_2}\sum_{j=1}^{n+n_1+n_2}\prod_{k=1}^{m+1}\bigtriangleup_{ijk}\nonumber\\
&=-\left(\frac{4}{3}\right)^{m+1}+\frac{1}{(n+n_1+n_2)^2}\left(\left(\frac{3}{2}\right)\sum_{i=1}^{n}\sum_{j=1}^{n}
\prod_{k=1}^{m}\bigtriangleup_{ijk}+\left(\frac{3}{2}\right)\sum_{i=n+1}^{n+n_1}\sum_{j=n+1}^{n+n_1}
\prod_{k=1}^{m}\bigtriangleup_{ijk}\right. \nonumber\\
&~~~~\left.+\left(\frac{3}{2}\right)\sum_{i=n+n_1+1}^{n+n_1+n_2}\sum_{j=n+n_1+1}^{n+n_1+n_2}
\prod_{k=1}^{m}\bigtriangleup_{ijk}
+2\left(\frac{23}{18}\right)\sum_{i=1}^{n}\sum_{j=n+1}^{n+n_1}
\prod_{k=1}^{m}\bigtriangleup_{ijk}\right.\nonumber\\
&~~~~\left.+2\left(\frac{23}{18}\right)\sum_{i=1}^{n}\sum_{j=n+n_1+1}^{n+n_1+n_2}
\prod_{k=1}^{m}\bigtriangleup_{ijk}+2\left(\frac{23}{18}\right)\sum_{i=n+1}^{n+n_1}\sum_{j=n+n_1+1}^{n+n_1+n_2}
\prod_{k=1}^{m}\bigtriangleup_{ijk}\right). \label{proof equi2
3}\end{align} It is easily known that the (\ref{proof equi2 2}) and
(\ref{proof equi2 3}) only have a  different positive coefficient in
each term. Hence, the row \aug design $D^{(3)}$ achieves its lower
bound then the corresponding column \aug design $D_{b}^{(3)}$
achieves its lower bound, and
vice versa. 

For a given initial design $\d_0$, let $r=0$ in Proposition 1, one
can get the WD-value of its two-stage row \aug design $D\in
\cC(n+n_1;2^{m_1}3^{m_2})$,
 \begin{align} &\mbox{WD}(D)=
  -\frac{n_1^2+2nn_1}{(n+n_1)^2}\left(\frac{4}{3}\right)^{m}
  +\frac {n_1}{(n+n_1)^2}\left(\frac{3}{2}\right)^{m}
  +\frac{n^2}{(n+n_1)^2}WD(\d_0) \nonumber\\
  & +\frac{1}{(n+n_1)^2}\left(\frac{5}{4}\right)^{m_1}\left(\frac{23}{18}\right)^{m_2}\left(
  \sum_{i=n+1}^{n+n_1}\sum_{j(\neq
i)=n+1}^{n+n_1}\left(\frac{6}{5}\right)^{\cf_{ij}}\left(\frac{27}{23}\right)^{\cv_{ij}}+2\sum_{i=1}^{n}\sum_{j=n+1}^{n+n_1}
\left(\frac{6}{5}\right)^{\cf_{ij}}\left(\frac{27}{23}\right)^{\cv_{ij}}\right).\label{two_stage_row}
\end{align}
For the $l$-stage row \aug design, one can regard the first $l-1$
stages as the initial design and the $l$-th stage as the follow-up
part. Replacing $n,n_1,\d_0$ in (\ref{two_stage_row}) by
$N_{l-1},n_{l-1}, D^{l-1}$
and using the similar technique in the proof of Theorem 1,  the
proof is complete. 

\vspace{2mm}

{For obtaining  the lower bounds of   $E(f_{NOD})$ of \caugs~$D_3$,
$D_{3b}$ and $D_{3B}$, some notations are given first.
 $Q_1=\zeta_1\psi_1^2+\eta_1(\psi_1+1)^2$,
~$Q_1'=\zeta_1(\psi_1+1)^2+\eta_1(\psi_1+2)^2$,
~$Q_1''=\zeta_1'(\psi_1'+1)^2+\eta_1'(\psi_1'+2)^2$,
~$Q_2=\zeta_2\psi_2^2+\eta_2(\psi_2+1)^2$,
~$Q_2'=\zeta_2'\psi_2'^2+\eta_2'(\psi_2'+1)^2$, and the parameters
$\zeta_i,\eta_i,$ and $\psi_i$ meet the following requirements:
$\zeta_1+\eta_1=n_1(n_1-1)$,~$\zeta_1\psi_1+\eta_1(\psi_1+1) =\frac
{m_1n_1(n_1-2)}{2}+\frac {m_2n_1(n_1-3)}{3}+\frac {rn_1(n_1-2)}{2}$,
~$\psi_1=\left\lfloor\left(\frac {m_1(n_1-2)}{2}+\frac
{m_2(n_1-3)}{3}+\frac {r(n_1-2)}{2}\right)/(n_1-1)\right\rfloor$;
~$\zeta_1'+\eta_1'=n_1(n_1-1)$,~$\zeta_1'\psi_1'+\eta_1'(\psi_1'+1)
=\frac {(m_1+1)n_1(n_1-2)}{2}+\frac {m_2n_1(n_1-3)}{3}+\frac
{rn_1(n_1-2)}{2}$, ~$\psi_1'=\left\lfloor\left(\frac
{(m_1+1)(n_1-2)}{2}+\frac {m_2(n_1-3)}{3}+\frac
{r(n_1-2)}{2}\right)/(n_1-1)\right\rfloor$;
~$\zeta_2+\eta_2=n_1n$,~$\zeta_2\psi_2+\eta_2(\psi_2+1) =\frac
{m_1n_1n}{2}+\frac {m_2n_1n}{3}$, ~$\psi_2=\left\lfloor\frac
{m_1}{2}+\frac {m_2}{3}\right\rfloor$;
~$\zeta_2'+\eta_2'=n_1n$,~$\zeta_2'\psi_2'+\eta_2'(\psi_2'+1) =\frac
{(m_1+1)n_1n}{2}+\frac {m_2n_1n}{3}$, ~$\psi_2'=\left\lfloor\frac
{m_1+1}{2}+\frac {m_2}{3}\right\rfloor$.  For the column augmented
designs, we have the following results under the criterion
$E(f_{NOD})$.

\begin{lemma}\label{thlbf}  Given  an initial design $\d_0\in\cu(n;2^{m_1}3^{m_2})$, the lower bounds of $E(f_{NOD})$
of \caugs~$D_3$, $D_{3b}$ and $D_{3B}$  respectively are,

\bea &LBf_{3} =
\frac{m(m-1)}{(m+r)(m+r-1)}E(f_{NOD}(\d_0))+\frac{1}{(m+r)(m+r-1)}(Q_1+2Q_2)
+\frac{(n+n_1)(m+r)}{m+r-1}\nonumber\\
&-\frac{1}{(m+r)(m+r-1)}\left[nm^2-n^2\left( \frac{m_1+m_1^2}{4}
+\frac{2m_2+m_2^2}{9}+\frac{m_1m_2}{3} \right) -rm_1n(n-2) \right.\nonumber\\
&-\frac{2rm_2n(n-3)}{3}-n(n-1)r^2+
 r\left(n^2+\frac{n_1^2}{2}\right)
 +\left(\frac{m_1}{2}+\frac{m_2}{3}+\frac{m_1(m_1-1)}{4}\right.\nonumber\\
 &\left.\left.+\frac{(m_2+r)(m_2+r-1)}{9}+\frac{m_1(m_2+r)}{3}\right)(n+n_1)^2\right],\label{LBf3}\eea \bea
& LBf_{3b} =  \frac{m(m-1)}{(m+r+1)(m+r)}E(f_{NOD}(\d_0))+\frac{1}{(m+r+1)(m+r)}(Q_1'+2Q_2)+\frac{(n+n_1)(m+r+1)}{m+r}\nonumber\\
 &-\frac{1}{(m+r+1)(m+r)}\left[nm^2-n^2\left( \frac{m_1+m_1^2}{4}
 +\frac{2m_2+m_2^2}{9}+\frac{m_1m_2}{3} \right)-(r+1)m_1n(n-2)\right. \nonumber\\
&  -\frac{2(r+1)m_2n(n-3)}{3}-n(n-1)(r+1)^2 +r\left(n^2+\frac{n_1^2}{2}\right)+n^2+n_1^2 \nonumber\\
 & +\left.\left(\frac{m_1}{2}+\frac{m_2}{3}+\frac{m_1(m_1+1)}{4}
 +\frac{(m_2+r)(m_2+r-1)}{9}+\frac{(m_1+1)(m_2+r)}{3}\right)(n+n_1)^2\right],\label{LBf3b} \eea \bea
& LBf_{3B} =  \frac{m(m-1)}{(m+r+2)(m+r+1)}E(f_{NOD}(\d_0))
+ \frac{1}{(m+r+2)(m+r+1)}(Q_1''+2Q_2')+\frac{(n+n_1)(m+r+2)}{m+r+1}\nonumber\\
 &-\frac{1}{(m+r+2)(m+r+1)}\left[nm^2-n^2\left( \frac{m_1+m_1^2}{4}
 +\frac{2m_2+m_2^2}{9}+\frac{m_1m_2}{3} \right) -(r+2)m_1n(n-2)\right. \nonumber\\
 &-\frac{2(r+2)m_2n(n-3)}{3}-n(n-1)(r+2)^2+(r+2)n^2+\frac{(3+r)n_1^2}{2}+nn_1\nonumber\\
&+\left.\left(\frac{m_1}{2}+\frac{m_2}{3}+\frac{(m_1+2)(m_1+1)}{4}
+\frac{(m_2+r)(m_2+r-1)}{9}+\frac{(m_1+2)(m_2+r)}{3}\right)(n+n_1)^2\right].\label{LBf3B2}
  \eea
Moreover,   $LBf_3$  can be achieved if and only if there are
 $\zeta_1$ number of $\cf_{ij}+\cv_{ij}+\ct_{ij}$
 take $\psi_1$, $\eta_1$ number of $\cf_{ij}+\cv_{ij}+\ct_{ij}$ take $\psi_1+1$$,i=n+1,\dots, n+n_1, j(\neq i)=n+1,\dots, n+n_1$, and $\zeta_2$ number of $\cf_{ij}+\cv_{ij}$
 take $\psi_2$, $\eta_2$ number of $\cf_{ij}+\cv_{ij}$ take $\psi_2+1$
$,i=1,\dots, n, j=n+1,\dots, n+n_1$.  Likewise, one can obtain the
conditions when the lower bounds $LBf_{3b}$  and $LBf_{3B}$ are
achievable.
\end{lemma}
}

\vspace{4mm} \noindent \textbf{Proof}. Similar to
\cite{FLL03b}, let $\N_{kl}(D_3)=(n_{uv}^{kl}(D_3))$ and
$\Z^j(D_3)=(z^j_{iu}(D_3))$ be an $(n+n_1)\times q_j$  matrix, where
\bea\label{proof11_1}
z^j_{iu}(D_3) = \left\{\ba{ll} 1, & if~ x_{ij}=u, \\
0, & otherwise, \ea\right. i=1,\dots,n+n_1,~u=0,\dots,q_j-1. \eea
Denote $\Z(D_3)=(\Z^1(D_3),\dots,\Z^{m+r}(D_3))$. From
(\ref{proof11_1}), one can obtain
$\Z(D_3)\Z(D_3)'=(\lambda_{ij}(D_3)), $ where $\lambda_{ij}(\d)$ is
the number of coincidences between the $i$-th and the $j$-th rows in
$\d$, moreover, $$ \Z(D_3)'\Z(D_3)= \left(\ba{ccccc}
\frac{n+n_1}{q_1}I_{q_1} & \dots &\N_{1m}(D_3)& \dots & \N_{1(m+r)}(D_3)\\
\dots & \dots &  \dots  & \dots &  \dots \\
\N_{m1}(D_3) & \dots  & \frac{n+n_1}{q_m}I_{q_m}  & \dots & \N_{m(m+r)}(D_3)\\
\dots & \dots &  \dots  & \dots &  \dots  \\
\N_{(m+r)1}(D_3)&\dots  & \dots &  \dots &  \Y_3 \ea\right), $$ and
$\Y_3= \left(\ba{ccc}
n& 0 &0\\
0 & \frac{n_1}{2} &  0   \\
0 & 0  & \frac{n_1}{2}  \ea\right).$ From the definition of
$f^{kl}_{NOD}$,
\begin{align*}
f^{kl}_{NOD}(D_3)&= \sum_{u=0}^{q_k-1}\sum_{v=0}^{q_l-1}\left((n_{uv}^{kl}(D_3))^2-2n_{uv}^{kl}(D_3)\frac{n+n_1}{q_kq_l}+\frac{(n+n_1)^2}{(q_kq_l)^2}\right)  \\
&=\sum_{u=0}^{q_k-1}\sum_{v=0}^{q_l-1}\left(n_{uv}^{kl}(D_3)\right)^2-\frac{(n+n_1)^2}{q_kq_l}
 =tr(\N_{kl}(D_3)\N_{lk}(D_3))-\frac{(n+n_1)^2}{q_kq_l},
\end{align*}
then
 \bea
&&E(f_{NOD}(D_3))= \sum_{1\leq k<l\leq m+r}f^{kl}_{NOD}(D_3)/\dbinom{m}{2}\nonumber\\
&=&\frac{1}{(m+r)(m+r-1)}\left(\sum_{k=1}^{m+r}\sum_{l(\neq k)=1}^{m+r}tr(\N_{kl}(D_3)\N_{lk}(D_3)-\sum_{k=1}^{m+r}\sum_{l(\neq k)=1}^{m+r}\frac{(n+n_1)^2}{q_kq_l}\right)\nonumber\\
&=&\frac{1}{(m+r)(m+r-1)}\left(tr(\Z(D_3)'\Z(D_3))^2-\sum_{k=1}^{m}\frac{(n+n_1)^2}{q_k}-r\left(n^2+\frac{n_1^2}{2}\right)\right.\nonumber\\
&&\left.-\sum_{k=1}^{m+r}\sum_{l(\neq k)=1}^{m+r}\frac{(n+n_1)^2}{q_kq_l} \right)\nonumber\\
&=&\frac{1}{(m+r)(m+r-1)}\left(\sum_{i=1}^{n+n_1}\sum_{j(\neq i)=1}^{n+n_1}\lambda_{ij}(D_3)^2-\sum_{k=1}^{m}\frac{(n+n_1)^2}{q_k}-r\left(n^2+\frac{n_1^2}{2}\right)\right.\nonumber\nonumber\\
&&\left.-\sum_{k=1}^{m+r}\sum_{l(\neq k)=1}^{m+r}\frac{(n+n_1)^2}{q_kq_l}\right)+\frac{(n+n_1)(m+r)}{m+r-1}\nonumber\\
&=&\frac{(n+n_1)(m+r)}{m+r-1}+\frac{1}{(m+r)(m+r-1)}\left(\sum_{i=1}^{n+n_1}\sum_{j(\neq i)=1}^{n+n_1}\lambda_{ij}(D_3)^2 -r\left(n^2+\frac{n_1^2}{2}\right)\right.\nonumber\\
&&\left. -\left(\frac{m_1}{2}+\frac{m_2}{3}+\frac{m_1(m_1-1)}{4}+\frac{(m_2+r)(m_2+r-1)}{9}+\frac{m_1(m_2+r)}{3}\right)(n+n_1)^2\right).\nonumber\\
\label{proof11_2} \eea
The value of $E(f_{NOD}(D_3))$ is a function of $\lambda_{ij}$, and
the other terms in (\ref{proof11_2}) are constant. For $D_3$,

\bea\label{sum}
&&\sum_{i=1}^{n+n_1}\sum_{j(\neq i)=1}^{n+n_1}\lambda_{ij}(D_3)^2=\sum_{i=1}^{n+n_1}\sum_{j(\neq i)=1}^{n+n_1}(\cf_{ij}+\cv_{ij}+\ct_{ij})^2=\sum_{i=1}^{n}\sum_{j(\neq i)=1}^{n}(\cf_{ij}+\cv_{ij}+r)^2\nonumber\\
&&+\sum_{i=n+1}^{n+n_1}\sum_{j(\neq
i)=n+1}^{n+n_1}(\cf_{ij}+\cv_{ij}+\ct_{ij})^2
+2\sum_{i=1}^{n}\sum_{j=n+1}^{n+n_1}(\cf_{ij}+\cv_{ij})^2. \eea
\noindent { In addition, for a given initial design $\d_0 \in
\cu(n;2^{m_1}3^{m_2})$, based on the Theorem 1 of Fang et al.
\cite{FLL03b}, we have $
E(f_{NOD}(\d_0))=\frac{1}{m(m-1)}\sum_{i=1}^{n}\sum_{j(\neq
i)=1}^{n}(\cf_{ij}+\cv_{ij})^2+\frac{nm}{m-1}-\frac{n^2}{m(m-1)}\left(
\frac{m_1+m_1^2}{4}+\frac{2m_2+m_2^2}{9}+\frac{m_1m_2}{3} \right). $
Then{\small \bea \label{sum3}  \sum_{i=1}^{n}\sum_{j(\neq
i)=1}^{n}(\cf_{ij}+\cv_{ij})^2&=&m(m-1)E(f_{NOD}(\d_0))-nm^2 +n^2
\left(
\frac{m_1+m_1^2}{4}+\frac{2m_2+m_2^2}{9}+\frac{m_1m_2}{3} \right).\nonumber\\
\eea } From (\ref{sum3}) and Lemma \ref{le1},  the first term of the
right side of (\ref{sum}) has, \bea &&\sum_{i=1}^{n}\sum_{j(\neq
i)=1}^{n}(\cf_{ij}+\cv_{ij}+r)^2 =\sum_{i=1}^{n}\sum_{j(\neq
i)=1}^{n}(\cf_{ij}+\cv_{ij})^2+2r\sum_{i=1}^{n}\sum_{j(\neq
i)=1}^{n}(\cf_{ij}+\cv_{ij})+n(n-1)r^2  \nonumber\\
&&=m(m-1)E(f_{NOD}(\d_0))-nm^2+n^2\left( \frac{m_1+m_1^2}{4} \right)
+rm_1n(n-2)
+\frac{2m_2+m_2^2}{9}+\frac{m_1m_2}{3} \nonumber\\
&&~~~~+\frac{2rm_2n(n-3)}{3}+n(n-1)r^2. \label{sum4} \eea }
Moreover, according to Lemma \ref{le1} and Lemma \ref{le2}(2), we
have \bea\label{sum1} \sum_{i=n+1}^{n+n_1}\sum_{j(\neq
i)=n+1}^{n+n_1}(\cf_{ij}+\cv_{ij}+\ct_{ij})^2 &\geq& Q_1,
\\ \label{sum2}
\sum_{i=1}^{n}\sum_{j=n+1}^{n+n_1}(\cf_{ij}+\cv_{ij})^2 &\geq& Q_2.
\eea \noindent {The equalities of (\ref{sum1}) and
(\ref{sum2}) hold if and only if there are $\zeta_1$ number of
$\cf_{ij}+\cv_{ij}+\ct_{ij}$ take $\psi_1$, $\eta_1$ number of
$\cf_{ij}+\cv_{ij}+\ct_{ij}$ take $\psi_1+1$$,i=n+1,\dots, n+n_1,
j(\neq i)=n+1,\dots, n+n_1$, and $\zeta_2$ number of
$\cf_{ij}+\cv_{ij}$
 take $\psi_2$, $\eta_2$ number of $\cf_{ij}+\cv_{ij}$ take $\psi_2+1$
$,i=1,\dots, n, j=n+1,\dots, n+n_1$.} Substituting (\ref{sum}) and
(\ref{sum4})-(\ref{sum2}) into (\ref{proof11_2}), we obtain
(\ref{LBf3}).

For $D_{3b}$, note that $ \Z(D_{3b})'\Z(D_{3b})=
\left(\ba{cccc}
\frac{n+n_1}{q_1}I_{q_1} & \dots &\N_{1(m+r)}(D_{3b})& \N_{1(m+r+1)}(D_{3b})\\
\dots & \dots &  \dots  &  \dots \\
\N_{(m+r)1}(D_{3b}) & \dots  & \Y_3 & \dots\\
\N_{(m+r+1)1}(D_{3b})&\dots   &  \dots &  \Y_{3b} \ea\right), $
where $\Y_{3b}= \left(\ba{cc}
n& 0 \\
0 & n_1  \ea\right),$ and \bea \sum_{i=1}^{n+n_1}\sum_{j(\neq
i)=1}^{n+n_1}\lambda_{ij}(D_{3b})^2 &&=\sum_{i=1}^{n}\sum_{j(\neq
i)=1}^{n}(\cf_{ij}+\cv_{ij}+r+1)^2 +\sum_{i=n+1}^{n+n_1}\sum_{j(\neq
i)=n+1}^{n+n_1}(\cf_{ij}+\cv_{ij}+\ct_{ij}+1)^2 \nonumber\\
&&+2\sum_{i=1}^{n}\sum_{j=n+1}^{n+n_1}(\cf_{ij}+\cv_{ij})^2.\label{lammdaij}
\eea

 For $D_{3B}$, we have
{\small $ \Z(D_{3B})'\Z(D_{3B})= \left(\ba{ccccc}
\frac{n+n_1}{q_1}I_{q_1} & \dots &\N_{1(m+r)}(D_{3B})& \dots & \N_{1(m+r+2)}(D_{3B})\\
\dots & \dots &  \dots  & \dots &  \dots \\
\N_{(m+r)1}(D_{3B}) & \dots  & \Y_3  & \dots & \dots\\
\dots & \dots &  \dots  & \Y_{3b} &  \dots  \\
\N_{(m+r+2)1}(D_{3B})&\dots  & \dots &  \dots &  \Y_{3B} \ea\right),
$} where $\Y_{3B}= \left(\ba{cc}
n+\frac{n_1}{2}& 0 \\
0 & \frac{n_1}{2}  \ea\right),$ and $ \sum_{i=1}^{n+n_1}\sum_{j(\neq
i)=1}^{n+n_1}\lambda_{ij}(D_{3B})^2 =\sum_{i=1}^{n}\sum_{j(\neq
i)=1}^{n}(\cf_{ij}+\cv_{ij}+r+2)^2+\sum_{i=n+1}^{n+n_1}\sum_{j(\neq
i)=n+1}^{n+n_1}(\cf_{ij}+\cv_{ij}+\ct_{ij}+1+\ca_{ij})^2
 +2\sum_{i=1}^{n}\sum_{j=n+1}^{n+n_1}(\cf_{ij}+\cv_{ij}+\ca_{ij})^2
$. The $E(f_{NOD}(D_{3b}))$ and $E(f_{NOD}(D_{3B}))$ can be
calculated similarly as $E(f_{NOD}(D_3))$, then, we obtain the
lower bounds in (\ref{LBf3b}) and (\ref{LBf3B2}).

\vspace{4mm} \noindent \textbf{Proof of Theorem 6}

(1) If the initial design is mixed-level, noting that all
$\cf_{ij}, \cv_{ij}$ and $\ct_{ij}$ are integers, and meet some
relations in Lemma \ref{le1}, we check the conditions that the lower
bound $LBW_3$ in Theorem 1 and the lower bound $LBf_3$ in Theorem 6
are achieved. Thus if  $LBW_3$ is reached then $LBf_3$ is also
reached, but, if $LBf_3$ is achieved, $LBW_3$ may not be achieved.
Likewise, the other situations can be proved.

 (2) When the initial design is symmetrical
 three-level, and $r>0$, the lower bound $LBW_3$ is in Theorem 3, and it can
  be reached if and only if $\cv_{ij}+\ct_{ij},i=n+1,\dots, n+n_1, j(\neq
  i)=n+1,\dots, n+n_1$ take
 $w_2$ or $w_2+1$, and $\cv_{ij},i=1,\dots, n, j=n+1,\dots, n+n_1$ take
 $w_3$ or $w_3+1$, in this time, the lower bound $LBf_3$ in Theorem 6 can also be reached. Actually, both of $LBW_3$ and $LBf_3$ use the tool of Lemma \ref{le2}(2).
 Other situations in (2)-(3) can be proved similarly, so  we omit  them.

\newpage

\section*{Appendix B}

\noindent In the Appendix B, the initial design in Example 3 and all the corresponding column augmented uniform designs in Example 1-3 are
given.

\vspace{4mm} \noindent\textbf{B1.} The
\caugs~in Example 1.

\bc

{\small  \tabcolsep=1pt \bt {c|c|c|c} \hline $D_3^*  \in
\cC_3(12+6;2^53^7\bullet 3^1)$ & $D_3^*  \in
\cC_3(12+6;2^53^7\bullet 3^2)$ & $D_3^*  \in
\cC_3(12+6;2^53^7\bullet 3^3)$ & $D_3^*  \in
\cC_3(12+6;2^53^7\bullet 3^4)$

\\\hline
 $ \ba{cc} \d_0^* & \zero_{12\times 1} \\
0   1   1   0   1   0   1   1   1   0   2   1   &1\\
1   0   1   0   0   1   2   2   0   1   0   0   &1\\
0   0   0   1   1   1   0   0   1   2   1   0   &1\\
0   1   1   0   0   2   0   0   2   1   2   2   &2\\
1   0   0   1   0   0   1   1   2   2   0   1   &2\\
1   1   0   1   1   2   2   2   0   0   1   2   &2

 \ea  $
  & $ \ba{cc}
\d_0^* & \zero_{12\times 2} \\
1   1   0   1   1   2   2   2   0   0   1   2   &1  1\\
0   0   0   1   1   1   0   0   1   2   1   0   &2  2\\
0   1   1   0   1   0   1   1   1   0   2   1   &1  2\\
0   1   1   0   0   2   0   0   2   1   2   2   &2  1\\
1   0   1   0   0   1   2   2   0   1   0   0   &1  2\\
1   0   0   1   0   0   1   1   2   2   0   1   &2  1

 \ea  $

 & $ \ba{cc}
\d_0^* & \zero_{12\times 3} \\
0   1   1   0   0   2   0   0   2   1   2   2&  1   1   2\\
1   1   0   1   1   2   2   2   0   0   1   2&  1   1   1\\
0   1   1   0   1   0   1   1   1   0   2   1&  2   2   1\\
1   0   0   1   0   0   1   1   2   2   0   1&  1   1   2\\
1   0   1   0   0   1   2   2   0   1   0   0&  2   2   1\\
0   0   0   1   1   1   0   0   1   2   1   0&  2   2   2

 \ea  $

  & $ \ba{cc}
\d_0^* & \zero_{12\times 4} \\
0   0   0   1   0   0   0   0   2   2   0   1&  1   2   2   2\\
0   1   1   0   1   0   0   0   1   0   2   2&  2   2   2   1\\
1   0   0   1   1   1   2   2   1   2   1   0&  2   2   1   1\\
1   1   1   0   0   2   2   2   2   0   2   1&  1   1   1   2\\
0   0   1   0   0   1   1   1   0   1   0   0&  2   1   1   1\\
1   1   0   1   1   2   1   1   0   1   1   2&  1   1   2   2

 \ea  $

  \\

 \hline $D_3^*  \in
\cC_3(12+12;2^53^7\bullet 3^1)$ & $D_3^*  \in
\cC_3(12+12;2^53^7\bullet 3^2)$ & $D_3^*  \in
\cC_3(12+12;2^53^7\bullet 3^3)$ & $D_3^*  \in
\cC_3(12+12;2^53^7\bullet 3^4)$

\\\hline
 $ \ba{cc} \d_0^* & \zero_{12\times 1} \\
1   0   0   0   0   1   2   2   2   2   0   0   &2\\
0   1   1   0   1   2   1   0   0   0   1   1   &2\\
0   0   0   1   1   2   1   1   2   2   1   1   &1\\
1   0   0   1   1   0   2   0   0   0   1   2   &2\\
1   0   1   0   0   0   1   1   0   1   0   1   &1\\
0   0   1   1   0   1   0   0   1   1   2   0   &2\\
0   0   1   0   1   0   0   2   1   0   0   0   &1\\
1   1   1   1   0   2   2   2   2   0   2   1   &1\\
1   1   1   0   1   1   2   1   1   2   2   2   &2\\
1   1   0   1   1   1   1   2   1   1   1   0   &1\\
0   1   0   0   0   2   0   0   2   1   2   2   &1\\
0   1   0   1   0   0   0   1   0   2   0   2   &2

 \ea  $
  & $ \ba{cc}
\d_0^* & \zero_{12\times 2} \\
0   0   1   1   1   0   0   1   1   0   2   0   &2  1\\
1   0   0   1   0   1   2   2   2   1   0   0   &2  2\\
1   0   1   0   0   0   2   1   2   2   0   1   &1  1\\
0   0   0   1   0   0   0   0   0   2   0   2   &1  2\\
0   1   1   0   0   2   1   1   0   1   0   2   &2  1\\
1   1   0   1   1   2   2   2   0   0   1   2   &1  1\\
0   1   1   0   0   2   0   2   2   0   2   0   &1  2\\
1   0   1   0   1   0   1   2   0   0   1   1   &2  2\\
1   1   1   0   1   1   2   0   1   2   2   2   &2  2\\
0   0   0   0   1   1   1   0   1   1   1   0   &1  1\\
1   1   0   1   0   1   1   1   1   1   2   1   &1  2\\
0   1   0   1   1   2   0   0   2   2   1   1   &2  1

 \ea  $

  & $ \ba{cc}
\d_0^* & \zero_{12\times 3} \\
1   0   1   0   0   0   2   2   2   0   0   0&  2   2   1\\
0   0   1   0   1   0   2   0   1   2   0   2&  1   1   2\\
0   0   0   1   0   0   0   1   2   2   0   1&  1   1   1\\
1   1   0   1   1   2   2   2   0   0   1   2&  1   1   1\\
1   0   1   1   1   0   1   1   0   0   1   1&  1   2   2\\
0   1   1   1   1   2   0   2   2   0   2   0&  2   1   2\\
1   1   0   1   0   1   2   1   1   2   2   2&  2   2   2\\
1   0   0   0   0   1   1   2   0   1   0   0&  1   1   2\\
1   1   1   0   1   1   1   1   1   1   2   1&  2   1   1\\
0   1   0   0   0   2   1   0   2   2   1   1&  2   2   2\\
0   1   1   0   0   2   0   0   0   1   2   2&  1   2   1\\
0   0   0   1   1   1   0   0   1   1   1   0&  2   2   1

 \ea  $

   & $ \ba{cc}
\d_0^* & \zero_{12\times 4} \\
 1  0   0   1   1   0   2   2   0   0   1   0&  1   1   1   1\\
1   0   1   0   0   1   1   2   0   1   0   0&  2   2   2   1\\
0   0   0   1   1   1   1   1   1   2   1   1&  2   1   2   2\\
1   0   0   1   0   0   2   2   2   2   0   2&  2   2   1   2\\
0   0   1   0   1   0   0   1   1   0   2   1&  2   2   1   1\\
1   1   1   1   1   2   1   2   2   0   2   1&  1   2   2   2\\
0   0   1   0   0   0   0   0   2   0   0   0&  1   1   2   2\\
1   1   1   1   0   1   2   1   1   1   2   2&  1   1   2   1\\
0   1   0   0   0   2   1   1   2   2   0   1&  1   1   1   1\\
0   1   0   1   0   2   0   0   0   1   1   2&  2   2   2   1\\
1   1   1   0   1   2   2   0   0   2   1   2&  2   1   1   2\\
0   1   0   0   1   1   0   0   1   1   2   0&  1   2   1   2
 \ea  $

   \\ \hline

 \et} \ec

\bc {\small  \tabcolsep=1pt \bt {c|c|c} \hline $D_{3B}^*  \in
\cC_{3B}(12+6;2^53^7\bullet 3^0 \bullet \B)$ & $D_{3B}^*  \in
\cC_{3B}(12+6;2^53^7\bullet 3^1\bullet \B)$ & $D_{3B}^*  \in
\cC_{3B}(12+6;2^53^7\bullet 3^2\bullet \B)$
\\\hline
 $ \ba{cc} \d_0^* & \zero_{12\times 2} \\
1   1   0   1   1   2   2   2   0   0   1   2&  1   0\\
0   1   1   0   0   2   0   0   2   1   2   2&  1   0\\
1   0   0   1   0   0   1   1   2   2   0   1&  1   0\\
1   0   1   0   0   1   2   2   0   1   0   0&  1   1\\
0   0   0   1   1   1   0   0   1   2   1   0&  1   1\\
0   1   1   0   1   0   1   1   1   0   2   1&  1   1

 \ea  $
  & $ \ba{cc}
\d_0^* & \zero_{12\times 3} \\
1   1   0   1   1   2   2   2   0   0   1   2&  2   1   0\\
1   0   0   1   0   0   1   1   2   2   0   1&  1   1   0\\
0   1   1   0   0   2   0   0   2   1   2   2&  1   1   0\\
1   0   1   0   0   1   2   2   0   1   0   0&  2   1   1\\
0   0   0   1   1   1   0   0   1   2   1   0&  1   1   1\\
0   1   1   0   1   0   1   1   1   0   2   1&  2   1   1

 \ea  $ & $ \ba{cc}
\d_0^* & \zero_{12\times 4}\\
0   0   0   1   1   1   0   0   1   2   1   0&  2   2   1   0\\
1   1   0   1   1   2   2   2   0   0   1   2&  1   1   1   0\\
0   1   1   0   0   2   0   0   2   1   2   2&  1   2   1   0\\
0   1   1   0   1   0   1   1   1   0   2   1&  2   1   1   1\\
1   0   1   0   0   1   2   2   0   1   0   0&  2   1   1   1\\
1   0   0   1   0   0   1   1   2   2   0   1&  1   2   1   1

 \ea  $  \\ \hline

 $D_{3B}^*  \in
\cC_{3B}(12+12;2^53^7\bullet 3^0\bullet \B)$ & $D_{3B}^*  \in
\cC_{3B}(12+12;2^53^7\bullet 3^1\bullet \B)$ & $D_{3B}^*  \in
\cC_{3B}(12+12;2^53^7\bullet 3^2\bullet \B)$

\\\hline
 $ \ba{cc} \d_0^* & \zero_{12\times 2} \\
1   0   1   0   0   1   2   2   2   2   0   0&  1   0\\
0   0   0   1   1   1   0   0   1   1   1   0&  1   0\\
1   0   0   1   0   0   1   1   0   1   0   1&  1   0\\
0   1   0   0   0   2   0   0   2   2   2   2&  1   0\\
1   1   1   1   1   2   2   2   0   0   1   2&  1   0\\
0   1   1   0   1   1   1   1   1   0   2   1&  1   0\\
1   1   0   1   0   1   2   2   1   1   2   2&  1   1\\
0   1   1   0   0   2   1   0   0   1   0   0&  1   1\\
1   0   1   0   1   0   0   2   1   0   2   0&  1   1\\
1   1   0   1   1   2   1   1   2   2   1   1&  1   1\\
0   0   0   0   1   0   2   1   0   2   0   2&  1   1\\
0   0   1   1   0   0   0   0   2   0   1   1&  1   1

 \ea  $
  & $ \ba{cc}
\d_0^* & \zero_{12\times 3} \\
1   1   0   1   1   1   2   2   1   2   2   2&  1   1   0\\
0   0   0   1   0   0   0   0   2   2   0   0&  2   1   0\\
0   0   0   1   1   0   1   1   0   0   1   1&  1   1   0\\
0   1   1   0   0   2   0   0   0   1   2   2&  1   1   0\\
1   0   1   0   0   0   2   2   2   0   0   1&  1   1   0\\
1   1   1   1   1   2   1   0   1   0   1   0&  2   1   0\\
0   1   0   0   1   2   0   0   2   2   1   1&  1   1   1\\
1   1   1   0   1   0   2   1   0   2   0   2&  2   1   1\\
1   0   0   1   0   2   2   2   0   1   1   2&  2   1   1\\
0   0   1   0   1   1   0   2   1   0   2   0&  2   1   1\\
1   0   0   0   0   1   1   1   1   1   0   0&  1   1   1\\
0   1   1   1   0   1   1   1   2   1   2   1&  2   1   1

 \ea  $ & $ \ba{cc}
\d_0^* & \zero_{12\times 4} \\
0   0   1   0   1   0   0   0   1   0   1   0&  2   2   1   0\\
0   1   1   0   0   2   1   1   0   1   0   2&  1   2   1   0\\
1   1   1   0   0   2   2   2   2   0   2   1&  2   1   1   0\\
1   0   0   1   1   1   1   1   0   1   1   0&  2   1   1   0\\
1   1   0   1   1   0   2   1   1   2   2   2&  2   2   1   0\\
0   0   0   1   0   0   0   0   2   2   0   1&  1   1   1   0\\
0   1   0   1   0   2   0   0   0   1   1   2&  2   1   1   1\\
0   0   0   0   0   1   2   2   2   2   0   0&  2   2   1   1\\
0   1   1   0   1   1   1   1   1   2   2   1&  1   1   1   1\\
1   0   1   1   0   1   0   2   1   1   2   0&  1   2   1   1\\
1   1   0   1   1   2   1   0   2   0   1   1&  1   2   1   1\\
1   0   1   0   1   0   2   2   0   0   0   2&  1   1   1   1
 \ea  $ \\ \hline

 \et} \ec

\vspace{10mm} \noindent\textbf{B2. }The
 \caugs~in Example 2.

\bc {\small  \tabcolsep=1pt \bt {c|c|c|c} \hline $D_3^*  \in
\cC_3(6+6;3^{10}\bullet 3^1)$ & $D_3^*  \in \cC_3(6+6;3^{10}\bullet
3^2)$ & $D_3^*  \in \cC_3(6+6;3^{10}\bullet 3^3)$ & $D_3^*  \in
\cC_3(6+6;3^{10}\bullet 3^4)$

\\\hline
 $ \ba{cc} \d_0^* & \zero_{6\times 1} \\
1   1   2   0   1   1   1   1   2   0 & 2\\
2   1   0   2   0   0   1   2   0   1 & 1\\
0   2   0   1   2   2   2   1   0   0 & 1\\
2   0   2   1   2   0   0   2   1   2 & 2\\
0   2   1   0   1   1   0   0   1   1 & 1\\
1   0   1   2   0   2   2   0   2   2 & 2

 \ea  $
  & $ \ba{cc}
\d_0^* & \zero_{6\times 2} \\
1   2   0   2   1   1   2   1   0   2 & 1   1\\
0   1   1   0   0   2   1   2   1   0 & 1   1\\
2   1   2   0   1   1   0   2   0   1 & 2   2\\
0   0   2   1   2   0   0   1   1   2 & 1   2\\
1   2   1   1   0   2   2   0   2   1 & 2   2\\
2   0   0   2   2   0   1   0   2   0 & 2   1\\

 \ea  $ & $ \ba{cc}
\d_0^* & \zero_{6\times 3} \\
1   1   2   2   1   1   0   0   0   1  & 1   2   1\\
0   2   0   0   0   2   2   1   1   1 & 1   2   2\\
2   1   0   0   2   0   1   2   2   0 & 1   1   1\\
2   2   1   1   1   0   0   2   1   2 & 2   2   2\\
0   0   2   1   2   1   1   1   0   0 & 2   1   2\\
1   0   1   2   0   2   2   0   2   2 & 2   1   1

 \ea  $
  &  $ \ba{cc}
\d_0^* & \zero_{6\times 4} \\
0   2   0   2   0   2   0   1   0   1 & 2   2   2   2\\
1   1   2   0   2   1   1   2   2   0 & 2   1   2   2\\
2   1   1   1   1   2   1   1   2   0 & 1   2   1   1\\
1   2   2   1   0   1   2   0   0   1 & 1   1   1   1\\
0   0   1   0   1   0   2   0   1   2 & 1   2   2   2\\
2   0   0   2   2   0   0   2   1   2 & 2   1   1   1

 \ea  $  \\ \hline

 \et} \ec

\bc {\small  \tabcolsep=1pt \bt {c|c|c} \hline $D_{3B}^*  \in
\cC_{3B}(6+6;3^{10}\bullet 3^0\bullet \B)$ & $D_{3B}^*  \in
\cC_{3B}(6+6;3^{10}\bullet 3^1\bullet \B)$ & $D_{3B}^*  \in
\cC_{3B}(6+6;3^{10}\bullet 3^2\bullet \B)$

\\\hline
 $ \ba{cc} \d_0^* & \zero_{6\times 2} \\
1   2   0   2   0   2   2   1   0   2&  1   0\\
1   1   1   0   1   0   1   0   2   1&  1   0\\
0   0   2   1   2   1   1   1   1   0&  1   0\\
0   0   1   2   2   2   0   2   2   2&  1   1\\
2   2   0   0   1   1   0   0   1   0&  1   1\\
2   1   2   1   0   0   2   2   0   1&  1   1

 \ea  $
  & $ \ba{cc}
\d_0^* & \zero_{6\times 3} \\
2   0   0   2   0   1   2   2   0   2&  1   1   0\\
1   2   1   1   1   1   1   1   2   2&  2   1   0\\
2   1   2   0   1   0   0   0   1   1&  1   1   0\\
0   0   2   0   0   2   2   1   2   1&  2   1   1\\
0   2   1   1   2   2   0   2   1   0&  1   1   1\\
1   1   0   2   2   0   1   0   0   0&  2   1   1

 \ea  $ & $ \ba{cc}
\d_0^* & \zero_{6\times 4}\\
0   0   1   1   2   0   2   0   0   1&  1   2   1   0\\
2   2   2   1   0   1   0   2   1   0&  1   1   1   0\\
1   1   1   2   1   2   1   1   2   2&  1   1   1   0\\
1   1   2   0   2   0   0   2   1   2&  2   2   1   1\\
0   2   0   0   0   2   2   1   2   1&  2   1   1   1\\
2   0   0   2   1   1   1   0   0   0&  2   2   1   1

 \ea  $  \\ \hline

 \et} \ec

\vspace{10mm} \noindent\textbf{B3.} The initial design and 
 \caugs~in Example 3.

\Bea
 \d_0^{*T}={\scriptsize \left[\ba{c}
1   0   0   1   1   0   1   0   1   1   1   1   1   0   0   0   0   0   1   1   0   0   0   1   0   1   0   1   1   0   0   1   0   0   0   1   1   1   1   0   0   1   0   1   0   1   1   1   1   0   0   1   0   1   0   0   1   1   0   1   1   0   0   1   0   0   1   0   1   0   0   1\\
0   0   0   1   1   0   0   1   0   0   1   0   1   1   0   0   1   1   0   0   1   0   1   1   1   0   1   0   0   0   1   0   0   0   0   1   1   1   1   1   0   0   0   1   1   1   1   1   0   1   1   0   0   0   1   0   0   1   0   1   1   0   1   1   1   1   0   0   0   0   1   1\\
1   1   1   0   1   0   1   1   0   1   0   1   1   0   0   1   0   0   1   1   0   0   0   1   1   0   1   1   1   1   1   0   1   0   0   1   1   0   0   1   1   0   0   0   1   0   0   1   1   0   1   0   0   1   1   0   0   0   1   1   0   0   0   1   1   1   0   0   0   0   1   0\\
1   1   0   0   1   0   0   1   0   0   1   1   0   1   1   0   0   0   1   1   1   1   1   0   0   0   0   1   1   1   0   0   1   0   0   0   0   1   1   0   0   1   1   1   1   1   1   0   0   1   0   1   1   1   1   0   0   1   0   1   1   1   0   0   1   0   0   0   0   1   0   0\\
2   0   0   2   0   0   0   1   0   2   1   0   2   1   2   2   2   1   1   0   0   1   0   2   0   0   1   2   2   0   1   2   2   2   1   1   1   1   0   2   1   1   2   1   0   0   2   0   1   2   1   1   2   1   1   1   0   2   2   0   1   2   0   2   2   0   0   2   1   0   1   0\\
1   0   2   2   0   2   2   2   0   1   2   2   1   2   1   2   2   2   2   0   0   0   1   1   0   2   2   2   0   2   1   2   1   2   0   0   0   0   0   2   1   1   2   1   1   1   2   2   0   0   0   1   1   0   0   1   0   2   1   1   2   1   0   1   0   0   0   1   1   0   1   1\\
0   1   2   2   2   2   1   0   0   0   2   0   1   0   2   0   0   1   1   1   1   1   0   2   0   0   1   0   2   2   2   0   1   0   1   2   2   0   1   2   2   0   2   2   2   0   1   1   0   2   1   1   2   0   2   0   1   0   1   1   1   0   2   2   0   0   2   2   1   1   1   1\\
0   1   0   0   0   1   0   2   2   1   2   2   1   0   2   0   2   2   2   0   2   0   1   0   0   0   0   2   1   1   2   2   0   1   0   2   1   0   1   2   1   2   1   0   2   1   0   1   1   1   1   1   0   1   0   1   2   1   2   1   0   0   1   0   2   0   2   1   2   2   2   2\\
2   0   1   1   1   2   0   2   2   1   1   1   0   0   0   2   1   0   2   1   1   0   0   1   1   0   2   0   0   2   1   2   1   0   1   2   0   0   2   0   1   0   1   1   0   1   2   0   0   2   2   1   2   1   0   0   1   1   1   2   2   2   0   2   1   2   0   2   2   2   2   0\\
2   2   0   1   0   2   1   2   2   0   2   0   0   1   2   0   0   2   1   2   2   0   1   2   0   2   2   2   0   1   0   1   1   0   0   0   2   1   2   0   1   1   2   0   1   2   0   2   2   1   1   0   2   1   1   2   1   1   1   0   0   1   1   1   2   1   0   0   0   0   2   1\\
2   0   0   0   2   1   1   2   1   1   1   1   1   2   1   2   2   2   0   0   1   0   1   1   0   1   0   0   0   2   0   0   1   0   2   0   2   1   2   1   2   2   2   0   1   2   1   0   1   1   1   0   0   2   1   0   2   0   1   2   1   0   2   2   2   0   2   2   1   2   0   0\\
2   2   2   1   0   1   0   1   1   2   0   1   1   2   2   2   1   0   1   1   2   1   0   0   2   1   1   2   1   0   0   2   2   0   0   2   2   0   1   0   1   1   2   2   1   0   2   0   0   2   1   1   0   0   1   0   1   0   2   0   2   0   2   1   1   2   0   1   2   0   0   2\\
2   0   0   1   1   2   0   1   0   0   1   1   2   0   0   2   0   0   0   1   2   0   2   0   2   0   2   2   1   1   0   1   1   1   2   0   0   1   0   2   1   2   1   1   2   0   1   1   2   0   1   0   2   2   0   1   2   2   0   0   2   0   2   2   2   1   1   2   1   1   1   2\\
0   0   2   1   2   1   0   2   2   1   0   1   2   1   2   0   0   2   1   1   1   2   2   1   1   2   0   1   0   0   0   2   0   2   1   0   2   2   1   0   2   0   1   1   2   0   0   1   1   0   0   2   1   2   1   0   0   2   2   2   2   0   1   2   0   1   0   1   2   1   1   0\\
1   1   1   1   0   0   2   0   1   0   2   2   1   2   2   0   0   0   2   0   0   0   1   2   2   2   2   1   1   1   0   2   2   2   1   0   1   0   2   2   2   0   0   0   1   1   2   0   0   0   1   0   2   2   1   2   1   1   0   2   1   1   1   0   2   0   0   2   1   1   1   2\\
0   2   0   0   1   1   1   0   2   0   0   0   0   1   0   1   1   2   2   2   2   0   1   2   2   2   1   1   1   2   0   0   1   2   1   1   1   0   1   2   2   0   1   2   1   2   1   0   0   2   0   1   2   1   0   0   1   2   2   0   2   0   0   2   0   2   2   2   1   0   1   1\\
2   2   2   0   1   1   2   0   0   0   1   1   2   1   2   0   2   2   1   0   2   1   1   0   1   0   2   1   0   0   1   2   2   1   0   2   0   2   2   0   1   1   0   2   0   1   1   0   1   1   1   2   1   2   2   0   1   0   1   2   0   0   2   1   0   0   2   2   0   2   1   0\\
2   1   1   2   1   0   2   0   1   0   2   2   2   1   2   0   2   1   2   0   2   2   2   1   1   0   0   2   1   2   2   1   1   0   1   0   2   0   1   1   0   1   1   1   1   0   0   2   1   2   1   1   0   2   0   2   2   0   0   1   2   1   0   0   0   2   0   2   0   0   1   0\\
0   2   0   1   1   1   1   2   2   1   0   0   2   1   0   0   0   1   1   0   2   0   2   0   1   1   0   2   1   2   2   2   1   2   2   2   1   2   0   2   2   0   0   2   1   0   2   2   1   1   0   2   1   0   2   1   1   1   1   0   1   2   0   0   2   0   0   2   0   1   0   1\\
2   2   0   1   1   2   2   1   1   1   0   1   0   0   2   0   2   2   2   0   0   2   1   2   1   1   1   0   2   1   0   0   1   1   2   2   1   2   1   1   2   1   0   0   0   2   2   2   0   0   0   1   1   2   1   1   0   0   0   0   1   2   2   2   0   1   0   0   0   1   2   2\\
0   0   2   1   2   0   2   0   1   0   1   0   1   0   2   1   2   1   2   2   1   0   1   0   1   0   0   1   1   0   0   1   1   1   2   2   2   2   1   2   2   2   2   0   0   1   0   0   0   1   0   2   2   1   1   1   0   0   2   0   2   2   2   1   2   2   1   0   1   0   1   2\\
0   1   0   2   0   2   1   0   1   1   0   0   2   1   2   2   1   0   1   1   0   2   1   0   2   0   1   0   0   2   1   1   2   0   0   2   1   2   1   1   1   1   1   0   2   1   2   2   2   1   0   2   2   2   1   0   2   1   0   0   2   0   0   0   2   1   2   0   1   2   2   0\\
0   0   1   1   2   0   1   2   1   0   1   0   2   1   2   0   2   0   2   1   0   1   1   0   2   2   2   2   1   1   0   1   1   1   2   1   1   1   2   2   0   1   1   2   0   2   0   2   2   2   2   0   0   0   0   0   0   2   2   1   0   0   0   1   1   0   2   1   2   2   1   0\\
2   0   1   2   1   0   2   0   0   1   1   0   0   1   1   1   0   2   0   1   1   2   2   0   0   2   2   2   2   1   2   2   0   1   1   2   0   0   1   2   1   2   0   2   0   0   1   2   1   2   1   0   2   0   1   1   2   0   2   2   1   0   2   1   1   0   0   0   1   2   0   1\\
0   0   2   2   1   0   2   1   1   1   0   2   2   1   2   2   1   2   2   1   1   1   2   0   0   0   0   1   0   0   0   0   1   0   2   2   1   0   2   0   2   0   0   2   2   2   2   1   2   1   1   0   1   1   0   2   0   2   0   1   0   1   1   0   2   2   2   1   1   0   1   1\\
0   1   1   1   2   0   2   2   1   0   1   2   2   0   2   1   0   1   0   0   0   1   0   2   2   0   2   1   0   2   0   0   1   2   0   0   2   0   2   1   1   2   1   2   0   1   1   1   1   2   0   0   0   1   1   0   2   2   1   1   1   2   2   0   2   1   0   1   2   0   2   2\\
1   0   2   2   0   2   1   2   1   2   1   0   0   2   0   0   1   2   2   1   1   1   0   1   2   0   2   1   0   1   0   1   2   2   0   0   2   2   1   1   1   0   0   2   1   0   1   0   1   0   0   1   0   2   1   1   2   0   1   2   2   2   1   0   2   0   2   2   2   1   0   0\\
0   2   0   2   1   2   0   0   0   1   0   2   1   0   1   0   2   1   1   0   1   2   0   0   0   1   1   2   1   1   2   1   2   2   1   1   2   2   2   2   0   0   2   2   1   1   1   0   2   0   2   0   2   2   1   1   1   0   0   0   0   1   2   2   0   1   0   1   2   0   1   2\\
1   1   1   1   2   1   1   1   2   1   1   0   1   0   2   0   2   1   0   1   0   0   2   2   1   0   2   2   2   0   2   1   2   1   1   1   2   0   1   0   2   1   0   2   1   0   2   0   2   1   1   0   1   1   0   2   0   2   0   2   0   2   0   0   0   2   2   0   0   2   0   2\\
1   1   2   1   0   1   0   0   0   2   0   2   0   0   1   0   2   1   1   1   1   2   1   2   1   1   0   2   0   2   1   0   2   2   2   2   1   0   2   1   0   2   1   2   2   0   1   0   0   0   1   0   0   0   2   2   0   2   2   1   2   1   0   1   1   0   1   0   1   2   2   2\\
0   0   1   0   2   0   2   0   1   0   1   1   0   0   2   2   2   1   1   2   0   0   0   0   1   1   0   2   1   0   2   0   2   2   0   1   0   2   2   2   0   2   1   0   1   0   1   1   2   1   2   1   2   1   2   1   2   2   0   2   2   1   1   1   1   2   0   2   1   1   0   0\\
2   0   2   2   2   0   0   1   2   2   2   1   2   2   1   1   0   1   1   1   2   0   1   1   0   0   2   0   1   0   1   0   1   1   0   2   2   1   1   2   2   2   2   0   0   0   1   2   0   0   2   0   2   2   1   0   1   1   1   0   2   2   1   0   0   0   1   1   0   2   1   0\\
0   1   1   0   2   2   0   1   0   1   2   2   0   2   1   2   1   2   0   2   1   1   2   2   1   2   0   2   1   1   0   1   0   0   0   2   0   0   2   1   2   1   0   0   0   0   1   2   1   2   0   0   0   1   2   0   2   1   0   1   1   2   0   1   0   2   1   2   2   2   1   1\\
1   0   1   2   1   0   2   2   1   1   1   0   1   2   2   0   1   1   1   2   0   1   1   0   2   2   1   2   2   1   2   0   1   1   2   0   0   0   0   0   0   0   2   1   0   2   1   0   1   2   0   2   0   2   2   0   1   0   2   0   2   0   1   1   1   1   0   2   0   2   2   2\\
0   2   1   1   2   0   0   0   0   1   2   2   1   2   0   1   1   2   1   1   2   0   0   1   0   1   2   1   0   2   2   2   2   1   2   0   2   1   0   1   1   0   0   0   1   1   0   0   2   1   0   2   2   0   2   0   1   2   0   1   2   1   1   2   0   0   2   0   2   1   1   2\\
1   0   1   2   0   0   0   2   0   0   0   2   0   1   1   0   2   0   0   0   2   0   2   1   2   1   2   1   2   1   1   0   1   1   1   2   1   2   2   0   2   1   2   0   1   1   1   1   2   1   1   1   2   0   0   2   2   0   2   2   2   2   1   2   0   1   1   0   2   0   0   0\\
1   2   1   2   0   2   0   0   2   2   2   0   2   2   2   0   0   0   0   1   0   0   1   0   0   1   1   2   0   2   2   0   0   1   2   2   0   2   1   1   0   0   0   1   1   0   1   1   1   1   0   1   2   1   1   0   2   2   2   2   0   1   2   1   2   2   1   1   2   0   1   1\\
0   1   2   1   0   0   2   1   2   0   1   2   1   1   1   1   1   0   1   0   2   2   0   0   0   2   2   0   1   0   2   2   0   1   2   0   2   0   1   0   2   0   2   1   1   0   2   1   1   2   0   1   2   2   2   0   0   2   1   2   0   1   2   1   0   1   2   0   0   1   2   1\\
0   1   2   0   2   2   1   0   0   1   1   0   2   0   0   2   1   0   2   0   0   0   1   2   1   1   2   2   1   1   2   2   0   0   0   0   1   1   2   0   1   2   2   1   0   0   1   2   0   0   1   2   1   1   2   2   0   2   1   0   2   2   2   0   2   1   2   1   1   1   1   0\\
0   1   0   0   0   0   0   1   2   2   2   1   1   0   1   2   1   2   2   1   0   0   1   1   2   1   2   2   0   2   0   1   1   0   1   1   1   2   0   0   1   2   2   1   2   0   2   1   0   2   2   0   1   1   2   2   1   0   2   2   1   0   2   1   1   0   2   0   0   0   0   2\\
1   1   0   2   1   2   0   0   1   0   0   1   1   1   2   1   2   0   2   2   0   1   1   2   2   1   2   0   0   1   1   2   0   2   2   1   0   1   1   1   2   2   1   0   2   0   2   2   2   2   0   0   0   1   0   1   2   0   1   2   2   1   2   0   0   0   1   0   0   0   2   1\\
2   1   2   1   1   2   0   0   2   0   0   1   0   1   1   2   2   2   0   1   2   0   0   0   0   2   1   0   2   2   2   2   2   0   0   1   2   2   2   1   2   1   0   2   1   1   1   1   2   2   0   0   0   1   0   1   0   1   1   0   2   0   1   1   2   0   0   1   1   1   2   0\\
1   0   2   0   0   1   0   2   2   0   1   0   1   1   2   2   1   0   2   1   1   0   2   2   1   1   0   0   0   2   0   2   1   1   0   1   1   0   0   2   0   1   2   2   0   2   1   2   2   0   2   1   2   2   2   1   2   1   1   1   0   0   1   0   0   1   0   0   2   2   1   2\\
1   2   1   2   2   2   0   0   1   0   0   1   1   1   1   2   1   1   2   0   1   0   2   2   1   0   1   1   1   1   0   0   0   2   2   0   2   1   0   0   1   2   0   2   0   0   2   1   0   2   1   1   1   2   2   2   0   1   2   2   0   0   0   2   2   0   0   1   2   2   0   1\\
1   1   1   0   1   0   0   0   0   1   1   0   2   0   0   0   1   2   1   0   2   2   0   2   1   2   2   0   0   2   1   1   2   2   2   2   1   1   1   0   2   2   2   0   1   1   1   1   2   2   0   0   1   2   1   0   1   2   1   0   0   2   0   0   2   2   0   1   2   2   0   1\\
2   2   1   2   2   0   0   1   0   0   0   0   1   1   0   0   1   1   2   1   0   1   2   1   0   2   0   2   1   1   0   1   0   2   0   2   1   0   1   1   2   0   2   0   0   0   1   0   2   1   2   2   1   1   2   1   1   2   0   2   2   0   1   0   2   1   2   2   1   0   2   2\\
1   1   2   0   1   2   1   2   1   1   0   0   2   1   0   0   1   0   2   0   2   0   0   0   0   2   0   2   2   0   2   0   2   1   2   0   0   2   2   2   2   2   1   1   1   2   2   1   0   1   0   0   1   2   1   1   1   0   0   2   1   0   0   1   1   2   1   0   1   2   1   2\\
0   2   2   2   2   1   1   0   0   2   0   1   1   0   2   0   1   1   0   1   2   0   0   2   0   2   2   1   0   1   0   2   0   0   0   0   1   2   0   2   0   1   1   1   0   1   1   1   1   2   2   2   1   2   1   2   2   0   1   1   0   1   0   1   2   2   1   0   0   2   2   2\\
0   1   1   1   1   2   1   1   2   0   0   2   1   2   0   1   1
0   1   2   1   2   1   1   0   1   2   0   2   1   1   2   0   2
0   2   0   0   0   0   0   1   0   2   1   2   1   0   1   0   2
0   1   2   2   0   0   0   2   2   2   2   1   0   2   2   2   1
1   0   0   2

 \ea\right]}. \label{d72}
 \Eea

\bc {\scriptsize  \tabcolsep=1pt \bt {c|c} \hline $D_{3}^* \in
\cC_{3}(72+24;2^{4}3^{45}\bullet 3^{1})$ & $D_{3}^* \in
\cC_{3}(72+24;2^{4}3^{45}\bullet 3^{2})$

\\\hline
$ \ba{cc}
\d_0^*  &\zero_{72\times 1}   \\
1   0   0   0   2   0   2   2   2   1   0   0   2   0   0   1   1   1   0   1   0   1   0   1   1   1   2   0   0   2   1   1   0   2   0   0   2   1   0   1   2   1   0   1   2   2   2   1   1   &   1   \\
0   0   1   1   1   0   0   0   2   2   2   1   1   2   1   2   0   2   1   2   1   0   1   2   2   1   2   2   0   2   2   2   1   0   2   0   1   1   2   0   2   0   1   1   1   0   0   2   1   &   2   \\
1   1   1   0   1   2   2   1   2   0   1   2   2   2   2   0   1   1   0   1   2   2   1   2   1   0   1   1   0   1   0   1   2   2   1   1   0   1   2   0   1   0   0   1   2   0   1   1   0   &   2   \\
1   0   0   1   0   2   2   2   0   1   1   2   2   1   2   0   0   0   2   2   0   0   2   0   2   1   0   0   1   1   2   1   0   1   2   2   0   0   1   0   1   0   1   0   2   1   1   2   2   &   2   \\
1   0   1   0   1   2   0   0   1   2   1   0   1   2   1   2   2   1   1   0   2   2   0   2   1   2   0   0   0   1   0   0   0   0   0   2   1   1   0   2   0   1   0   0   2   1   2   0   0   &   1   \\
1   1   0   0   0   0   2   0   0   0   2   1   1   2   2   1   1   2   2   2   0   1   0   1   1   1   0   0   1   2   0   0   0   2   1   2   2   1   2   2   0   2   2   0   2   0   0   2   0   &   2   \\
1   0   1   1   1   1   2   0   2   0   2   0   0   0   0   0   0   2   1   0   1   1   2   0   0   2   1   2   1   2   0   1   2   1   0   2   2   2   1   1   2   1   2   1   0   2   2   2   2   &   2   \\
0   0   0   1   2   0   1   2   1   0   0   2   0   2   2   2   0   0   0   1   0   0   1   0   2   2   1   1   1   0   0   2   2   2   2   2   1   0   1   2   2   2   2   2   0   2   2   1   0   &   1   \\
0   0   1   1   0   2   0   1   0   1   1   1   0   0   0   2   2   1   0   0   1   2   0   2   1   2   0   1   2   2   1   2   0   2   2   2   2   0   1   0   2   0   1   2   1   1   0   0   0   &   2   \\
0   1   1   1   0   1   0   2   1   1   2   2   1   2   0   1   2   2   1   1   0   0   2   1   0   1   0   1   0   1   2   0   0   0   1   0   1   2   1   1   2   0   0   1   0   1   0   0   2   &   1   \\
1   1   1   1   2   2   1   1   1   2   2   2   0   0   1   0   1   2   2   1   2   0   0   2   2   0   2   2   2   0   1   1   1   0   2   0   1   0   1   0   2   1   2   0   2   0   1   0   2   &   1   \\
1   0   0   0   0   1   2   1   1   2   2   2   1   0   2   0   1   1   2   2   2   2   2   0   0   0   2   1   0   2   2   2   2   1   2   0   1   2   0   2   0   2   0   2   0   1   0   0   1   &   2   \\
0   1   1   0   1   2   1   1   0   0   2   2   2   1   1   2   1   0   0   0   1   1   1   0   0   0   1   0   2   2   1   0   0   0   0   0   1   0   0   1   0   0   0   0   1   1   2   2   1   &   2   \\
1   1   1   1   2   1   0   1   0   0   0   1   1   1   0   1   2   0   1   2   1   0   2   0   2   2   1   0   0   1   2   2   2   2   0   1   2   2   0   2   1   2   2   0   2   0   0   2   1   &   1   \\
0   0   0   0   1   0   1   0   1   2   1   0   2   0   0   1   1   0   2   1   1   0   1   0   2   2   1   1   2   1   2   2   1   2   1   1   2   1   2   1   0   1   2   1   1   1   1   0   2   &   2   \\
0   1   1   1   2   0   1   2   2   0   1   1   0   1   2   2   0   1   1   2   1   0   0   1   0   0   2   2   2   0   0   0   0   1   1   1   2   2   2   0   1   1   2   0   0   2   0   1   0   &   1   \\
0   0   0   1   0   2   1   0   1   1   0   0   2   2   0   0   0   2   0   2   1   1   2   1   0   0   2   2   2   0   0   1   2   0   0   1   0   1   0   2   0   2   1   0   1   2   1   1   1   &   1   \\
0   0   1   1   0   1   2   2   1   0   1   0   2   1   1   1   1   0   2   0   2   1   2   2   1   1   2   1   1   0   1   0   1   0   0   1   0   1   2   1   1   2   1   2   0   2   0   2   0   &   1   \\
1   1   0   1   1   0   0   1   2   1   0   0   0   1   2   1   0   0   0   0   2   2   0   2   0   1   0   2   2   1   0   2   1   1   1   0   0   0   0   1   0   2   1   1   0   0   1   0   1   &   2   \\
1   0   1   0   2   1   1   2   0   1   2   0   0   1   1   2   2   1   2   1   0   2   1   1   2   0   0   2   2   1   2   0   2   1   0   2   0   2   2   2   2   1   1   1   1   1   1   2   2   &   1   \\
0   1   0   0   0   1   0   0   2   2   0   1   0   1   1   1   2   2   2   0   2   2   1   1   0   2   1   2   1   0   1   1   1   2   2   0   0   0   0   2   1   1   0   2   2   2   2   1   0   &   2   \\
0   1   0   0   1   1   0   1   0   1   1   2   1   0   2   2   0   2   0   1   0   2   2   2   1   0   1   0   1   0   1   1   1   1   1   2   2   2   2   0   0   2   0   2   1   2   2   1   2   &   1   \\
0   1   0   0   2   2   2   2   2   2   0   1   1   2   1   0   2   1   1   2   0   1   0   1   2   2   0   0   0   0   2   2   1   1   1   1   1   0   1   1   1   0   1   2   1   0   2   1   2   &   2   \\
1   1   0   0   2   0   1   0   0   2   0   1   2   0   0   0   2   0   1   0   2   1   1   0   1   1   2   1   1   2   1   0   2   0   2   1   0   2   1   0   1   0   2   2   0   0   1   0   1   &   1   \\

 \ea  $ &
$ \ba{cc}
\d_0^*  &\zero_{72\times 2}   \\
0   0   0   1   2   1   1   0   1   2   0   0   2   2   1   1   1   2   2   2   0   1   2   1   0   2   2   1   1   2   2   2   2   0   2   1   0   1   1   2   1   2   1   0   0   1   0   2   2   &   2   2   \\
1   1   1   1   2   0   1   2   0   0   2   0   1   1   1   2   1   0   1   0   2   1   0   2   0   1   2   1   1   2   0   0   2   1   0   2   1   2   2   1   0   2   1   1   0   0   0   2   1   &   1   1   \\
1   0   0   0   0   0   0   1   1   1   1   2   0   1   2   0   1   2   0   2   2   2   2   2   0   1   0   1   2   2   0   2   0   0   2   0   1   0   0   2   0   2   0   2   0   1   1   0   1   &   1   1   \\
1   1   0   0   2   1   2   2   0   1   2   2   0   1   1   1   1   2   2   1   0   0   2   1   2   1   2   2   1   0   1   1   1   0   0   0   0   2   2   0   2   0   1   0   2   0   1   2   2   &   2   1   \\
0   1   0   0   1   2   1   1   2   1   0   1   0   1   0   0   2   0   1   0   1   0   1   0   0   2   1   1   2   1   2   2   1   1   0   1   1   0   1   1   1   1   0   0   0   0   1   2   1   &   2   2   \\
1   1   0   0   2   1   0   0   0   2   0   1   1   0   1   2   2   1   0   0   2   2   0   1   0   0   2   2   2   0   2   0   2   2   0   0   0   2   0   2   0   1   2   2   2   2   2   0   0   &   2   2   \\
1   0   1   1   1   2   2   2   1   1   2   2   2   2   1   2   2   0   0   1   1   0   0   0   1   0   0   0   2   1   2   0   0   1   0   0   2   1   0   0   2   2   0   0   2   1   0   1   1   &   2   1   \\
1   1   1   0   0   2   1   0   0   2   2   1   2   0   2   0   1   1   2   0   2   1   0   0   1   0   1   1   1   2   1   1   1   2   2   2   2   0   1   0   1   2   2   2   2   0   0   0   0   &   2   1   \\
0   0   1   1   1   0   1   0   2   2   2   2   1   2   2   1   0   2   2   0   2   1   2   2   0   0   0   2   2   0   1   1   1   0   1   0   1   0   0   0   2   0   1   1   1   0   1   1   1   &   2   2   \\
0   1   0   1   1   0   2   1   2   1   1   1   1   1   2   1   0   0   2   2   1   0   0   2   0   1   2   0   1   2   2   0   1   1   1   1   0   2   0   0   1   1   2   1   1   2   0   1   0   &   2   1   \\
1   0   1   1   2   2   2   1   2   1   0   0   2   2   0   0   0   1   0   1   2   2   2   1   2   0   1   2   2   1   0   0   2   2   1   1   0   2   1   1   1   1   1   1   1   0   1   0   2   &   1   1   \\
0   0   0   0   1   1   1   0   0   1   1   0   1   0   0   1   2   2   0   1   0   0   2   0   2   1   0   0   0   1   1   2   1   2   1   2   2   2   2   1   0   2   0   2   1   1   0   0   2   &   1   2   \\
0   0   1   0   1   2   0   0   1   0   1   0   2   1   1   2   2   0   0   0   2   1   1   2   1   2   1   1   0   0   1   1   1   0   0   1   0   1   0   1   1   0   0   2   1   2   2   1   0   &   1   1   \\
0   0   1   1   2   2   1   1   0   0   2   1   0   1   0   2   0   2   0   2   1   2   0   2   2   2   2   0   2   2   1   0   0   0   2   2   2   0   1   0   2   1   2   0   1   2   2   2   2   &   1   2   \\
1   0   1   1   0   0   2   2   2   2   1   0   0   0   0   0   1   0   2   0   1   0   0   0   2   2   2   2   0   2   2   2   0   2   2   0   1   1   1   1   2   1   1   1   2   1   2   0   0   &   1   2   \\
1   0   0   1   0   1   2   0   0   0   0   0   2   0   2   0   0   1   1   2   1   1   2   0   2   1   1   0   0   0   0   1   2   1   0   2   2   1   0   2   0   2   1   1   2   2   1   2   1   &   2   2   \\
1   1   1   0   0   1   0   1   2   0   1   2   0   1   1   1   0   0   2   2   0   2   1   2   1   0   1   0   1   1   2   1   0   2   2   2   0   1   2   2   1   1   0   1   2   1   1   2   2   &   1   2   \\
0   1   0   0   0   2   2   0   1   2   2   2   1   2   0   0   1   1   1   1   1   2   1   2   1   2   1   1   0   1   2   0   0   0   1   0   1   1   2   1   0   0   2   1   2   2   2   0   2   &   2   1   \\
0   0   1   1   0   2   0   2   2   0   1   1   1   2   1   2   2   1   1   2   0   1   1   1   2   2   0   2   2   0   0   2   0   2   1   2   1   0   2   2   0   0   1   2   1   1   0   2   0   &   2   1   \\
0   0   0   0   1   0   1   2   1   1   2   2   0   0   2   2   0   2   2   1   0   2   1   0   1   0   0   2   1   0   0   0   2   1   1   2   1   1   2   2   0   1   2   0   0   2   2   1   0   &   1   2   \\
0   1   1   1   0   1   0   2   1   2   0   2   0   0   0   1   0   1   1   2   1   2   2   1   0   0   2   2   0   0   1   2   1   1   0   1   2   2   1   1   2   0   0   2   0   2   0   1   1   &   1   1   \\
0   1   0   0   2   0   0   1   0   0   0   1   1   2   2   1   1   0   0   1   0   0   1   0   2   2   1   0   0   1   0   2   2   2   2   0   2   0   0   2   2   2   2   0   1   0   2   1   1   &   1   1   \\
1   1   0   0   2   0   2   2   2   2   0   0   2   2   2   0   2   1   1   0   0   1   0   1   1   1   0   0   0   2   1   1   0   1   1   1   0   0   1   0   1   0   0   2   0   0   2   1   2   &   1   2   \\
1   1   1   1   1   1   0   1   1   0   1   1   2   0   0   2   2   2   1   1   2   0   1   1   1   1   0   1   1   1   0   1   2   0   2   1   2   2   2   0   2   0   2   0   0   1   1   0   0   &   2   2   \\

 \ea  $
 \\ \hline

$D_{3}^* \in \cC_{3}(72+24;2^{4}3^{45}\bullet 3^{3})$ &$D_{3}^* \in
\cC_{3}(72+24;2^{4}3^{45}\bullet 3^{4})$
\\\hline
$ \ba{cc}
\d_0^*  &\zero_{72\times 3}   \\

0   0   1   1   0   1   0   1   0   0   0   1   1   1   0   1   0   1   0   2   1   2   2   2   2   2   1   0   0   1   2   2   0   1   2   2   2   2   0   2   0   2   1   2   1   1   2   2   2   &   2   1   2   \\
0   1   1   0   0   2   0   0   2   0   2   1   0   2   2   0   0   1   1   0   1   2   0   2   1   0   1   2   2   0   0   1   1   2   0   1   2   0   1   1   0   1   0   2   2   2   1   1   0   &   1   2   2   \\
1   0   1   1   2   2   1   2   1   1   1   0   0   1   1   0   1   1   0   1   2   2   2   1   0   0   2   1   2   0   0   0   2   0   0   0   1   2   0   2   0   2   1   1   2   1   1   0   1   &   1   2   1   \\
1   0   0   0   1   2   2   0   0   2   2   2   2   0   0   0   1   1   2   0   1   2   2   0   0   2   0   0   2   2   1   1   2   1   1   2   0   2   1   1   0   1   0   0   1   2   2   2   2   &   2   2   1   \\
1   1   1   0   0   1   0   2   0   2   2   2   1   0   1   2   1   0   2   0   2   2   0   0   0   0   2   2   0   2   2   2   0   2   0   0   1   2   0   0   2   0   0   2   2   1   0   0   1   &   2   1   2   \\
0   0   0   0   1   2   2   0   0   2   1   1   1   1   1   2   2   2   0   0   2   0   0   2   2   2   0   2   2   0   1   2   1   2   1   0   1   0   0   1   1   0   1   2   1   0   1   0   2   &   1   1   1   \\
0   0   0   0   1   0   0   1   1   1   1   2   0   0   2   1   0   2   0   1   0   2   1   1   0   1   0   0   2   1   1   1   0   0   2   0   1   0   0   0   2   0   0   1   1   1   1   1   1   &   2   2   2   \\
1   1   1   1   2   2   1   0   1   2   2   1   2   2   0   0   2   0   1   1   1   1   0   1   0   1   0   1   2   1   1   0   0   0   1   1   2   2   1   0   1   0   2   1   0   1   0   0   2   &   2   1   2   \\
1   0   1   1   0   0   2   0   2   0   1   0   0   0   2   2   0   0   2   0   1   1   0   2   1   1   0   2   1   2   0   0   0   1   0   2   0   1   2   0   2   0   1   1   0   0   0   2   0   &   2   2   1   \\
1   0   1   0   1   0   0   2   2   1   1   0   2   2   0   2   2   1   2   1   0   0   1   1   2   2   0   2   0   1   0   2   2   2   1   1   0   1   2   2   0   1   1   0   1   1   1   0   0   &   2   2   2   \\
0   1   0   1   0   1   2   2   2   1   0   1   0   1   1   1   0   2   2   2   0   0   2   1   0   1   2   0   1   0   2   0   0   0   1   0   0   2   1   1   1   1   1   0   2   2   0   1   0   &   2   2   2   \\
0   1   1   1   2   1   1   1   0   0   2   2   2   0   0   2   0   2   1   1   0   0   1   0   2   0   1   2   1   0   1   1   2   2   1   2   2   2   2   0   2   0   2   2   1   0   0   2   2   &   1   2   1   \\
1   1   0   1   0   2   0   1   2   0   1   1   2   0   0   0   1   1   0   2   2   0   1   2   1   1   2   1   2   2   0   1   0   2   2   0   2   1   1   0   2   1   2   0   2   1   1   2   0   &   1   1   1   \\
0   0   1   0   0   1   2   0   1   2   0   0   2   2   1   1   2   2   1   1   0   1   2   1   1   0   1   2   0   0   1   1   1   0   0   0   2   1   2   2   2   2   0   0   2   2   2   1   2   &   2   1   1   \\
1   1   0   1   1   1   0   1   2   0   0   0   2   1   2   1   2   0   1   0   2   1   2   0   0   1   1   1   1   1   1   1   1   1   1   1   0   0   0   1   1   2   0   2   0   0   1   2   1   &   2   1   1   \\
0   1   0   0   1   1   0   2   1   1   1   2   1   0   1   2   2   0   0   2   0   2   1   2   1   2   1   1   1   0   2   0   2   1   2   2   0   1   2   0   1   0   0   2   0   2   2   1   0   &   1   1   2   \\
0   1   0   0   2   0   1   1   2   2   1   0   1   1   2   1   0   1   2   1   2   1   0   2   2   0   2   0   0   2   2   0   1   1   1   1   0   1   0   2   0   1   0   1   1   2   0   1   2   &   1   1   1   \\
1   1   1   1   2   0   2   1   0   2   0   1   1   2   2   2   2   2   0   0   2   2   2   1   1   0   2   0   1   2   0   0   2   1   2   1   0   2   1   2   1   1   2   0   0   0   2   0   1   &   1   2   2   \\
1   1   0   0   2   0   2   2   2   2   0   2   1   2   1   0   1   1   1   2   0   1   1   1   2   1   1   0   0   1   2   1   0   2   2   2   1   0   2   1   1   0   1   1   2   0   2   2   1   &   1   2   1   \\
1   0   1   1   1   2   2   2   0   1   0   2   0   1   2   1   2   0   2   1   1   0   0   0   2   0   0   0   0   1   2   0   1   1   0   2   2   0   1   0   2   2   2   0   2   0   1   1   0   &   1   1   2   \\
1   0   0   0   0   1   1   0   1   1   2   0   0   1   2   0   1   2   2   2   0   1   2   0   1   1   0   1   1   2   1   2   1   2   2   2   1   0   2   2   0   2   1   2   0   1   0   0   2   &   1   1   2   \\
0   0   0   1   2   0   1   2   1   0   0   1   0   0   0   0   0   0   1   2   1   0   1   0   2   2   2   2   0   0   0   2   2   0   2   1   1   0   1   2   2   1   2   1   0   2   2   1   1   &   2   1   1   \\
0   1   0   0   2   0   1   0   0   1   2   0   1   2   0   1   1   2   0   0   2   0   0   0   0   2   1   1   1   1   0   2   1   0   0   0   2   1   2   1   0   2   2   1   0   2   2   0   0   &   2   2   1   \\
0   0   1   1   1   2   1   1   1   0   2   2   2   2   1   2   1   0   1   2   1   1   1   2   1   2   2   1   2   2   2   2   2   0   0   1   1   1   0   1   1   2   2   0   1   0   0   2   1   &   1   2   2   \\

 \ea  $

 &
$ \ba{cc}
\d_0^*  &\zero_{72\times 4}   \\

0   1   0   0   2   0   0   0   2   2   2   0   1   2   0   1   1   2   1   0   0   1   0   2   2   2   1   1   0   1   1   2   1   2   1   0   1   1   2   1   0   0   0   1   1   2   2   0   2   &   1   1   2   2   \\
1   0   1   1   1   2   2   2   1   0   2   0   0   0   1   2   2   0   1   1   0   1   1   1   1   1   0   1   0   1   2   0   0   0   0   2   1   1   2   1   2   0   1   1   2   0   0   2   1   &   1   2   2   2   \\
0   0   1   1   2   2   1   1   1   1   1   0   1   1   1   1   1   0   0   1   2   0   2   2   2   0   2   0   2   0   1   0   1   0   0   0   2   1   0   2   1   2   1   0   1   2   1   2   1   &   2   1   2   2   \\
1   0   1   1   2   0   2   2   2   2   1   0   2   1   0   1   2   1   2   1   1   0   0   1   2   2   0   0   0   1   2   0   0   1   1   1   0   2   1   0   2   1   2   0   1   1   0   0   2   &   2   2   1   1   \\
1   1   1   0   1   2   0   1   2   1   0   1   2   1   2   1   2   0   2   0   2   0   1   2   2   1   1   1   1   1   1   1   1   2   2   1   0   0   1   0   1   0   2   2   2   0   1   1   0   &   2   2   2   2   \\
0   1   1   1   0   1   0   2   0   1   0   2   0   0   0   1   0   2   0   1   0   0   2   1   0   0   2   2   1   0   1   1   0   0   1   0   2   2   1   0   2   2   2   1   2   2   0   1   0   &   1   1   1   2   \\
1   1   0   0   1   0   2   1   0   1   0   2   1   0   2   1   0   1   0   1   1   2   0   0   0   0   1   0   0   1   2   0   2   1   0   0   2   2   0   1   2   0   0   1   2   0   2   1   1   &   2   2   2   1   \\
1   0   1   1   0   1   1   0   0   1   2   0   2   0   2   0   1   0   2   0   2   1   2   0   0   1   0   1   1   2   1   2   2   1   0   2   2   2   0   0   0   2   0   2   0   1   0   0   2   &   1   1   2   1   \\
1   0   0   1   2   0   2   2   2   1   0   0   0   2   2   0   0   0   0   0   0   1   0   2   2   1   2   2   0   2   0   2   2   2   2   1   0   0   1   2   2   2   1   1   0   0   1   2   1   &   1   1   1   1   \\
1   1   0   0   0   0   1   2   1   2   1   2   0   0   1   0   0   2   2   1   0   0   1   0   2   0   0   2   1   0   2   0   1   2   2   2   1   0   2   2   0   1   1   2   2   1   1   0   2   &   2   1   2   1   \\
1   0   0   0   0   2   2   0   1   2   0   2   2   2   1   0   1   1   2   2   1   1   0   2   1   0   2   1   0   0   2   1   1   1   0   1   0   1   1   2   1   1   0   2   2   2   0   1   0   &   1   2   1   2   \\
1   1   1   1   0   2   0   1   0   0   1   1   1   2   2   2   1   1   0   2   2   2   0   2   1   2   1   1   2   1   0   2   0   2   2   2   2   1   2   2   0   2   2   0   2   1   1   0   0   &   1   1   1   1   \\
0   0   0   0   0   2   2   1   0   0   0   0   2   0   0   0   2   1   0   2   0   1   2   1   1   2   1   0   0   0   0   1   0   2   1   2   2   0   0   1   0   2   1   2   1   2   2   2   2   &   2   2   2   2   \\
0   1   0   0   0   1   0   2   1   2   0   1   1   2   2   1   2   2   1   2   2   2   2   1   0   1   0   0   0   2   2   2   0   0   0   0   1   0   0   2   0   2   0   2   0   1   0   1   1   &   2   2   1   1   \\
0   0   1   1   1   0   0   1   2   2   2   2   0   2   2   2   0   2   2   2   1   2   2   2   0   0   2   2   2   2   0   1   2   1   1   2   1   2   2   0   0   1   1   0   0   0   1   1   1   &   2   2   2   2   \\
0   1   0   1   2   0   1   0   0   0   1   1   1   0   2   2   2   0   1   0   1   0   0   0   1   0   0   0   1   2   0   2   1   1   0   2   2   0   2   0   2   1   2   0   1   2   0   1   0   &   1   2   1   2   \\
0   1   1   1   0   1   0   2   2   0   1   0   0   1   1   2   0   0   1   0   2   2   1   1   0   2   1   2   1   0   0   0   2   1   0   1   0   2   2   1   1   1   0   1   0   2   2   0   0   &   2   2   2   1   \\
1   1   0   0   2   0   1   1   1   2   1   1   2   0   0   0   2   1   1   1   2   2   1   1   0   1   2   1   2   2   0   1   2   0   1   1   1   2   1   1   1   0   2   1   0   0   1   0   1   &   1   1   1   2   \\
0   0   1   0   2   2   1   0   0   1   2   2   0   1   0   0   1   2   0   0   1   1   0   0   1   2   2   1   2   2   1   0   0   0   2   0   1   0   1   1   2   1   2   0   0   2   2   0   0   &   2   2   2   1   \\
1   1   0   0   2   1   2   2   0   1   2   1   1   1   1   2   1   1   2   2   0   2   2   1   1   1   2   0   1   2   1   0   0   1   2   2   0   2   2   0   1   0   1   0   1   0   2   2   2   &   2   1   1   2   \\
0   1   1   1   1   2   1   1   2   0   2   2   1   1   0   0   0   1   1   2   1   1   1   0   2   0   1   2   2   0   2   2   1   2   1   1   1   0   0   1   1   0   0   2   1   0   0   2   1   &   1   1   1   1   \\
1   0   0   0   1   1   2   0   2   0   2   1   0   1   1   1   1   2   2   0   2   0   2   0   0   1   0   2   2   1   0   1   1   0   1   0   0   1   0   2   0   1   1   1   2   1   1   2   2   &   1   2   1   1   \\
0   0   1   0   1   1   0   0   1   2   1   2   2   2   1   2   2   0   0   1   0   2   1   2   1   2   0   2   2   0   1   1   2   2   2   0   0   1   0   0   2   0   0   2   1   1   2   1   0   &   1   1   1   1   \\
0   0   0   1   1   1   1   0   1   0   0   1   2   2   0   2   0   2   1   2   1   0   1   0   2   2   1   0   1   1   2   2   2   0   2   1   2   1   1   2   1   2   2   0   0   1   2   2   2   &   2   1   2   2   \\

 \ea  $    \\ \hline
 \et} \ec

\bc {\scriptsize  \tabcolsep=1pt \bt {c|c} \hline $D_{3}^* \in
\cC_{3}(72+48;2^{4}3^{45}\bullet 3^{1})$ & $D_{3}^* \in
\cC_{3}(72+48;2^{4}3^{45}\bullet 3^{2})$

\\\hline
$ \ba{cc}
\d_0^*  &\zero_{72\times 1}   \\

0   0   0   0   0   1   2   1   2   2   0   2   1   2   1   0   0   1   2   1   2   2   2   2   2   2   1   0   0   0   2   1   2   1   2   0   0   0   0   2   2   0   1   2   1   0   2   1   2       1   \\
1   1   1   0   2   2   1   1   0   2   2   2   2   2   1   2   1   1   0   0   2   2   1   1   2   0   2   2   2   0   1   1   2   2   1   1   0   0   1   0   1   0   0   0   1   0   1   2   2   &   1   \\
1   0   0   0   0   2   2   2   0   1   2   2   1   1   2   0   2   1   1   2   1   1   0   0   0   1   0   2   2   1   2   2   2   1   1   1   1   0   1   1   1   0   0   1   2   0   1   0   2   &   2   \\
0   0   1   0   0   0   0   2   2   1   1   0   0   2   0   0   0   2   0   2   0   0   1   1   2   2   2   2   2   0   0   2   2   2   2   0   1   0   2   2   1   1   1   0   2   2   2   0   0   &   1   \\
0   1   0   0   1   2   1   1   2   0   2   1   1   1   1   2   1   0   1   0   2   0   1   0   0   2   1   0   2   0   1   0   1   0   2   0   1   0   0   1   1   0   1   0   1   0   0   2   1   &   2   \\
1   0   0   1   0   1   1   1   2   1   0   0   2   1   0   0   0   0   0   0   2   1   2   2   0   1   2   1   1   2   1   1   0   1   2   1   0   1   1   1   1   2   0   1   0   2   0   2   1   &   1   \\
1   1   0   0   2   1   0   0   0   2   1   1   0   1   2   1   2   2   2   0   2   2   2   2   1   1   0   2   1   0   1   1   1   0   2   0   0   0   2   2   0   1   0   2   0   0   2   0   2   &   2   \\
0   0   1   1   2   2   1   0   0   1   2   0   2   2   0   0   0   1   0   1   1   1   0   1   1   2   1   2   0   1   1   1   0   2   2   0   2   1   1   0   2   0   1   1   1   1   2   2   2   &   1   \\
1   0   1   1   1   1   0   1   2   0   1   1   1   1   0   2   0   1   0   2   2   2   2   1   0   1   0   0   2   2   0   0   0   2   1   0   2   2   2   0   1   1   1   1   1   1   1   0   2   &   2   \\
0   1   0   0   0   1   0   0   1   0   1   2   2   0   2   0   0   1   2   2   1   1   0   1   1   0   2   2   2   0   1   1   0   2   0   0   0   1   0   1   0   1   0   0   2   2   0   1   1   &   2   \\
1   0   1   1   0   1   2   0   0   2   1   1   1   0   0   2   2   0   2   0   1   2   0   1   1   2   2   2   1   2   1   2   1   1   1   2   0   2   0   0   2   0   2   2   1   2   0   0   0   &   2   \\
0   0   1   1   0   0   0   0   1   2   0   1   0   2   2   0   1   0   2   0   1   0   2   0   0   2   1   1   2   1   1   2   1   0   0   1   2   0   1   0   2   1   2   0   0   1   1   1   1   &   1   \\
1   1   0   0   2   0   1   1   1   1   1   0   1   2   2   0   1   1   0   1   2   0   0   2   1   0   2   1   0   2   1   0   2   1   0   0   1   2   2   0   2   2   2   1   0   0   2   0   0   &   2   \\
0   0   0   0   1   1   0   0   2   2   1   0   1   2   1   1   2   0   1   0   0   1   1   0   1   1   1   1   0   1   2   2   1   2   0   1   0   1   2   1   2   0   0   2   1   2   2   2   0   &   2   \\
0   1   1   0   2   1   1   2   0   2   0   0   0   1   1   1   1   0   0   1   1   2   0   0   0   0   2   0   0   0   2   0   1   1   0   1   2   2   0   1   1   1   0   2   2   1   2   0   1   &   2   \\
0   0   1   1   1   0   2   1   1   0   1   0   0   1   1   1   1   0   1   2   2   1   1   2   1   1   2   0   1   1   0   1   1   1   2   2   1   0   2   0   2   1   2   1   2   0   1   1   2   &   1   \\
0   1   0   0   1   0   1   0   2   1   1   2   0   0   2   0   0   2   0   0   1   2   2   0   0   0   1   0   2   2   0   1   2   0   1   2   2   2   0   0   0   1   0   1   1   2   2   1   2   &   1   \\
1   0   0   0   1   2   2   0   1   1   0   0   2   0   0   1   2   0   0   1   0   0   2   0   1   1   1   1   2   0   1   1   1   2   1   1   2   0   0   0   1   2   2   2   1   0   1   0   0   &   1   \\
1   1   1   1   0   1   2   0   1   0   0   1   0   0   1   0   0   2   1   2   1   1   2   1   2   1   2   1   0   0   2   1   1   0   0   1   1   2   1   2   1   0   1   1   2   0   1   2   0   &   2   \\
0   0   0   1   1   2   2   1   1   2   1   0   2   1   2   1   2   1   2   2   0   0   2   2   2   2   2   0   0   2   2   0   2   1   1   2   0   1   1   2   0   1   1   0   1   1   1   2   2   &   2   \\
0   1   1   1   0   2   0   2   2   0   2   2   0   2   2   1   0   1   1   0   2   2   0   2   2   0   0   2   0   0   0   1   1   1   1   1   2   0   1   0   2   0   2   2   2   2   0   1   0   &   2   \\
1   1   0   1   2   0   2   2   2   2   1   2   2   0   0   0   1   2   1   1   0   1   1   0   2   0   2   1   2   1   2   0   1   2   1   0   0   1   2   2   1   1   2   1   2   1   0   1   2   &   2   \\
0   0   0   1   1   2   1   0   1   2   0   2   0   0   0   0   0   2   1   0   1   0   1   0   2   2   0   2   1   0   2   0   0   1   2   2   1   0   1   1   2   1   2   1   0   2   2   0   1   &   2   \\
1   0   1   0   1   0   1   0   2   2   2   0   1   0   2   2   1   1   2   1   1   1   0   2   2   1   2   2   0   2   0   2   1   2   1   1   2   1   2   1   0   1   1   1   0   0   0   0   1   &   1   \\
0   1   0   0   1   1   0   1   1   0   0   1   0   0   1   2   2   0   0   1   0   2   0   2   2   2   1   2   1   1   1   2   2   2   2   1   2   1   1   2   1   2   2   0   0   2   1   1   2   &   2   \\
0   0   0   1   2   0   2   2   0   0   1   1   1   1   2   2   0   2   0   0   2   1   0   2   1   2   2   1   1   0   0   0   1   0   1   2   1   1   1   1   1   2   1   2   1   2   2   1   0   &   1   \\
1   1   1   1   1   0   0   1   0   1   1   2   1   1   2   1   0   0   0   2   0   2   1   2   2   2   1   1   2   1   2   0   2   0   2   0   0   1   0   1   2   0   0   1   2   0   1   0   1   &   1   \\
0   1   0   1   2   2   1   2   1   2   1   1   1   2   1   2   2   1   0   2   0   2   0   1   0   0   0   0   2   0   2   2   0   1   0   2   1   0   2   0   1   0   0   2   0   1   0   1   1   &   1   \\
1   0   1   1   1   2   2   1   1   0   2   2   0   0   1   2   1   0   0   1   1   0   2   0   1   0   0   0   0   2   2   2   0   0   0   0   2   1   0   2   2   2   0   0   2   2   2   2   1   &   1   \\
0   0   1   0   0   2   1   1   1   1   2   1   0   1   2   0   0   0   2   2   1   2   2   2   1   0   0   1   2   1   2   2   1   2   2   2   1   1   2   2   0   2   1   2   1   1   0   0   0   &   2   \\
1   0   0   1   2   0   1   2   1   1   2   0   0   2   1   2   2   0   2   1   2   2   1   1   2   1   0   2   1   1   0   0   2   0   0   1   0   1   0   2   0   0   1   1   1   1   1   2   1   &   2   \\
0   0   1   0   1   1   0   2   1   1   2   0   1   1   1   1   2   2   1   1   0   1   2   1   0   0   0   2   1   0   1   1   1   0   1   0   1   2   0   2   0   2   1   1   0   2   0   1   1   &   1   \\
1   1   1   1   1   1   2   0   0   0   0   1   2   2   1   0   2   0   1   0   2   1   0   0   0   1   0   2   2   1   0   0   2   1   0   0   0   0   0   2   2   2   1   0   0   1   0   0   0   &   1   \\
1   0   0   0   2   0   2   0   0   1   0   0   0   0   1   1   1   2   2   0   0   0   0   0   2   0   0   0   0   2   2   2   0   2   2   0   1   0   1   1   2   2   1   2   2   1   0   1   2   &   1   \\
0   1   1   0   2   2   2   1   1   2   2   0   2   2   0   1   2   2   1   0   1   1   0   2   0   2   0   1   2   2   1   0   0   0   0   1   1   1   0   1   0   0   2   0   0   2   2   0   2   &   1   \\
0   1   1   1   2   0   2   2   0   0   2   2   1   2   2   0   1   2   0   2   0   1   2   0   2   1   1   1   1   2   1   2   0   2   2   2   2   0   0   0   0   2   1   2   0   0   1   2   1   &   2   \\
1   0   0   1   2   0   1   1   2   0   0   1   2   2   2   1   0   2   2   2   0   1   1   1   0   2   0   0   0   1   0   2   2   1   1   1   2   2   0   2   0   2   0   0   0   0   2   1   1   &   2   \\
0   0   1   1   0   1   1   2   2   0   1   0   2   1   0   2   2   1   1   2   1   0   1   2   2   2   1   2   0   0   1   2   2   0   0   1   2   2   2   1   0   2   0   2   1   1   0   2   2   &   1   \\
0   1   0   1   0   2   2   0   0   1   0   1   2   2   0   1   1   1   0   1   1   2   2   1   1   1   1   0   1   1   0   0   0   0   1   0   2   0   2   1   0   2   2   0   2   2   0   2   0   &   2   \\
1   1   1   0   2   2   0   1   2   2   0   1   2   0   0   0   2   1   1   2   2   1   1   0   1   0   1   0   0   2   0   0   2   2   2   2   1   2   1   2   0   1   2   2   2   2   1   0   1   &   1   \\
0   1   1   1   2   0   0   1   0   1   0   1   1   0   0   2   2   2   1   1   2   0   2   1   0   0   2   2   2   2   0   2   2   2   0   2   2   2   2   1   2   0   2   0   1   0   2   1   1   &   2   \\
1   1   1   1   1   1   1   2   0   1   2   2   2   0   2   2   0   0   2   1   0   0   2   0   1   0   1   1   1   0   2   0   2   1   0   2   0   2   2   0   1   2   1   2   0   2   0   1   2   &   2   \\
1   0   0   0   0   1   1   0   0   2   2   2   1   2   0   2   1   2   2   2   0   0   1   0   2   0   0   0   1   2   2   1   1   0   1   2   2   2   2   2   0   1   2   0   2   1   1   2   0   &   1   \\
1   1   0   0   2   0   0   2   2   0   2   1   0   1   0   1   0   2   2   2   1   0   0   1   0   1   2   0   1   1   2   0   0   1   2   1   0   2   1   0   2   1   2   0   1   0   0   2   0   &   1   \\
1   0   0   0   0   2   0   2   0   0   0   2   0   0   1   2   1   1   1   0   0   2   1   1   1   2   1   1   1   1   0   1   0   1   0   2   1   2   2   1   0   1   0   2   0   0   2   2   0   &   1   \\
0   1   0   1   0   1   0   2   1   2   2   2   2   0   0   1   1   0   2   1   2   0   1   2   0   1   2   1   0   2   0   2   0   0   1   0   1   2   1   0   2   0   0   1   0   1   1   0   0   &   1   \\
1   1   1   0   0   2   0   0   2   0   1   0   1   1   1   1   1   2   2   0   2   0   0   2   0   2   1   1   1   1   0   1   0   2   0   2   0   1   0   2   1   2   2   1   2   1   1   2   2   &   2   \\
1   1   1   0   1   0   2   2   2   1   0   2   2   1   2   2   2   1   1   1   0   2   1   1   1   1   0   0   0   2   1   1   0   0   2   2   0   1   1   0   1   0   0   0   2   1   2   1   0   &   2   \\
 \ea  $
 &
 $ \ba{cc}
\d_0^*  &\zero_{72\times 2}   \\
1   1   0   0   1   0   1   2   1   2   2   2   1   0   2   2   2   0   1   1   0   0   1   0   2   1   0   2   2   1   1   2   2   2   1   0   1   0   2   0   2   0   1   1   0   1   1   0   1   &   1   1   \\
0   0   1   0   1   1   0   0   1   2   2   1   0   1   1   1   1   0   1   0   2   2   1   2   1   1   1   1   2   0   1   1   1   0   2   1   1   0   0   2   1   0   0   2   0   1   2   2   1   &   2   1   \\
0   0   1   1   2   0   1   1   2   0   2   0   0   2   2   0   0   0   1   2   1   1   0   2   2   0   0   0   0   2   2   2   0   2   2   2   1   0   0   2   0   2   1   0   1   1   0   2   2   &   1   1   \\
0   0   0   1   0   2   2   0   1   0   1   0   2   0   1   2   1   2   1   2   0   1   1   1   1   1   0   0   2   0   0   0   0   0   0   2   1   1   2   1   1   0   1   0   0   1   2   2   0   &   1   2   \\
0   0   1   1   1   2   2   0   0   0   2   2   2   1   0   0   1   2   0   0   1   1   0   0   2   2   2   1   1   1   1   1   0   0   1   2   1   0   1   0   0   1   1   0   2   2   2   2   0   &   2   1   \\
0   1   0   0   0   0   0   0   2   0   1   0   0   1   0   0   0   2   0   0   1   0   1   2   2   2   0   2   1   2   1   1   1   2   2   0   2   0   2   1   2   1   2   1   1   2   2   1   0   &   1   2   \\
1   0   1   0   1   2   0   2   2   1   2   2   2   2   0   2   1   1   0   1   1   2   0   2   1   0   1   2   2   1   0   2   2   2   1   1   1   1   1   1   0   1   0   0   1   2   0   0   2   &   1   1   \\
1   0   1   0   1   1   0   0   0   0   0   1   2   0   0   2   1   0   2   2   0   0   1   0   2   2   1   1   1   1   2   2   2   2   2   2   2   2   1   2   1   2   2   2   0   1   0   2   0   &   1   2   \\
1   1   0   0   0   1   2   0   1   2   1   0   1   2   1   1   2   0   2   0   1   0   0   1   1   0   0   0   1   2   2   2   2   2   0   1   0   1   1   2   2   2   0   0   2   1   2   1   2   &   2   1   \\
1   0   1   0   2   0   1   0   0   2   1   0   0   0   2   1   0   2   0   0   2   2   0   0   0   2   2   2   0   2   0   2   2   1   0   2   1   1   1   2   0   1   0   1   0   2   2   0   1   &   2   2   \\
1   0   0   0   2   0   2   0   2   1   1   0   1   1   1   1   2   0   0   1   0   1   2   2   1   1   2   1   0   2   1   0   1   0   0   2   0   1   0   1   1   0   1   1   1   2   0   2   1   &   1   1   \\
1   0   0   0   0   2   2   0   1   2   0   1   0   2   1   0   2   1   1   2   0   2   2   2   1   2   0   2   0   1   2   1   1   0   1   2   1   0   1   2   0   2   2   2   2   0   1   0   2   &   2   2   \\
0   1   1   1   0   2   0   1   2   2   1   2   1   2   1   2   0   2   0   2   0   2   1   2   2   0   1   2   1   0   2   0   0   2   2   0   0   1   2   0   2   0   1   2   2   0   2   2   0   &   2   1   \\
1   1   0   1   0   1   1   0   0   0   2   1   1   0   1   2   1   2   2   1   2   1   0   2   2   0   2   2   2   0   1   0   1   0   1   0   2   1   0   2   0   2   2   1   2   2   0   0   0   &   1   1   \\
1   0   0   0   0   0   0   1   2   1   1   1   0   0   2   0   1   1   2   0   2   2   2   2   0   1   2   1   2   2   0   2   2   2   2   0   0   1   2   1   0   1   1   2   2   1   1   1   0   &   2   1   \\
1   1   1   1   0   2   2   2   2   1   1   2   0   1   1   0   1   0   2   0   2   0   1   0   2   0   2   2   0   0   0   2   1   0   0   1   1   0   1   0   2   0   0   2   2   1   0   0   0   &   1   2   \\
0   1   0   1   2   0   2   2   0   2   0   2   0   0   1   2   1   1   2   2   0   2   1   0   1   0   2   0   1   2   2   0   2   1   2   2   0   2   2   2   2   1   2   0   1   2   2   1   2   &   1   1   \\
0   1   1   0   2   2   1   0   1   1   0   0   2   2   0   2   2   1   2   1   0   1   0   0   2   0   1   2   0   0   1   1   1   2   0   1   0   0   0   2   1   0   2   2   2   2   2   1   2   &   1   2   \\
0   1   1   0   0   1   0   1   0   1   0   0   1   1   0   2   2   2   1   0   0   1   1   1   0   2   0   1   0   2   1   2   1   1   1   2   0   1   2   1   0   2   0   0   0   2   2   0   2   &   2   1   \\
0   0   0   0   2   2   1   2   0   0   1   1   2   1   1   0   2   0   1   1   2   0   1   0   2   2   0   0   2   0   0   0   2   1   2   1   0   0   1   2   1   1   2   0   1   2   1   1   0   &   2   1   \\
1   0   1   1   1   1   0   1   1   1   1   1   1   1   1   2   0   1   2   2   1   0   0   2   2   2   2   0   0   1   2   0   1   1   2   1   1   2   2   2   1   1   1   0   0   0   1   2   1   &   1   1   \\
0   0   0   1   1   0   1   1   1   1   0   2   2   2   2   2   0   0   0   0   2   1   1   2   0   0   2   1   1   2   1   1   2   0   0   1   0   1   1   0   2   2   2   1   0   0   1   1   1   &   1   2   \\
1   1   1   1   2   2   1   1   0   0   0   2   0   1   1   0   1   0   0   0   1   2   0   1   0   0   2   2   2   0   1   0   2   2   0   0   0   2   0   1   1   2   1   0   1   0   0   2   1   &   2   2   \\
0   1   0   0   0   1   1   2   1   1   2   2   0   1   2   0   0   2   0   1   0   2   2   0   0   0   1   1   1   0   2   2   0   1   2   0   1   2   0   1   1   2   0   2   0   2   0   1   0   &   1   1   \\
1   0   0   1   1   0   2   2   2   1   0   0   0   1   2   1   0   0   2   1   0   2   0   1   1   1   0   0   1   1   0   0   1   1   1   1   0   0   2   0   0   2   0   2   0   0   0   0   0   &   2   2   \\
0   1   0   0   2   0   2   2   0   0   0   1   1   2   1   1   0   1   0   0   0   1   0   1   1   1   1   0   1   0   0   1   0   2   1   2   1   0   1   1   1   0   0   1   2   0   1   1   0   &   1   1   \\
0   0   1   1   0   1   2   2   2   2   2   1   2   2   0   1   2   1   2   2   1   1   1   1   0   2   2   2   2   2   1   2   0   2   1   0   1   0   0   1   1   0   1   1   1   0   0   1   1   &   2   2   \\
1   1   0   0   2   2   1   1   2   2   0   1   2   1   0   0   2   2   1   2   1   1   1   1   0   1   2   0   2   2   2   2   0   0   2   1   2   2   1   1   1   0   2   2   2   0   2   0   1   &   1   1   \\
1   1   0   1   2   1   0   2   0   2   1   1   2   1   0   1   1   1   0   1   0   1   0   0   0   1   0   0   0   1   1   1   0   1   2   0   2   1   1   0   2   2   0   2   1   1   0   2   2   &   2   1   \\
1   1   0   1   0   1   0   2   2   0   0   0   2   2   2   0   0   1   1   2   2   1   2   1   2   1   2   2   1   0   0   1   2   2   0   1   2   2   1   2   0   2   1   2   0   2   1   2   1   &   2   1   \\
1   1   0   1   2   0   1   0   2   2   2   0   1   2   0   0   1   1   1   0   2   1   0   0   0   2   1   1   1   1   0   0   1   1   1   1   1   2   2   1   1   1   2   1   0   0   2   0   2   &   1   2   \\
0   1   0   0   1   0   0   1   2   0   0   2   1   2   1   1   1   2   1   1   2   0   1   0   0   0   1   0   0   1   2   1   1   0   0   0   1   1   0   2   0   1   0   0   2   0   1   1   1   &   2   2   \\
0   0   1   0   0   2   1   0   2   1   2   2   1   2   2   0   1   2   2   2   1   0   2   0   0   1   0   0   2   2   1   0   0   1   2   2   2   0   0   0   2   2   2   0   0   2   1   2   1   &   2   2   \\
0   1   0   0   1   1   2   1   2   2   2   1   2   0   2   1   0   2   2   0   2   2   2   2   2   2   1   0   2   2   1   1   1   1   1   1   0   2   1   0   0   0   1   0   1   0   0   2   2   &   1   2   \\
0   1   1   1   1   1   2   2   1   0   1   2   1   1   2   0   0   0   1   2   1   2   2   2   0   1   0   2   1   1   2   0   2   1   1   1   0   1   2   1   1   0   0   1   2   1   1   1   2   &   1   2   \\
0   1   1   1   0   1   2   2   1   2   0   1   2   0   0   1   0   0   1   1   1   0   0   2   1   1   2   1   0   2   2   0   0   0   0   0   0   2   1   0   2   1   2   1   2   2   0   0   0   &   2   2   \\
0   0   1   1   2   0   1   2   0   2   1   0   1   2   0   1   1   2   0   1   0   0   2   2   2   0   2   1   2   0   2   2   1   0   1   1   2   2   2   1   0   2   1   2   1   1   1   0   2   &   2   2   \\
0   1   1   1   2   0   0   0   1   0   2   1   0   0   2   1   2   0   1   2   1   0   2   1   0   2   1   1   0   0   0   2   1   0   0   1   2   1   2   0   2   0   2   0   0   2   0   1   0   &   2   1   \\
1   0   1   1   0   0   2   1   1   1   0   2   2   2   2   1   2   2   0   1   0   0   2   1   0   0   0   0   2   1   0   1   0   1   1   0   2   2   0   2   0   1   2   0   2   1   2   1   0   &   1   1   \\
1   0   0   0   2   2   1   2   1   0   2   2   0   1   0   2   0   1   2   2   0   2   1   1   0   2   1   1   2   0   2   0   0   0   0   0   2   1   2   2   0   1   0   1   2   1   2   2   2   &   2   2   \\
1   0   0   0   1   0   1   0   1   2   1   2   2   0   2   0   0   1   2   2   1   1   2   0   1   0   0   0   0   0   1   1   1   2   1   2   2   2   2   0   1   1   0   2   1   1   0   1   1   &   2   1   \\
1   0   1   0   1   2   0   2   0   2   2   2   1   0   1   2   2   1   1   0   2   1   0   1   1   0   0   1   0   2   2   0   2   1   0   2   2   2   0   0   2   0   1   2   0   0   0   1   0   &   1   2   \\
1   1   1   1   2   1   0   0   0   1   2   1   0   2   0   2   2   2   2   0   1   2   2   1   1   1   0   2   1   2   0   1   0   0   2   0   0   0   1   2   2   1   1   0   0   1   1   0   1   &   1   2   \\
0   0   0   1   1   2   2   1   0   1   0   0   0   0   0   1   0   0   1   0   1   0   2   0   2   2   1   2   0   1   2   2   2   1   0   0   2   2   0   1   2   0   0   1   1   0   2   2   1   &   1   1   \\
1   1   1   0   0   2   2   1   0   1   2   0   2   0   2   1   1   0   0   1   2   0   2   0   1   1   1   1   1   1   1   1   0   2   2   2   2   1   2   1   1   2   1   1   2   0   1   0   2   &   2   2   \\
0   0   1   0   2   1   0   1   1   1   0   0   1   0   1   0   2   1   0   1   2   2   2   1   2   2   2   0   0   1   0   0   2   0   2   0   1   0   0   0   2   0   0   1   1   1   2   0   2   &   2   2   \\
1   0   1   1   1   1   1   1   2   0   1   0   0   0   0   2   2   2   2   1   2   0   1   1   1   1   1   2   1   0   0   1   2   1   1   2   0   2   0   0   2   1   2   1   1   0   1   2   2   &   2   1   \\
0   1   0   1   1   2   0   1   0   0   1   1   1   1   2   2   2   1   0   2   2   2   0   2   1   2   1   1   2   1   0   2   0   2   0   2   2   0   0   0   0   2   2   2   1   2   1   0   1   &   1   2   \\
 \ea  $
  \\ \hline
 \et} \ec

\bc {\scriptsize  \tabcolsep=1pt \bt {c|c} \hline $D_{3}^* \in
\cC_{3}(72+48;2^{4}3^{45}\bullet 3^{3})$ & $D_{3}^* \in
\cC_{3}(72+48;2^{4}3^{45}\bullet 3^{4})$

\\\hline
$ \ba{cc}
\d_0^*  &\zero_{72\times 3}   \\
1   0   1   1   0   2   2   0   0   0   2   0   2   1   2   0   1   1   2   2   1   2   2   1   0   1   1   1   1   1   2   2   0   0   1   2   0   1   1   1   0   2   1   0   0   1   1   0   2   &   2   2   1   \\
0   0   0   0   0   2   0   1   2   0   1   0   2   2   0   2   1   2   0   2   2   0   1   2   2   2   2   0   2   2   0   2   0   0   2   0   2   1   1   1   1   1   0   0   1   1   2   2   0   &   1   2   1   \\
0   1   1   0   2   1   1   1   1   2   1   1   1   1   0   1   0   2   2   2   2   0   1   2   2   0   2   2   1   0   2   0   1   1   2   2   2   2   2   2   1   1   2   0   0   2   1   0   0   &   2   2   2   \\
1   0   0   1   2   0   2   1   1   2   0   0   0   1   2   1   2   0   1   0   2   2   2   2   0   1   0   0   0   2   0   2   1   1   0   0   1   0   0   2   2   2   0   0   0   1   0   0   1   &   1   2   1   \\
0   1   1   1   2   0   2   2   0   0   1   2   2   2   2   2   1   0   0   2   0   0   1   0   2   2   1   1   2   0   1   2   0   2   2   0   1   0   2   0   0   0   2   2   0   1   0   1   0   &   2   1   2   \\
0   1   0   1   1   2   1   1   1   0   1   2   0   1   1   0   0   1   0   1   1   2   1   2   0   0   0   2   2   0   2   0   2   0   1   0   1   0   0   1   0   1   1   0   1   2   1   1   1   &   1   1   2   \\
0   1   1   1   2   0   1   2   2   1   0   2   2   2   2   0   0   0   0   1   1   0   2   1   2   0   2   2   0   0   0   0   2   0   0   1   0   2   1   0   2   2   0   0   0   2   0   1   0   &   1   1   1   \\
0   1   0   0   2   1   1   0   1   0   0   1   2   2   1   2   2   0   2   2   0   1   1   1   2   1   0   2   2   2   1   1   2   2   0   1   0   0   1   2   1   2   1   0   2   1   1   2   2   &   1   2   1   \\
1   0   1   1   0   1   0   2   2   0   1   1   0   0   1   1   0   0   1   2   1   0   1   1   0   1   0   2   1   0   2   1   0   1   1   1   1   2   2   1   2   0   0   2   1   0   0   2   0   &   2   1   1   \\
1   0   1   0   2   1   2   1   2   2   0   0   2   2   1   1   0   2   0   1   2   2   2   2   1   1   1   1   1   1   1   1   2   1   2   1   0   2   0   2   1   1   1   2   1   0   2   0   2   &   1   1   2   \\
0   1   0   0   0   2   2   1   0   1   0   1   1   0   1   1   1   1   0   1   2   2   0   1   1   0   2   1   0   2   2   2   0   2   0   0   0   1   0   1   1   2   0   2   2   0   0   1   0   &   2   1   2   \\
0   0   0   1   1   2   2   0   0   0   2   1   0   0   1   2   0   2   2   0   2   0   2   0   2   0   1   0   1   0   2   1   1   1   2   2   2   0   1   0   2   0   1   2   1   0   0   2   2   &   1   2   2   \\
1   1   1   0   0   2   0   0   1   2   1   0   0   0   1   0   1   0   1   1   2   2   1   0   1   0   1   2   2   1   0   1   1   0   0   0   1   1   2   2   0   1   0   1   2   1   1   0   0   &   1   1   2   \\
1   1   0   1   0   0   2   2   2   1   0   1   0   1   2   0   1   2   2   2   1   0   2   0   0   1   2   0   2   2   0   0   0   2   2   2   1   0   2   1   1   1   1   0   2   0   1   2   1   &   2   1   2   \\
1   1   1   1   0   2   1   0   0   1   1   0   1   0   0   2   2   1   0   0   1   0   0   0   2   0   0   1   0   2   0   2   2   2   0   2   2   0   2   1   2   2   2   2   1   1   2   0   0   &   1   1   2   \\
1   1   1   1   1   1   0   1   0   0   1   0   1   2   1   2   0   1   0   1   2   0   2   1   0   0   1   2   1   1   0   1   2   2   0   0   2   2   2   0   2   2   1   1   2   0   1   1   1   &   2   2   1   \\
0   0   1   1   1   1   2   0   2   1   1   0   1   1   1   1   2   0   0   0   2   1   2   2   2   1   2   1   2   1   1   0   1   0   0   0   0   0   0   1   1   0   1   1   1   1   1   2   1   &   2   1   1   \\
1   0   0   0   1   1   2   0   1   1   1   2   2   0   0   1   1   1   2   1   0   0   2   0   0   1   2   0   1   0   1   1   2   0   0   2   2   2   0   0   0   1   0   1   0   1   2   2   2   &   2   1   1   \\
0   0   1   0   1   1   0   0   1   0   0   1   1   1   1   0   2   2   1   2   1   1   2   0   0   2   1   0   0   1   2   2   2   0   2   1   1   2   0   2   1   2   2   2   0   1   0   2   1   &   1   2   2   \\
0   1   0   1   2   0   0   0   1   2   2   1   1   2   1   1   1   2   1   0   2   0   0   2   0   1   0   1   1   1   2   0   1   0   1   0   1   1   1   2   2   1   2   0   2   0   2   1   0   &   1   1   1   \\
0   0   1   0   2   0   0   2   0   2   2   0   1   2   2   0   2   1   1   2   0   1   1   2   1   1   0   0   0   2   0   1   2   1   2   2   2   1   2   0   0   0   1   0   0   0   2   2   1   &   2   1   2   \\
0   0   1   0   1   0   1   1   1   1   2   2   0   0   2   2   1   0   1   0   1   0   1   0   0   2   1   1   2   1   1   2   1   1   0   1   2   1   2   1   2   1   0   0   0   2   2   1   1   &   2   2   1   \\
0   1   0   0   0   0   0   0   2   1   0   1   1   2   0   1   1   1   0   2   1   2   0   2   1   1   1   0   1   0   1   1   1   2   2   1   2   0   2   1   0   0   0   1   2   2   1   1   1   &   1   2   1   \\
1   1   1   0   1   2   1   2   1   2   2   2   1   2   1   2   2   0   0   1   0   2   0   2   1   0   0   0   1   2   2   0   0   1   2   1   1   0   1   0   1   0   0   2   0   1   2   0   1   &   2   2   1   \\
0   1   1   1   1   1   2   2   1   1   2   0   2   0   0   2   2   2   1   1   1   2   2   1   1   1   0   2   2   2   1   0   0   0   1   0   2   1   2   0   2   0   1   1   0   2   0   1   2   &   1   2   2   \\
1   1   1   0   2   2   0   1   0   1   2   2   0   1   2   2   2   2   1   0   2   2   0   1   0   0   2   0   2   0   1   1   0   0   0   2   1   2   0   0   0   0   2   0   2   0   2   0   2   &   1   2   1   \\
0   0   0   1   0   0   0   0   1   2   1   2   0   0   2   1   0   2   0   0   0   2   1   1   1   0   2   1   1   2   0   1   0   2   1   2   0   2   1   0   1   1   2   2   0   2   1   1   1   &   1   2   1   \\
0   1   1   0   0   1   1   2   0   2   0   0   0   1   0   1   2   2   1   0   0   0   2   1   0   0   1   2   0   0   1   2   1   1   1   1   2   0   1   1   0   0   0   2   2   1   1   0   2   &   2   1   2   \\
1   1   1   1   2   0   1   0   0   0   2   2   1   1   0   1   1   2   2   1   0   1   1   0   1   0   1   0   1   2   1   1   0   0   2   0   2   1   0   2   2   2   2   1   2   2   0   2   1   &   1   1   1   \\
1   0   0   0   0   1   1   0   0   2   2   2   1   2   0   2   0   1   2   2   1   1   2   1   0   2   2   1   0   2   2   2   1   0   0   0   2   2   0   2   0   0   0   0   1   2   0   1   1   &   1   1   2   \\
1   1   0   0   1   2   2   2   2   0   2   0   2   1   0   0   0   1   1   0   2   1   0   1   0   2   1   0   1   2   0   1   0   1   1   1   1   2   1   1   1   1   2   1   0   2   2   1   2   &   1   1   2   \\
1   1   0   0   1   1   1   2   0   1   0   0   0   1   2   2   0   0   0   0   0   1   1   0   1   1   1   1   1   0   2   0   1   1   0   2   0   1   2   0   1   2   0   2   2   0   2   2   0   &   1   2   2   \\
1   0   0   0   0   1   2   2   1   2   1   2   1   0   1   0   1   2   2   2   0   1   0   2   2   2   2   2   0   0   1   2   1   2   2   0   0   1   2   1   2   0   1   2   2   2   2   0   2   &   1   1   1   \\
0   0   1   1   2   2   2   0   2   2   2   0   0   2   1   0   1   1   1   0   1   1   0   0   2   0   2   2   0   2   2   0   2   1   0   2   1   2   1   1   2   1   1   1   2   2   0   2   2   &   2   2   1   \\
1   0   1   0   1   0   0   2   2   1   0   0   0   1   1   1   2   1   2   1   0   0   0   1   2   2   0   0   0   1   2   0   0   2   1   2   0   1   2   2   0   1   1   0   1   1   0   1   0   &   1   2   2   \\
1   1   0   0   2   0   2   1   2   2   0   2   1   0   0   0   2   0   1   1   0   0   1   0   2   2   1   1   0   1   2   0   2   2   1   0   0   2   0   1   1   0   2   1   1   0   1   0   1   &   2   2   1   \\
1   1   0   0   0   1   0   1   2   2   2   2   2   0   2   0   0   1   2   0   2   1   2   0   0   0   0   2   2   2   1   1   2   2   2   0   0   1   1   2   2   0   0   2   0   1   0   1   0   &   2   2   2   \\
0   1   0   1   1   1   0   2   2   0   1   2   1   1   2   0   1   1   2   1   2   0   2   0   1   1   0   1   0   0   0   2   1   1   1   1   1   1   2   2   1   2   2   2   1   2   1   0   2   &   1   2   1   \\
1   0   0   1   2   0   0   2   2   0   0   1   0   0   0   2   2   2   1   1   0   1   1   1   2   1   2   1   0   0   0   2   2   0   2   1   1   0   1   2   2   1   1   1   0   0   2   0   0   &   2   2   2   \\
0   0   1   1   0   2   0   1   1   1   0   2   2   2   0   0   0   0   2   2   0   1   2   2   2   1   1   1   0   1   2   1   0   2   2   2   1   0   0   0   2   0   0   1   2   1   1   2   2   &   1   1   2   \\
0   1   0   1   1   2   0   1   0   0   0   1   2   0   1   2   2   1   1   2   0   2   0   1   1   2   1   2   1   1   0   2   2   1   0   2   0   2   0   2   0   1   2   0   1   2   2   2   0   &   2   1   1   \\
1   1   1   1   2   2   1   1   0   2   1   1   2   0   0   1   1   0   2   1   1   1   0   2   1   2   2   0   0   1   1   0   1   0   1   1   2   1   1   0   1   0   2   0   1   0   0   2   2   &   2   2   2   \\
0   0   0   0   0   1   1   0   0   1   2   1   0   1   0   0   0   0   2   0   1   2   1   1   1   2   0   0   2   1   1   0   0   1   2   0   0   0   1   0   2   1   2   1   1   1   2   0   0   &   2   1   1   \\
1   0   0   1   2   2   1   2   1   1   1   1   2   1   0   2   2   1   0   2   2   2   1   0   1   2   0   1   2   1   2   0   2   2   1   1   2   2   1   2   0   2   0   2   2   2   0   0   2   &   2   1   1   \\
0   0   0   0   2   0   1   1   0   1   2   2   1   2   2   1   0   2   0   0   0   2   0   2   2   2   0   2   1   1   0   2   1   0   1   2   1   0   2   1   0   2   1   1   1   0   0   0   2   &   2   2   2   \\
0   0   0   1   0   0   2   2   1   0   1   0   2   2   2   1   0   0   0   2   1   1   0   2   0   0   2   0   2   0   0   0   0   1   0   1   0   1   0   2   0   2   2   1   1   2   2   1   1   &   2   2   2   \\
1   0   1   0   1   2   2   0   2   0   0   1   2   2   2   0   2   2   1   1   1   1   0   0   1   2   0   2   2   0   0   1   1   2   1   2   2   0   0   0   1   2   2   2   2   0   1   1   0   &   2   1   1   \\
1   0   1   1   1   0   1   1   2   2   2   1   0   0   2   2   1   0   2   0   2   2   0   2   2   2   2   2   2   2   1   2   2   2   1   1   0   2   0   0   0   2   1   1   0   0   1   2   1   &   1   1   2   \\

 \ea  $
 &
 $ \ba{cc}
\d_0^*  &\zero_{72\times 4}   \\
1   1   0   1   0   1   2   1   0   2   0   1   1   0   1   0   1   0   2   1   2   1   2   1   1   0   2   0   0   2   2   1   1   1   1   2   2   2   0   0   1   2   2   2   2   1   0   2   0   &   2   1   2   1   \\
1   1   1   1   2   1   1   2   0   1   0   2   0   1   1   2   1   0   0   1   0   2   2   1   0   0   2   1   2   0   1   0   2   0   0   0   0   2   2   1   1   2   1   0   0   1   1   0   1   &   2   1   2   1   \\
0   1   0   0   1   1   0   0   1   2   2   1   0   0   1   2   1   2   2   0   2   2   1   2   0   0   0   2   1   1   2   0   1   0   1   0   1   1   2   2   0   1   2   1   2   1   2   0   2   &   2   2   1   1   \\
1   1   0   0   1   0   2   2   1   1   1   1   1   0   0   1   0   2   1   1   1   0   2   0   1   1   1   0   0   2   0   0   0   1   0   1   1   0   2   1   2   0   0   1   2   2   2   1   1   &   2   1   2   2   \\
1   1   0   1   2   0   1   1   2   0   2   1   0   1   2   0   0   2   1   2   2   1   1   2   0   1   0   2   2   2   0   1   1   1   2   1   1   0   1   2   1   2   1   0   2   0   1   2   2   &   2   1   1   2   \\
0   0   1   1   1   1   2   2   0   0   0   0   2   2   1   2   0   2   1   2   0   1   1   1   2   2   2   2   1   0   0   1   2   1   2   2   1   0   2   2   0   0   1   1   0   0   2   2   2   &   1   2   2   2   \\
0   0   0   1   1   2   1   2   1   2   0   2   2   2   1   0   0   2   0   2   0   0   0   1   1   2   0   0   2   0   1   1   0   0   1   0   2   0   1   0   2   2   0   2   1   1   2   1   2   &   2   1   1   1   \\
1   0   1   1   1   2   2   0   2   2   2   2   2   0   0   0   1   1   2   1   1   1   1   0   2   0   0   0   2   2   2   1   0   2   2   0   1   1   1   1   2   0   0   1   2   1   0   2   2   &   2   2   1   2   \\
1   1   0   0   2   0   1   0   0   0   1   2   2   2   2   2   1   2   0   1   2   1   1   0   1   0   1   1   1   0   0   1   0   0   2   1   2   1   2   0   0   1   2   2   0   2   2   0   0   &   1   1   1   2   \\
0   0   0   1   0   1   2   1   2   1   1   2   2   2   2   0   0   1   0   1   2   2   2   2   0   1   2   1   0   1   2   0   2   1   2   0   0   1   0   0   0   0   0   1   1   0   2   1   2   &   2   1   2   1   \\
0   1   1   1   0   0   0   2   0   0   2   2   0   2   2   2   1   0   0   0   1   0   0   2   2   0   2   2   2   2   2   2   0   2   0   0   1   0   2   0   2   2   2   1   2   1   0   1   0   &   1   1   2   2   \\
1   1   1   1   2   2   2   0   0   0   2   2   1   1   0   1   0   2   1   2   1   0   2   0   0   0   1   0   1   1   2   0   0   0   0   0   2   2   0   0   2   1   2   1   2   0   1   2   2   &   1   2   2   1   \\
0   0   1   1   0   2   2   1   2   0   1   2   0   1   0   0   1   0   1   2   1   0   1   2   1   2   1   1   0   0   0   0   0   0   2   2   1   0   2   1   1   0   2   0   1   2   2   0   0   &   2   2   1   1   \\
0   1   1   0   0   2   1   0   1   1   0   0   0   2   1   1   0   2   0   1   1   2   0   2   1   0   2   0   1   2   0   1   2   0   0   0   0   2   0   2   1   1   1   1   1   2   2   2   0   &   2   2   1   2   \\
1   1   1   1   2   2   0   1   2   1   1   1   0   1   1   2   2   0   2   0   1   0   0   1   1   1   2   0   1   0   2   0   0   1   2   1   0   1   1   0   1   1   0   2   0   0   2   2   0   &   2   1   2   1   \\
0   0   0   0   1   1   1   0   0   2   1   1   2   1   0   0   2   2   1   0   1   1   0   0   0   2   2   0   2   1   1   2   1   1   1   1   0   1   1   1   1   1   2   0   1   2   1   1   2   &   1   2   1   1   \\
1   0   1   1   1   2   0   1   0   1   0   0   2   1   0   1   2   0   0   1   2   0   1   0   0   1   1   1   2   1   0   2   2   2   1   0   2   1   0   2   0   2   0   2   2   1   1   0   1   &   1   2   1   2   \\
0   0   1   0   2   1   1   0   0   0   2   0   1   0   1   2   0   1   0   0   0   1   1   0   0   2   1   2   0   0   1   2   2   2   0   1   2   1   0   1   2   2   0   1   1   2   2   2   0   &   2   2   2   1   \\
0   1   0   0   1   0   0   1   1   1   1   2   1   1   2   2   2   0   2   2   0   1   1   2   0   2   2   1   2   2   1   2   1   0   0   0   0   1   0   0   2   0   2   1   0   2   1   1   1   &   2   2   1   2   \\
0   0   1   1   0   0   1   1   1   0   1   0   2   0   0   2   0   1   2   2   1   0   2   0   0   0   2   0   0   2   2   2   2   1   1   2   1   1   1   2   0   1   1   0   0   1   1   0   1   &   2   1   2   2   \\
1   0   0   0   2   0   2   1   2   2   0   0   0   0   1   1   1   0   0   0   0   1   0   2   2   1   2   0   0   2   1   2   0   0   2   0   0   1   0   1   2   0   1   0   2   0   2   2   1   &   1   1   1   1   \\
0   1   1   1   2   2   1   0   2   2   0   0   1   2   2   1   1   2   1   0   0   1   2   0   2   1   1   1   0   2   1   2   1   0   1   2   1   0   1   1   0   2   0   2   0   0   0   0   2   &   2   2   1   1   \\
1   0   1   0   0   1   0   2   1   1   1   0   1   1   2   1   1   2   2   2   2   2   2   1   0   1   0   0   2   1   1   1   0   0   2   1   1   0   0   0   0   0   0   0   2   1   0   0   2   &   1   2   2   1   \\
1   0   1   0   1   1   1   0   0   2   2   1   1   2   1   2   2   0   1   1   1   0   0   1   2   2   1   2   0   1   2   2   1   0   0   1   2   2   2   0   0   0   0   0   1   0   0   0   2   &   2   1   1   2   \\
0   0   0   0   0   1   2   0   1   1   0   2   0   2   0   0   2   2   2   0   0   0   2   0   2   1   1   2   1   0   1   1   1   2   2   2   0   0   0   2   2   1   1   2   2   2   1   0   2   &   2   1   2   2   \\
0   1   1   0   1   2   0   1   2   1   1   1   1   0   1   2   0   1   0   2   2   0   2   2   2   0   1   2   2   0   1   1   1   2   2   0   2   0   2   1   1   0   1   2   1   0   1   1   0   &   1   1   1   1   \\
0   0   1   1   1   2   0   0   1   0   1   0   2   1   1   0   1   0   1   0   2   1   0   2   1   1   0   1   1   1   0   1   1   1   0   2   0   1   2   1   2   2   1   0   0   1   1   2   0   &   1   1   2   2   \\
0   1   1   1   1   0   0   2   2   0   0   2   1   0   0   1   2   1   1   1   2   0   1   1   2   1   2   2   0   1   0   1   2   0   1   1   1   2   1   1   2   1   2   1   1   0   1   0   1   &   1   1   1   1   \\
1   1   1   0   2   2   2   2   0   1   2   0   1   1   1   0   2   1   2   1   0   1   0   0   1   0   0   2   0   1   1   0   2   1   1   2   1   0   1   0   2   0   1   2   1   0   2   0   2   &   1   1   1   2   \\
0   0   0   0   1   0   2   2   2   1   0   0   0   1   2   1   0   0   0   0   1   2   1   0   0   2   1   0   2   0   2   0   1   1   0   2   1   2   0   1   0   1   0   2   1   0   0   2   1   &   1   1   1   2   \\
1   0   0   0   1   2   1   1   2   1   2   2   1   2   1   2   1   2   2   0   2   0   0   0   2   2   0   1   1   2   2   1   2   2   0   2   1   2   0   2   1   2   2   0   0   0   2   1   1   &   1   2   2   1   \\
0   1   0   1   2   0   2   1   1   1   0   1   2   1   2   1   2   1   0   1   0   0   0   2   2   0   0   0   0   0   0   0   1   0   1   1   2   1   0   2   0   2   2   0   1   2   0   1   2   &   1   2   2   2   \\
1   0   1   0   0   0   2   0   2   1   0   0   0   2   2   0   2   1   1   2   1   1   2   1   1   0   0   2   0   1   2   1   1   2   0   1   2   1   2   2   1   2   1   1   2   1   1   1   1   &   1   2   1   1   \\
1   1   0   0   0   2   0   1   0   0   0   1   2   0   2   2   2   1   0   0   1   2   0   1   1   1   1   1   2   1   0   1   0   2   0   2   2   0   1   1   0   0   2   0   2   0   1   0   1   &   2   2   2   1   \\
1   0   1   0   2   1   1   0   0   0   2   1   2   0   0   0   1   1   2   2   1   2   2   2   1   2   2   1   2   2   1   0   2   2   0   1   0   2   0   2   0   2   2   2   1   2   0   2   1   &   1   1   1   2   \\
0   0   0   1   0   2   0   0   2   0   2   1   1   2   0   1   2   2   2   1   1   1   0   1   0   2   0   2   2   0   0   0   0   0   1   0   2   0   0   2   0   1   1   0   0   2   0   2   0   &   1   1   2   1   \\
0   1   0   0   2   0   0   0   1   2   0   1   2   2   0   0   1   1   1   2   0   2   1   1   2   2   1   0   0   2   0   2   2   2   2   0   0   2   1   2   1   1   0   0   0   1   2   1   1   &   1   1   2   1   \\
1   1   0   1   0   1   0   0   2   0   1   0   0   0   2   0   0   0   1   0   2   0   2   0   0   0   0   2   0   0   0   2   2   1   0   0   0   2   0   2   2   0   0   2   0   1   0   0   0   &   1   2   1   2   \\
1   0   0   1   1   0   2   2   0   1   2   2   0   0   2   2   0   1   2   0   0   2   2   1   1   1   0   1   1   2   0   2   0   1   1   2   0   2   2   0   0   1   1   2   1   2   1   1   0   &   1   2   1   1   \\
0   0   1   1   2   0   0   1   0   2   2   1   1   2   2   2   2   2   1   0   2   2   1   0   1   0   0   0   2   2   1   0   2   2   2   2   2   2   2   0   2   1   1   0   1   0   0   1   1   &   2   2   2   2   \\
1   1   1   0   0   1   0   1   2   2   2   0   0   1   0   1   0   1   1   0   2   2   0   1   0   2   1   1   1   2   1   1   0   1   1   0   2   2   1   1   0   1   2   1   0   2   2   1   2   &   2   1   2   2   \\
0   1   0   0   0   2   1   2   1   2   2   2   2   0   1   0   2   0   1   2   2   1   1   1   0   1   2   1   1   0   2   0   2   0   0   1   1   2   1   1   1   0   1   2   2   2   0   1   0   &   1   2   2   2   \\
1   1   0   0   2   0   1   2   2   2   1   0   2   2   0   1   1   2   2   1   0   0   1   1   2   1   2   2   2   1   1   0   0   2   1   2   0   1   2   0   1   0   1   1   1   1   0   2   0   &   2   2   2   2   \\
0   1   0   0   0   1   0   2   0   0   1   2   1   1   1   1   2   1   0   2   0   2   2   2   1   1   0   0   1   1   2   2   0   2   2   2   0   2   2   1   1   2   0   2   0   1   2   2   2   &   2   2   1   2   \\
1   0   0   1   2   1   1   2   1   2   1   0   1   1   0   2   0   1   0   2   0   1   0   2   2   0   2   1   1   0   2   0   1   2   1   2   2   0   2   2   1   1   0   1   2   2   0   0   1   &   1   1   1   1   \\
1   0   1   0   2   0   2   2   1   2   1   1   0   0   0   0   2   0   0   1   2   2   0   2   2   2   0   1   0   0   0   2   1   2   2   1   1   0   1   0   2   2   2   2   0   1   1   0   0   &   2   2   2   1   \\
1   1   1   1   0   1   2   2   1   0   0   1   2   2   2   1   0   0   2   2   0   2   1   0   1   2   1   2   1   1   2   2   2   1   2   1   0   1   1   2   2   2   2   0   2   2   0   1   0   &   2   1   1   1   \\
0   0   0   1   1   2   1   1   1   2   2   2   0   0   2   1   1   0   2   1   1   2   2   2   2   2   1   2   1   1   1   2   1   1   1   1   2   0   1   0   1   0   0   1   0   0   1   2   1   &   1   2   2   2   \\
 \ea  $
  \\ \hline
 \et} \ec

\bc {\scriptsize  \tabcolsep=1pt \bt {c} \hline $D_{3}^{*T}$, where
$D_{3}^{*} \in \cC_{3}(72+72;2^{4}3^{45}\bullet 3^{1})$

\\\hline
$ \ba{cc}
\d_0^{*T}  & \ba{c}0    1   0   1   0   0   1   0   0   0   1   1   0   1   1   1   0   1   0   0   0   0   1   0   0   0   0   1   0   0   1   0   1   1   0   0   1   1   1   1   0   1   0   1   0   0   1   1   1   0   1   1   1   0   0   1   1   1   0   0   0   1   1   0   0   0   1   1   0   1   1   1   \\
1   1   1   1   0   0   0   1   0   0   0   0   1   1   0   0   0   1   0   0   0   1   1   1   0   0   1   1   1   1   1   0   0   0   0   0   0   1   0   0   1   0   1   0   0   0   1   1   0   1   0   1   1   1   0   0   1   1   0   1   0   1   1   1   1   0   0   0   1   1   1   1   \\
1   0   0   1   0   1   1   0   0   0   0   0   0   1   0   1   1   0   1   1   0   0   1   1   1   1   0   1   1   0   1   1   0   1   0   0   0   0   1   0   1   0   1   1   1   1   0   1   1   1   0   0   1   1   0   0   0   0   0   0   1   1   0   0   0   1   0   1   1   1   0   1   \\
1   1   1   1   0   1   0   0   0   1   0   1   1   0   1   0   1   0   0   0   1   0   1   1   1   1   1   0   1   1   0   1   0   1   1   0   1   0   0   0   0   1   0   1   0   0   0   1   1   0   0   0   1   1   1   0   1   0   1   0   0   1   0   1   0   1   0   0   0   1   0   1   \\
1   1   1   2   1   0   0   1   2   2   0   2   2   2   0   2   0   1   1   0   0   2   1   1   1   1   0   1   0   0   2   1   1   2   1   0   0   0   2   2   2   0   2   2   1   2   0   0   0   2   1   0   0   2   1   1   2   2   1   2   0   1   2   2   0   0   2   1   1   0   1   0   \\
2   1   0   0   2   0   1   2   0   0   1   2   0   2   0   0   0   0   1   2   1   0   0   2   2   2   1   1   1   1   1   2   2   0   1   0   1   0   2   2   0   1   1   1   0   0   2   2   2   0   1   2   2   2   2   0   2   0   2   1   1   1   1   1   1   2   0   0   1   2   0   1   \\
0   0   2   0   2   1   0   1   0   1   2   2   1   1   0   0   2   2   0   0   1   1   0   0   2   1   2   2   0   2   0   0   1   1   2   2   1   2   2   2   0   2   1   0   0   1   0   1   2   2   0   1   2   0   1   2   1   2   0   1   1   2   0   2   0   1   1   1   1   2   1   0   \\
2   1   2   2   0   0   2   1   1   2   0   2   1   2   1   2   0   0   0   2   2   2   2   0   2   1   2   1   1   1   1   1   1   2   2   2   1   1   1   0   1   1   0   0   0   2   2   1   0   1   2   0   1   2   0   0   1   2   0   1   1   0   0   0   0   0   0   2   1   2   0   0   \\
0   0   1   0   1   2   2   1   2   1   0   2   0   1   1   0   1   0   0   2   0   1   2   2   2   0   0   2   1   2   0   1   0   1   0   1   2   2   1   2   2   0   0   0   1   1   1   0   2   2   1   2   1   1   0   1   2   2   2   2   0   2   2   0   1   1   0   2   1   0   0   1   \\
2   0   1   0   2   0   0   0   1   2   2   0   0   1   1   1   1   2   0   0   0   2   0   0   1   2   1   2   2   0   1   1   1   2   0   2   2   0   0   2   2   2   2   1   1   2   1   0   1   1   0   1   2   0   1   0   2   1   2   1   0   2   0   1   2   1   0   2   1   2   1   0   \\
0   2   1   1   1   1   1   1   2   0   0   2   1   2   0   1   2   1   1   2   0   1   0   2   1   2   1   0   1   1   2   0   0   2   0   1   0   0   1   1   0   2   0   1   0   2   2   1   0   0   0   0   2   0   2   2   2   1   1   0   0   1   2   2   2   2   1   2   2   0   2   0   \\
2   1   1   2   0   1   0   0   1   0   2   0   0   1   1   2   1   2   1   2   0   0   0   0   0   2   2   0   2   1   2   2   0   0   1   0   1   1   1   1   1   2   1   0   2   2   2   1   0   0   2   2   1   1   1   2   2   0   1   2   0   0   1   0   1   2   0   2   2   1   0   2   \\
1   1   2   0   2   0   0   0   0   1   2   2   2   1   0   1   2   2   2   1   2   1   0   0   2   0   1   2   2   1   1   2   1   2   1   0   1   2   0   1   1   0   0   0   2   0   0   2   0   1   0   2   2   0   2   0   2   2   1   0   2   1   2   0   1   1   1   1   1   0   1   0   \\
2   1   0   0   0   1   1   1   0   2   0   2   0   1   0   2   0   0   1   2   0   2   0   2   1   1   1   0   0   0   0   1   0   1   1   2   2   2   0   0   2   0   1   1   2   0   1   0   2   1   1   1   2   1   2   2   1   2   0   2   2   1   2   1   2   1   1   0   0   2   2   2   \\
0   2   0   2   0   0   1   2   0   1   1   0   1   0   2   0   2   2   1   1   0   1   2   1   2   1   1   2   0   2   1   1   0   2   1   2   0   0   1   1   2   0   0   0   0   1   2   0   0   2   1   1   1   1   2   0   0   2   1   0   2   0   2   0   1   2   2   2   2   1   1   2   \\
0   2   2   2   1   1   0   0   1   2   1   2   0   1   1   1   0   0   2   0   1   2   1   1   2   2   1   1   0   0   2   1   1   0   0   1   2   0   2   0   1   0   2   2   0   0   0   0   0   0   2   0   2   1   2   0   2   1   2   1   2   1   0   1   1   0   1   2   2   1   2   2   \\
1   0   2   2   2   0   1   2   0   0   1   2   0   2   1   0   1   1   2   1   0   1   0   1   0   0   1   0   0   0   2   1   2   2   2   0   2   1   1   2   2   0   0   2   2   1   2   1   0   1   0   0   1   2   0   0   1   0   2   1   2   2   0   1   2   0   1   1   2   2   1   1   \\
1   2   1   0   0   0   1   1   2   2   1   1   1   1   0   2   2   2   2   2   0   0   1   2   0   0   1   2   2   1   1   0   1   0   2   0   2   1   2   0   1   0   1   1   1   1   2   0   0   0   0   1   1   2   1   2   2   2   1   2   2   0   0   2   0   0   2   1   2   0   0   0   \\
0   2   2   0   2   1   1   1   2   1   2   1   0   1   0   0   0   1   0   1   1   0   0   2   0   0   2   2   1   0   2   2   0   2   1   2   1   0   0   0   2   2   2   0   1   2   2   1   2   1   0   0   1   1   0   1   2   1   1   0   0   1   1   0   1   2   2   2   1   2   0   2   \\
1   2   2   1   1   2   2   1   1   2   0   1   0   1   2   2   1   0   0   2   1   1   0   2   2   0   1   1   2   2   2   2   1   0   0   0   2   2   1   1   0   2   2   0   1   1   0   0   0   2   1   0   0   2   2   2   1   2   0   0   1   1   0   1   0   2   0   0   0   0   1   1   \\
1   2   0   2   1   2   2   1   1   0   0   0   0   0   0   0   1   2   0   2   1   2   1   0   1   1   2   1   2   2   1   2   2   2   2   0   1   0   1   2   2   0   1   2   0   1   0   1   1   0   0   0   2   0   1   1   2   1   2   2   2   0   2   0   1   1   1   0   0   2   0   1   \\
1   0   2   0   0   1   2   2   0   1   1   1   2   1   0   2   1   0   0   2   0   1   0   1   1   2   0   2   0   1   1   2   1   1   1   1   0   2   1   2   2   2   2   2   0   2   2   0   0   0   0   1   2   2   2   1   0   0   0   0   1   1   1   2   2   0   1   2   1   0   0   2   \\
0   2   1   1   1   1   0   1   0   0   1   0   1   0   0   1   2   2   1   1   1   1   2   1   1   0   2   0   2   0   2   2   2   0   2   2   1   1   2   2   0   2   2   0   2   1   1   0   2   1   1   0   1   2   1   0   1   2   0   0   2   1   2   2   0   2   0   2   0   0   0   0   \\
1   2   0   0   1   2   1   0   2   1   2   1   2   1   1   1   1   0   1   2   0   2   2   1   2   1   1   1   1   2   1   0   1   2   0   2   1   1   2   0   2   1   1   1   2   0   0   0   0   0   2   0   2   1   0   0   2   1   2   0   0   0   0   2   2   2   0   0   2   0   0   2   \\
2   1   0   1   1   0   2   0   2   1   2   1   2   0   1   2   0   0   0   1   2   0   0   0   1   0   0   0   0   2   0   0   1   2   1   2   0   2   2   1   1   1   1   1   2   1   0   1   1   2   2   1   1   2   2   0   2   1   2   2   0   2   0   1   0   2   0   0   1   2   2   1   \\
0   0   1   1   0   2   1   2   1   2   1   2   0   1   1   0   0   1   2   0   1   1   1   1   0   2   2   1   0   0   0   2   1   2   1   1   1   2   2   0   2   2   0   1   2   0   0   2   2   2   2   0   1   1   2   0   0   0   2   2   0   0   1   1   2   1   2   2   0   0   1   0   \\
1   0   0   0   0   0   2   0   2   0   2   0   2   1   0   1   2   0   1   0   1   2   1   1   2   0   2   0   2   2   2   1   2   2   1   0   2   1   1   0   0   1   2   0   0   1   0   2   1   1   2   1   0   2   1   2   1   1   1   2   0   0   2   2   1   0   0   2   1   2   1   1   \\
2   0   1   2   0   2   2   1   0   2   0   0   2   0   0   1   2   0   2   1   1   1   0   1   0   0   0   2   1   2   0   1   1   1   0   2   2   0   1   0   1   1   0   1   2   0   2   1   2   0   1   1   1   0   1   1   2   2   2   0   0   2   1   1   2   2   0   2   0   2   2   1   \\
0   1   2   1   2   2   0   2   0   0   0   0   1   1   2   2   0   2   2   2   1   2   0   1   2   2   1   1   0   2   2   1   2   2   1   1   2   1   2   0   0   1   1   1   0   0   0   2   2   0   0   2   1   2   0   1   2   0   0   1   0   0   1   0   2   1   1   1   1   0   1   0   \\
1   1   1   1   2   2   1   1   2   2   1   1   0   1   1   1   0   0   1   1   1   2   2   0   0   0   2   1   0   0   2   1   0   1   0   0   0   2   1   2   2   2   0   2   0   0   2   2   1   1   0   1   2   0   2   2   0   0   1   1   2   2   2   1   0   2   0   2   2   0   0   0   \\
0   1   2   0   1   0   0   0   0   2   1   0   2   1   2   0   0   2   1   2   2   1   2   1   2   1   2   0   0   1   2   1   0   0   0   1   1   0   2   0   2   2   2   1   1   1   1   1   0   2   2   2   0   0   2   2   1   0   0   0   1   2   1   1   2   0   0   1   1   2   1   2   \\
2   0   2   1   1   0   1   0   2   2   2   0   0   0   2   2   0   0   0   1   2   0   1   1   0   1   0   1   1   2   1   2   1   2   1   1   0   1   2   0   0   1   2   0   1   0   2   2   2   0   2   1   1   0   2   0   1   2   1   2   2   2   1   1   2   0   1   2   1   0   0   0   \\
0   1   2   2   1   1   2   2   0   2   0   0   2   1   0   2   0   0   0   0   2   0   0   0   0   0   0   2   1   1   2   2   1   2   1   1   2   2   1   2   0   0   1   0   1   2   0   0   2   1   1   1   1   2   1   0   1   0   2   1   2   1   2   1   1   0   0   1   2   1   2   2   \\
2   1   1   2   2   1   2   0   2   0   0   1   2   2   1   1   1   0   2   2   1   0   0   0   0   1   1   1   0   2   2   1   1   2   0   1   0   2   0   1   2   2   1   0   0   1   0   2   2   1   2   2   1   0   0   1   0   2   1   0   2   2   1   0   0   0   0   2   1   1   2   0   \\
1   1   1   0   1   1   0   0   2   0   2   1   1   1   2   1   1   1   0   2   1   0   2   1   0   2   0   1   2   2   0   2   0   1   1   2   1   2   0   2   1   1   2   2   1   0   1   0   2   2   1   0   0   2   2   0   2   0   0   1   2   0   2   0   0   0   2   1   2   0   2   0   \\
1   0   0   0   0   1   1   1   0   0   0   2   2   1   2   2   0   2   0   0   1   1   1   2   1   0   0   0   0   0   2   1   2   1   1   2   1   1   0   1   2   2   1   0   1   2   0   2   2   0   2   0   1   2   2   0   1   1   2   0   2   1   1   2   1   2   2   1   2   2   0   0   \\
1   2   0   1   2   1   1   2   1   1   0   2   0   2   2   2   2   1   0   1   2   1   2   2   0   1   0   0   1   2   2   0   1   0   2   0   0   0   2   0   0   0   2   0   0   1   1   2   1   1   1   0   1   1   1   1   0   1   2   0   2   2   0   0   1   2   0   2   1   0   2   2   \\
0   2   1   1   1   0   1   2   2   2   1   0   2   2   1   2   2   2   1   2   1   2   2   1   1   0   2   2   1   0   1   2   0   1   0   0   2   0   0   1   0   2   2   0   1   0   0   2   0   1   2   0   1   0   1   1   0   0   2   0   0   2   2   1   0   0   1   1   1   2   0   1   \\
2   1   2   2   1   2   1   2   1   1   1   2   2   1   2   0   1   1   2   2   1   0   1   2   0   1   0   0   1   0   0   0   1   1   0   2   2   1   0   2   0   1   2   1   0   1   0   0   0   0   2   1   2   1   1   0   0   2   0   2   0   0   1   2   0   2   0   2   2   1   2   0   \\
1   2   2   0   0   1   1   2   1   2   0   1   1   0   2   0   0   1   2   2   1   0   1   1   1   0   1   2   0   1   0   1   1   2   0   0   1   2   1   0   0   1   2   0   2   2   1   0   2   1   2   1   2   2   0   2   2   0   0   0   2   2   0   2   1   0   2   0   1   2   1   0   \\
1   2   2   2   2   0   0   0   0   2   2   2   0   1   0   0   1   1   2   1   2   0   2   0   1   2   1   0   2   0   2   0   0   1   1   1   2   1   2   2   1   0   1   1   0   2   0   1   2   1   1   1   1   1   0   1   1   0   0   2   0   2   2   0   1   2   0   0   1   2   2   0   \\
0   1   0   0   1   1   2   1   1   2   0   1   1   0   2   1   0   1   2   0   2   0   1   2   1   0   2   1   1   2   2   1   0   2   2   0   0   1   0   2   0   0   2   1   2   1   0   2   1   2   2   2   2   0   0   1   0   0   0   1   2   0   0   2   1   2   1   1   0   1   2   2   \\
0   0   1   2   2   2   0   2   1   2   1   1   1   2   0   2   2   2   2   1   0   0   2   1   0   0   1   0   0   2   1   2   1   1   1   1   0   2   1   0   2   0   1   1   0   0   0   2   1   1   1   0   0   0   1   2   2   1   1   2   0   0   2   2   0   1   0   0   2   2   1   2   \\
2   1   2   1   0   1   2   0   0   0   1   1   1   2   0   2   1   2   0   1   1   1   1   2   2   1   0   1   2   0   0   0   1   0   2   2   0   1   0   2   2   2   2   0   0   0   2   2   2   1   2   2   0   0   0   1   2   0   2   1   2   1   0   1   1   1   0   1   0   0   1   2   \\
2   1   0   2   2   0   1   2   1   0   2   1   1   2   0   1   0   0   1   2   1   1   2   0   1   1   1   2   0   2   2   0   0   0   2   1   2   0   1   1   1   2   0   0   1   0   0   1   1   2   2   2   2   1   1   1   1   2   0   0   0   0   0   2   0   2   2   1   0   2   0   2   \\
1   0   1   0   2   2   1   0   0   0   0   1   0   0   1   2   0   1   1   1   2   2   2   0   2   0   1   1   1   2   0   0   0   2   2   1   0   2   2   0   0   2   2   2   2   1   2   1   0   0   1   0   0   1   1   2   0   0   2   1   1   1   1   2   0   1   2   1   2   2   2   1   \\
0   0   1   0   1   1   1   1   2   2   0   2   1   0   1   2   1   2   2   1   2   2   1   0   0   2   0   1   0   2   1   0   0   0   1   0   0   2   1   0   1   0   2   2   1   2   2   0   2   2   2   1   2   0   1   0   1   0   2   1   0   2   1   2   0   1   2   0   2   1   1   0   \\
2   2   0   1   1   0   0   1   2   0   2   1   2   0   1   1   2   0   0   2   1   0   2   1   2   1   1   0   1   1   2   2   0   0   1   2   2   2   0   0   1   0   0   0   1   1   0   2   0   0   1   1   2   2   2   1   0   1   2   0   2   2   1   2   1   0   2   1   1   0   2   0   \\
2   2   0   1   2   0   2   1   0   2   1   2   2   2   1   2   0   2   2   0   0   1   1   0   0   2   0   1   1   1   1   1   1   1   1   2   1   0   2   0   1   2   1   0   2   1   1   0   1   2   0   0   2   0   1   1   2   0   0   0   2   1   2   2   2   2   0   0   0   0   1   0   \ea  \\
 & \\
\zero_{1\times 72} & \ba{c} 2   1   2   2   2   1   2   2   1   2
2   1   1   1   1   2   2   1   1   1   1   1   2   2   2   1   2
1   1   2   1   1   2   1   1   1   1   1   1   1   2   2   2   2
1   1   1   1   2   2   2   1   2   2   1   2   2   1   2   2   2
1   2   1   2   2   2   2   1   1   2   1   \\ \ea

 \ea  $

  \\ \hline
 \et} \ec

\bc {\scriptsize  \tabcolsep=1pt \bt {c} \hline $D_{3}^{*T}$, where
$D_{3}^{*}  \in \cC_{3}(72+72;2^{4}3^{45}\bullet 3^{2})$

\\\hline
$ \ba{cc}
\d_0^{*T}  & \ba{c}1    1   1   0   0   1   0   1   0   0   1   1   1   0   0   1   1   1   1   0   1   1   0   0   1   0   0   0   0   1   0   0   0   0   1   0   0   1   1   0   1   0   0   0   1   0   1   1   0   0   1   0   0   1   1   1   1   1   0   1   0   1   1   0   0   1   1   0   0   1   1   0   \\
1   0   0   0   1   0   0   0   0   1   0   1   0   0   1   0   1   1   0   0   1   0   0   1   0   1   0   1   1   0   0   1   0   0   1   0   0   1   1   1   0   1   0   1   1   1   0   1   0   1   1   0   1   0   1   1   0   0   1   1   1   0   0   1   1   0   1   0   1   1   1   0   \\
0   0   1   1   0   1   1   1   0   1   0   0   1   1   1   0   1   1   0   0   1   0   0   0   0   0   1   0   1   0   1   1   1   1   1   0   1   0   0   1   1   0   1   1   1   0   1   1   1   0   0   0   1   0   0   0   0   0   1   1   0   1   1   0   0   1   0   0   0   0   1   1   \\
0   1   0   0   1   0   0   1   0   1   0   0   0   0   1   0   1   0   1   0   1   0   1   1   0   1   1   0   1   1   0   1   1   0   0   0   1   1   0   0   1   1   1   0   0   0   0   1   1   0   1   0   1   0   1   1   1   1   0   0   0   0   1   1   0   1   0   1   1   1   0   1   \\
1   0   0   1   2   2   2   1   2   0   2   0   0   2   0   1   2   1   0   1   1   0   1   2   1   2   0   0   0   2   2   1   0   0   0   0   0   2   1   2   1   0   1   1   2   0   0   2   2   0   2   1   1   0   2   2   1   1   2   2   1   0   2   1   2   2   0   1   1   0   1   1   \\
1   1   2   0   1   2   0   1   0   1   0   0   1   0   2   1   0   0   0   2   2   2   1   2   2   2   2   1   0   0   0   2   1   1   2   1   2   2   0   0   2   1   2   2   1   1   2   1   1   2   0   0   1   1   0   2   1   0   0   2   2   0   1   2   0   0   1   0   0   1   1   2   \\
2   0   2   1   0   1   0   2   1   0   0   0   2   1   0   2   1   0   0   1   0   2   0   1   2   2   2   1   2   2   2   2   2   0   0   1   1   1   2   0   2   1   0   1   1   1   2   1   0   2   2   1   0   1   2   1   0   1   0   1   0   0   1   2   0   2   0   1   2   2   0   1   \\
0   1   2   1   2   2   1   2   0   1   0   2   0   2   2   2   2   1   2   0   2   1   0   0   1   0   2   0   1   0   1   0   0   2   1   2   1   1   0   1   1   0   0   1   2   2   0   1   1   0   1   0   0   1   2   1   1   2   2   0   1   2   0   2   2   0   0   1   2   1   0   2   \\
2   2   0   2   1   1   1   2   0   1   0   1   2   2   0   1   1   0   2   2   1   1   1   1   2   0   2   1   1   0   1   0   0   0   2   0   1   0   1   0   0   0   1   1   0   0   2   2   2   0   2   1   2   2   2   1   0   0   2   0   2   2   0   0   2   0   2   1   2   1   1   1   \\
2   0   0   1   0   1   1   2   2   2   2   1   2   0   1   2   1   0   2   0   2   1   1   2   1   2   0   1   0   0   0   0   2   2   1   1   0   0   1   2   1   0   2   0   0   2   2   2   1   1   2   2   0   0   1   1   0   0   1   2   2   0   1   1   0   1   0   1   2   0   2   1   \\
2   0   0   1   1   2   0   0   0   0   2   0   1   2   2   1   2   2   1   1   1   1   1   0   2   2   0   2   1   1   0   2   1   1   1   1   2   0   1   0   1   1   1   2   2   2   0   1   1   0   0   1   2   2   0   1   0   2   0   0   2   2   0   0   1   2   0   0   0   2   2   2   \\
1   2   1   2   1   0   1   0   2   2   0   2   1   0   2   2   0   2   0   2   2   1   1   1   2   1   2   0   0   2   0   1   0   0   1   2   1   0   2   0   1   1   0   2   1   1   0   2   1   2   2   0   2   0   0   2   0   2   1   0   0   1   2   0   1   1   0   1   1   1   2   0   \\
1   1   2   1   0   2   2   2   0   2   0   0   2   1   2   1   0   0   0   2   2   0   2   1   0   1   0   1   1   1   0   0   2   1   1   1   2   2   2   0   2   0   1   1   1   1   0   2   0   2   2   0   0   0   0   1   1   2   1   0   2   1   1   1   2   1   2   2   2   0   0   0   \\
2   1   1   0   2   1   2   0   2   1   2   0   1   2   0   0   1   1   2   2   0   1   0   2   1   0   2   0   2   0   0   1   1   0   0   2   2   0   0   0   1   0   2   2   1   1   2   2   1   2   2   1   1   0   1   1   1   1   1   1   1   0   0   2   2   0   2   2   0   0   2   0   \\
2   2   0   2   2   2   0   1   1   2   0   2   1   2   0   1   1   2   0   2   2   1   0   1   2   0   2   0   2   1   2   0   2   1   0   0   1   0   2   1   0   0   1   1   2   2   0   0   0   1   2   0   1   1   1   0   1   1   1   1   0   0   0   2   1   2   0   2   0   1   2   1   \\
1   2   0   0   0   2   1   0   1   0   0   2   0   0   0   1   2   2   1   0   2   1   1   0   2   2   1   0   1   2   1   1   1   2   1   1   2   1   0   1   2   2   2   2   1   2   2   0   2   1   0   0   1   0   0   0   1   2   1   0   1   0   2   1   2   1   2   2   2   0   0   0   \\
0   1   2   0   2   2   2   2   1   0   1   1   1   0   1   1   0   1   0   1   2   0   0   1   1   2   0   2   2   1   1   1   0   1   2   0   1   0   2   2   0   0   2   0   2   1   2   2   0   2   1   0   2   1   2   1   2   0   1   1   2   0   2   0   2   1   1   2   0   0   0   0   \\
2   0   0   1   1   0   2   0   2   0   2   2   2   2   1   1   2   0   0   1   1   1   2   0   1   2   2   1   0   0   0   2   2   1   1   2   1   1   1   1   0   0   2   2   1   2   2   0   1   0   1   0   0   0   2   2   2   1   0   0   2   1   1   0   2   0   1   1   1   2   0   0   \\
1   2   0   0   1   2   2   2   0   1   2   0   1   1   0   2   0   2   0   0   1   2   1   2   0   2   2   1   1   1   0   2   2   2   2   2   2   0   0   1   0   0   0   0   1   0   1   1   0   0   2   2   1   0   0   2   0   1   1   1   1   1   2   2   0   0   2   1   1   1   2   1   \\
2   2   1   1   2   1   0   2   1   2   0   1   2   2   1   1   1   0   1   1   2   2   0   1   1   2   0   1   2   2   1   1   0   0   0   2   2   2   0   1   0   1   0   2   1   0   0   1   2   2   1   2   0   0   0   2   1   0   2   0   1   2   1   0   2   0   0   0   2   0   1   0   \\
1   1   0   2   0   0   0   1   1   0   0   0   2   2   1   0   1   2   2   2   2   1   2   2   0   2   0   1   0   1   1   1   1   2   2   0   2   1   1   0   2   0   0   1   2   2   1   1   1   0   2   1   1   2   2   0   2   0   1   2   2   1   2   0   0   0   0   1   0   2   0   1   \\
1   2   1   1   0   2   1   1   2   1   1   2   0   2   2   1   0   0   0   2   0   2   1   0   2   2   1   1   1   1   0   0   0   2   1   0   1   1   0   2   2   2   2   1   1   0   0   0   2   0   2   1   2   0   2   1   0   0   0   1   0   2   0   2   1   1   2   1   0   2   2   0   \\
0   1   1   1   2   0   0   2   0   1   2   1   1   2   2   1   0   0   1   2   1   2   2   0   0   2   2   0   2   0   2   2   0   0   0   2   1   0   2   1   0   1   1   1   0   2   0   1   2   1   0   1   0   0   2   1   2   1   1   2   1   2   2   0   0   1   1   2   0   1   2   0   \\
2   1   1   2   2   2   1   0   0   2   1   2   1   0   1   0   2   0   2   0   0   1   0   2   1   1   1   0   1   1   0   0   2   2   0   0   2   1   2   1   2   2   2   2   1   1   0   0   2   0   2   2   2   0   1   0   2   0   1   0   1   2   1   0   2   0   1   1   0   0   1   1   \\
1   1   1   2   2   1   1   1   0   0   1   1   0   0   0   2   2   0   0   1   1   0   0   0   2   1   2   1   2   1   1   1   2   0   0   2   1   0   1   1   1   0   1   2   0   0   0   2   1   0   2   2   0   2   0   0   2   2   2   2   0   2   2   1   2   2   2   0   2   1   0   1   \\
0   2   1   0   2   2   2   2   0   1   0   1   1   2   0   1   0   0   1   2   0   0   2   1   2   1   0   0   1   2   0   1   0   2   2   1   0   1   1   1   1   2   0   2   1   0   2   0   0   0   2   2   2   1   1   2   2   1   1   0   2   1   0   1   0   2   1   1   2   0   0   2   \\
1   0   0   2   0   0   0   2   2   1   0   1   2   0   0   2   2   2   2   1   0   0   0   0   0   2   2   2   2   1   1   1   2   2   1   1   2   2   0   1   2   2   0   1   0   0   0   0   2   1   1   0   1   0   0   2   1   0   2   1   1   1   1   0   2   1   1   2   1   1   2   1   \\
0   0   2   0   2   0   1   0   0   2   1   1   1   2   1   0   2   1   2   2   1   0   2   1   1   0   0   0   1   2   1   0   2   1   2   1   1   1   1   0   1   1   1   2   0   0   2   2   0   0   0   0   1   2   0   1   0   0   0   1   1   1   2   2   2   2   2   2   0   2   2   2   \\
1   2   2   2   1   2   1   1   2   0   0   2   2   1   2   0   0   2   2   1   0   2   0   2   1   2   2   0   0   0   0   1   0   1   0   1   2   1   0   1   1   1   1   0   0   1   0   2   0   1   0   2   2   1   0   1   0   1   1   0   0   0   1   0   2   2   1   2   2   2   1   2   \\
2   1   0   1   0   1   2   2   0   2   1   2   1   1   1   1   0   2   2   0   0   2   2   1   0   2   0   2   0   0   1   1   0   0   1   1   2   0   1   1   1   0   1   1   0   2   2   0   2   2   2   0   1   1   0   2   1   2   0   2   0   2   0   1   2   2   2   1   1   0   0   0   \\
2   2   0   2   0   1   1   2   2   1   1   1   2   0   2   2   2   2   0   1   0   1   0   0   0   1   0   2   2   0   0   1   1   1   0   2   0   2   0   0   0   0   2   1   1   0   0   1   1   2   0   2   1   1   0   2   0   1   0   2   1   2   1   2   1   1   1   1   2   2   2   0   \\
1   0   1   0   2   2   0   0   1   1   0   1   1   0   1   0   0   2   0   0   0   1   1   0   0   0   1   0   0   2   2   1   2   1   1   2   2   0   1   0   0   2   1   1   2   0   2   0   2   2   1   2   0   1   1   2   2   1   2   2   1   2   1   2   1   2   2   2   2   0   0   1   \\
1   0   1   2   2   0   0   2   2   2   2   0   0   2   0   1   2   1   0   0   2   1   2   1   1   2   0   1   0   0   2   0   1   1   1   1   2   0   0   2   2   1   1   2   0   0   0   1   1   2   2   1   0   2   0   0   1   2   1   1   1   0   2   2   0   1   2   1   0   1   2   0   \\
2   1   2   2   2   1   2   1   2   0   0   0   0   0   2   2   0   2   1   0   0   0   1   0   2   2   1   0   2   0   1   0   2   1   2   0   0   1   1   2   0   0   1   1   1   1   0   2   1   1   1   1   0   2   0   2   2   1   0   2   0   2   1   1   0   0   2   1   2   1   1   2   \\
1   1   0   2   2   0   2   1   2   1   0   2   2   2   1   0   1   0   1   1   0   0   0   1   1   2   1   0   2   0   0   2   2   0   1   1   0   0   1   1   2   2   1   2   1   1   2   0   2   2   2   1   0   2   0   2   2   0   0   0   1   2   2   0   0   1   1   1   1   2   0   0   \\
1   2   1   2   2   2   0   2   1   2   2   1   0   1   0   0   1   2   0   0   1   0   0   0   2   0   1   0   0   0   0   2   0   2   1   2   0   2   2   2   1   1   2   0   2   0   2   1   1   2   1   1   1   0   0   1   0   1   1   1   1   2   2   1   2   0   1   1   0   2   0   2   \\
2   1   2   2   2   0   0   0   1   2   2   1   0   1   2   2   0   2   1   2   1   0   1   2   1   2   0   1   1   0   2   1   1   2   0   1   1   2   0   1   0   0   0   1   1   0   2   0   2   0   0   0   0   0   1   1   0   2   2   1   2   1   2   0   1   2   0   2   2   1   0   1   \\
1   0   0   1   0   1   0   2   0   1   1   1   2   1   1   1   1   1   1   0   0   2   2   0   2   0   0   2   1   1   2   1   0   2   0   0   1   2   1   2   2   1   1   0   0   2   1   2   0   1   0   0   2   1   0   2   2   0   2   0   1   2   2   2   2   0   2   1   2   0   2   0   \\
2   2   2   0   1   0   1   1   0   0   1   1   1   2   1   2   1   2   2   0   2   0   2   1   2   2   0   1   1   2   0   1   2   2   0   2   0   0   2   2   1   0   2   2   1   2   0   2   1   1   1   1   0   1   0   0   0   1   2   1   0   0   0   0   0   0   1   2   1   2   0   1   \\
2   1   1   0   2   2   0   0   2   1   2   0   2   2   0   1   1   0   1   1   2   0   0   2   2   0   1   1   2   1   2   1   2   1   2   0   2   2   1   1   0   0   0   0   0   0   0   0   1   1   0   2   1   1   2   1   2   0   0   1   1   1   2   2   2   0   2   1   1   2   0   0   \\
0   0   1   1   2   1   2   1   1   2   2   0   1   0   2   2   1   2   2   1   0   1   0   0   0   2   0   2   1   2   2   1   2   0   0   0   0   0   1   0   2   2   1   2   2   1   2   2   2   2   1   0   1   2   0   1   0   1   1   1   0   0   1   2   0   0   1   0   1   1   0   2   \\
2   2   2   2   2   2   2   0   1   2   0   2   1   1   0   0   0   1   1   2   0   0   2   2   0   0   1   1   0   1   2   1   2   1   1   2   1   2   1   1   0   2   1   0   1   2   0   1   1   0   0   1   2   2   2   0   1   2   0   2   0   0   0   2   0   0   1   0   1   0   1   1   \\
2   2   1   2   1   0   2   1   0   1   0   0   2   1   0   0   1   2   0   0   0   1   2   2   1   1   2   2   1   1   2   1   1   1   0   0   1   0   2   2   2   2   1   0   2   0   2   2   0   0   0   1   0   1   1   2   0   2   0   0   2   1   2   1   2   1   0   0   1   2   1   0   \\
1   0   2   1   0   0   2   2   2   1   1   2   0   0   0   1   1   0   1   2   2   2   1   0   1   0   0   0   1   0   1   0   0   1   2   2   2   0   2   0   1   2   2   1   2   1   0   1   2   1   0   1   1   2   1   0   0   0   0   2   2   1   2   2   0   1   1   0   2   2   2   1   \\
2   1   0   1   2   0   0   0   1   2   1   0   2   2   2   2   2   0   1   1   0   1   0   0   2   1   2   2   0   2   2   2   1   1   2   0   1   0   1   0   2   1   2   1   2   0   0   1   0   0   0   0   0   0   2   0   1   1   1   2   1   1   2   1   2   1   2   0   1   2   0   1   \\
2   1   2   0   2   1   0   0   1   1   1   2   2   2   1   1   2   0   0   2   1   0   1   0   0   2   2   2   1   0   2   0   0   2   2   1   1   2   2   0   1   2   0   0   1   0   1   2   1   2   0   2   2   1   1   1   1   2   1   0   0   0   0   0   0   2   1   0   1   0   1   2   \\
0   0   1   1   1   2   2   0   2   2   0   1   2   0   0   2   1   1   0   2   0   1   1   1   0   2   1   0   1   2   1   1   0   2   2   1   0   2   2   2   1   0   1   2   0   1   2   1   2   2   1   2   0   0   2   0   1   2   0   0   1   1   0   0   2   0   1   0   0   2   2   1   \\
1   2   2   1   2   2   1   2   2   0   0   0   2   0   1   2   2   2   0   1   1   1   2   0   0   2   1   0   0   2   1   2   1   0   0   0   2   2   0   1   0   1   1   2   1   0   0   0   0   1   1   0   2   2   1   1   0   1   1   0   2   2   0   2   1   2   2   1   0   0   1   1   \\
2   0   0   0   2   2   1   2   1   1   2   2   2   1   2   0   2   1   1   2   1   1   2   0   2   1   2   1   0   2   1   0   0   2   0   2   0   1   0   2   0   0   1   2   2   0   0   2   2   0   1   1   1   1   1   1   1   0   0   2   2   1   1   2   0   1   0   1   0   0   0   0   \\  \ea  \\
 & \\
\zero_{2\times 72} & \ba{c} 2   2   2   1   1   1   2   1   1   2   1   1   2   2   2   2   2   2   1   2   2   2   2   1   1   1   1   1   2   1   1   2   1   2   1   1   1   2   2   2   1   2   1   2   1   2   2   1   2   1   2   2   1   1   1   1   2   1   1   2   1   1   2   2   2   2   2   1   1   1   1   2   \\
1   2   2   1   1   1   1   1   1   1   2   2   2   2   1   2   2
2   1   2   1   1   2   1   2   2   1   1   1   1   2   1   2   1
2   1   2   1   1   2   1   2   2   1   2   2   1   2   1   2   2
2   2   1   2   2   1   1   1   1   1   2   2   2   2   1   1   2
2   1   1   2       \\ \ea

 \ea  $

  \\ \hline
 \et} \ec

\bc {\scriptsize  \tabcolsep=1pt \bt {c} \hline $D_{3}^{*T}$, where
$D_{3}^{*}  \in \cC_{3}(72+72;2^{4}3^{45}\bullet 3^{3})$

\\\hline
$ \ba{cc}
\d_0^{*T}  & \ba{c}1    0   0   0   1   1   0   0   1   1   1   1   1   0   1   0   0   0   1   1   1   0   0   1   1   1   1   0   1   0   1   1   1   0   0   0   0   0   1   1   1   0   0   1   0   1   1   0   1   0   0   1   1   1   1   1   0   0   0   1   0   0   1   0   0   0   0   0   0   0   1   1   \\
1   0   1   0   0   1   1   1   1   1   0   0   0   0   1   0   0   1   1   1   0   0   0   0   0   1   1   1   1   0   0   0   1   0   1   0   0   1   0   0   1   1   0   0   1   1   0   1   0   1   1   0   1   0   0   1   0   0   1   1   0   1   0   0   1   1   0   1   1   1   0   1   \\
0   0   1   1   0   0   0   1   0   0   0   1   0   0   0   1   0   1   1   1   1   0   0   1   1   0   1   0   1   0   0   0   1   1   0   1   1   1   0   1   0   0   1   1   0   0   1   1   0   0   0   1   0   1   1   0   1   0   1   0   1   1   1   0   1   1   0   0   0   1   1   1   \\
0   0   0   1   1   0   1   0   0   0   1   0   0   0   0   1   1   1   1   1   0   0   1   1   0   0   1   0   1   1   1   1   1   1   1   0   1   1   0   1   1   0   1   0   0   1   0   0   1   1   0   1   0   0   1   0   1   0   0   1   1   1   0   0   0   0   1   0   1   1   0   0   \\
2   2   1   2   2   1   0   2   0   1   0   0   1   1   0   2   1   1   2   2   0   1   1   2   2   1   0   0   1   2   0   0   0   0   2   0   0   0   0   1   1   2   1   1   0   2   2   0   2   2   1   1   0   0   1   2   1   0   2   2   0   1   2   1   1   1   2   0   0   2   1   2   \\
1   2   0   2   0   2   2   0   1   1   0   2   1   0   2   1   0   2   1   1   1   1   2   1   2   1   0   0   1   0   2   0   2   1   1   1   0   2   1   0   2   0   1   1   2   0   0   1   2   1   0   2   2   0   0   0   1   0   0   0   2   1   0   2   2   2   0   1   1   2   2   2   \\
1   1   1   2   1   0   0   0   0   0   1   2   2   2   2   2   2   0   2   2   0   2   0   0   2   1   0   2   0   2   2   2   1   1   1   2   0   0   0   2   2   1   0   1   0   0   1   0   1   2   0   0   2   2   1   0   0   1   1   1   2   1   2   1   2   1   0   1   1   1   0   1   \\
0   0   2   1   0   1   2   2   2   1   1   2   1   1   0   2   1   1   0   1   2   0   1   2   1   1   2   2   0   2   1   0   0   0   0   0   0   2   0   0   2   0   1   2   0   1   0   1   1   0   2   0   2   0   2   2   1   2   1   2   2   1   2   0   0   0   2   1   2   1   1   1   \\
2   0   0   2   0   1   1   0   0   1   0   2   1   0   1   1   2   2   1   0   2   1   2   2   1   2   1   2   2   1   0   1   0   1   0   0   0   0   2   0   2   1   2   1   2   0   2   0   1   2   1   1   0   0   2   2   1   1   2   0   2   1   0   0   1   0   2   1   2   0   1   2   \\
1   0   0   0   2   2   0   0   0   0   0   2   1   1   1   0   2   0   2   2   2   1   0   0   0   2   1   2   0   1   0   1   1   0   1   1   2   1   0   0   0   2   1   1   1   1   1   1   2   2   1   2   2   0   2   2   2   0   2   0   1   1   2   2   2   0   1   2   1   0   1   2   \\
0   1   1   0   0   0   2   2   2   0   2   1   1   0   0   0   1   1   2   1   0   2   1   1   1   2   1   0   0   2   1   2   0   2   2   0   1   2   1   2   1   1   0   1   2   0   1   1   1   1   1   2   2   0   0   2   2   0   2   0   0   0   2   2   0   2   1   1   0   0   2   2   \\
0   0   2   0   1   2   2   2   1   2   2   2   1   1   2   2   2   1   1   2   0   0   1   1   0   2   1   0   0   2   0   0   1   0   2   0   1   0   1   2   0   0   0   2   1   1   1   2   1   0   0   2   2   2   0   0   1   1   1   0   2   0   0   2   1   1   1   2   1   2   0   1   \\
2   1   1   2   2   2   1   1   1   0   1   2   2   0   1   2   1   2   1   0   0   0   0   0   1   2   1   2   0   0   2   0   1   0   2   2   0   0   1   2   2   0   1   1   0   2   2   0   0   1   0   2   2   2   1   2   1   2   1   0   0   2   1   1   0   0   2   1   1   1   0   0   \\
1   0   2   1   1   2   2   1   1   2   0   2   0   2   0   0   2   0   0   1   0   1   1   0   2   0   1   0   2   1   0   1   2   2   2   2   0   1   1   1   1   0   1   1   0   1   0   0   2   2   1   0   0   2   2   2   0   1   2   0   1   0   2   0   2   1   2   2   1   0   2   1   \\
0   1   0   0   1   2   1   2   0   1   0   2   2   2   0   0   2   2   1   2   0   0   0   1   1   1   0   1   1   1   1   2   2   1   2   0   2   1   1   0   2   1   0   1   0   1   0   0   0   0   2   0   0   2   2   2   2   2   2   2   1   0   1   1   1   2   1   2   2   0   0   1   \\
1   0   1   2   0   1   2   2   2   2   2   1   2   1   0   1   0   2   1   0   2   1   1   1   2   0   2   1   0   2   0   1   2   0   2   2   1   1   1   0   0   0   1   2   0   2   0   2   0   1   2   1   0   0   2   1   1   0   0   1   0   0   1   2   1   2   2   0   1   0   0   2   \\
1   1   0   2   2   2   2   2   1   0   0   1   2   0   2   0   0   2   2   1   2   1   0   2   1   0   2   2   1   0   0   0   0   0   1   2   1   0   1   1   2   1   2   1   1   1   0   0   2   0   1   0   1   1   0   2   1   0   1   0   1   0   2   1   2   2   2   2   0   1   2   0   \\
2   0   2   2   0   0   2   1   1   0   1   1   1   2   0   1   2   0   2   2   2   0   1   0   0   1   0   0   0   0   1   0   0   2   2   1   0   0   2   0   1   1   2   0   1   1   2   2   2   0   1   1   0   1   1   2   1   2   0   2   1   0   0   2   2   2   1   1   1   2   2   1   \\
2   1   0   0   2   0   1   1   1   2   1   2   0   0   2   1   2   1   1   2   2   1   0   1   0   2   1   0   1   2   0   2   0   0   2   2   2   0   2   0   1   0   1   0   2   0   2   0   1   0   2   2   0   0   1   1   2   1   1   0   1   1   2   1   1   2   0   1   0   2   0   2   \\
1   0   2   0   1   2   2   0   2   0   1   0   1   0   2   1   1   2   1   2   0   1   1   1   2   2   2   1   0   2   0   0   0   2   1   2   2   0   2   0   2   0   2   1   0   1   0   1   0   2   1   0   1   2   1   0   0   2   2   2   2   1   1   0   1   0   2   0   0   1   1   1   \\
2   2   0   1   2   0   2   1   0   1   0   0   0   0   1   0   2   2   1   0   2   1   1   2   2   1   2   2   2   0   2   1   1   1   0   0   1   0   1   2   1   0   0   0   1   0   0   2   0   0   2   1   1   1   1   2   2   1   0   0   1   1   0   2   0   1   1   2   2   2   2   2   \\
0   1   0   1   0   2   2   2   1   0   0   0   2   2   1   1   0   0   2   1   2   0   2   1   0   2   0   1   2   2   2   0   2   1   1   2   0   2   1   1   1   1   0   2   2   1   2   0   0   1   2   1   2   1   1   0   2   1   1   2   0   0   0   2   1   0   0   1   0   0   1   2   \\
2   0   1   0   1   2   2   0   1   1   1   1   0   2   2   2   2   1   2   1   1   2   1   0   2   1   0   1   1   1   1   0   0   2   0   2   2   0   0   2   2   1   1   1   0   2   2   2   1   2   2   0   0   0   0   0   2   1   1   2   2   1   0   0   0   1   1   0   2   0   0   0   \\
1   0   1   1   1   1   1   2   1   2   0   2   1   1   0   0   2   0   1   0   1   2   2   1   2   1   2   1   0   2   0   0   2   1   0   1   0   1   2   1   0   2   1   0   1   2   0   2   1   2   0   2   0   2   2   0   2   2   0   0   1   0   0   0   1   0   1   0   2   2   2   2   \\
0   1   0   0   2   0   1   0   1   2   2   2   1   0   0   1   1   0   1   0   1   2   2   0   1   0   2   2   2   1   1   0   1   0   2   2   1   2   0   0   1   1   1   1   0   1   2   2   2   2   1   1   2   1   2   0   2   2   2   0   0   0   1   0   2   1   2   0   0   2   0   1   \\
2   1   2   2   1   1   2   0   1   2   1   1   2   0   2   1   1   1   0   2   2   2   2   1   0   1   0   0   1   2   2   0   2   0   0   2   2   1   1   2   0   0   1   0   1   1   0   0   0   1   1   2   0   0   2   0   0   1   2   0   0   0   1   2   1   0   2   2   1   2   1   0   \\
2   0   1   2   0   1   1   1   2   1   0   0   0   0   1   1   2   0   1   0   2   1   2   2   2   2   2   1   1   0   1   1   2   2   1   0   0   0   0   1   1   2   1   0   1   1   1   1   0   2   2   0   2   2   0   0   2   2   1   2   0   2   2   1   0   0   2   0   0   2   0   1   \\
0   0   0   1   2   1   2   2   1   0   0   2   0   0   2   1   1   1   0   1   2   1   2   0   0   2   2   0   2   0   1   2   1   2   2   2   0   2   0   1   0   1   0   1   1   1   2   2   2   2   2   0   0   0   2   0   1   1   0   1   2   1   1   2   2   0   1   1   0   0   1   1   \\
2   2   2   2   2   1   1   0   0   0   1   2   1   2   1   1   1   2   0   1   1   0   2   1   1   0   2   0   1   2   0   1   1   2   1   0   0   1   0   0   2   0   1   2   2   0   0   0   2   1   2   2   1   2   0   0   0   1   2   2   0   2   0   1   0   1   1   2   1   0   0   2   \\
1   0   0   2   0   1   1   0   2   0   1   0   1   2   2   2   2   2   0   0   0   1   0   1   2   2   2   2   1   0   1   0   2   0   2   1   2   1   1   1   1   2   1   0   1   1   2   0   2   0   2   2   2   0   1   2   1   1   2   0   2   0   1   0   0   1   0   1   0   1   2   0   \\
1   0   2   1   1   2   0   1   1   2   1   0   0   1   1   0   0   1   0   1   1   2   0   2   2   2   0   0   0   2   2   1   0   0   2   2   2   1   2   1   2   1   1   1   0   0   0   1   2   1   1   1   2   2   0   0   2   2   1   0   0   2   2   2   1   0   0   1   2   0   0   2   \\
1   0   0   2   0   2   1   1   2   2   1   1   1   1   1   1   1   1   0   0   2   2   2   1   0   0   0   2   1   0   0   2   2   1   2   0   2   2   2   0   0   1   1   2   0   1   2   1   2   0   0   1   0   1   0   1   0   2   2   2   2   2   2   1   1   0   0   2   0   0   0   1   \\
1   1   0   0   2   2   2   2   0   1   0   0   1   0   2   2   1   0   0   1   0   1   1   2   0   2   1   1   2   0   2   1   1   0   0   2   2   0   0   1   0   1   1   2   0   0   2   0   0   2   2   1   0   2   0   0   1   0   1   2   1   2   1   1   1   2   2   2   1   2   2   1   \\
0   1   2   0   0   0   0   1   2   2   0   1   2   1   1   1   2   0   0   1   2   1   0   2   1   1   2   1   2   0   2   1   0   0   2   2   1   1   0   0   1   0   1   0   2   2   2   0   2   2   0   0   2   0   1   0   0   2   2   2   2   1   1   1   0   1   1   1   1   0   2   2   \\
0   0   1   1   2   1   1   0   1   0   1   1   1   2   2   0   1   0   0   2   2   0   2   0   2   0   0   1   2   2   0   1   0   1   0   0   1   2   2   1   0   0   2   1   1   2   2   2   1   1   2   2   2   0   1   1   2   2   2   0   0   1   2   0   0   2   0   2   1   1   0   1   \\
0   1   0   2   1   1   0   1   1   2   0   0   2   0   2   2   1   0   0   2   0   1   1   1   2   1   1   2   1   0   1   1   1   1   0   1   2   2   0   0   1   0   2   1   2   0   0   2   2   0   2   1   0   2   2   1   0   1   1   2   0   0   2   2   2   2   2   0   0   1   0   2   \\
2   1   1   2   1   0   1   0   2   1   2   0   2   0   0   1   0   2   2   1   0   2   2   0   0   0   2   0   0   1   2   0   1   1   0   0   2   0   0   0   1   0   1   1   1   0   2   1   1   2   1   2   1   2   1   1   1   1   2   2   1   0   2   2   2   2   0   2   2   0   1   0   \\
2   0   2   0   0   2   0   1   2   0   0   1   1   0   2   2   0   1   2   0   0   0   0   2   1   2   2   1   1   1   2   0   1   1   1   0   1   2   1   1   0   1   1   1   1   0   2   0   2   1   2   1   0   2   2   2   2   1   0   0   0   2   1   2   0   0   2   1   2   2   1   0   \\
0   0   2   0   1   0   2   0   1   1   0   0   2   1   0   1   2   2   0   2   2   2   0   0   0   2   2   2   1   0   0   0   1   1   0   0   1   1   0   1   1   2   2   2   2   2   1   1   2   2   2   1   1   2   2   1   0   1   0   0   2   0   1   2   0   2   1   1   1   1   0   1   \\
2   1   0   1   2   2   2   1   1   2   0   0   0   1   0   0   0   1   2   2   1   1   2   0   2   0   0   2   1   0   1   0   0   2   2   2   2   0   1   0   0   1   2   0   1   1   2   0   1   1   2   1   1   2   2   0   0   2   0   2   1   1   1   0   0   0   1   2   1   1   2   2   \\
1   2   2   2   2   0   1   1   0   2   2   1   0   1   1   2   2   0   2   0   2   0   0   1   1   2   2   0   2   2   1   0   1   0   2   0   2   2   0   1   1   0   0   1   0   1   0   2   0   1   0   2   1   0   0   2   0   1   1   0   2   1   2   0   1   1   1   2   1   0   2   1   \\
2   2   1   1   2   0   0   0   2   1   0   1   2   1   2   1   2   2   1   2   0   2   0   0   2   0   1   0   1   0   2   1   0   1   2   2   2   2   1   0   0   1   1   2   0   2   1   0   0   0   1   2   1   0   1   2   0   2   0   2   0   0   0   1   0   1   1   1   2   1   2   1   \\
2   1   2   0   1   2   0   2   1   0   2   0   1   2   1   2   1   2   2   2   1   2   2   0   2   0   2   0   2   1   2   0   0   1   2   1   0   1   1   1   0   2   0   0   2   0   1   0   2   1   0   0   1   1   1   2   1   2   0   0   1   1   2   0   1   2   1   0   0   0   0   1   \\
1   1   2   2   0   0   1   2   0   0   2   1   1   0   2   1   0   0   0   2   2   0   1   2   0   1   1   2   1   0   2   1   2   0   0   2   2   0   1   1   0   1   0   2   1   2   0   1   2   1   2   0   1   0   1   2   2   1   1   2   0   1   1   2   2   0   0   0   2   0   1   2   \\
1   0   2   1   2   0   1   0   2   1   1   2   2   1   2   0   1   1   2   0   0   0   2   2   2   0   0   1   2   0   1   1   0   0   2   1   0   0   1   0   2   2   2   1   2   0   0   2   1   1   0   1   1   2   1   0   1   0   1   2   2   0   2   2   0   1   1   0   2   1   0   2   \\
1   0   2   1   2   2   0   2   0   1   1   1   2   0   0   2   0   0   0   1   2   2   0   1   1   0   1   2   0   0   2   2   0   2   2   1   1   1   2   0   1   1   0   1   1   0   1   2   0   2   2   2   1   0   0   1   1   1   2   0   1   1   0   1   0   2   2   2   0   2   2   0   \\
0   0   0   2   1   1   0   2   0   2   2   1   1   2   1   2   1   1   1   0   2   2   0   0   2   0   1   0   1   2   2   0   2   0   0   1   0   1   2   1   2   1   2   0   2   0   2   1   0   2   2   0   0   1   0   2   0   1   1   1   0   1   0   2   2   0   1   2   1   2   1   1   \\
1   2   0   0   0   0   0   1   1   2   0   2   1   1   0   2   1   2   2   0   0   0   2   2   1   2   0   1   1   1   0   1   2   0   2   2   2   2   2   2   0   0   2   0   0   2   0   1   2   2   2   1   2   2   0   1   1   1   1   1   1   1   0   1   1   2   1   1   0   0   0   0   \\
1   0   0   2   1   0   2   1   1   1   0   0   2   1   1   2   2   1   2   2   0   0   1   1   0   1   0   0   0   1   1   2   0   1   0   2   1   2   2   2   2   1   1   2   0   0   0   2   2   2   1   0   1   0   1   2   0   2   1   1   1   0   2   2   0   0   0   2   0   2   1   2   \\  \ea  \\
 & \\
\zero_{3\times 72} & \ba{c} 1   1   1   1   2   2   1   1   1   2   1   1   1   1   1   2   2   1   1   1   2   1   2   2   2   1   2   1   1   2   2   1   1   2   2   2   1   1   2   2   1   2   2   2   2   2   1   2   1   1   1   2   2   1   2   2   1   1   2   2   1   2   2   2   2   2   1   2   2   1   1   1   \\
2   1   2   2   2   1   2   2   2   1   1   1   1   1   1   2   2   2   1   1   1   2   1   1   2   2   1   2   2   2   2   2   1   1   1   1   2   2   2   1   1   1   1   2   2   1   2   2   2   1   1   1   2   2   2   1   2   1   1   2   1   1   1   1   2   2   1   2   2   1   2   1   \\
1   2   2   1   2   1   2   1   1   2   1   1   2   2   2   2   1
2   1   1   2   1   2   1   2   2   2   1   1   1   2   2   1   1
2   1   2   2   1   2   1   1   2   2   2   2   1   1   1   2   1
1   1   2   1   2   2   2   1   2   2   1   1   1   2   1   1   2
1   2   1   2   \\ \ea

 \ea  $

  \\ \hline
 \et} \ec

\bc {\scriptsize  \tabcolsep=1pt \bt {c} \hline $D_{3}^{*T}$, where
$D_{3}^{*}  \in \cC_{3}(72+72;2^{4}3^{45}\bullet 3^{4})$

\\\hline
$ \ba{cc}
\d_0^{*T}  & \ba{c}1    0   1   1   0   0   0   0   0   0   0   0   1   1   1   1   1   0   1   0   0   0   1   1   0   1   0   0   1   0   1   1   1   1   1   1   1   1   0   0   1   0   1   1   0   1   1   0   1   0   1   0   1   0   1   0   1   0   0   0   1   1   0   1   0   1   0   0   0   1   0   0   \\
1   0   1   0   1   0   1   1   0   1   0   1   1   0   0   0   0   1   0   1   1   1   0   0   0   1   0   1   1   0   0   0   1   1   0   1   1   1   0   0   0   1   1   1   1   1   1   1   0   1   0   0   0   0   0   0   1   0   0   0   1   1   1   1   1   1   0   0   0   0   0   1   \\
0   1   1   0   0   0   0   1   0   0   1   1   1   0   0   0   1   1   1   0   1   1   0   1   1   0   1   0   1   1   0   0   0   0   1   0   1   1   1   0   0   0   1   1   0   0   0   0   0   0   1   0   1   1   1   1   1   1   1   0   1   0   0   0   1   1   0   0   1   1   1   0   \\
1   1   0   1   0   0   1   0   0   1   0   1   1   0   0   1   0   1   0   0   1   0   0   0   1   0   1   0   1   0   0   1   1   0   0   0   0   0   0   1   1   0   1   1   1   0   1   1   1   0   1   0   0   1   1   1   1   1   1   0   1   0   0   0   0   1   1   0   1   0   1   1   \\
2   1   0   0   2   0   2   1   0   1   1   2   2   2   1   0   1   0   1   1   0   2   1   2   1   2   2   2   1   0   1   2   1   2   2   0   2   2   1   2   2   1   0   2   2   0   2   1   0   1   0   0   2   1   1   0   1   0   0   0   0   1   2   0   0   0   0   1   2   1   1   0   \\
1   2   1   1   1   2   0   2   1   2   1   1   0   0   2   2   0   0   2   1   2   2   2   0   1   1   0   0   2   0   2   0   2   1   1   0   0   0   0   1   0   0   0   2   1   2   0   1   1   0   1   0   0   1   2   1   2   1   0   2   1   1   0   1   2   2   2   2   2   1   0   2   \\
2   1   0   2   1   1   0   1   1   2   2   1   0   0   0   2   0   0   1   0   0   2   0   1   0   1   0   1   1   2   0   2   2   2   2   1   2   1   2   0   1   2   0   2   0   2   1   1   1   0   0   2   1   1   2   2   0   2   0   1   1   2   1   0   0   2   0   1   1   0   2   2   \\
0   1   0   1   1   2   2   0   0   0   2   2   1   0   1   0   1   0   1   2   2   0   2   2   1   1   0   1   1   1   0   1   0   1   2   2   1   2   2   1   0   0   1   2   2   1   2   1   0   0   1   2   0   0   1   0   2   2   2   1   2   2   0   2   0   1   2   1   2   0   0   0   \\
0   1   1   2   0   1   2   2   2   2   1   2   1   2   2   1   0   2   2   1   1   0   0   1   2   2   2   1   0   1   1   1   0   0   2   0   2   1   0   1   0   1   2   2   1   2   0   1   1   0   0   2   0   0   2   2   1   2   2   1   0   0   1   1   0   1   0   0   0   2   0   0   \\
1   1   2   1   2   1   0   2   1   1   2   0   2   0   1   2   1   1   1   1   1   1   0   1   0   1   2   0   1   1   2   1   0   1   2   0   2   0   1   2   2   0   0   1   1   2   2   2   2   2   0   2   1   0   2   0   0   0   2   0   2   0   2   0   0   0   1   0   2   2   0   0   \\
2   1   1   1   1   1   2   2   2   0   0   2   0   1   1   2   2   1   2   0   2   1   0   0   0   0   2   0   2   0   2   1   2   0   1   0   2   2   0   0   0   1   1   1   0   0   1   1   1   0   1   1   2   2   1   0   0   0   0   2   0   2   2   1   1   2   2   2   1   2   1   0   \\
1   1   2   2   1   0   1   0   1   1   1   2   2   1   1   0   2   0   2   1   1   2   0   0   0   1   1   2   2   2   2   2   0   0   0   0   0   0   0   0   1   2   1   0   1   2   1   0   1   2   1   2   0   2   0   0   2   0   2   1   2   2   0   0   1   1   2   2   1   1   0   2   \\
0   0   2   1   0   2   1   2   1   0   2   0   2   1   2   2   0   2   1   0   2   1   2   2   2   2   1   1   1   0   0   2   2   0   1   1   1   0   1   1   2   2   2   0   0   1   1   1   2   1   1   0   0   1   0   2   2   0   0   2   0   1   1   0   1   0   1   2   0   0   2   0   \\
1   0   2   2   1   1   2   2   0   2   1   2   2   1   0   2   1   0   2   1   2   1   2   2   1   0   1   1   1   0   0   1   0   2   1   0   0   1   2   0   1   0   2   0   2   0   0   2   0   0   0   2   1   0   0   0   1   1   0   2   2   0   2   1   0   1   1   1   1   2   2   2   \\
2   2   1   0   1   2   2   2   1   1   1   1   1   0   2   0   0   0   2   2   0   1   1   0   1   0   1   2   2   2   1   2   2   0   2   2   0   0   1   1   1   0   2   1   1   1   0   0   0   2   0   2   1   0   0   1   1   0   2   0   2   1   1   0   0   2   0   2   1   2   0   2   \\
2   0   2   0   0   2   2   0   1   1   0   2   1   2   2   1   1   0   1   1   0   1   1   0   0   0   0   0   0   2   0   0   1   1   1   2   1   0   2   2   1   0   1   0   2   0   2   2   0   1   2   0   2   1   2   0   2   1   1   2   1   2   1   1   2   2   1   0   1   2   2   0   \\
1   0   2   1   0   2   1   1   0   2   2   0   2   2   2   0   0   0   1   0   1   1   2   1   0   0   1   2   2   1   1   0   0   1   2   0   2   0   0   2   2   1   0   2   2   1   1   1   0   1   1   0   2   2   0   1   1   1   0   0   0   0   2   2   2   2   2   1   1   1   2   0   \\
2   0   0   0   0   0   2   2   2   0   2   2   0   2   1   1   2   1   1   0   0   1   1   0   2   1   1   0   0   2   2   0   0   2   1   1   1   1   0   1   2   1   2   2   1   2   1   1   2   2   0   1   0   2   1   0   1   0   0   2   2   0   0   2   0   0   2   1   1   1   2   1   \\
2   1   2   1   0   0   1   2   1   2   1   0   0   0   0   2   0   0   2   2   0   0   0   2   0   2   2   1   0   0   2   2   0   0   1   0   1   1   1   0   1   0   2   1   1   2   2   1   0   1   2   2   2   2   0   1   1   2   1   2   2   0   0   1   1   1   2   1   0   1   1   2   \\
2   2   2   0   0   0   2   0   2   2   1   1   1   2   0   2   2   0   1   0   2   1   0   1   2   2   0   2   2   1   2   1   1   0   2   2   1   0   2   0   0   1   0   1   1   1   1   2   1   0   0   2   1   1   0   0   2   2   1   1   0   1   0   1   0   2   1   1   0   0   2   0   \\
0   2   0   1   2   0   2   1   1   1   1   0   2   1   2   1   1   2   0   1   1   0   0   0   1   0   2   0   1   0   2   2   0   2   2   2   0   1   1   0   1   1   0   2   1   2   2   2   0   0   1   1   0   1   1   2   2   0   0   2   2   0   2   0   2   1   2   2   0   1   0   1   \\
1   2   2   1   0   0   2   1   0   0   1   2   1   1   0   0   0   1   2   2   2   2   1   0   1   1   2   1   0   2   2   2   1   2   1   2   0   1   0   2   0   0   0   0   0   0   0   2   1   2   1   1   2   0   2   0   1   2   0   2   1   2   1   0   0   2   1   0   1   1   1   2   \\
0   2   1   0   1   1   1   0   2   2   1   2   0   0   0   2   1   0   0   2   2   2   0   1   0   2   2   1   2   0   2   1   0   2   2   2   2   0   0   2   2   1   1   0   0   2   0   1   2   0   2   2   1   1   1   2   1   1   1   1   0   1   0   1   1   0   0   0   1   1   2   0   \\
2   2   1   0   0   0   1   0   0   2   1   1   2   2   1   0   0   2   2   0   1   2   2   0   1   1   0   0   2   0   1   2   0   1   0   1   1   2   2   2   1   1   2   1   0   0   0   2   1   0   2   1   1   2   2   0   2   1   2   2   1   0   2   2   0   1   1   0   1   1   0   0   \\
0   2   1   0   1   0   1   2   0   1   1   2   2   2   1   2   2   2   1   0   1   2   2   2   0   2   0   2   0   1   0   0   1   0   1   0   0   1   1   1   2   0   0   1   0   1   2   2   1   1   1   0   1   1   2   2   0   2   2   1   0   2   0   0   0   2   2   0   0   2   1   1   \\
0   1   0   2   0   2   1   0   1   1   0   0   1   0   2   0   2   1   0   0   1   2   1   2   2   2   1   1   0   0   2   1   0   1   1   0   1   2   1   2   2   1   0   0   0   2   2   0   1   0   1   0   1   0   0   0   2   2   2   2   1   1   2   1   2   2   2   1   0   1   2   2   \\
1   2   1   1   2   0   0   2   1   1   0   1   2   0   0   2   1   0   0   0   0   2   1   0   2   1   1   1   2   2   0   0   1   2   2   2   1   1   2   1   0   2   2   2   2   1   2   1   2   0   2   0   0   0   0   1   1   2   1   1   0   1   0   0   1   0   2   0   2   2   0   1   \\
1   2   2   2   0   2   1   2   2   0   0   2   1   2   1   1   0   2   0   0   1   0   1   2   0   0   1   1   1   1   2   2   1   1   1   2   2   1   0   2   0   1   2   0   2   1   0   0   0   0   1   0   1   2   1   2   0   0   1   1   2   2   0   0   1   2   2   0   1   2   0   0   \\
2   2   2   2   2   1   0   2   1   2   1   2   1   2   2   0   0   1   0   1   1   0   0   2   2   0   1   0   1   2   1   2   1   1   1   2   1   0   1   0   0   0   0   0   0   2   1   1   2   0   0   0   2   1   2   0   2   1   0   1   1   1   2   0   1   0   0   2   2   0   2   1   \\
0   1   1   1   2   0   2   2   0   2   0   0   0   0   0   0   2   2   2   1   1   0   1   0   0   2   0   1   2   2   1   2   1   2   2   0   1   1   1   1   1   0   1   0   2   1   1   0   2   2   2   2   1   1   2   0   1   1   0   1   0   0   2   2   2   2   0   1   0   1   2   0   \\
2   2   2   0   1   2   0   1   1   1   1   1   1   2   0   2   1   0   1   0   2   0   0   1   2   2   2   0   2   2   0   1   1   0   1   0   1   0   2   1   0   1   0   0   2   2   2   1   0   2   1   2   1   0   0   0   1   2   2   2   1   2   1   2   1   0   1   0   0   2   0   0   \\
2   2   0   1   0   2   0   2   1   0   0   1   1   2   1   0   2   1   1   1   2   0   1   2   0   0   0   0   0   0   1   0   2   2   1   2   1   0   1   2   2   1   1   0   2   2   0   2   1   1   2   2   1   2   0   1   1   2   2   1   0   0   0   0   2   1   0   2   1   2   0   1   \\
1   2   2   0   1   1   0   1   1   0   2   2   2   1   0   1   2   1   1   1   0   1   0   0   0   2   2   2   0   0   2   1   2   2   2   2   0   1   1   1   1   1   2   1   2   1   2   1   0   0   1   0   0   2   2   2   0   0   1   2   1   0   2   0   0   1   0   2   0   2   0   0   \\
2   0   2   1   2   1   0   2   0   0   2   1   0   2   2   0   1   2   2   1   0   0   2   2   0   1   1   2   2   2   1   0   1   2   1   0   2   1   1   2   0   0   1   1   1   0   2   1   2   1   1   1   0   0   1   2   0   2   0   2   1   1   0   0   0   2   2   1   0   0   1   0   \\
0   2   0   1   2   0   1   1   0   2   1   2   2   2   2   0   2   2   1   1   1   2   1   2   1   1   2   2   2   1   2   1   1   0   2   0   1   0   0   2   0   1   2   0   0   1   1   1   0   0   2   2   0   1   0   0   0   2   1   2   1   0   0   2   0   1   0   0   1   0   1   2   \\
0   0   2   0   2   1   1   1   2   1   1   0   1   1   2   2   2   2   0   1   0   0   1   0   0   2   2   2   1   2   0   1   2   2   1   1   0   1   0   2   0   0   1   2   1   0   1   0   1   0   0   2   1   2   0   0   2   1   1   1   2   0   0   0   1   2   0   2   2   1   2   2   \\
2   1   2   0   1   0   1   0   2   0   0   1   2   1   2   1   2   1   2   1   2   2   0   1   2   0   0   2   0   1   1   0   1   0   1   0   2   1   1   0   0   0   0   0   2   1   2   2   2   0   2   1   0   1   0   2   1   0   1   0   0   0   1   1   2   2   1   2   1   2   2   2   \\
1   0   2   2   2   1   0   2   2   0   2   0   1   0   0   0   1   0   0   0   2   2   1   0   1   2   2   0   0   1   1   1   2   0   1   2   2   1   0   0   0   1   1   1   2   0   2   1   1   1   2   1   2   1   2   1   2   1   2   1   0   0   1   2   1   2   2   2   0   0   0   0   \\
2   1   2   0   1   2   2   0   2   2   2   0   0   1   0   0   2   2   0   0   1   1   1   1   0   0   0   2   1   2   2   0   2   0   1   0   2   1   2   0   1   1   2   2   1   0   2   2   1   1   2   2   0   2   1   1   1   1   1   1   1   2   0   0   0   1   0   0   2   0   0   1   \\
2   0   0   0   2   1   2   2   1   1   0   0   2   0   1   2   0   1   0   0   0   1   1   2   1   2   0   1   0   1   2   1   1   0   2   2   1   1   2   2   2   1   2   0   1   1   0   1   0   1   1   2   0   2   0   2   2   1   0   2   0   2   0   2   0   2   1   0   1   1   2   0   \\
2   2   0   1   0   1   2   1   0   1   2   1   2   0   1   1   2   0   0   1   1   1   0   2   0   0   1   1   1   1   0   0   0   0   1   0   2   2   0   2   1   1   0   1   2   2   1   0   2   2   2   2   2   2   2   2   1   0   0   1   2   2   1   1   0   0   2   0   0   1   0   2   \\
2   0   0   0   1   2   0   1   1   0   0   2   2   0   2   0   2   1   0   2   2   0   0   0   1   2   1   2   1   1   1   2   0   0   2   2   1   1   2   0   1   0   1   0   1   2   0   1   2   1   2   0   1   2   1   2   0   2   1   1   0   0   2   1   0   2   2   1   2   0   1   1   \\
1   1   0   0   2   0   0   2   1   1   1   1   2   2   0   1   0   1   2   0   2   0   1   1   0   2   1   2   0   2   2   2   0   0   0   1   2   0   2   0   2   1   1   2   2   0   2   0   0   0   1   1   0   2   2   0   2   2   0   1   0   1   1   2   2   1   1   1   1   0   1   2   \\
0   1   2   1   0   1   1   2   1   0   1   1   1   2   2   1   2   1   1   2   0   1   0   0   1   0   2   2   0   0   1   0   1   1   2   2   0   1   2   2   0   2   0   2   0   2   1   0   2   1   1   0   0   0   1   0   0   1   2   2   2   1   1   2   2   2   0   0   2   0   2   0   \\
2   0   1   2   2   0   1   1   2   0   2   1   2   2   1   2   1   1   2   1   2   1   1   0   1   0   0   0   2   0   0   1   0   0   0   2   1   2   1   1   1   0   2   2   0   2   1   0   0   2   0   1   0   1   1   2   0   2   0   1   0   2   0   2   1   0   1   1   2   2   0   2   \\
2   1   1   2   0   0   0   2   2   0   2   0   0   0   2   1   2   1   1   1   2   2   0   0   1   0   0   1   1   2   2   2   0   1   2   2   0   2   0   1   2   1   2   0   1   0   1   2   1   0   1   1   1   1   0   1   1   1   2   0   0   2   1   1   2   0   2   0   0   0   2   2   \\
0   1   0   1   2   1   0   1   0   1   1   2   1   0   2   0   1   1   1   0   0   2   0   0   0   0   1   2   0   1   2   0   2   2   0   1   1   2   2   2   1   2   2   0   2   2   0   1   1   0   0   2   2   0   2   1   1   2   0   1   0   1   0   1   2   1   2   0   1   2   2   2   \\
1   0   0   2   0   0   1   2   1   2   1   1   1   0   1   2   2   2   1   1   0   0   0   1   1   2   0   1   1   0   0   0   1   0   2   0   0   0   0   2   1   1   1   2   1   0   0   1   1   2   2   1   2   2   1   0   2   2   1   2   0   2   2   0   1   2   0   2   2   2   0   2   \\
2   2   2   2   2   2   2   1   1   0   2   2   1   0   1   2   2   2   2   0   0   0   2   1   0   0   1   0   1   0   0   1   1   1   2   1   2   1   1   1   2   2   0   0   2   2   1   1   1   1   0   1   0   1   0   0   1   2   0   0   0   0   0   1   0   1   0   2   0   1   2   2   \\  \ea  \\
 & \\
\zero_{4\times 72} & \ba{c} 2   1   1   2   2   1   2   1   1   1   2   2   1   1   2   2   2   2   1   1   1   2   2   2   1   1   2   2   2   1   2   2   2   1   1   2   2   1   2   1   1   1   1   2   2   1   1   1   2   2   2   1   1   1   1   1   2   2   2   2   2   1   2   1   2   1   1   1   1   2   2   1   \\
2   1   2   2   2   2   1   1   1   2   1   1   2   1   2   1   1   2   2   2   2   1   1   2   2   1   2   2   1   1   1   1   1   2   2   2   2   2   1   1   2   2   1   2   1   1   1   1   1   2   2   2   1   2   2   1   2   1   2   2   1   2   1   1   1   1   2   1   1   2   2   1   \\
1   2   2   1   1   1   2   2   2   1   1   1   2   2   2   1   2   1   2   2   1   1   1   1   1   2   2   1   2   1   1   2   2   2   1   1   2   2   1   1   2   2   1   2   2   1   1   1   1   1   2   2   1   1   1   2   2   2   2   2   1   2   2   1   1   1   2   1   2   1   2   2   \\
2   1   2   2   1   1   2   1   2   2   1   1   2   1   1   1   1
2   2   2   1   2   1   2   1   2   1   1   2   2   2   1   2   1
2   2   1   2   1   2   2   2   1   1   1   1   2   1   1   1   1
1   1   1   2   2   1   2   1   2   2   1   1   2   2   1   2   2
2   2   2   1    \ea

 \ea  $

  \\ \hline
 \et} \ec

\bc {\scriptsize  \tabcolsep=1pt \bt {c|c} \hline $D_{3B}^* \in
\cC_{3B}(72+24;2^{4}3^{45}\bullet 3^{0}\bullet \B)$ & $D_{3}^* \in
\cC_{3B}(72+24;2^{4}3^{45}\bullet 3^{1}\bullet \B)$

\\\hline
$ \ba{cc}
\d_0^*  &\zero_{72\times 2}   \\
0   0   0   0   1   0   0   0   1   2   1   0   1   0   0   1   2   2   1   1   0   0   1   2   2   2   0   0   0   1   2   0   1   0   1   0   1   1   2   1   0   1   0   1   1   1   2   0   2   &   1   0   \\
1   0   1   0   0   1   2   0   0   2   1   0   1   0   1   2   2   0   2   0   2   1   0   0   1   0   2   2   0   2   2   2   1   2   0   0   0   1   0   2   2   2   0   2   2   1   0   0   0   &   1   0   \\
1   0   1   1   0   2   2   2   2   0   2   0   2   2   0   0   0   0   1   2   1   1   2   2   1   1   1   1   0   1   2   1   0   0   2   1   0   1   1   1   2   0   1   1   2   0   1   2   2   &   1   0   \\
0   1   1   1   2   0   0   2   0   0   0   1   1   2   1   2   2   1   1   2   0   0   1   1   2   2   1   2   0   1   0   2   2   2   2   1   2   0   2   0   2   0   2   0   2   0   2   1   2   &   1   0   \\
1   1   1   1   1   2   0   1   1   1   2   2   1   1   2   2   2   0   0   0   1   2   0   2   0   0   0   1   2   2   1   2   0   0   0   0   2   0   0   0   2   0   0   0   0   1   1   1   1   &   1   0   \\
0   0   0   0   1   2   2   0   0   1   0   2   2   0   1   2   0   1   0   1   0   2   2   0   1   0   1   0   1   0   1   1   1   1   1   2   2   0   0   1   0   0   0   2   1   2   2   1   2   &   1   0   \\
0   1   0   0   2   0   1   2   1   1   1   0   2   2   2   1   1   0   0   1   0   1   0   2   1   1   2   1   1   2   1   0   1   2   0   1   2   1   2   0   1   2   2   1   0   2   0   1   0   &   1   0   \\
0   0   1   1   1   0   1   0   1   0   2   2   2   2   2   0   0   2   1   1   1   0   1   0   0   0   0   2   2   0   0   1   2   0   0   1   1   1   0   2   2   2   2   0   1   2   0   1   1   &   1   0   \\
1   0   1   1   2   2   1   2   0   2   2   1   1   0   0   2   2   1   1   0   2   1   0   1   1   2   0   1   1   2   0   0   0   1   1   2   1   2   1   0   1   1   2   2   0   0   0   0   0   &   1   0   \\
1   1   0   0   0   0   0   0   2   0   0   1   0   2   2   1   1   2   2   0   2   2   0   1   0   1   0   0   1   1   0   1   0   0   2   0   0   1   1   2   0   1   2   0   2   1   1   1   0   &   1   0   \\
1   1   0   0   2   1   2   2   2   2   2   2   2   2   0   0   2   1   2   1   0   2   1   1   2   1   0   2   2   0   1   0   2   2   1   1   0   2   1   0   1   0   1   2   1   1   1   0   2   &   1   0   \\
0   1   0   0   0   1   1   1   1   2   0   2   0   1   0   0   0   0   1   0   1   2   1   0   0   2   1   1   1   2   2   2   2   1   2   1   0   2   1   1   1   1   0   2   0   2   2   0   2   &   1   0   \\
1   1   1   0   1   1   2   1   0   2   2   2   1   1   2   1   1   2   2   1   2   0   2   0   2   0   1   0   1   1   2   1   1   1   1   2   2   2   2   0   1   1   2   1   2   0   1   2   2   &   1   1   \\
0   0   0   1   1   0   1   1   2   1   1   1   0   1   1   2   0   0   0   0   2   0   1   2   2   1   2   2   2   0   0   0   1   1   1   1   0   0   2   1   1   1   1   1   1   0   1   2   0   &   1   1   \\
0   0   1   1   0   1   0   0   2   1   1   0   0   0   0   1   0   1   2   0   1   0   2   1   0   2   2   2   0   0   1   1   0   1   0   0   2   2   1   0   2   1   1   1   0   2   2   1   0   &   1   1   \\
1   0   0   1   2   1   2   2   0   1   0   0   0   1   1   1   2   0   2   1   0   0   2   1   1   1   0   0   1   1   0   0   2   1   0   2   0   2   2   2   2   2   0   0   1   1   0   2   1   &   1   1   \\
0   1   0   1   0   1   0   2   2   2   0   2   1   1   2   1   2   2   1   2   2   1   2   1   0   1   2   0   0   0   1   2   0   0   1   1   1   0   0   1   1   2   0   2   2   2   0   1   1   &   1   1   \\
0   0   0   0   2   0   2   0   1   2   2   1   1   2   1   0   1   2   2   2   1   1   2   2   1   2   2   0   2   2   1   0   1   0   2   0   1   0   1   2   1   0   1   0   0   0   2   2   1   &   1   1   \\
1   1   1   0   2   2   2   2   1   0   1   2   0   1   1   2   1   2   0   2   0   2   1   1   0   0   2   1   2   0   2   0   0   0   0   2   1   1   2   2   0   1   2   0   2   2   2   0   0   &   1   1   \\
0   1   1   1   0   1   0   1   0   0   1   1   2   0   2   2   0   2   0   2   2   2   1   2   1   0   2   1   1   2   1   1   2   2   2   2   0   2   2   0   0   0   1   2   1   0   0   2   1   &   1   1   \\
1   0   1   0   2   0   1   1   1   1   0   0   0   2   2   0   0   0   0   2   0   2   2   0   2   2   1   2   0   1   0   2   2   2   2   2   1   0   0   2   0   2   1   1   0   1   1   0   1   &   1   1   \\
0   0   1   1   1   2   1   1   2   0   1   1   2   1   0   1   1   1   2   2   1   1   0   2   2   2   0   0   2   1   2   2   0   2   1   2   2   1   0   2   0   2   1   0   1   1   0   2   2   &   1   1   \\
1   1   0   0   1   2   0   1   2   0   2   0   2   0   1   0   1   1   1   0   2   1   0   0   0   1   1   2   0   2   0   1   2   1   0   0   1   2   0   1   0   0   0   1   0   0   2   2   1   &   1   1   \\
1   1   0   1   0   2   1   0   0   1   0   1   0   0   0   0   1   1   0   1   1   0   0   0   2   0   1   1   2   0   2   2   1   2   2   0   2   0   1   1   2   2   2   2   2   2   1   0   0   &   1   1   \\

 \ea  $ &
$ \ba{cc}
\d_0^*  &\zero_{72\times 3}   \\
1   1   1   1   2   2   1   1   1   2   1   2   1   1   1   0   1   2   1   2   2   0   1   0   0   0   2   1   2   0   2   0   2   0   0   0   1   2   0   2   1   1   2   0   2   1   1   0   1   &   2   1   0   \\
1   1   0   0   2   0   1   0   0   2   2   2   1   0   2   0   1   2   2   1   2   1   0   0   1   0   0   0   1   2   1   1   1   2   1   2   1   0   2   0   0   2   2   2   2   0   0   1   2   &   1   1   0   \\
1   0   0   0   1   0   2   1   2   1   1   0   0   1   0   1   0   2   0   1   2   0   0   2   2   2   1   0   0   1   2   0   1   1   1   0   1   2   2   1   0   1   0   1   1   0   2   0   2   &   2   1   0   \\
0   1   1   1   0   1   0   1   2   0   0   1   0   0   0   1   0   2   1   0   2   0   1   1   0   1   1   2   1   0   0   1   2   1   2   1   0   2   2   1   2   1   2   0   0   0   2   2   0   &   1   1   0   \\
0   0   0   1   1   0   2   0   2   1   1   0   2   0   2   1   1   2   0   0   0   1   2   1   0   1   2   1   2   2   1   1   0   0   1   0   0   1   0   1   1   1   1   1   1   2   1   1   0   &   1   1   0   \\
1   0   1   1   1   0   0   2   1   1   1   2   0   0   2   2   0   0   2   2   1   0   2   2   2   1   0   2   1   2   0   1   0   0   2   2   1   0   2   0   2   0   1   1   0   1   1   2   1   &   2   1   0   \\
0   0   1   0   0   2   0   1   0   1   0   2   0   1   0   0   1   0   0   0   2   2   2   0   0   2   1   1   2   2   1   2   0   2   2   2   1   0   0   1   0   2   0   2   0   2   2   0   1   &   1   1   0   \\
0   1   1   0   2   0   0   0   1   2   1   1   1   2   0   2   2   2   1   2   0   1   1   2   1   2   2   1   0   2   1   0   1   0   2   1   2   1   2   0   1   0   2   2   0   2   2   0   0   &   2   1   0   \\
1   0   1   0   2   1   2   2   0   1   2   0   1   1   1   1   2   1   2   1   0   2   2   1   1   1   0   0   0   1   1   0   2   1   0   2   0   2   0   0   0   0   1   2   1   1   0   2   2   &   1   1   0   \\
1   0   0   0   1   2   1   0   1   0   2   0   2   1   1   2   1   0   2   0   0   1   1   0   2   2   1   1   1   0   2   2   1   0   0   1   0   1   1   2   1   2   0   2   0   1   0   2   0   &   1   1   0   \\
1   1   0   0   0   2   2   2   2   2   0   1   2   0   1   0   2   1   2   0   2   2   0   1   2   1   2   2   2   0   0   1   1   2   1   1   0   0   1   0   1   0   1   2   2   0   1   0   0   &   2   1   0   \\
1   0   1   1   0   2   2   0   2   2   0   0   2   2   0   0   2   1   1   0   1   1   0   1   1   2   0   0   0   2   0   1   0   1   0   2   2   1   1   2   2   1   0   0   2   2   2   2   2   &   2   1   0   \\
0   1   0   0   1   1   0   1   2   2   0   1   1   1   2   1   2   0   1   0   1   0   0   0   0   1   0   0   0   2   2   2   1   1   1   1   2   0   0   1   1   0   0   2   1   1   0   1   1   &   1   1   1   \\
0   0   1   1   2   0   1   2   2   0   2   0   1   2   0   0   0   1   1   2   1   1   0   2   2   2   2   2   2   0   0   0   2   2   0   1   1   0   0   1   1   0   1   1   1   0   0   2   1   &   1   1   1   \\
1   1   1   1   1   2   2   1   0   1   2   2   0   2   2   2   0   0   0   0   1   2   0   2   1   0   2   1   1   2   2   0   0   1   0   2   0   2   1   0   1   1   2   1   0   0   0   1   0   &   2   1   1   \\
1   0   1   1   0   1   2   0   0   0   1   1   2   0   0   1   1   0   2   2   1   0   2   0   2   2   2   0   0   1   2   2   0   2   2   0   2   2   1   0   2   2   2   0   2   1   0   0   0   &   1   1   1   \\
0   0   0   0   2   0   2   2   1   2   0   2   0   2   1   0   0   1   2   2   0   2   0   0   2   0   2   0   0   0   2   2   2   1   2   0   0   1   1   2   0   1   1   0   0   2   2   1   1   &   1   1   1   \\
1   1   0   0   2   0   0   1   1   0   0   1   2   2   2   1   1   1   0   1   2   2   2   2   1   1   1   1   1   1   0   1   2   0   2   1   2   1   0   2   0   2   2   0   2   0   1   2   1   &   2   1   1   \\
0   1   0   1   1   2   0   0   0   0   2   1   0   1   1   2   2   2   0   2   1   2   1   1   0   2   0   0   2   1   1   2   0   0   2   0   1   0   1   2   2   0   1   0   1   0   2   2   2   &   2   1   1   \\
0   1   1   0   0   1   1   0   1   1   2   2   0   2   1   1   0   2   0   1   1   0   2   2   0   0   1   2   2   0   1   1   1   0   0   0   2   0   0   0   2   2   0   1   2   2   1   1   2   &   1   1   1   \\
1   1   0   0   0   1   0   2   2   0   2   2   1   0   1   2   1   1   1   1   0   1   1   1   1   0   1   2   0   1   2   2   0   2   0   0   1   1   2   1   2   0   0   1   2   1   2   1   1   &   2   1   1   \\
0   0   1   1   1   2   1   1   1   2   1   0   2   2   2   2   2   0   2   1   2   2   1   2   1   0   0   2   2   1   0   0   2   2   1   1   0   1   2   2   0   1   0   0   1   1   1   0   2   &   1   1   1   \\
0   1   0   1   2   1   1   2   0   1   0   0   2   0   0   0   2   0   1   1   0   0   1   0   2   1   1   1   1   1   1   0   2   2   1   2   2   2   1   1   2   2   2   1   0   2   1   0   2   &   2   1   1   \\
0   0   0   1   0   1   1   2   0   0   1   1   1   1   2   2   0   1   0   2   0   1   2   1   0   0   0   2   1   0   0   2   1   1   1   2   2   2   2   2   0   2   1   2   1   2   0   1   0   &   2   1   1   \\

 \ea  $
 \\ \hline

$D_{3B}^* \in \cC_{3B}(72+24;2^{4}3^{45}\bullet 3^{2}\bullet \B)$ \\
\cline{1-1} $ \ba{cc}
\d_0^*  &\zero_{72\times 4}   \\

0   1   0   1   1   1   1   0   1   0   2   1   2   1   0   2   0   2   1   2   1   1   2   1   0   2   1   2   1   0   1   1   1   0   1   1   2   1   1   2   1   2   2   0   0   2   0   2   2   &   2   1   1   0   \\
0   0   1   0   2   0   2   0   2   2   1   1   1   2   2   1   1   2   2   2   1   1   0   2   1   0   2   0   0   2   1   0   0   0   2   1   2   1   0   0   1   1   2   0   1   2   2   1   0   &   2   2   1   0   \\
0   0   0   0   2   0   2   2   1   0   0   0   2   2   1   0   2   1   0   1   0   1   0   1   2   2   0   0   0   0   0   0   2   2   0   2   1   0   0   2   0   2   1   0   1   2   2   1   1   &   1   1   1   0   \\
0   0   1   1   0   0   1   2   0   1   2   1   0   1   2   0   0   0   2   0   1   1   0   2   2   2   0   2   2   0   1   2   1   1   2   2   1   0   1   0   2   0   1   2   0   0   1   0   0   &   2   1   1   0   \\
1   0   0   0   0   2   1   0   1   2   2   2   1   0   0   0   1   2   2   1   2   0   0   0   0   0   2   2   2   0   2   0   1   0   1   0   1   2   0   2   0   1   2   2   2   2   0   0   0   &   1   1   1   0   \\
1   0   1   0   1   2   2   1   2   1   1   0   2   1   0   0   2   1   0   1   2   2   2   1   0   2   1   1   2   1   1   1   0   2   1   0   0   0   0   1   1   0   1   1   1   1   1   0   2   &   2   2   1   0   \\
1   1   1   1   2   1   0   1   2   2   0   0   1   1   1   1   2   0   1   0   2   1   0   2   0   1   2   1   0   2   2   0   1   1   0   1   0   2   0   1   1   0   0   2   0   0   0   0   1   &   1   1   1   0   \\
1   1   1   1   0   2   0   1   2   0   0   1   2   2   0   0   1   1   2   2   1   2   1   0   2   0   2   2   2   2   0   2   2   2   2   1   0   2   1   2   0   1   1   0   2   1   1   2   1   &   2   1   1   0   \\
0   0   0   0   1   1   2   0   1   2   1   0   0   1   1   1   0   0   2   0   1   2   1   2   2   2   1   2   0   1   2   2   1   1   0   1   0   1   2   2   2   1   0   1   1   1   2   2   2   &   1   2   1   0   \\
1   0   1   1   1   2   1   0   0   0   1   0   2   0   2   2   2   0   1   0   2   0   1   0   1   1   0   1   2   1   0   0   2   0   0   2   2   1   2   0   2   2   2   0   1   1   1   2   0   &   2   1   1   0   \\
0   0   0   1   0   0   0   2   2   0   1   0   0   1   0   1   0   2   1   2   0   0   1   1   0   1   0   1   1   0   0   1   0   1   1   2   1   2   2   1   0   1   1   1   0   2   1   1   0   &   1   2   1   0   \\
1   1   1   0   0   1   0   2   1   1   1   2   0   0   1   2   2   1   2   1   0   2   1   1   1   0   0   2   0   0   1   1   1   2   0   0   0   1   2   0   2   0   0   2   2   1   0   1   0   &   2   2   1   0   \\
1   1   0   0   2   0   0   1   1   1   1   1   1   2   2   1   1   2   0   1   2   2   0   2   1   1   1   1   1   1   0   1   1   2   1   1   2   0   2   0   0   2   2   1   2   0   1   0   2   &   1   1   1   1   \\
0   1   0   0   1   0   0   1   1   2   2   2   1   2   2   2   2   2   0   2   0   2   1   2   0   2   0   1   2   1   2   2   0   0   2   0   1   1   0   2   0   0   0   2   0   1   2   0   1   &   2   1   1   1   \\
0   1   0   0   0   1   2   2   2   2   0   2   2   0   2   0   0   1   1   0   2   2   2   0   1   1   1   0   1   1   2   1   2   1   1   2   0   2   1   0   1   1   0   2   2   2   1   1   2   &   2   1   1   1   \\
0   1   1   1   2   1   1   1   0   2   2   2   1   0   0   2   1   1   1   1   0   0   1   0   2   2   1   2   0   2   1   0   2   1   1   2   2   2   2   0   2   0   2   1   1   0   2   0   2   &   1   2   1   1   \\
0   1   0   1   2   2   1   1   0   0   0   1   0   1   1   1   1   0   0   1   2   0   2   0   0   0   1   0   2   0   1   2   1   0   0   0   2   0   0   1   0   2   0   0   2   0   1   1   1   &   2   2   1   1   \\
1   1   1   1   0   2   2   2   0   0   2   2   1   2   1   1   1   0   2   2   0   1   0   1   1   1   2   1   1   1   2   0   0   2   0   2   1   1   2   2   1   2   1   1   2   0   0   2   2   &   2   2   1   1   \\
0   0   1   1   1   1   0   2   0   0   1   1   1   0   1   2   0   1   0   2   1   0   2   2   2   0   2   0   1   0   2   2   2   2   2   0   2   2   2   1   1   0   1   2   1   2   0   2   0   &   1   1   1   1   \\
1   1   0   0   2   0   0   0   2   2   2   0   0   0   1   0   1   2   1   0   2   1   2   1   2   1   2   2   0   2   0   1   2   0   2   0   1   0   1   1   2   0   1   1   0   0   2   2   1   &   2   2   1   1   \\
1   0   0   1   1   2   2   1   2   1   2   2   0   0   2   2   0   1   0   0   1   2   0   2   1   0   0   0   2   2   0   1   0   1   1   2   1   1   1   1   0   1   0   0   1   0   0   2   1   &   1   2   1   1   \\
1   0   1   1   2   1   2   0   0   1   0   1   0   1   0   1   2   2   2   1   0   0   2   1   2   1   0   0   1   1   0   0   0   0   2   0   0   0   1   2   2   1   2   0   2   1   2   0   2   &   1   1   1   1   \\
0   1   0   0   0   2   1   0   0   1   0   0   2   2   0   2   2   0   0   0   0   1   1   0   1   2   1   1   1   2   1   2   0   2   2   1   0   0   1   1   1   2   2   2   0   2   2   1   0   &   1   2   1   1   \\
1   0   1   0   1   0   1   2   1   1   0   2   2   2   2   0   0   0   1   2   1   0   2   0   0   0   2   0   0   2   2   2   2   1   0   1   2   2   0   0   2   2   0   1   0   1   0   1   1   &   1   2   1   1   \\

 \ea  $\\
     \cline{1-1}
 \et} \ec

\bc {\scriptsize  \tabcolsep=1pt \bt {c|c} \hline $D_{3B}^* \in
\cC_{3B}(72+48;2^{4}3^{45}\bullet 3^{0}\bullet \B)$ & $D_{3B}^* \in
\cC_{3B}(72+48;2^{4}3^{45}\bullet 3^{1}\bullet \B)$

\\\hline
$ \ba{cc}
\d_0^*  &\zero_{72\times 2}   \\

0   0   1   0   1   2   0   1   0   0   1   2   1   2   1   2   2   1   1   2   2   2   1   0   1   2   1   0   2   0   2   0   2   2   1   1   1   0   2   1   0   1   0   2   1   0   0   2   2   &   1   0   \\
1   1   1   1   2   1   0   0   2   2   2   1   1   0   0   2   2   1   1   2   1   2   2   1   1   2   1   2   0   2   0   2   1   2   2   1   2   2   2   2   1   0   1   2   0   0   0   0   2   &   1   0   \\
1   1   0   0   2   0   1   2   2   2   0   0   1   1   1   0   2   2   2   1   2   1   0   0   2   0   0   2   1   0   1   1   1   2   2   1   1   0   0   2   0   2   2   2   0   0   0   0   1   &   1   0   \\
0   1   0   1   0   2   2   2   0   0   0   1   1   0   0   0   2   1   1   1   1   1   0   1   1   0   2   0   1   1   2   0   0   1   1   1   2   2   1   1   1   1   2   1   2   2   0   1   0   &   1   0   \\
0   0   0   1   1   0   0   1   2   0   2   2   0   0   2   1   1   0   2   0   2   0   0   0   0   2   1   1   1   2   0   1   0   0   1   0   1   2   0   1   0   1   2   1   0   1   1   1   0   &   1   0   \\
0   1   1   1   1   1   0   1   1   0   0   1   1   1   1   1   2   0   0   2   2   2   2   2   0   1   0   0   0   1   2   1   2   1   2   0   2   1   0   2   1   0   1   1   2   1   1   2   1   &   1   0   \\
1   1   1   0   1   2   1   1   1   0   2   2   2   1   2   0   0   2   1   2   1   0   1   0   0   0   1   1   2   0   2   1   1   0   1   1   2   1   0   2   1   2   2   0   2   2   1   1   2   &   1   0   \\
1   1   1   1   0   2   2   1   1   2   1   2   0   0   1   1   1   0   2   1   2   0   0   2   2   0   2   2   0   0   2   2   1   2   1   0   1   1   1   1   2   0   1   2   2   0   0   0   0   &   1   0   \\
0   1   1   1   2   1   1   2   0   2   2   2   0   1   1   2   1   0   2   1   0   0   1   0   2   0   0   0   1   0   2   0   1   1   0   2   0   2   2   2   2   1   0   0   2   1   1   2   2   &   1   0   \\
1   1   0   0   0   0   2   2   1   0   2   2   0   0   1   2   1   0   0   1   0   2   1   2   1   1   1   1   1   1   2   0   0   0   2   2   0   1   2   1   1   0   0   1   0   0   2   2   2   &   1   0   \\
0   0   1   0   2   1   2   1   2   2   1   0   1   2   0   1   1   1   2   0   0   1   0   2   2   0   2   0   0   2   1   0   1   2   2   2   1   0   0   1   1   0   0   1   1   1   2   2   2   &   1   0   \\
0   0   0   0   1   0   2   0   1   1   1   2   2   0   2   2   2   2   0   0   0   1   2   0   1   0   0   1   1   2   1   2   2   1   2   0   0   1   2   0   1   2   1   2   0   2   2   1   0   &   1   0   \\
0   0   0   1   1   2   2   0   2   1   2   2   0   2   0   0   0   1   0   1   1   0   2   2   2   1   1   2   2   0   1   1   1   0   2   0   2   0   1   1   2   0   1   1   1   2   1   2   1   &   1   0   \\
0   0   1   1   1   2   2   2   0   0   2   1   0   1   1   2   0   2   1   0   1   1   0   1   0   2   0   2   2   2   0   2   0   1   0   2   1   0   0   0   2   0   1   0   1   0   2   2   2   &   1   0   \\
0   1   1   0   0   2   2   2   2   1   1   2   2   2   2   0   1   1   2   1   1   2   2   1   1   0   0   2   2   2   1   1   0   2   0   0   2   1   1   0   0   1   0   0   2   1   0   1   0   &   1   0   \\
0   1   1   1   2   2   1   0   1   1   1   0   0   1   1   0   0   1   0   0   2   2   0   2   0   2   1   1   1   0   0   0   2   1   2   2   0   0   1   0   1   1   2   2   0   2   2   0   0   &   1   0   \\
1   0   0   1   2   0   2   2   1   2   1   0   2   2   0   0   2   1   1   1   2   0   1   0   1   1   0   1   0   1   0   0   2   0   1   1   2   1   1   0   2   0   2   1   0   1   1   0   2   &   1   0   \\
0   0   0   1   0   1   1   0   0   0   0   0   2   1   2   1   1   0   1   0   1   1   2   0   2   1   1   2   0   0   1   2   1   1   1   1   2   0   1   1   1   2   0   0   0   2   0   2   2   &   1   0   \\
1   0   0   1   2   0   2   1   0   2   0   1   2   0   2   0   0   0   2   2   1   2   2   0   0   1   2   0   2   2   2   2   2   1   2   1   0   2   1   0   2   1   2   0   1   1   0   2   1   &   1   0   \\
0   1   0   0   2   0   1   1   0   1   0   0   1   2   2   1   0   1   0   1   2   2   2   1   0   0   0   0   0   1   0   2   2   1   1   0   0   2   0   2   0   2   0   2   1   0   2   1   2   &   1   0   \\
0   0   1   0   0   1   1   2   1   1   1   1   2   1   0   2   0   0   0   2   0   0   1   2   2   2   2   2   2   0   1   2   0   0   0   1   0   2   0   0   1   2   2   2   1   2   0   0   0   &   1   0   \\
1   0   1   0   0   1   0   0   1   0   2   0   2   2   1   2   0   1   2   0   1   1   0   2   1   0   2   2   0   1   0   1   2   2   0   1   0   1   2   2   2   1   0   1   2   1   1   1   1   &   1   0   \\
1   1   0   1   2   1   0   1   0   1   0   0   0   1   0   1   2   0   0   1   0   0   1   1   0   2   1   1   2   1   1   2   0   2   0   0   2   0   1   1   2   2   0   0   0   1   2   1   1   &   1   0   \\
1   1   0   0   2   0   2   0   2   2   2   0   0   1   0   1   1   2   1   2   1   1   0   1   0   1   2   0   0   2   2   0   0   0   0   0   0   1   0   1   0   1   1   1   2   2   2   0   1   &   1   0   \\
1   1   1   1   0   1   0   0   2   2   0   1   2   2   0   1   1   2   2   0   2   2   2   2   0   1   1   1   1   2   0   1   0   0   1   1   0   2   1   0   0   2   0   2   2   2   0   2   2   &   1   1   \\
1   0   0   0   1   1   2   1   2   1   1   2   1   1   1   1   2   2   2   1   2   0   2   1   2   1   2   0   1   1   2   1   1   1   1   2   0   2   0   0   1   1   1   0   1   0   1   1   0   &   1   1   \\
1   0   1   0   1   2   0   1   2   1   0   1   2   0   2   2   2   1   1   0   2   0   0   2   2   2   0   0   0   1   0   1   1   1   0   1   2   1   1   1   0   1   1   0   0   2   2   2   0   &   1   1   \\
0   1   0   0   0   2   1   0   2   2   1   1   1   0   1   0   2   0   1   0   2   0   1   0   0   1   0   2   2   2   0   2   2   0   0   0   0   0   0   1   2   0   0   1   2   2   1   1   1   &   1   1   \\
0   0   1   0   1   1   1   0   0   2   2   0   0   0   0   1   0   2   2   1   0   2   1   1   1   2   1   2   1   1   1   1   2   1   1   2   2   1   2   0   0   1   2   1   1   2   2   0   2   &   1   1   \\
1   0   1   1   2   2   1   2   0   1   2   0   2   1   1   1   1   0   0   1   2   1   1   1   1   1   2   1   1   2   0   0   2   0   0   2   1   1   0   2   0   2   1   2   1   2   0   2   1   &   1   1   \\
1   1   0   0   0   1   2   0   0   1   2   1   1   0   0   2   1   1   0   0   2   2   0   0   2   0   1   1   2   1   2   2   1   2   0   0   2   0   0   0   0   2   2   2   1   1   1   0   0   &   1   1   \\
1   0   1   1   0   2   1   0   2   0   1   0   0   0   0   2   2   2   2   2   1   0   2   0   1   0   0   0   0   2   0   0   0   2   2   2   2   0   0   2   2   2   2   0   1   1   0   2   0   &   1   1   \\
0   1   0   1   1   1   0   2   1   0   1   0   2   0   1   2   0   1   1   2   0   1   2   1   0   1   2   2   1   0   1   1   2   0   0   2   2   2   0   1   0   2   1   2   0   1   0   1   2   &   1   1   \\
0   1   1   1   2   0   1   2   2   1   0   0   2   0   0   2   0   2   0   1   1   0   1   0   2   0   2   2   2   0   0   0   2   2   1   0   1   2   2   1   1   1   1   1   1   0   2   0   1   &   1   1   \\
0   0   0   0   2   2   0   2   1   0   2   1   1   1   1   1   2   2   2   1   1   2   0   2   1   2   0   0   2   1   1   0   0   0   2   0   1   2   1   2   2   2   0   0   1   2   2   0   2   &   1   1   \\
0   1   0   0   0   1   0   1   2   0   0   2   2   1   0   0   0   2   1   2   2   1   1   1   1   1   1   0   0   1   0   1   0   1   2   2   1   0   2   0   2   0   2   2   2   0   2   2   0   &   1   1   \\
0   0   0   1   0   0   0   2   0   2   1   2   1   2   2   1   0   1   0   2   0   0   0   2   2   2   0   0   1   0   0   2   0   2   1   2   2   1   2   0   0   2   1   0   1   2   0   1   1   &   1   1   \\
1   0   0   0   2   2   1   0   0   2   0   2   2   0   1   0   1   1   0   0   0   1   2   0   1   2   1   0   0   0   2   1   0   0   0   2   1   2   1   1   1   0   0   2   2   1   2   0   0   &   1   1   \\
1   0   1   1   1   0   2   0   2   0   0   0   1   2   2   0   0   0   1   2   0   1   0   0   2   2   1   2   1   1   2   2   0   2   2   0   0   1   0   2   2   2   1   1   2   0   1   2   2   &   1   1   \\
1   1   1   1   1   0   1   0   0   1   2   2   1   2   2   2   2   0   1   0   0   0   0   1   0   1   0   0   2   1   1   0   0   0   1   0   1   0   2   0   1   0   2   0   0   1   0   0   0   &   1   1   \\
1   1   0   0   0   2   1   1   1   2   2   1   2   2   0   0   1   1   2   2   0   1   1   1   2   2   2   1   2   2   1   1   2   2   1   0   0   2   1   2   1   0   1   0   0   0   1   2   1   &   1   1   \\
0   1   0   0   1   0   2   0   1   2   0   1   0   2   2   1   1   0   2   2   1   0   1   2   1   2   2   1   0   0   2   0   1   1   0   1   1   0   2   2   1   2   2   1   0   0   2   1   0   &   1   1   \\
1   0   0   0   0   1   2   2   1   1   0   1   0   1   2   0   2   2   2   2   0   2   2   1   0   1   0   2   2   0   0   0   1   2   1   2   1   0   2   2   0   1   1   2   2   1   1   0   2   &   1   1   \\
1   1   1   0   2   0   0   1   0   0   1   1   1   0   2   2   1   2   0   0   2   1   0   2   0   1   2   1   1   2   1   0   1   2   2   2   2   2   2   0   2   1   2   1   2   0   1   1   1   &   1   1   \\
0   0   1   1   0   2   1   1   1   1   1   2   1   2   0   0   1   2   1   2   1   1   2   2   0   2   2   1   1   1   2   2   1   0   2   1   1   1   2   2   0   1   0   0   0   1   2   0   1   &   1   1   \\
0   0   1   1   2   0   0   2   0   0   0   1   0   2   1   1   2   2   0   0   0   2   1   1   2   0   2   2   0   0   1   2   2   0   2   1   0   0   1   2   2   0   2   0   2   0   1   1   0   &   1   1   \\
1   1   1   0   1   1   0   2   1   1   2   2   0   2   2   0   0   0   1   1   1   2   2   0   2   0   1   1   0   2   2   2   2   1   0   2   1   2   1   0   2   0   0   1   0   2   0   1   1   &   1   1   \\
1   0   0   1   1   0   0   1   2   2   1   0   0   1   2   2   0   0   0   0   0   2   1   2   1   0   0   1   2   2   1   1   1   1   0   1   0   1   2   1   0   0   0   2   1   0   1   0   1   &   1   1   \\
 \ea  $
 &
 $ \ba{cc}
\d_0^*  &\zero_{72\times 3}   \\
0   1   1   1   2   1   0   2   2   0   0   1   1   0   1   2   1   1   1   2   2   2   2   1   0   2   2   2   0   0   0   2   2   0   0   0   2   2   0   2   0   0   1   1   1   1   2   1   2   &   2   1   0   \\
0   0   1   1   1   0   2   2   2   0   0   0   2   2   2   0   0   1   0   2   0   1   2   0   2   0   2   0   0   0   2   1   0   1   1   2   1   0   1   2   0   2   1   2   0   1   0   1   0   &   1   1   0   \\
1   0   0   1   2   1   1   0   1   1   2   0   2   1   0   1   2   1   1   1   0   1   2   2   0   1   0   1   1   1   1   1   1   0   1   1   2   1   2   1   0   2   2   1   0   2   0   0   2   &   2   1   0   \\
1   0   0   0   1   2   2   1   1   0   0   2   2   0   1   0   1   2   0   1   2   2   2   0   1   0   1   1   2   1   0   1   2   2   1   0   1   2   0   2   0   1   2   1   2   0   1   0   1   &   2   1   0   \\
1   1   1   1   1   2   1   1   0   1   1   1   2   0   0   2   2   1   0   0   2   0   0   0   0   2   2   1   2   1   2   0   2   2   0   2   2   1   0   0   1   2   2   2   2   2   1   0   0   &   1   1   0   \\
0   1   0   1   0   2   1   1   0   0   2   0   0   0   1   0   0   0   1   0   1   0   0   1   2   0   1   2   1   0   0   1   1   2   0   2   1   0   1   1   2   2   1   2   2   0   0   2   2   &   2   1   0   \\
0   1   0   1   1   0   2   0   0   0   0   2   1   0   2   1   0   0   2   1   1   0   1   0   1   0   0   0   0   1   2   0   2   0   1   0   2   2   2   0   2   2   2   0   1   2   0   1   0   &   2   1   0   \\
0   0   1   1   2   2   1   1   2   2   1   0   0   1   1   1   2   0   2   0   1   2   2   1   1   0   0   0   2   0   1   0   0   1   0   2   0   2   0   1   1   1   0   0   1   2   0   2   2   &   1   1   0   \\
1   0   0   1   1   0   0   1   1   2   2   2   0   2   2   2   1   0   2   0   0   2   0   2   2   2   0   1   0   2   2   2   0   0   1   0   0   1   1   0   1   0   0   2   0   1   0   1   1   &   2   1   0   \\
0   0   1   0   0   0   0   1   1   2   2   2   1   2   0   1   0   1   1   2   1   2   1   2   0   2   1   0   1   1   2   0   2   1   2   2   1   2   2   2   0   1   0   1   0   2   2   2   1   &   1   1   0   \\
0   0   0   0   0   1   2   2   0   2   1   2   1   0   0   0   2   1   1   2   1   1   2   0   1   1   0   0   0   2   2   2   0   2   2   2   1   0   1   0   2   0   0   2   1   1   0   0   2   &   2   1   0   \\
0   1   0   1   2   2   0   1   0   1   0   1   0   2   0   1   1   1   0   1   0   0   0   1   2   2   1   0   0   1   1   2   0   1   2   0   2   0   1   1   2   1   2   0   1   0   2   2   0   &   1   1   0   \\
1   1   1   0   0   2   0   2   0   1   2   0   1   0   0   2   2   0   0   1   2   2   0   1   0   0   0   2   2   2   1   2   0   2   1   0   1   2   2   1   0   0   0   2   1   1   0   0   0   &   1   1   0   \\
1   1   0   0   0   2   0   0   2   0   0   2   2   0   1   2   1   1   0   2   0   1   1   1   1   1   1   0   1   0   1   1   0   0   2   2   0   1   0   1   0   0   0   2   2   2   2   1   0   &   2   1   0   \\
0   1   0   0   0   1   1   0   2   2   1   1   1   1   1   0   1   1   2   0   2   1   0   2   2   0   2   2   2   0   2   0   1   2   1   0   0   0   2   1   1   1   1   1   2   2   0   0   0   &   2   1   0   \\
1   0   0   1   0   1   0   1   2   0   1   1   2   0   2   0   2   0   2   0   2   0   2   0   0   1   0   2   1   2   0   2   2   1   2   1   2   1   1   2   2   2   0   0   0   1   1   2   0   &   2   1   0   \\
1   1   1   0   2   1   2   2   1   2   0   1   0   0   1   1   2   2   2   1   0   2   0   1   2   0   2   0   1   0   2   0   2   1   2   1   0   2   1   0   2   1   2   2   0   1   1   0   2   &   2   1   0   \\
0   1   0   0   2   0   1   2   0   2   2   1   2   1   0   2   1   2   2   2   0   1   1   1   1   1   2   2   2   2   1   1   2   2   1   1   0   1   1   2   1   2   2   0   0   0   2   2   1   &   1   1   0   \\
0   1   1   1   0   1   0   0   0   0   1   0   2   2   0   1   1   0   2   0   2   2   1   2   1   1   1   1   1   1   1   1   0   2   2   0   0   0   2   0   2   0   2   1   2   1   1   2   2   &   1   1   0   \\
1   0   0   0   1   2   2   0   2   1   2   0   0   2   0   0   0   2   1   0   1   2   0   0   1   2   0   2   1   2   1   1   2   1   1   1   0   0   0   0   2   1   1   1   1   2   2   1   2   &   1   1   0   \\
0   1   1   0   1   0   1   1   2   0   0   0   1   1   2   1   2   0   1   0   1   1   0   0   0   1   1   1   2   1   1   2   1   1   1   1   2   0   0   1   1   0   0   2   0   0   1   1   1   &   2   1   0   \\
1   0   1   0   1   2   2   0   2   1   2   1   0   1   1   0   0   0   1   2   1   0   2   2   2   1   1   0   0   1   2   0   1   0   2   0   2   1   1   1   1   0   1   1   2   1   1   0   0   &   1   1   0   \\
0   0   0   0   2   0   1   0   1   1   1   2   1   2   1   2   2   2   0   1   0   0   1   2   2   0   0   2   2   0   1   2   1   0   0   0   1   0   0   2   2   0   1   0   1   0   1   1   1   &   2   1   0   \\
0   1   1   1   0   0   2   0   1   2   1   0   0   2   2   1   1   2   1   2   1   0   0   2   0   0   2   2   0   2   0   2   1   0   0   1   1   1   2   2   0   1   2   0   2   2   1   0   0   &   1   1   0   \\
0   0   0   0   1   0   0   2   2   1   1   1   0   0   2   1   0   1   0   2   1   0   0   2   2   1   2   0   2   2   0   1   0   2   0   1   2   2   2   1   0   1   0   1   1   2   1   1   1   &   2   1   1   \\
1   1   0   0   2   0   2   0   2   2   1   0   0   0   0   0   2   0   1   0   0   1   1   2   1   2   1   1   0   2   0   0   2   1   0   2   0   2   2   0   1   0   1   1   0   0   2   2   0   &   2   1   1   \\
1   1   0   1   2   0   0   1   2   0   1   1   1   2   0   1   0   2   1   2   2   0   2   1   0   1   2   0   1   1   2   0   0   0   1   1   0   2   0   0   1   0   2   1   1   0   0   0   1   &   1   1   1   \\
1   1   1   0   2   0   2   1   2   2   2   2   1   2   2   1   2   2   2   1   0   0   2   0   2   0   1   0   0   2   1   1   1   2   1   2   2   1   2   0   2   0   0   1   2   0   1   2   2   &   1   1   1   \\
1   1   0   1   2   0   1   2   1   0   2   2   0   1   2   0   0   0   0   1   0   1   1   0   0   1   0   2   2   1   0   0   0   1   2   0   1   0   2   2   2   1   0   1   2   1   2   0   2   &   1   1   1   \\
0   1   1   1   0   1   0   1   1   1   0   2   2   0   2   0   0   0   0   2   1   1   2   2   0   0   2   1   2   2   1   2   2   2   2   0   0   1   0   0   2   2   1   0   0   2   0   1   1   &   1   1   1   \\
0   1   0   0   1   2   0   1   1   0   1   1   2   1   2   0   2   1   1   1   2   2   0   2   1   2   0   0   0   0   0   0   1   0   1   1   1   1   1   2   1   2   0   0   1   1   2   1   2   &   1   1   1   \\
0   0   0   0   0   1   0   0   0   1   1   0   1   1   1   1   0   2   0   0   0   0   2   0   0   2   1   0   0   1   2   1   1   1   0   0   0   1   0   2   0   2   0   2   1   1   2   2   2   &   2   1   1   \\
1   1   1   1   0   1   1   2   0   1   0   2   0   2   2   2   0   1   0   1   1   2   2   1   1   0   1   2   0   0   2   2   1   1   1   1   0   1   1   1   1   1   0   1   2   0   1   1   1   &   2   1   1   \\
1   0   1   0   2   0   1   0   0   2   0   0   2   2   2   0   1   0   0   0   2   2   1   0   2   2   1   1   0   0   0   2   2   2   2   1   2   0   1   2   0   1   2   0   1   2   2   0   2   &   2   1   1   \\
1   0   1   0   1   1   0   0   1   2   0   1   0   0   0   2   2   2   2   1   1   1   2   1   1   2   1   2   0   0   0   1   1   0   0   0   2   0   0   2   2   2   0   2   0   2   0   0   1   &   1   1   1   \\
0   1   1   1   1   1   1   2   1   2   1   2   2   1   1   2   0   0   1   1   2   0   1   0   2   0   1   2   1   0   1   1   2   1   0   1   2   2   1   0   1   0   0   2   0   2   2   1   0   &   1   1   1   \\
1   1   1   1   2   1   1   0   0   0   2   1   1   1   2   2   0   2   2   2   2   1   0   2   2   2   0   0   1   2   1   2   1   1   2   2   1   2   0   2   0   2   2   0   2   2   0   2   1   &   2   1   1   \\
1   1   0   0   1   1   0   0   0   0   2   2   0   1   1   2   1   1   1   0   1   0   1   0   0   0   2   1   2   2   2   2   0   0   0   2   1   2   0   1   2   1   2   0   0   1   2   2   1   &   1   1   1   \\
1   0   1   1   1   2   2   1   1   0   1   0   2   1   1   2   2   1   2   2   1   1   1   2   2   2   2   2   2   1   2   2   2   2   2   1   0   1   1   2   2   1   1   1   1   0   1   2   2   &   2   1   1   \\
0   1   0   0   2   2   1   2   1   1   2   0   1   2   2   0   1   1   2   2   0   2   0   0   1   1   1   1   1   2   2   0   2   0   0   2   2   2   0   1   0   2   1   2   0   0   1   0   2   &   1   1   1   \\
1   0   0   0   0   2   2   2   0   2   0   1   2   2   1   1   2   1   2   0   2   2   0   1   0   1   0   0   2   1   0   1   1   1   0   2   0   0   1   0   0   0   1   0   2   0   0   1   1   &   1   1   1   \\
1   0   1   1   1   1   1   1   2   1   0   0   2   0   0   2   1   2   1   1   0   0   1   1   2   1   0   2   0   1   0   0   2   2   1   2   0   2   2   1   0   2   1   0   1   1   1   2   1   &   1   1   1   \\
1   0   0   0   0   1   2   2   2   1   0   2   0   1   0   1   1   0   2   2   0   0   1   0   2   2   2   1   1   0   2   2   1   2   2   1   1   0   2   2   1   2   2   2   2   1   0   1   0   &   1   1   1   \\
1   0   1   1   2   0   2   2   0   1   1   0   1   1   1   1   1   0   0   1   2   1   0   1   1   1   2   1   1   0   0   0   0   0   0   0   1   1   2   0   2   2   1   1   0   2   2   2   1   &   2   1   1   \\
0   0   0   1   2   0   2   1   1   2   0   1   1   0   1   2   0   2   0   0   2   2   2   2   1   1   2   1   1   2   1   1   1   1   2   2   1   0   2   1   1   1   1   2   1   0   2   0   0   &   1   1   1   \\
0   0   1   0   0   2   0   2   1   0   2   1   1   2   2   2   2   2   1   2   0   1   1   1   1   0   0   1   1   1   0   0   0   2   0   2   2   1   2   0   1   1   2   0   2   1   0   1   0   &   2   1   1   \\
0   0   1   1   1   2   1   2   2   1   2   2   0   1   0   0   0   2   2   0   2   2   1   2   0   2   0   2   2   0   0   1   0   0   2   1   1   0   0   0   0   0   2   2   0   0   1   2   0   &   2   1   1   \\
0   0   1   1   0   2   2   0   0   2   2   2   2   2   0   0   1   2   0   1   1   1   2   1   0   2   2   1   2   2   1   0   0   0   1   0   2   2   1   1   1   0   1   0   2   0   2   0   2   &   2   1   1   \\
 \ea  $
  \\ \hline
 \et} \ec

\bc {\scriptsize  \tabcolsep=1pt \bt {c} \hline $D_{3B}^* \in
\cC_{3B}(72+48;2^{4}3^{45}\bullet 3^{2}\bullet \B)$
\\\hline
$ \ba{cc}
\d_0^*  &\zero_{72\times 4}   \\
1   0   0   0   1   2   1   1   1   1   2   2   2   0   2   0   1   2   0   1   2   2   1   0   1   0   1   1   2   1   0   1   2   0   0   0   1   1   0   2   0   2   0   0   2   0   1   2   1   &   2   1   1   0   \\
1   1   0   0   1   0   0   0   2   0   0   1   2   0   2   0   0   1   1   2   0   1   1   0   1   1   2   0   1   2   0   1   2   1   2   1   1   0   1   2   1   1   2   2   0   2   2   1   0   &   2   1   1   0   \\
0   0   0   0   2   0   0   2   1   1   0   2   0   2   1   2   2   0   2   1   0   2   0   1   2   0   2   2   2   0   0   0   2   1   0   1   0   2   0   2   2   1   2   0   0   2   2   1   0   &   1   2   1   0   \\
1   1   0   0   1   2   0   2   2   0   2   1   0   1   1   0   0   0   1   0   2   2   1   2   0   2   0   2   2   0   1   1   1   0   0   1   1   0   0   0   0   0   1   2   1   0   0   2   0   &   1   2   1   0   \\
0   1   0   1   2   1   1   2   0   0   1   1   1   0   0   0   2   1   1   1   2   2   0   0   2   2   1   0   0   1   0   2   2   0   1   1   2   2   1   0   2   2   2   2   1   1   0   0   2   &   2   2   1   0   \\
1   1   1   1   1   2   2   1   2   0   0   2   1   2   2   1   0   0   1   2   1   1   0   2   2   0   1   2   0   1   2   1   2   1   1   1   2   2   0   1   1   0   2   1   0   0   1   2   0   &   2   1   1   0   \\
1   1   1   1   0   1   0   1   1   2   0   2   0   0   2   1   1   2   2   2   2   2   2   0   0   1   0   0   0   2   2   2   1   1   2   2   1   2   0   2   0   0   0   2   0   1   0   0   1   &   1   2   1   0   \\
1   0   0   0   1   0   2   2   2   1   0   0   0   1   0   1   0   2   0   1   0   2   2   1   0   1   0   0   1   1   0   0   2   1   1   2   0   2   2   1   0   1   1   1   1   0   1   0   2   &   2   1   1   0   \\
1   0   1   0   1   0   2   0   1   2   1   2   1   0   1   0   2   2   1   2   0   1   2   1   1   2   1   1   0   0   2   0   1   0   1   2   1   0   2   1   1   0   1   2   0   0   2   0   2   &   1   1   1   0   \\
1   1   1   0   2   0   0   1   1   2   1   1   0   1   1   1   1   1   0   1   2   0   0   2   2   1   2   1   0   1   2   0   1   0   2   0   0   1   1   0   2   1   2   1   2   1   1   0   2   &   2   1   1   0   \\
0   1   1   0   2   0   1   1   1   1   1   0   2   1   2   0   0   0   0   2   1   1   2   2   1   0   2   0   1   2   1   0   0   0   2   2   1   2   0   0   1   2   1   0   1   2   0   1   1   &   1   2   1   0   \\
1   1   0   1   2   0   1   1   0   1   2   2   1   2   0   2   0   0   0   1   2   0   1   0   0   0   1   1   1   2   1   0   2   1   2   0   2   0   2   0   2   0   2   1   0   0   2   0   1   &   1   2   1   0   \\
1   1   0   0   0   2   2   0   2   2   2   0   2   1   0   0   1   2   2   0   2   1   0   1   2   1   1   2   2   2   1   1   0   2   1   0   0   0   1   1   1   1   1   1   2   2   1   0   1   &   1   2   1   0   \\
0   0   1   1   0   2   2   2   1   1   1   0   0   0   2   0   2   0   1   2   1   0   2   0   1   1   0   2   2   0   0   2   1   2   0   0   2   0   2   0   2   2   1   1   0   2   1   0   1   &   2   2   1   0   \\
0   0   1   1   1   1   0   2   2   0   1   2   0   0   0   1   0   0   1   1   1   0   1   1   0   2   1   1   1   0   2   0   0   0   0   0   1   2   1   1   2   1   0   1   2   2   2   1   0   &   1   2   1   0   \\
0   0   1   1   0   1   0   1   2   2   1   1   1   2   0   2   1   0   0   0   1   2   1   2   1   1   2   1   2   2   0   2   2   2   2   1   0   1   2   1   1   0   1   2   1   2   0   2   0   &   2   1   1   0   \\
0   0   0   1   1   2   0   0   0   1   1   2   2   1   1   1   2   0   0   0   0   0   2   0   2   1   0   0   2   1   1   2   0   2   2   0   0   0   0   1   2   0   0   1   1   1   2   2   2   &   1   1   1   0   \\
1   1   1   1   2   2   1   0   0   0   1   0   0   0   1   2   1   2   2   0   0   0   1   0   2   0   2   2   1   0   0   0   1   2   1   2   0   1   2   2   1   1   2   0   2   0   0   2   0   &   2   2   1   0   \\
0   0   1   0   1   0   1   0   1   2   2   0   2   2   2   2   0   1   2   2   1   2   1   2   2   2   1   2   2   1   2   2   2   2   1   1   2   1   2   2   0   1   1   1   1   1   1   0   2   &   1   2   1   0   \\
1   0   1   0   0   1   2   0   0   1   2   0   2   1   0   1   1   2   2   2   0   1   2   2   0   2   2   1   0   2   1   2   0   0   0   1   2   1   0   2   0   2   2   1   0   2   2   2   2   &   2   2   1   0   \\
0   1   0   0   2   1   1   0   0   2   2   2   1   2   2   1   0   2   2   0   1   2   0   2   0   2   0   0   1   0   1   1   1   1   2   2   0   1   1   0   1   1   0   0   2   2   2   1   2   &   2   1   1   0   \\
0   1   0   1   0   1   1   2   1   0   2   2   1   1   1   2   0   2   2   2   1   1   2   1   0   0   2   2   2   0   0   2   0   0   0   0   1   0   2   2   2   0   0   0   1   1   0   1   1   &   2   1   1   0   \\
1   0   1   1   2   0   2   2   0   0   0   0   2   2   1   1   2   1   1   1   0   1   0   1   1   1   0   0   0   1   0   1   0   1   0   2   1   1   1   2   2   2   1   0   2   0   0   2   1   &   1   1   1   0   \\
0   1   1   0   0   1   0   1   2   1   0   1   1   2   2   1   2   2   1   0   2   0   2   1   1   0   1   0   0   1   1   1   1   2   1   1   2   0   2   0   0   2   0   0   2   0   1   1   0   &   1   2   1   0   \\
1   1   1   1   2   1   0   2   2   1   2   2   0   1   1   2   0   1   2   2   0   2   2   1   2   1   1   1   1   1   1   1   0   2   2   2   2   0   2   0   1   2   1   2   0   1   2   2   0   &   2   1   1   1   \\
0   0   0   1   2   2   2   0   1   0   0   1   0   2   1   1   1   2   0   2   1   1   2   0   0   1   1   1   1   0   0   1   2   0   2   2   2   0   1   2   0   1   1   0   2   2   1   2   2   &   1   2   1   1   \\
0   0   0   1   1   0   0   0   2   2   0   0   0   0   0   1   0   2   0   1   0   0   0   2   2   2   2   2   0   0   1   1   1   2   1   0   2   1   0   2   0   2   1   2   1   2   0   1   1   &   2   1   1   1   \\
0   0   1   0   1   2   2   2   0   0   2   2   1   1   1   2   1   1   0   1   2   0   0   2   1   0   1   0   1   0   2   0   1   1   0   2   2   0   0   1   1   2   0   2   1   2   0   2   2   &   2   2   1   1   \\
0   1   0   1   0   2   2   2   1   1   0   1   0   0   2   2   0   1   0   0   1   2   0   2   2   0   2   0   0   2   2   0   0   0   1   0   1   0   1   1   1   0   1   1   2   0   2   1   2   &   1   1   1   1   \\
0   1   0   1   1   1   2   2   2   2   1   0   1   1   1   1   2   0   1   2   2   0   1   2   0   1   2   0   0   2   2   0   2   1   0   1   0   1   2   2   1   1   0   0   0   1   1   2   1   &   1   2   1   1   \\
0   0   1   0   2   0   1   1   2   0   2   1   0   1   0   1   1   2   1   2   1   0   1   1   2   2   0   0   2   1   1   2   1   0   2   1   1   0   1   1   1   0   2   0   0   0   2   2   1   &   2   2   1   1   \\
0   0   1   1   1   0   2   1   2   0   1   2   2   2   2   0   1   1   2   0   2   1   0   0   0   2   0   1   2   2   1   0   0   1   1   0   0   0   0   0   0   1   1   0   0   1   0   1   0   &   2   1   1   1   \\
0   1   1   0   0   0   0   0   0   2   0   0   0   2   1   0   2   1   0   0   2   2   1   1   0   0   1   1   0   1   0   1   0   0   0   0   0   1   0   1   2   0   0   1   0   2   2   1   1   &   2   2   1   1   \\
1   1   0   1   0   0   0   2   0   0   2   2   1   2   2   0   1   1   1   2   0   1   0   2   0   1   2   1   1   1   2   2   0   2   1   1   1   1   2   1   0   2   2   1   2   2   0   0   1   &   1   2   1   1   \\
1   1   0   0   2   2   0   1   2   1   1   0   1   0   0   2   2   1   1   0   1   1   0   0   0   2   1   2   0   2   0   2   2   2   0   2   2   1   0   1   0   0   0   0   1   2   2   1   1   &   1   1   1   1   \\
1   0   0   0   2   0   2   0   0   2   2   1   0   0   2   1   1   0   2   0   1   0   2   0   1   0   1   0   0   2   2   2   0   1   0   0   0   2   1   0   2   1   2   2   1   0   0   2   0   &   1   1   1   1   \\
1   0   1   1   2   2   1   2   0   2   0   1   2   1   2   2   1   1   0   1   2   2   2   1   0   2   2   2   2   0   0   2   2   2   1   1   0   2   1   0   1   2   0   2   0   0   1   1   1   &   1   1   1   1   \\
1   0   1   1   0   1   2   1   2   1   2   2   2   0   0   0   0   1   2   1   1   0   2   2   2   0   2   2   0   0   0   1   2   1   2   2   0   2   1   0   2   0   1   0   2   1   1   2   2   &   1   2   1   1   \\
0   1   0   1   1   1   1   1   0   1   0   1   2   0   0   2   2   1   2   0   2   0   2   1   0   1   2   1   2   0   1   1   1   0   0   2   2   2   0   1   1   1   2   2   1   0   1   0   0   &   1   1   1   1   \\
0   0   1   1   0   1   0   0   0   0   1   0   1   1   0   2   0   2   1   0   2   2   2   2   1   2   0   2   1   1   0   1   1   1   2   0   1   2   2   2   0   1   2   1   0   1   2   0   2   &   1   1   1   1   \\
1   0   1   0   0   1   1   0   2   2   2   1   1   2   1   0   2   0   2   0   2   0   0   0   2   1   0   2   1   2   2   0   2   2   2   2   1   2   0   2   2   2   1   1   1   0   1   1   1   &   2   1   1   1   \\
1   1   1   1   1   2   0   0   1   2   2   1   2   2   0   2   2   2   1   1   1   0   1   1   1   0   0   2   2   2   1   0   0   0   1   1   2   2   1   0   2   0   2   0   1   2   0   0   2   &   2   1   1   1   \\
1   0   0   1   0   2   1   0   1   1   0   0   2   1   0   0   0   0   2   0   1   1   0   0   1   2   0   1   1   1   2   0   1   1   0   1   0   1   1   1   1   2   0   2   0   1   0   0   0   &   2   2   1   1   \\
1   0   0   0   1   0   1   2   1   2   1   2   1   0   2   2   2   0   2   1   0   1   1   2   1   1   0   0   1   2   1   1   1   0   1   2   2   1   2   0   2   0   0   1   2   1   1   1   0   &   2   2   1   1   \\
0   1   1   0   0   2   1   1   0   0   1   0   0   1   2   1   2   0   0   1   0   2   1   0   2   0   0   1   2   0   2   2   0   1   1   2   2   0   1   2   0   2   2   2   2   1   1   0   0   &   1   1   1   1   \\
0   1   0   0   0   1   2   1   1   0   0   1   2   0   1   2   2   2   0   2   0   1   1   1   1   2   0   1   1   1   2   2   2   2   2   0   0   1   0   0   2   2   2   2   2   0   0   1   2   &   1   2   1   1   \\
0   1   0   0   2   1   1   2   1   2   0   0   2   2   1   0   1   0   1   1   0   1   1   0   2   0   1   2   0   0   1   2   1   2   0   1   1   2   2   1   0   0   0   2   2   1   2   2   2   &   2   1   1   1   \\
1   0   0   0   2   2   2   1   0   1   1   1   1   2   0   0   1   1   0   2   0   2   0   1   1   2   2   0   2   2   2   0   0   2   2   0   1   2   2   2   0   1   0   0   1   1   2   0   0   &   2   2   1   1   \\

 \ea  $
   \\ \hline
 \et} \ec

\bc {\scriptsize  \tabcolsep=1pt \bt {c} \hline $D_{3B}^{*T}$, where
$D_{3B}^{*}  \in \cC_{3B}(72+72;2^{4}3^{45}\bullet 3^{0}\bullet \B)$

\\\hline
$ \ba{cc}
\d_0^{*T}  & \ba{c}1    1   1   0   1   0   1   0   0   0   0   1   0   0   0   1   0   0   0   0   1   0   1   1   0   1   1   1   1   0   0   1   0   1   1   1   0   0   1   1   1   1   1   1   1   1   0   1   1   0   1   0   0   0   0   1   0   0   1   1   0   0   1   0   0   1   0   1   0   1   0   0   \\
0   1   1   1   0   1   1   0   0   0   0   0   1   0   0   0   0   0   1   1   0   1   1   0   0   1   1   0   1   1   1   1   1   0   1   0   1   0   0   0   1   0   0   1   1   0   1   1   1   0   1   1   1   0   1   0   0   0   1   0   1   0   1   0   0   0   0   0   1   1   1   1   \\
0   0   1   0   1   0   0   1   0   1   0   0   0   1   1   1   1   0   0   0   0   1   0   1   0   1   1   1   1   1   0   0   1   0   1   0   0   1   0   0   1   1   0   1   0   0   1   1   1   1   0   1   1   0   1   1   0   1   0   1   1   1   1   0   1   0   1   0   0   0   0   0   \\
1   1   1   1   0   0   1   1   0   0   0   0   0   0   1   0   1   1   0   1   0   1   1   0   1   1   1   0   0   1   0   0   0   1   1   0   0   1   0   1   0   1   0   1   0   1   1   0   1   1   1   1   0   1   1   0   1   1   0   0   0   1   1   0   1   1   0   0   0   0   0   0   \\
0   0   1   1   1   1   2   2   2   1   2   1   2   0   2   0   2   0   0   0   1   1   1   0   1   2   0   2   1   2   0   0   1   1   2   2   0   1   0   0   0   2   1   0   2   2   0   2   2   1   2   1   2   0   2   0   1   0   2   0   1   1   1   2   0   1   0   2   1   2   1   0   \\
0   0   2   1   0   0   1   0   0   2   2   2   1   1   1   1   0   2   1   0   2   0   0   2   1   2   2   2   2   0   2   0   1   0   1   0   2   0   1   2   1   0   0   1   1   1   2   0   0   1   0   2   0   2   1   0   2   2   1   1   1   2   1   0   2   1   2   1   1   2   0   2   \\
2   2   0   1   0   2   2   1   2   1   2   1   0   2   0   1   1   1   0   0   2   0   2   2   0   1   0   0   0   1   2   2   1   2   1   0   0   0   0   1   0   2   1   2   0   2   1   2   1   2   0   2   0   0   1   0   0   2   1   2   1   2   0   0   2   1   1   1   2   1   1   1   \\
0   2   2   1   0   2   2   1   0   0   0   1   2   2   2   0   2   2   1   1   1   0   1   1   1   1   0   2   0   1   1   2   2   0   0   2   2   2   0   2   1   1   2   0   1   0   0   2   2   0   0   0   2   1   0   0   0   1   1   2   1   1   2   1   2   2   1   0   1   1   0   0   \\
0   1   1   2   1   1   0   0   1   1   1   2   0   0   2   2   1   2   2   2   2   2   0   1   1   1   0   0   2   1   1   2   0   1   0   0   1   2   2   0   1   2   1   2   2   0   0   2   0   2   1   0   0   0   0   0   1   2   2   2   1   2   2   0   1   1   0   0   2   1   0   1   \\
1   2   0   0   0   1   2   1   1   2   2   2   1   2   0   1   0   1   2   2   0   2   0   1   0   1   2   2   1   0   2   1   1   2   0   1   0   2   0   0   2   1   2   1   0   2   0   2   1   0   0   2   0   1   2   0   1   0   2   2   1   1   0   2   2   1   1   0   0   1   0   1   \\
2   0   0   2   1   1   2   2   1   1   2   0   0   2   0   1   1   0   1   1   0   2   1   0   1   2   0   1   2   2   0   1   2   0   2   2   2   0   0   1   1   1   2   0   1   0   1   0   0   2   0   2   2   0   1   0   1   1   2   2   0   2   1   1   1   2   2   2   0   1   0   0   \\
1   2   1   2   0   1   0   1   0   2   1   1   2   0   1   0   0   0   0   1   0   0   2   0   1   0   1   1   2   2   0   1   2   2   2   2   1   0   2   0   2   2   0   0   1   1   1   0   0   1   0   2   0   1   0   1   2   2   1   2   1   0   2   2   1   2   2   0   2   1   1   2   \\
1   2   0   0   0   0   0   0   0   1   0   2   0   1   1   0   2   2   1   1   0   0   2   2   1   2   0   1   2   1   2   1   2   1   1   2   1   2   0   0   1   2   2   2   1   2   2   1   0   2   2   1   1   1   1   0   1   2   1   0   2   0   2   2   0   0   1   2   0   0   1   0   \\
0   0   1   1   0   1   1   2   2   0   0   2   0   1   1   1   2   2   0   1   1   0   0   2   0   1   2   1   0   2   2   2   1   2   1   1   2   2   1   0   2   2   0   0   1   0   0   2   0   2   1   2   2   2   2   0   0   1   0   0   0   1   1   0   1   2   0   1   2   2   1   1   \\
0   1   1   0   1   2   1   2   0   1   1   0   2   0   1   1   0   0   0   2   1   0   2   2   0   0   0   1   0   2   1   2   2   2   0   1   2   0   0   2   0   2   2   1   1   0   0   2   0   1   1   2   1   1   1   2   1   1   2   2   1   0   2   0   2   2   1   2   1   0   2   0   \\
2   0   0   2   2   0   2   1   0   2   1   2   1   0   1   1   1   0   1   0   1   2   2   0   0   1   1   2   2   0   0   1   2   1   0   0   2   1   0   1   1   2   0   0   0   1   0   1   2   1   2   0   2   2   2   2   0   2   2   0   1   1   2   1   1   2   2   0   0   0   1   1   \\
0   1   0   1   2   0   2   2   1   1   1   2   2   2   2   1   0   0   2   2   1   2   1   2   0   1   1   2   0   0   0   0   0   2   1   1   2   0   1   0   0   0   1   2   0   2   0   1   1   1   2   2   1   1   0   0   2   1   1   2   2   0   2   1   1   0   2   0   0   2   1   0   \\
1   1   0   2   0   0   0   0   1   0   2   1   2   2   1   2   0   2   1   0   0   2   1   1   2   2   0   1   1   1   1   2   1   2   0   0   0   1   0   2   2   0   0   0   1   0   1   0   1   1   2   2   2   2   2   2   1   0   2   2   0   0   1   2   1   1   0   2   2   1   1   0   \\
1   1   0   1   1   0   2   2   0   1   2   0   0   1   1   2   0   0   2   1   0   0   0   2   2   1   2   0   2   1   1   1   2   2   2   0   0   0   1   2   2   0   2   0   2   1   2   1   0   1   0   1   1   0   0   2   1   2   2   2   2   1   1   0   2   2   0   1   0   1   1   0   \\
2   2   2   2   0   0   1   1   1   1   1   1   1   0   2   0   2   2   0   0   1   0   1   0   1   1   2   2   0   1   0   2   1   1   0   2   2   1   1   0   2   0   2   1   0   0   2   0   1   0   1   2   1   0   0   0   2   2   0   1   1   0   2   2   1   2   0   2   2   1   0   2   \\
1   0   1   0   0   1   0   1   0   1   1   2   0   1   0   0   2   1   2   2   2   2   1   0   2   2   1   2   0   1   1   2   2   0   1   2   0   0   1   2   1   0   0   2   2   0   1   2   0   1   1   2   0   2   2   2   2   2   2   1   0   1   1   1   0   1   1   0   0   2   0   0   \\
0   0   2   1   1   2   2   0   1   2   2   0   2   1   1   0   1   1   2   0   0   2   1   2   0   0   0   2   0   0   1   1   2   1   1   2   1   1   2   0   2   2   1   1   1   2   1   0   2   1   0   1   0   2   0   2   0   2   1   0   0   0   2   0   1   0   1   1   2   2   2   0   \\
2   1   2   1   1   0   1   0   2   1   2   0   2   2   2   0   1   1   0   1   2   0   0   2   2   0   0   0   1   2   0   2   2   0   0   1   0   1   1   1   1   0   2   2   0   2   0   1   0   2   1   2   0   2   2   2   0   2   0   0   1   0   1   1   1   2   1   1   0   1   1   2   \\
1   1   1   1   2   0   1   2   0   0   2   0   1   0   1   1   2   1   2   0   1   2   2   1   2   2   0   2   0   0   1   2   2   1   0   0   0   0   2   0   1   2   0   0   2   0   1   2   1   2   2   2   1   2   1   1   1   1   0   1   0   0   1   0   2   2   2   0   0   1   0   2   \\
2   2   0   0   1   0   2   2   1   2   1   2   0   1   2   0   2   1   0   0   0   2   1   1   2   0   2   1   1   0   0   0   1   0   2   0   1   0   1   1   0   2   0   0   2   1   1   1   2   2   1   2   1   2   0   0   1   1   1   2   2   0   1   2   1   0   0   2   1   2   0   2   \\
1   0   0   2   1   1   0   0   0   0   2   1   1   2   1   2   2   1   1   0   2   2   0   0   2   0   2   2   2   1   2   1   0   2   1   1   2   0   2   0   0   2   1   1   1   1   2   0   2   2   1   1   1   0   0   1   2   0   0   0   1   1   1   2   1   2   0   0   0   2   2   0   \\
2   2   2   1   0   0   0   0   2   1   0   0   1   2   1   0   2   0   2   0   1   1   1   0   2   2   1   0   1   1   1   1   2   2   0   2   0   0   1   1   0   1   0   2   1   1   2   2   0   1   2   0   1   2   2   2   0   0   2   2   1   1   1   2   2   0   1   0   0   2   0   1   \\
2   0   2   2   2   0   2   0   1   1   0   2   2   0   0   1   2   2   0   2   0   1   1   1   1   1   2   0   2   2   1   1   2   0   1   1   1   2   2   2   1   1   1   0   2   1   0   0   0   2   1   2   0   1   2   1   0   0   0   2   0   1   1   0   1   0   2   0   2   1   0   0   \\
1   1   0   1   0   1   1   2   1   1   2   2   0   0   2   1   1   2   1   2   1   2   2   0   0   1   0   2   0   0   0   2   2   2   2   2   1   0   2   1   1   0   2   1   1   1   2   0   0   0   0   1   1   0   2   1   0   2   0   0   0   2   1   0   2   1   1   0   2   2   0   2   \\
2   0   0   0   2   1   0   0   1   0   2   0   0   0   1   0   0   1   1   2   1   2   1   2   0   2   1   1   1   0   0   2   2   1   2   2   2   2   1   0   2   2   1   2   0   1   2   1   0   0   0   2   1   1   0   0   1   2   2   2   1   1   1   2   0   2   1   1   0   1   2   0   \\
2   2   0   1   0   0   1   0   1   2   1   0   1   2   0   1   0   0   2   0   2   1   0   1   1   0   2   2   0   2   2   0   1   1   1   1   1   2   0   0   2   2   2   1   1   0   1   1   0   0   1   1   0   2   1   0   2   0   2   0   2   1   2   0   2   2   1   2   2   0   1   2   \\
2   1   1   0   0   0   0   0   1   2   1   1   2   0   2   2   2   1   0   2   0   1   2   1   1   0   2   0   1   0   2   1   2   1   2   0   2   0   1   1   1   0   2   1   1   0   0   2   2   1   0   0   1   2   0   2   0   1   0   2   0   2   1   2   0   2   1   2   2   0   1   1   \\
2   0   2   0   0   0   2   1   1   1   0   1   1   0   0   1   2   0   0   1   0   1   2   2   1   2   1   1   0   2   1   1   2   2   0   0   0   2   2   0   0   0   2   2   2   2   0   1   0   0   1   1   2   0   1   2   2   1   1   0   2   1   0   2   1   1   2   1   2   2   1   0   \\
1   0   1   1   0   1   2   1   1   2   0   2   0   0   2   2   0   0   2   1   0   0   2   2   2   0   0   1   2   2   1   1   1   1   2   2   2   2   2   2   1   2   0   1   1   1   2   2   1   0   0   0   1   1   2   0   1   0   0   1   1   1   0   2   0   2   0   1   2   0   0   0   \\
1   2   0   1   0   1   0   1   2   1   2   1   1   2   2   0   2   0   1   2   0   2   0   0   1   1   0   1   1   2   0   0   0   1   2   2   0   1   1   2   1   2   1   0   2   1   1   1   1   0   0   1   0   2   0   2   0   1   0   2   2   0   2   0   0   1   2   2   2   0   2   2   \\
1   2   2   2   1   2   0   0   2   1   2   2   1   0   2   2   1   1   2   0   0   0   0   2   0   1   1   2   0   1   1   2   0   1   2   0   1   1   0   2   0   2   0   0   1   2   2   1   0   0   1   2   2   0   0   1   0   1   1   1   2   1   1   2   0   2   2   1   0   1   0   0   \\
2   1   1   2   1   0   0   2   2   1   1   2   2   1   2   0   1   0   1   1   2   0   2   2   2   0   2   0   0   2   0   1   0   0   1   1   2   1   1   1   2   0   0   2   0   0   2   0   1   0   2   1   1   1   2   0   1   0   2   0   0   1   2   0   1   1   2   2   0   1   2   0   \\
2   2   2   0   1   1   2   0   1   0   0   0   2   0   0   2   1   0   1   2   2   2   1   0   2   1   0   2   1   2   1   0   2   1   0   0   0   2   1   0   2   2   2   1   0   2   2   0   0   2   0   0   2   1   1   0   0   1   2   1   1   2   1   1   1   1   1   1   0   2   1   0   \\
2   1   2   2   2   1   0   2   2   2   1   1   1   0   2   2   1   2   2   0   0   1   2   0   0   1   1   0   1   1   1   0   2   0   0   1   2   0   2   2   2   2   2   1   1   1   1   2   0   0   0   2   2   1   2   0   1   0   0   1   0   0   1   0   1   1   0   1   0   0   2   0   \\
1   0   2   1   0   0   2   0   1   1   2   2   1   1   2   1   2   1   1   2   2   1   1   2   0   2   0   0   2   1   0   0   0   2   0   1   2   0   2   0   2   0   0   0   1   1   0   0   1   2   1   0   1   0   2   2   1   1   2   1   2   1   0   0   2   2   0   2   2   0   1   1   \\
0   1   1   1   2   1   2   2   2   1   0   1   0   1   2   0   0   0   0   2   1   2   0   0   2   1   2   0   2   2   2   0   1   1   2   1   1   0   0   2   2   1   2   2   1   1   1   1   0   0   1   0   2   1   1   0   2   1   0   2   0   0   0   2   2   0   2   1   2   1   0   0   \\
0   1   0   0   2   2   2   1   1   1   0   0   0   0   2   2   0   1   1   1   1   0   1   2   2   2   2   1   1   2   0   0   0   2   0   2   2   2   2   2   1   0   0   2   0   1   1   0   2   1   1   2   2   2   1   1   0   0   2   1   0   1   0   0   0   1   0   2   1   1   1   2   \\
1   2   1   2   0   2   1   2   1   0   2   1   0   1   1   1   1   0   1   0   0   2   2   0   1   2   2   2   0   2   0   1   0   0   0   0   2   2   1   2   2   1   1   2   0   0   1   0   1   2   1   1   0   1   0   2   2   1   0   1   2   1   0   0   0   0   2   2   2   2   1   0   \\
1   2   2   0   1   1   1   0   1   1   0   2   2   2   0   2   0   1   0   2   1   1   0   0   2   1   0   2   0   1   2   1   2   0   1   0   2   2   1   2   1   2   1   0   1   2   1   1   1   1   2   2   0   0   0   2   1   0   0   0   2   0   2   1   2   0   1   2   0   0   0   2   \\
1   2   2   0   0   2   1   1   0   1   0   0   2   1   0   2   0   1   2   0   1   1   2   2   1   2   2   1   2   0   0   2   0   0   2   1   1   1   2   1   0   0   0   0   1   0   0   2   0   2   2   0   1   2   1   0   2   1   0   2   1   0   1   1   2   1   2   2   2   1   1   0   \\
0   1   0   0   1   1   2   0   2   2   0   0   2   2   1   0   2   1   1   1   1   0   0   1   1   2   1   2   2   0   0   2   1   2   1   0   2   0   1   2   2   1   1   2   2   0   2   0   1   2   2   0   0   0   1   1   1   1   0   0   1   2   0   2   1   2   0   1   0   0   2   2   \\
0   1   0   2   0   1   2   0   2   1   2   1   1   2   1   0   2   2   0   2   2   1   2   1   2   0   1   2   2   1   2   0   1   0   0   0   0   0   1   1   2   2   1   2   1   0   0   2   0   1   1   1   2   1   0   2   2   1   0   2   1   0   0   1   0   1   0   0   0   1   2   2   \\
2   1   2   2   2   1   2   1   0   1   2   2   0   0   1   0   1   2   2   0   1   0   0   1   1   0   0   0   1   2   1   1   1   1   2   0   2   1   0   0   0   2   0   2   2   0   1   0   0   2   1   0   2   2   0   0   1   2   1   1   0   2   1   2   2   0   1   2   1   0   1   2   \\
2   0   1   2   0   2   1   0   2   0   2   1   2   0   1   0   0   1   1   0   0   0   1   2   2   2   0   2   0   1   2   0   1   1   2   1   0   2   0   2   1   0   2   0   0   0   0   1   1   2   1   1   2   0   2   1   1   2   0   2   0   1   2   1   1   1   1   2   0   2   1   2   \\  \ea  \\
 & \\
\zero_{2\times 72} & \ba{c} 1   1   1   1   1   1   1   1   1   1   1   1   1   1   1   1   1   1   1   1   1   1   1   1   1   1   1   1   1   1   1   1   1   1   1   1   1   1   1   1   1   1   1   1   1   1   1   1   1   1   1   1   1   1   1   1   1   1   1   1   1   1   1   1   1   1   1   1   1   1   1   1   \\
0   0   0   0   0   0   0   0   0   0   0   0   0   0   0   0   0
0   0   0   0   0   0   0   0   0   0   0   0   0   0   0   0   0
0   0   1   1   1   1   1   1   1   1   1   1   1   1   1   1   1
1   1   1   1   1   1   1   1   1   1   1   1   1   1   1   1   1
1   1   1   1   \\\ea

 \ea  $

  \\ \hline
 \et} \ec

\bc {\scriptsize  \tabcolsep=1pt \bt {c} \hline $D_{3B}^{*T}$, where
$D_{3B}^{*}  \in \cC_{3B}(72+72;2^{4}3^{45}\bullet 3^{1}\bullet \B)$

\\\hline
$ \ba{cc}
\d_0^{*T}  & \ba{c}1    0   1   1   0   1   0   0   0   1   0   0   1   0   1   1   1   1   0   1   1   1   1   1   1   0   1   1   0   0   0   0   0   1   0   0   1   1   1   0   0   1   1   0   1   0   1   0   1   1   1   1   0   1   0   0   0   0   1   1   0   0   0   0   0   1   1   0   0   0   0   1   \\
1   1   1   1   0   0   1   1   1   1   1   0   1   0   1   1   1   0   0   1   0   1   1   0   0   0   0   0   1   0   0   0   0   1   1   0   0   0   0   1   0   0   0   0   1   1   1   1   0   0   1   1   1   0   1   1   1   1   0   0   1   0   0   1   1   0   0   0   1   0   0   1   \\
1   0   0   0   1   0   1   0   0   1   1   1   1   0   1   0   0   1   1   0   1   1   0   1   0   1   1   0   0   0   0   1   0   1   0   1   0   0   0   0   1   0   0   0   0   0   1   1   1   1   1   1   1   1   0   1   1   1   0   1   0   0   1   0   0   0   1   1   0   1   1   0   \\
1   0   1   0   0   0   1   0   0   1   0   0   0   1   1   0   0   1   1   0   0   1   0   1   0   1   0   1   0   1   1   1   0   1   1   0   0   0   0   1   1   1   1   1   1   1   1   1   0   0   1   0   1   1   0   1   0   0   1   1   1   0   0   1   0   1   0   0   0   1   0   0   \\
2   2   1   2   2   1   2   1   1   2   1   1   1   0   0   2   0   2   0   0   2   2   0   1   1   1   2   1   2   0   0   1   0   0   0   0   0   2   2   2   1   1   2   0   2   2   1   2   0   1   1   0   0   2   1   1   0   2   0   2   0   1   0   1   2   0   0   0   1   2   1   1   \\
0   1   2   2   0   0   2   0   1   1   2   2   1   0   0   0   1   0   2   2   0   1   1   0   0   1   1   2   0   1   1   2   2   2   0   2   1   2   0   2   2   2   1   1   2   1   2   1   2   1   2   1   0   0   0   1   1   0   2   0   0   2   1   2   0   1   2   1   0   0   0   1   \\
1   1   0   1   0   1   0   1   2   2   1   2   0   2   2   2   1   2   1   2   1   1   0   1   2   0   0   2   1   0   1   1   2   2   0   0   2   1   0   1   1   0   1   1   1   1   0   0   1   2   2   0   0   2   2   2   1   1   2   2   2   2   0   0   2   0   0   0   0   2   0   0   \\
2   1   0   0   1   1   1   1   2   2   1   2   1   0   2   1   0   0   0   1   0   0   0   2   0   1   2   2   0   0   2   1   0   2   2   2   0   2   2   1   0   1   0   0   2   1   2   1   1   2   1   0   2   1   0   2   0   0   1   0   0   2   1   2   1   1   0   2   1   1   2   2   \\
1   1   0   1   0   0   1   1   0   2   1   0   2   2   0   0   1   2   2   2   2   0   0   1   0   1   1   1   2   0   1   0   1   2   2   2   2   0   1   0   1   2   0   0   0   0   2   2   1   0   2   2   0   1   1   0   0   1   1   2   1   2   2   2   0   1   1   1   2   2   0   2   \\
1   2   1   1   1   1   2   0   2   0   1   0   0   2   0   2   0   2   0   1   2   1   1   2   1   1   0   1   2   0   0   2   2   1   1   0   1   0   1   0   1   0   2   1   0   0   0   2   2   2   1   2   0   0   2   2   2   1   2   0   0   1   1   2   0   1   2   1   0   2   0   2   \\
0   1   1   2   1   2   2   1   1   1   1   0   2   0   2   0   1   2   1   0   2   0   2   1   0   0   0   2   0   2   2   2   1   2   1   1   1   0   0   0   2   1   2   1   2   0   0   2   0   2   2   2   0   1   0   1   2   0   1   1   0   0   0   0   2   0   1   2   2   1   1   1   \\
1   1   1   2   0   2   0   0   2   2   0   0   2   0   1   2   0   0   2   1   1   1   0   2   0   1   0   2   0   1   1   1   2   2   2   1   0   1   2   0   0   2   0   1   0   2   0   1   0   2   2   0   0   2   1   1   2   0   1   1   1   1   1   2   2   2   1   2   2   0   1   0   \\
2   2   0   0   0   1   2   0   0   1   1   2   1   1   0   0   2   2   2   1   1   2   1   1   2   2   0   2   2   1   0   0   1   0   0   1   0   1   0   1   0   2   0   2   2   2   2   0   0   0   1   2   0   2   2   1   1   1   1   1   0   2   0   0   1   1   2   1   2   0   1   2   \\
0   2   1   1   2   0   2   0   1   1   1   1   2   0   2   0   0   1   2   0   2   1   1   2   2   0   1   0   0   1   1   0   0   2   2   1   1   0   0   1   1   0   2   2   1   2   2   0   2   2   1   0   0   1   2   0   0   2   2   2   0   2   1   0   2   2   0   0   1   1   1   1   \\
0   1   2   0   0   1   2   2   0   2   1   1   0   0   2   2   2   0   0   0   2   1   1   1   2   0   1   2   0   0   1   0   1   0   1   2   0   0   2   1   0   2   1   0   2   2   1   1   0   1   2   1   2   1   2   0   2   2   1   0   1   1   2   1   1   0   1   2   0   2   2   0   \\
2   2   1   0   2   1   1   1   1   0   0   0   2   1   2   0   0   1   1   0   0   1   2   0   2   2   1   2   2   2   0   2   0   0   1   2   0   1   0   2   0   1   2   0   2   0   0   2   1   2   1   1   0   2   2   1   0   1   2   1   2   0   1   0   1   1   1   2   0   0   1   2   \\
1   2   2   1   1   2   2   1   0   0   0   1   0   0   0   2   1   2   2   1   2   1   1   1   2   2   2   0   0   1   0   1   2   1   0   2   2   2   0   2   1   1   0   0   0   0   2   0   0   1   2   1   2   2   0   0   2   0   2   1   1   1   0   0   1   0   1   1   2   1   0   2   \\
1   0   0   2   1   0   2   0   0   0   2   0   2   2   2   0   2   2   2   1   2   1   1   0   1   1   2   1   1   0   0   1   1   2   0   1   1   2   1   0   0   1   2   2   1   1   1   0   0   1   2   0   0   0   2   2   2   0   1   1   2   0   1   2   1   0   0   2   1   0   2   2   \\
2   2   0   2   0   2   1   0   2   1   0   1   2   1   2   1   2   2   1   1   1   0   2   0   0   1   0   2   0   1   1   0   0   0   2   0   1   2   0   1   2   2   2   0   0   0   1   1   2   0   0   2   0   1   2   1   2   1   1   2   0   1   2   1   1   0   0   1   1   2   2   0   \\
1   0   2   2   0   1   2   1   1   1   2   2   2   2   2   0   0   0   1   1   0   0   0   0   1   1   2   0   2   2   2   1   1   1   0   2   2   1   2   0   1   2   2   1   0   2   2   0   1   0   0   1   0   1   2   0   2   0   2   2   1   0   0   0   1   2   0   0   1   1   1   1   \\
0   2   1   2   2   0   0   2   0   2   1   0   0   0   2   1   2   0   1   2   1   2   0   2   0   1   2   1   1   1   1   0   0   1   0   1   1   0   1   2   1   2   0   2   1   0   0   2   1   0   2   1   1   0   1   1   2   0   2   1   0   2   2   2   1   0   1   0   0   2   2   2   \\
0   2   2   1   1   1   2   2   2   0   0   1   0   1   0   0   1   1   0   2   0   1   0   2   0   0   0   2   2   0   1   2   1   2   0   1   2   0   2   2   1   0   2   0   1   1   1   0   0   2   2   1   2   1   1   1   1   2   2   1   0   1   2   0   2   1   0   2   0   0   2   1   \\
1   1   2   1   2   0   0   2   1   2   2   1   2   2   2   2   1   2   1   0   0   0   0   2   1   2   0   0   1   0   1   0   2   0   1   1   1   1   2   1   0   2   2   2   0   0   2   1   0   1   1   1   0   0   1   0   2   0   0   0   1   2   2   0   2   1   0   2   1   1   0   2   \\
1   2   2   1   0   1   2   0   0   0   2   0   2   1   1   0   0   0   1   2   0   2   1   2   1   1   1   0   2   0   2   1   2   1   2   2   1   0   2   2   0   2   1   0   0   1   1   1   0   2   0   2   1   2   1   2   0   2   1   0   0   1   0   0   1   2   0   1   0   2   2   1   \\
2   2   1   0   1   2   2   0   2   1   2   1   2   2   0   0   1   0   2   2   1   1   0   0   2   0   0   1   1   0   2   1   1   1   0   1   1   2   2   0   2   2   1   1   2   0   0   1   0   2   0   1   0   0   1   2   0   0   1   2   2   2   0   1   1   1   0   2   0   2   1   0   \\
2   0   1   2   0   0   2   1   2   0   0   1   2   1   1   0   0   1   0   1   2   1   1   1   0   0   2   1   2   1   2   2   2   0   1   0   2   1   1   2   2   2   0   1   0   0   1   1   1   0   2   1   2   2   2   1   0   0   0   2   0   2   2   0   1   2   0   1   1   0   2   0   \\
2   2   0   0   2   0   1   0   1   1   2   0   1   2   0   2   1   1   1   1   0   2   1   0   2   1   1   0   1   2   1   2   2   2   0   0   0   0   2   2   0   2   2   1   0   0   2   0   2   1   1   2   2   0   2   2   2   1   0   1   1   1   0   1   1   1   1   0   0   0   2   0   \\
0   2   0   2   1   2   1   1   0   2   0   0   1   0   2   0   1   1   2   2   0   1   0   1   2   2   0   1   2   1   2   0   1   1   2   0   1   0   1   1   2   0   2   1   2   0   2   2   1   1   1   2   2   1   1   0   0   0   2   0   1   1   2   0   0   0   2   2   1   0   2   0   \\
1   2   2   2   0   0   0   1   0   0   1   1   1   0   1   0   2   1   1   1   0   1   0   2   0   2   0   1   2   2   0   2   0   2   1   2   1   1   0   0   0   2   1   2   2   0   0   1   2   0   2   0   2   1   2   0   1   2   2   0   0   2   0   1   1   1   1   2   2   2   1   1   \\
2   1   2   0   2   1   1   2   2   0   0   1   2   2   1   1   2   2   0   2   2   1   1   0   0   1   0   1   1   2   0   0   0   2   0   1   1   1   2   2   1   0   0   0   0   2   1   2   2   2   1   0   0   0   2   0   2   0   2   1   1   0   0   1   0   1   2   1   0   2   1   1   \\
1   2   1   1   0   2   0   1   1   2   2   0   2   2   0   2   1   1   0   0   0   0   0   0   1   1   2   2   0   2   0   2   1   1   2   2   2   1   0   1   2   0   2   0   1   1   2   1   0   2   0   1   0   1   2   0   1   1   0   0   2   1   2   2   2   1   2   1   0   0   1   0   \\
2   0   2   1   1   0   2   0   2   0   0   1   2   1   1   2   1   1   1   1   0   0   0   2   1   1   2   1   2   2   1   0   0   0   0   2   0   0   2   2   2   1   2   0   0   1   0   0   0   1   0   1   2   0   0   0   2   1   2   2   2   2   1   1   2   1   2   2   1   0   1   1   \\
2   2   0   1   2   2   2   1   1   2   1   2   0   0   2   0   0   1   2   1   0   1   0   2   1   2   1   0   2   1   0   0   2   0   1   0   0   0   2   2   1   0   1   2   2   0   2   0   2   0   2   0   0   0   1   1   2   1   1   2   0   0   1   2   1   1   1   1   1   1   2   0   \\
0   2   0   0   1   1   2   1   1   2   0   1   2   0   1   1   2   1   0   2   0   0   1   1   2   0   0   0   2   2   0   2   0   2   2   1   1   2   2   0   1   1   1   0   1   0   2   1   0   0   2   0   2   0   1   1   2   1   1   2   0   2   1   1   0   1   2   2   0   2   0   2   \\
1   2   2   0   2   1   1   0   1   0   1   0   2   2   2   0   2   1   0   1   1   2   0   1   0   1   2   0   0   0   2   1   0   0   1   2   1   2   2   0   0   1   0   1   1   1   0   0   0   0   2   2   1   2   1   0   2   2   1   0   1   0   0   2   2   2   2   1   2   2   1   1   \\
0   1   1   2   2   0   2   0   2   1   2   2   0   2   1   0   1   1   1   1   1   0   2   1   1   1   1   2   0   0   0   1   0   0   2   2   2   1   1   1   1   2   2   0   2   0   0   2   2   0   2   0   2   2   1   1   2   0   0   0   2   1   1   0   2   1   2   0   1   0   0   0   \\
1   0   0   1   1   0   1   2   0   1   2   0   0   1   0   1   1   0   0   0   2   1   0   0   2   2   2   2   2   2   1   1   0   2   1   2   0   2   1   0   1   2   0   2   0   0   2   2   1   2   1   0   0   1   0   2   2   2   1   2   2   1   0   1   1   2   0   1   2   0   1   1   \\
2   1   0   0   0   2   1   1   0   2   0   2   2   0   0   0   1   2   1   0   2   1   1   2   0   2   1   1   2   1   0   0   1   2   0   1   2   0   1   2   1   2   2   0   0   2   1   1   2   1   0   0   1   0   1   2   2   0   0   0   0   0   2   2   1   1   2   2   2   1   1   1   \\
1   1   0   2   2   0   2   2   1   2   0   0   0   0   1   1   0   1   2   2   1   1   2   0   0   0   1   0   2   0   1   1   0   1   2   2   1   2   1   1   1   1   2   1   1   0   2   0   0   2   0   1   0   2   2   2   2   0   2   0   2   0   0   1   0   2   1   1   1   2   2   0   \\
0   2   0   0   1   0   0   2   0   0   1   2   2   2   0   1   1   1   2   1   0   2   2   2   0   1   2   0   1   1   2   2   1   1   1   0   1   1   2   0   1   0   0   0   2   1   2   0   0   1   1   1   1   2   2   1   2   0   2   2   2   0   2   1   1   1   0   0   2   0   2   0   \\
1   1   2   0   2   2   0   0   2   1   1   1   2   0   2   2   1   1   2   0   0   0   0   0   2   0   2   2   0   1   1   1   1   2   0   2   1   1   1   2   2   0   1   2   0   0   2   2   0   0   1   2   1   2   1   1   0   1   1   2   0   0   0   2   0   0   1   2   1   2   0   1   \\
2   1   2   0   1   2   0   1   0   2   2   1   1   2   2   0   2   0   1   0   1   2   2   0   1   0   1   0   2   1   0   0   0   0   1   0   1   2   1   0   2   1   1   2   1   1   0   0   1   2   1   1   2   2   0   1   2   0   2   0   1   2   2   1   0   2   0   1   2   0   2   0   \\
1   2   0   2   0   1   0   0   0   0   1   2   2   1   1   2   0   1   2   2   2   2   0   1   2   0   0   1   1   1   1   1   0   0   1   1   1   2   1   0   2   2   1   2   0   1   1   2   2   2   1   1   2   2   2   0   2   2   0   0   2   1   0   2   0   0   0   0   0   1   1   0   \\
1   0   2   2   1   1   2   0   1   2   0   2   2   0   0   0   2   2   1   1   0   0   0   1   2   0   0   1   2   1   1   2   2   0   0   2   1   2   1   1   0   0   0   2   1   0   1   1   2   2   0   1   2   0   0   1   1   1   2   0   1   0   2   2   2   1   2   0   2   0   1   1   \\
0   0   0   2   0   0   0   0   1   2   2   2   1   1   2   2   2   0   2   0   2   2   1   1   1   1   1   0   1   2   0   1   0   2   1   1   1   1   2   1   1   2   0   1   2   1   2   1   0   2   1   2   0   0   0   0   0   2   2   0   2   1   2   0   2   0   1   0   0   1   1   2   \\
1   1   2   1   0   0   1   1   2   0   0   1   1   2   0   1   2   0   2   0   0   0   1   2   0   2   2   1   2   1   0   0   2   2   2   0   2   1   1   0   2   0   2   1   2   1   0   0   0   2   0   1   1   2   0   1   2   2   1   2   0   1   1   2   0   0   1   0   2   1   2   1   \\
1   1   2   0   2   0   1   0   2   0   1   2   0   0   1   1   1   0   2   1   2   0   2   1   0   0   2   1   2   1   2   1   2   0   0   0   1   2   2   0   0   0   2   0   1   2   1   2   1   2   2   1   0   2   2   0   0   2   0   2   1   1   1   2   0   1   0   2   1   1   1   0   \\
1   2   2   0   0   1   1   2   1   0   2   2   0   2   1   0   1   2   1   2   1   0   2   0   1   2   0   0   0   1   2   0   0   1   1   1   0   0   1   0   2   2   2   2   2   1   2   2   0   2   2   0   0   1   1   1   0   1   0   2   0   1   1   1   1   2   2   0   0   0   2   1   \\
1   2   1   0   1   0   2   1   2   2   2   0   0   1   2   1   0   2   2   0   0   0   1   1   1   2   2   2   1   1   1   0   0   1   0   0   2   2   0   0   2   2   0   0   1   2   0   1   1   2   1   2   0   1   1   0   2   2   2   0   0   2   0   0   2   1   1   2   1   0   1   1   \\  \ea  \\
 & \\
\zero_{3\times 72} & \ba{c} 1   1   1   2   1   2   2   2   2   2   2   1   1   2   2   1   1   1   2   1   2   2   1   2   1   2   1   1   1   2   1   2   1   1   1   2   2   2   2   2   1   1   2   1   2   1   1   1   1   2   1   2   2   2   2   1   1   1   1   2   1   2   1   2   1   2   2   2   2   1   1   2   \\
1   1   1   1   1   1   1   1   1   1   1   1   1   1   1   1   1   1   1   1   1   1   1   1   1   1   1   1   1   1   1   1   1   1   1   1   1   1   1   1   1   1   1   1   1   1   1   1   1   1   1   1   1   1   1   1   1   1   1   1   1   1   1   1   1   1   1   1   1   1   1   1   \\
0   0   0   0   0   0   0   0   0   0   0   0   0   0   0   0   0
0   0   0   0   0   0   0   0   0   0   0   0   0   0   0   0   0
0   0   1   1   1   1   1   1   1   1   1   1   1   1   1   1   1
1   1   1   1   1   1   1   1   1   1   1   1   1   1   1   1   1
1   1   1   1       \\\ea

 \ea  $

  \\ \hline
 \et} \ec

\bc {\scriptsize  \tabcolsep=1pt \bt {c} \hline $D_{3B}^{*T}$, where
$D_{3B}^{*}  \in \cC_{3B}(72+72;2^{4}3^{45}\bullet 3^{2}\bullet \B)$

\\\hline
$ \ba{cc}
\d_0^{*T}  & \ba{c}0    1   1   1   1   1   1   1   0   1   0   1   0   1   0   0   0   1   0   1   0   1   1   1   0   0   0   0   1   0   0   0   0   1   1   1   0   1   0   0   0   1   1   0   1   1   0   0   1   0   1   1   1   1   0   1   0   1   1   1   0   0   0   0   0   1   1   0   0   0   1   0   \\
0   0   0   1   1   0   0   0   0   0   1   0   1   0   0   1   1   1   1   1   0   0   0   0   1   1   0   1   1   1   0   1   1   1   1   1   1   0   0   1   1   1   1   0   1   1   1   1   1   0   1   1   0   0   1   0   0   1   0   0   0   0   0   0   0   0   1   0   0   0   1   1   \\
1   0   1   0   0   1   1   0   1   0   0   0   1   1   0   1   1   1   1   0   0   1   0   0   0   1   1   1   1   0   1   0   0   1   0   0   0   1   1   1   0   1   1   0   0   0   1   1   0   1   0   1   0   1   0   1   0   0   1   0   1   1   0   1   1   1   0   0   0   0   1   0   \\
1   1   0   0   0   0   0   0   1   1   0   0   1   1   1   1   0   1   0   0   0   1   1   1   1   0   1   0   1   0   1   0   1   0   1   1   0   0   0   0   1   1   1   0   0   0   1   1   1   1   1   0   0   0   1   1   1   0   0   0   1   0   0   1   0   0   1   0   1   1   1   0   \\
1   1   1   2   0   1   1   0   0   0   0   2   0   2   1   0   2   2   1   2   2   1   2   1   1   0   1   0   0   1   2   2   1   1   2   0   2   0   0   2   2   1   2   1   1   2   2   2   0   0   1   0   0   0   0   2   0   0   1   1   0   1   0   2   2   2   2   1   2   1   0   1   \\
0   2   1   1   1   0   2   0   1   2   0   2   2   0   1   2   0   0   0   1   0   1   2   1   2   2   1   1   1   1   2   0   1   2   2   2   0   2   1   0   1   0   0   1   1   2   0   2   0   1   2   1   2   1   2   0   0   0   0   0   2   0   0   1   2   2   0   2   1   1   1   2   \\
0   1   0   0   2   2   1   1   1   2   0   1   0   0   0   2   2   1   0   1   1   2   0   2   2   2   1   1   2   0   1   1   2   1   1   0   0   2   1   1   2   1   0   2   0   2   1   1   2   0   1   0   0   1   2   0   2   0   0   1   0   1   2   2   2   2   2   1   0   2   0   0   \\
0   1   0   1   1   2   1   1   0   2   2   0   0   0   0   0   2   0   0   0   0   0   1   0   2   1   1   1   2   0   2   2   1   2   1   2   2   2   2   0   1   1   2   0   2   0   0   1   2   2   1   1   0   1   1   1   0   0   1   1   1   2   2   2   0   2   1   2   2   1   2   0   \\
0   0   1   2   0   1   2   0   0   1   1   1   2   0   2   1   2   2   1   0   0   2   2   1   1   2   0   1   0   1   0   1   1   1   2   2   1   2   0   2   0   1   0   0   1   0   1   1   1   2   2   0   1   2   0   1   2   2   2   1   1   2   2   0   0   0   0   2   2   2   0   0   \\
2   1   0   0   2   1   2   0   2   0   2   2   0   1   0   1   2   2   0   1   1   2   2   1   0   1   0   2   1   0   0   1   2   1   2   2   2   2   0   1   1   0   0   2   1   0   1   1   2   1   0   1   1   0   0   0   1   0   1   2   2   0   0   2   1   0   2   1   1   0   2   2   \\
2   0   2   0   0   2   2   1   0   1   0   2   2   0   2   1   0   1   1   2   1   2   1   0   1   2   1   1   0   2   0   2   1   0   1   0   0   1   2   2   0   0   1   1   1   2   0   2   2   1   2   1   1   0   0   2   2   1   0   1   0   1   0   2   0   1   2   1   0   0   2   2   \\
0   1   1   1   0   1   2   2   1   1   2   0   2   2   0   1   1   1   0   1   0   0   0   2   1   2   0   2   2   0   0   1   2   0   0   0   0   0   2   2   1   2   1   2   2   1   1   0   2   1   2   0   2   0   2   0   0   1   2   2   1   2   1   2   0   0   1   1   1   0   2   1   \\
2   0   1   2   0   2   1   1   2   1   2   2   0   2   0   0   1   1   0   1   0   1   2   2   2   1   1   1   0   2   1   0   1   2   0   0   1   2   0   1   0   0   0   1   1   2   0   2   0   2   2   2   0   2   1   1   0   0   0   2   1   2   0   0   2   2   1   0   1   1   2   1   \\
2   2   1   0   0   1   2   1   0   0   2   0   0   0   1   2   2   1   0   1   0   2   2   2   1   2   1   2   0   2   1   2   0   0   1   1   2   0   1   1   1   2   0   0   1   2   0   1   1   2   1   0   2   1   0   0   2   1   0   0   0   2   2   1   1   2   2   1   0   2   1   0   \\
2   0   2   1   2   0   0   2   0   0   1   1   0   0   0   2   0   1   1   1   0   2   0   1   1   1   1   0   2   2   2   2   1   2   2   1   1   1   1   1   0   1   1   0   0   2   2   0   0   2   0   0   2   2   2   2   2   1   1   2   1   2   0   0   1   0   2   1   0   1   2   0   \\
2   2   2   2   0   1   2   1   0   0   1   2   0   2   1   1   1   1   1   0   2   2   1   0   2   0   2   0   2   0   0   2   1   1   1   0   1   2   2   2   1   2   2   1   2   0   0   0   1   1   0   2   0   1   2   2   1   0   0   0   1   0   0   0   0   1   2   1   1   0   1   2   \\
0   0   2   0   1   2   2   2   1   0   1   1   1   2   0   2   1   1   2   1   2   0   2   2   2   0   1   2   1   1   0   0   0   0   0   2   2   2   1   0   0   1   2   0   1   0   0   2   0   0   1   0   2   2   2   0   1   1   1   1   2   0   0   0   2   1   1   1   2   1   2   1   \\
2   1   2   2   0   0   0   1   0   0   0   1   2   1   1   0   1   2   0   2   2   2   2   1   0   1   0   2   2   1   1   2   0   1   1   0   2   1   1   0   0   0   1   2   1   2   0   2   2   1   1   1   0   0   2   0   2   1   2   0   1   2   1   0   1   2   1   2   0   2   0   0   \\
0   0   1   0   1   2   1   0   2   2   0   2   1   1   1   0   1   2   2   1   0   2   0   2   1   2   0   2   0   1   0   1   0   0   2   1   2   1   2   0   0   1   0   1   0   1   2   0   2   1   2   2   0   1   1   2   2   2   0   0   1   1   0   1   2   0   2   2   1   1   2   0   \\
0   0   2   2   0   1   0   2   0   2   2   1   2   1   2   0   2   1   0   0   0   0   1   2   2   2   2   1   1   1   1   2   1   1   0   0   1   0   1   1   0   0   1   1   2   0   2   2   1   2   1   0   1   1   1   0   2   0   2   1   0   0   1   1   2   2   2   0   2   0   2   0   \\
1   1   1   2   1   1   0   1   2   1   0   1   1   0   1   2   1   2   0   2   0   0   2   1   0   0   2   1   0   0   0   0   1   2   2   1   2   2   2   0   1   1   1   0   0   0   1   1   1   1   2   1   2   0   2   0   0   2   2   0   2   2   2   2   1   0   2   0   0   2   0   2   \\
2   2   2   0   1   0   2   1   2   1   1   1   0   1   0   0   0   0   1   1   2   1   2   0   1   2   0   1   2   2   1   2   2   0   0   0   2   1   2   0   2   2   0   0   2   0   2   1   2   1   2   0   1   2   0   1   0   2   1   0   0   0   1   1   2   1   1   0   0   1   1   2   \\
2   1   0   1   1   2   1   0   1   0   0   0   2   0   1   0   0   0   1   2   2   0   1   2   1   2   0   1   2   1   2   1   2   0   2   0   0   1   2   1   0   0   1   2   1   0   2   1   2   2   0   1   2   2   1   2   1   0   2   1   2   0   1   0   2   0   2   1   0   2   1   0   \\
2   0   1   0   2   0   1   2   2   2   1   2   0   1   2   2   0   2   0   1   0   2   1   1   0   2   2   0   1   1   1   1   2   0   0   1   1   1   0   2   1   2   1   1   2   0   2   1   2   1   0   1   0   1   0   0   0   2   1   0   2   2   2   0   0   2   1   1   2   0   0   0   \\
1   0   0   2   1   2   2   0   2   1   2   1   0   2   0   0   0   2   1   0   0   1   1   0   1   2   2   0   1   2   1   1   1   0   2   2   0   0   1   0   0   1   2   2   2   1   2   0   0   1   1   1   1   0   2   0   2   1   0   1   1   0   2   2   1   2   0   2   1   2   2   0   \\
2   2   2   0   1   2   0   1   1   2   0   0   1   2   0   0   1   2   1   0   0   0   1   1   2   0   0   2   1   2   0   1   0   0   2   2   1   1   2   1   2   2   1   2   1   1   0   0   0   1   1   2   1   1   0   1   2   0   2   0   2   0   2   1   0   2   0   2   1   1   0   2   \\
0   0   0   1   1   1   1   2   2   0   2   0   1   1   2   2   2   1   0   1   2   2   1   0   1   2   2   2   2   1   0   0   1   1   2   0   2   0   2   1   1   1   2   2   0   0   2   2   0   1   2   2   0   0   0   1   1   0   2   0   1   1   2   0   1   2   0   0   0   1   0   1   \\
0   0   0   2   0   0   2   1   1   2   0   0   2   2   0   1   0   0   2   1   1   2   1   1   1   2   2   1   2   2   1   0   1   0   2   1   2   2   2   2   2   1   0   0   1   1   2   1   2   1   0   1   2   0   0   1   1   0   0   1   2   2   1   1   0   0   2   2   0   0   0   1   \\
2   2   1   1   0   1   2   2   1   1   1   2   1   0   2   2   0   2   0   2   1   1   1   0   2   2   0   0   0   2   1   1   0   0   1   0   2   2   0   1   2   1   2   0   1   0   0   2   1   1   0   1   2   1   1   0   1   0   2   2   0   2   0   0   2   0   0   0   2   1   2   1   \\
2   1   1   1   1   2   2   2   0   1   0   1   0   0   0   1   2   1   1   0   2   0   1   0   0   2   2   2   1   1   0   2   1   0   0   2   0   2   0   0   0   0   0   0   1   2   2   1   2   1   1   2   0   1   0   2   1   2   2   0   1   1   1   1   1   0   2   0   2   2   2   2   \\
1   2   2   2   0   0   2   2   2   0   2   1   1   1   0   1   2   0   1   1   0   2   0   2   2   1   0   1   0   1   0   1   2   0   1   0   0   0   2   0   1   2   1   2   2   0   2   0   1   0   0   1   0   1   2   2   2   2   0   1   0   2   0   1   2   1   0   1   1   1   2   1   \\
2   2   1   2   2   0   0   2   1   0   1   1   1   2   1   0   2   0   2   1   0   0   1   0   0   1   2   2   1   0   2   0   0   1   0   2   0   1   0   1   2   1   0   2   2   1   2   0   1   2   1   1   2   0   0   0   2   0   0   1   1   2   1   0   1   0   2   2   1   1   2   2   \\
0   2   1   2   1   1   0   2   1   0   0   2   0   0   1   0   0   1   1   1   2   2   2   2   2   1   1   2   2   2   0   0   0   2   1   1   1   0   2   1   1   1   2   1   0   0   2   0   1   1   0   0   0   0   1   0   0   0   2   1   2   2   2   2   1   1   1   0   0   2   2   2   \\
0   1   2   2   1   2   2   2   2   1   2   0   0   2   0   0   0   0   1   0   1   1   1   0   0   2   0   2   0   2   2   1   1   1   1   1   0   2   1   1   2   0   2   0   0   2   2   0   1   2   0   0   0   0   1   0   2   2   1   2   2   1   0   1   1   1   1   0   1   2   1   2   \\
1   1   2   2   1   2   0   2   2   0   0   2   2   1   0   0   2   1   1   0   0   1   2   1   0   2   2   0   0   1   0   2   2   0   1   0   0   1   0   1   2   0   2   1   1   0   1   1   1   0   1   0   2   1   1   0   2   2   2   2   2   0   1   2   0   2   0   1   1   0   2   1   \\
1   2   2   1   2   0   1   1   1   1   2   0   2   2   0   0   1   2   1   1   2   0   1   1   2   2   2   0   0   0   0   2   0   2   2   1   1   2   1   0   0   1   2   2   0   1   1   1   1   1   0   0   2   1   0   2   2   0   0   1   0   0   2   2   2   1   0   0   0   1   2   0   \\
2   1   2   0   0   1   0   2   0   1   1   0   2   2   1   0   2   2   0   1   0   2   1   0   0   2   1   1   1   1   0   1   2   2   0   2   0   1   0   0   0   1   2   0   1   0   1   1   0   0   2   2   1   2   1   1   2   0   0   2   1   1   2   1   2   0   2   1   2   2   0   2   \\
0   2   1   2   1   0   0   0   2   2   1   1   0   2   2   1   1   2   0   0   1   1   1   0   1   2   1   1   0   2   0   0   2   2   0   2   2   0   2   1   0   0   2   2   2   2   1   0   1   1   1   2   0   2   0   1   2   0   1   1   1   2   0   2   0   0   2   1   0   0   1   1   \\
0   1   2   0   2   1   2   0   1   2   1   1   1   0   0   0   2   2   2   1   1   2   0   1   0   0   2   0   2   2   0   2   1   0   1   1   0   1   2   2   0   0   1   1   1   1   1   2   2   2   1   2   2   0   0   0   1   0   0   2   1   1   2   0   2   1   0   0   2   0   1   2   \\
0   1   2   2   1   0   2   0   0   0   1   2   1   1   2   0   1   2   0   0   2   0   1   1   2   0   1   2   2   1   2   0   1   0   1   1   0   0   1   1   1   0   0   0   0   2   2   1   2   0   0   1   2   2   0   1   2   1   0   1   2   1   2   0   2   2   2   2   1   1   2   0   \\
0   0   0   1   0   1   2   1   1   2   0   2   2   0   0   1   1   0   2   1   2   1   0   2   1   0   2   2   0   1   2   1   1   0   1   1   0   1   1   1   2   2   2   2   1   0   1   0   1   0   2   2   2   1   0   0   2   0   1   0   2   2   0   2   1   0   2   0   2   1   2   0   \\
2   1   1   2   2   2   0   2   1   1   2   0   0   1   1   0   1   1   0   2   1   2   0   2   0   0   2   1   0   1   2   0   1   2   1   0   1   1   2   2   0   1   0   0   0   1   0   2   2   0   0   2   2   2   1   0   2   2   1   0   2   0   1   2   1   0   0   1   1   0   1   2   \\
0   0   1   2   1   1   0   2   2   2   0   1   0   2   2   2   2   2   0   0   2   1   0   1   1   1   1   0   1   2   1   2   0   0   1   0   1   1   0   0   1   2   2   2   0   2   1   2   2   0   2   1   2   0   0   1   1   1   0   0   2   1   0   1   0   1   0   1   0   2   2   2   \\
0   2   0   2   1   0   1   1   2   1   0   0   0   2   1   1   1   0   1   2   0   2   1   0   2   2   2   2   1   1   1   2   0   2   1   2   2   1   0   1   0   1   0   1   2   0   1   0   2   2   1   0   2   0   0   2   1   2   1   2   2   0   0   1   2   2   0   0   1   0   0   1   \\
1   1   2   2   2   1   1   0   0   1   2   1   2   0   2   2   1   2   0   2   0   2   1   0   0   1   1   2   0   1   2   0   0   1   0   0   1   2   2   1   2   0   1   1   2   2   2   0   0   1   2   1   0   0   2   1   0   1   0   2   2   0   1   0   2   2   0   1   0   1   1   0   \\
2   0   1   1   2   0   1   0   1   2   0   0   2   2   2   1   2   0   0   2   2   0   0   1   2   0   1   1   1   1   0   0   2   1   2   0   2   1   0   2   0   0   1   0   0   0   2   0   1   1   1   2   2   2   0   0   1   1   0   1   0   2   2   1   1   2   1   1   1   2   2   2   \\
0   0   2   1   1   1   2   0   1   0   0   2   0   2   2   1   2   1   0   0   2   0   1   1   0   1   0   2   1   0   1   2   2   2   1   2   1   0   1   2   2   2   0   0   1   1   0   0   0   0   0   2   2   0   0   1   2   0   2   1   2   1   1   2   1   2   0   2   1   2   1   1   \\
0   2   1   0   2   2   2   1   1   0   1   2   2   1   2   2   0   0   1   1   0   2   2   1   2   2   2   0   0   0   1   1   1   1   0   0   1   0   2   1   0   0   2   0   0   2   1   0   2   1   1   1   0   2   1   0   2   0   2   0   0   1   1   0   0   2   2   1   2   1   1   2   \\
2   1   2   0   2   1   1   1   0   0   1   2   1   2   1   0   0   1   1   0   0   0   2   1   2   2   2   2   1   0   1   1   0   2   2   0   1   2   2   0   2   2   0   2   2   2   1   2   2   0   0   1   0   0   0   0   2   1   0   1   1   1   0   1   0   1   1   0   2   2   0   1   \\  \ea  \\
 & \\
\zero_{4\times 72} & \ba{c} 1   2   1   2   1   2   1   2   1   2   2   2   2   2   2   2   1   1   1   1   1   2   1   1   1   1   2   2   2   1   1   2   2   1   2   1   2   2   1   2   1   2   1   2   1   2   1   2   1   2   1   1   1   2   1   1   1   2   1   2   1   1   1   2   2   1   2   2   1   2   2   2   \\
1   2   2   2   2   2   1   2   1   1   1   1   2   1   1   2   1   2   1   1   2   2   1   1   1   2   1   1   1   2   2   1   2   1   1   2   2   2   2   2   1   2   1   2   1   1   2   2   1   2   2   2   2   1   1   1   1   1   1   2   1   1   1   2   1   2   2   2   2   2   2   1   \\
1   1   1   1   1   1   1   1   1   1   1   1   1   1   1   1   1   1   1   1   1   1   1   1   1   1   1   1   1   1   1   1   1   1   1   1   1   1   1   1   1   1   1   1   1   1   1   1   1   1   1   1   1   1   1   1   1   1   1   1   1   1   1   1   1   1   1   1   1   1   1   1   \\
0   0   0   0   0   0   0   0   0   0   0   0   0   0   0   0   0
0   0   0   0   0   0   0   0   0   0   0   0   0   0   0   0   0
0   0   1   1   1   1   1   1   1   1   1   1   1   1   1   1   1
1   1   1   1   1   1   1   1   1   1   1   1   1   1   1   1   1
1   1   1   1   \\\ea

 \ea  $

  \\ \hline
 \et} \ec

\end{document}